\documentclass[12pt,oneside,letterpaper]{book}

\usepackage[matrix,arrow,curve,tips,frame]{xy}
\usepackage{epsfig}
\usepackage{doublespace}
\usepackage{vmargin} 

\usepackage[leqno]{amsmath}
\usepackage{amssymb}
\usepackage{enumerate} 
\usepackage{proof}
\usepackage{fancybox} 

\usepackage{theorem}
\usepackage{QED}

\newtheorem{theorem}{Theorem}[section]
\newtheorem{lemma}[theorem]{Lemma}
\newtheorem{proposition}[theorem]{Proposition}
\newtheorem{corollary}[theorem]{Corollary}
{\theorembodyfont{\rmfamily} \newtheorem{definition}[theorem]{Definition}}
{\theorembodyfont{\rmfamily} \newtheorem{example}[theorem]{Example}}
{\theorembodyfont{\rmfamily} \newtheorem{remark}[theorem]{Remark}}

\usepackage{euscript}
\usepackage{stmaryrd}

\newenvironment{player}{\( \left\{\begin{array}{c}}{\end{array}\right\} \)}
\newenvironment{opponent}{\( \left(\begin{array}{c}}{\end{array}\right) \)}

\newcommand{\<}{\ensuremath{\langle}}
\renewcommand{\>}{\ensuremath{\rangle}}

\newcommand{\cat}[1]{\ensuremath{\mathbf{#1}}}
\newcommand{\category}[1]{\ensuremath{\mathbf{#1}}}

\newcommand{\script}[1]{\ensuremath{\EuScript {#1}}}
\newcommand{\gothic}[1]{\ensuremath{\mathfrak{#1}}}
\newcommand{\Set}{\ensuremath{\mathbf{Set}}}
\newcommand{\id}{\ensuremath{\mathrm{id}}}
\newcommand{\Hom}{\ensuremath{\mathrm{Hom}}}
\newcommand{\bag}{\ensuremath{\mathrm{bag}}}
\newcommand{\cut}{\ensuremath{\mathrm{cut}}}
\newcommand{\op}{\ensuremath{\mathrm{op}}}
\newcommand{\ob}{\ensuremath{\mathrm{ob}}}
\newcommand{\dom}{\ensuremath{\mathrm{dom}}}
\newcommand{\cod}{\ensuremath{\mathrm{cod}}}

\newcommand{\oi}{\ensuremath{\gothic{I}}}

\newcommand{\x}{\ensuremath{\times}}

\newcommand{\vd}{\ensuremath{\vdash}}

\newcommand{\meet}{\ensuremath{\wedge}}
\newcommand{\join}{\ensuremath{\vee}}

\newcommand{\bigjoin}{\ensuremath{\bigvee}}
\newcommand{\prefix}{\ensuremath{\sqsubseteq}}
\newcommand{\iso}{\ensuremath{\cong}}

\newcommand{\Ra}{\ensuremath{\Rightarrow}}

\newcommand{\ra}{\ensuremath{\rightarrow}}

\newcommand{\Lra}{\ensuremath{\Longrightarrow}}

\newcommand{\Lla}{\ensuremath{\Longleftarrow}}

\newcommand{\pc}{\ensuremath{\Mapstochar=\!=\hspace{-1.2ex}=\Mapsfromchar}}
\newcommand{\ora}[1]{\ensuremath{\overrightarrow {#1}}}
\newcommand{\ola}[1]{\ensuremath{\overleftarrow {#1}}}
\newcommand{\ol}[1]{\ensuremath{\overline {#1}}}

\newcommand{\dual}[1]{\ensuremath{#1^{*}}}

\newcommand{\hgt}{\ensuremath{\mathrm{hgt}}}
\newcommand{\cuthgt}{\ensuremath{\mathrm{cuthgt}}}

\newcommand{\bs}{\ensuremath{\backslash}}

\newcommand{\qqquad}{\quad \quad \quad}
\newcommand{\qqqquad}{\quad \quad \quad \quad}

\newcommand{\centro}[1]{\begin{center} #1 \end{center}}

\newcommand{\senta}[1]{\begin{center} #1 \end{center}}

\makeatletter

\newdimen\w@dth

\def\setw@dth#1#2{\setbox\z@\hbox{\scriptsize $#1$}\w@dth=\wd\z@
\setbox\@ne\hbox{\scriptsize $#2$}\ifnum\w@dth<\wd\@ne \w@dth=\wd\@ne \fi
\advance\w@dth by 1.2em}

\def\t@^#1_#2{\allowbreak\def\n@one{#1}\def\n@two{#2}\mathrel
{\setw@dth{#1}{#2}
\mathop{\hbox to \w@dth{\rightarrowfill}}\limits
\ifx\n@one\empty\else ^{\box\z@}\fi
\ifx\n@two\empty\else _{\box\@ne}\fi}}
\def\t@@^#1{\@ifnextchar_ {\t@^{#1}}{\t@^{#1}_{}}}

\def\t@left^#1_#2{\def\n@one{#1}\def\n@two{#2}\mathrel{\setw@dth{#1}{#2}
\mathop{\hbox to \w@dth{\leftarrowfill}}\limits
\ifx\n@one\empty\else ^{\box\z@}\fi
\ifx\n@two\empty\else _{\box\@ne}\fi}}
\def\t@@left^#1{\@ifnextchar_ {\t@left^{#1}}{\t@left^{#1}_{}}}

\def\two@^#1_#2{\def\n@one{#1}\def\n@two{#2}\mathrel{\setw@dth{#1}{#2}
\mathop{\vcenter{\hbox to \w@dth{\rightarrowfill}\kern-4ex
                 \hbox to \w@dth{\rightarrowfill}}%
       }\limits
\ifx\n@one\empty\else ^{\box\z@}\fi
\ifx\n@two\empty\else _{\box\@ne}\fi}}
\def\tw@@^#1{\@ifnextchar_ {\two@^{#1}}{\two@^{#1}_{}}}

\def\tofr@^#1_#2{\def\n@one{#1}\def\n@two{#2}\mathrel{\setw@dth{#1}{#2}
\mathop{\vcenter{\hbox to \w@dth{\rightarrowfill}\kern-4ex
                 \hbox to \w@dth{\leftarrowfill}}%
       }\limits
\ifx\n@one\empty\else ^{\box\z@}\fi
\ifx\n@two\empty\else _{\box\@ne}\fi}}
\def\t@fr@^#1{\@ifnextchar_ {\tofr@^{#1}}{\tofr@^{#1}_{}}}

\newdimen\W@dth
\def\setW@dth#1#2{\setbox\z@\hbox{$#1$}\W@dth=\wd\z@
\setbox\@ne\hbox{$#2$}\ifnum\W@dth<\wd\@ne \W@dth=\wd\@ne \fi
\advance\W@dth by 1.2em}

\def\T@^#1_#2{\allowbreak\def\N@one{#1}\def\N@two{#2}\mathrel
{\setW@dth{#1}{#2}
\mathop{\hbox to \W@dth{\rightarrowfill}}\limits
\ifx\N@one\empty\else ^{\box\z@}\fi
\ifx\N@two\empty\else _{\box\@ne}\fi}}
\def\T@@^#1{\@ifnextchar_ {\T@^{#1}}{\T@^{#1}_{}}}

\def\T@left^#1_#2{\def\N@one{#1}\def\N@two{#2}\mathrel{\setW@dth{#1}{#2}
\mathop{\hbox to \W@dth{\leftarrowfill}}\limits
\ifx\N@one\empty\else ^{\box\z@}\fi
\ifx\N@two\empty\else _{\box\@ne}\fi}}
\def\T@@left^#1{\@ifnextchar_ {\T@left^{#1}}{\T@left^{#1}_{}}}

\def\Tofr@^#1_#2{\def\N@one{#1}\def\N@two{#2}\mathrel{\setW@dth{#1}{#2}
\mathop{\vcenter{\hbox to \W@dth{\rightarrowfill}\kern-4ex
                 \hbox to \W@dth{\leftarrowfill}}%
       }\limits
\ifx\N@one\empty\else ^{\box\z@}\fi
\ifx\N@two\empty\else _{\box\@ne}\fi}}
\def\T@fr@^#1{\@ifnextchar_ {\Tofr@^{#1}}{\Tofr@^{#1}_{}}}

\def\Two@^#1_#2{\def\N@one{#1}\def\N@two{#2}\mathrel{\setW@dth{#1}{#2}
\mathop{\vcenter{\hbox to \W@dth{\rightarrowfill}\kern-4ex
                 \hbox to \W@dth{\rightarrowfill}}%
       }\limits
\ifx\N@one\empty\else ^{\box\z@}\fi
\ifx\N@two\empty\else _{\box\@ne}\fi}}
\def\Tw@@^#1{\@ifnextchar_ {\Two@^{#1}}{\Two@^{#1}_{}}}

\def\to{\@ifnextchar^ {\t@@}{\t@@^{}}}
\def\from{\@ifnextchar^ {\t@@left}{\t@@left^{}}}
\def\two{\@ifnextchar^ {\tw@@}{\tw@@^{}}}
\def\tofro{\@ifnextchar^ {\t@fr@}{\t@fr@^{}}}
\def\To{\@ifnextchar^ {\T@@}{\T@@^{}}}
\def\From{\@ifnextchar^ {\T@@left}{\T@@left^{}}}
\def\Two{\@ifnextchar^ {\Tw@@}{\Tw@@^{}}}
\def\Tofro{\@ifnextchar^ {\T@fr@}{\T@fr@^{}}}

\makeatother

\setpapersize{USletter}
\setmarginsrb{1.5in}{1.1in}{1.1in}{1.1in}{12pt}{25pt}{0pt}{30pt}

\newcommand{\gp}{\gothic{p}}
\newcommand{\gq}{\gothic{q}}
\newcommand{\gr}{\gothic{r}}

\newcommand{\gu}{\gothic{u}}

\newcommand{\alphatri}{
  \begin{picture}(16,15)(2,0)
  \put(0,0){\line(1,0){20}}
  \put(0,0){\line(2,3){10}}
  \put(20,0){\line(-2,3){10}}
  \put(7,3){$\alpha$}
  \end{picture}}

\newcommand{\betatri}{
  \begin{picture}(16,15)(2,0)
  \put(0,0){\line(1,0){20}} 
  \put(0,0){\line(2,3){10}}
  \put(20,0){\line(-2,3){10}}
  \put(6,3){$\beta$}
  \end{picture}}

\newcommand{\deltatri}{
  \begin{picture}(16,15)(2,0)
  \put(0,0){\line(1,0){20}} 
  \put(0,0){\line(2,3){10}}
  \put(20,0){\line(-2,3){10}}
  \put(7,3){$\delta$}
  \end{picture}}

\newcommand{\alphabox}{
  \mbox{\begin{picture}(12,12)(2,0)
  \put(0,0){\line(1,0){12}}
  \put(0,0){\line(0,1){12}}
  \put(12,12){\line(-1,0){12}}
  \put(12,12){\line(0,-1){12}}
  \put(3,4){$\alpha$}
  \end{picture}}}

\newcommand{\betabox}{
  \mbox{\begin{picture}(12,12)(2,0)
  \put(0,0){\line(1,0){12}}
  \put(0,0){\line(0,1){12}}
  \put(12,12){\line(-1,0){12}}
  \put(12,12){\line(0,-1){12}}
  \put(2,3){$\beta$}
  \end{picture}}}

\newcommand{\gammabox}{
  \mbox{\begin{picture}(12,12)(2,0)
  \put(0,0){\line(1,0){12}}
  \put(0,0){\line(0,1){12}}
  \put(12,12){\line(-1,0){12}}
  \put(12,12){\line(0,-1){12}}
  \put(3,4){$\gamma$}
  \end{picture}}}

\newcommand{\deltabox}{
  \mbox{\begin{picture}(12,12)(2,0)
  \put(0,0){\line(1,0){12}}
  \put(0,0){\line(0,1){12}}
  \put(12,12){\line(-1,0){12}}
  \put(12,12){\line(0,-1){12}}
  \put(3,3){$\delta$}
  \end{picture}}}

\begin{document}
\frontmatter

\thispagestyle{empty}
\begin{spacing}{1.66}
\begin{center}
UNIVERSITY OF CALGARY
\vfill
$\Sigma\Pi$-Polycategories, Additive Linear Logic, and Process Semantics
\bigskip\bigskip\\
by
\bigskip\bigskip\\
Craig Antonio Pastro
\vfill
A THESIS \\
SUBMITTED TO THE FACULTY OF GRADUATE STUDIES \\
IN PARTIAL FULFILMENT OF THE REQUIREMENTS FOR THE \\
DEGREE OF MASTER OF SCIENCE
\vfill
DEPARTMENT OF COMPUTER SCIENCE
\vfill
CALGARY, ALBERTA \\
MARCH, 2004
\vfill
\copyright\ Craig Antonio Pastro 2004
\end{center}
\end{spacing}
\clearpage

\pagestyle{plain} 
\begin{spacing}{1.66}
\begin{center}
UNIVERSITY OF CALGARY \\
FACULTY OF GRADUATE STUDIES
\end{center}

\noindent The undersigned certify that they have read, and recommend to the
Faculty of Graduate Studies for acceptance, a thesis entitled
``$\Sigma\Pi$-Polycategories, Additive Linear Logic, and Process Semantics''
submitted by Craig Antonio Pastro in partial fulfillment of the requirements
for the degree of Master of Science.
\end{spacing}
\vfill

\begin{spacing}{1.1}
\begin{flushright}
\begin{minipage}{8cm}
\bigskip\bigskip
--------------------------------------------------------- \\
Supervisor, Dr. James Robin B. Cockett \\
Department of Computer Science
\bigskip\bigskip\\

--------------------------------------------------------- \\
Dr. Robert J. Walker \\
Department of Computer Science
\bigskip\bigskip\\

--------------------------------------------------------- \\
Dr. Richard Zach \\
Department of Philosophy
\bigskip\bigskip\\
\end{minipage}
\end{flushright}
\vfill
\begin{flushleft}
-------------------------------------------- \\
Date
\end{flushleft}
\end{spacing}
\clearpage

\chapter*{Abstract}
\begin{spacing}{1.24}
We present a process semantics for the purely additive fragment of linear
logic in which formulas denote protocols and (equivalence classes of) proofs
denote multi-channel concurrent processes. The polycategorical model induced
by this process semantics is shown to be equivalent to the free polycategory
based on the syntax (i.e., it is full and faithfully complete). This
establishes that the additive fragment of linear logic provides a semantics
of concurrent processes. Another property of this semantics is that it gives a
canonical representation of proofs in additive linear logic. \\

\noindent 
This arXived version omits Section~\ref{sec-pol-cd}: ``Circuit diagrams for
polycategories'' as the \Xy-pic diagrams would not compile due to lack of
memory. For a complete version see
\senta{\texttt{<http://www.cpsc.ucalgary.ca/$\sim$pastroc/>}.}
\end{spacing}
\clearpage

\chapter*{Acknowledgements}
\begin{spacing}{1.24}
Without the support, inspiration, and guidance of innumerable friends and
colleagues this work would never have been accomplished. It gives me great
pleasure to thank them all here.

First of all I would like to thank my supervisor, Professor Robin Cockett.
His guidance, support, and encouragement has been invaluable for this thesis,
as well as for myself. I consider myself extremely fortunate to have worked
under his supervision.

A thank you (in alphabetical order) to Robin Cockett, Rob Walker, and
Richard Zach for serving on my thesis committee.

To the wonderful staff in the Department of Computer Science at the University
of Calgary. I am also grateful to the Department of Computer Science at the
University of Calgary for financial support.

To my many mentors (academic or otherwise) who have helped me to find myself
(or, indeed, to lose myself when that was necessary) --- for their guidance
and encouragement I will always be indebted.

My friends\ldots what can I say really? You all mean more to me than I could 
write on a few lines here. I love you all. Cheers!

To both of my parents whose love and support (emotionally and financially)
for me is absolute. To both of them my deepest and most heartfelt thanks.

This thesis was typeset using \LaTeX\ with diagrams constructed 
using the \Xy-pic package of K. Rose and R. Moore and with inferences 
constructed using the \texttt{proof.sty} package of M. Tatsuya.
\end{spacing}
\clearpage

\chapter*{}
\begin{spacing}{1.24}
\begin{center}
\emph{In memory of my Mother, Lydia Javier Gladman} \\ \emph{1945 - 2000}
\end{center}
\end{spacing}
\clearpage

\begin{spacing}{1.24}
\tableofcontents
\listoftables
\end{spacing}

\mainmatter
\markright{}\pagestyle{myheadings}
\begin{spacing}{1.24}

\nocite{lambek86:intro}
\nocite{girard89:proofs}
\nocite{blass92:game}

\setcounter{chapter}{-1} \chapter{Introduction}

Since the introduction of linear logic by Girard~\cite{girard87:linear},
people have believed that it should somehow provide a logic of concurrent
communication. This intuition, however, has turned out to be rather
difficult to substantiate in practice. In a seminal paper, Abramsky and
Jagadeesan~\cite{abramsky94:games} describe a game model for the
multiplicative fragment of linear logic (MLL) with the MIX rule: formulas
denote games and proofs denote winning strategies. As they point out, a
game can be seen as a process; in fact, they provide a ``dictionary''
translating between the terminology of the two areas. The basic idea behind
game semantics is to interpret a formula as a two-player game between
``player'' and ``opponent'', and proofs of formulas as ``winning
strategies'' for the player. Using this idea, a proof can then be seen as
an interaction between player and opponent. In terms of processes, one
thinks of the player as the ``system'', the opponent as the ``environment'',
and winning strategies for the player as ``deadlock free processes''. In
this view, therefore, a proof can be seen as a process or system
interacting with its environment.  

Game semantics has turned out to be a remarkably effective tool for providing
fully abstract semantics for programming languages. Its initial success was
in providing the first syntax-free, fully abstract model for Scott's language
PCF~\cite{scott93:type}, given independently by Abramsky, Jagadeesan, and
Malacaria~\cite{abramsky00:full}, Hyland and Ong~\cite{hyland00:full}, and
Nickau~\cite{nickau96:hereditarily}. Since then game semantics has provided
fully abstract models for other programming languages with various other
features, such as richer type structures, different evaluation strategies,
non-determinism, etc. (see, e.g., \cite{abramsky98:games,
danos00:probabilistic,harmer99:fully, laird98:semantic, mccusker96:games}).

As pointed out by Abramsky and Melli\`es~\cite{abramsky99:concurrent},
these types of game models do not, however, provide a model of concurrent
communication. A play in these ``sequential'' games is a fixed interleaving
of player and opponent moves. In~\cite{abramsky99:concurrent}, Abramsky and
Melli\`es generalize games to ``concurrent games'', in which they abandon
this interleaving of moves, and construct games in which both the player and
opponent act in a distributed, asynchronous fashion. They announce
in~\cite{abramsky99:concurrent} that the concurrent games model is a
``good'' (i.e., fully complete; see below) model of multiplicative-additive
linear logic.

In~\cite{abramsky94:games}, Abramsky and Jagadeesan introduce the
notion of a fully complete model of a logic. Given a logic \script{L} and a
categorical model of that logic \script{M}, full completeness says that for
any formulas $A$ and $B$ of the logic, any morphism $f: \llbracket A
\rrbracket \ra \llbracket B \rrbracket$ in the \script{M}-interpretation of
\script{L} is the interpretation of a proof of $A \vd B$. One may also ask
for a stronger notion: that any map in the model is the denotation of a
unique proof. The term ``full'' is derived from category theory; that a
model be fully complete is equivalent to requiring that the functor from
the free category based on the syntax to the model is full. That any map
in the model is the denotation of a unique proof is equivalent to requiring
that this functor also be faithful. The construction of fully complete and
full and faithfully complete models is of interest in its own right as it
provides insights into the proofs of a logic. For a recent survey
see~\cite{blute03:softness}.

In this thesis we present a process semantics for the additive fragment of
linear logic. As indicated above, this semantics could alternatively
have been described as a game theoretic semantics, and indeed, many of the
ideas and terminology are derived from this view. The games (processes) which
we introduce are quite different from the Abramsky-Jagadeesan or Hyland-Ong
style of games, and are perhaps more in the spirit of the money games of
Joyal~\cite{joyal95:lattices} and the games of Luigi
Santocanale~\cite{santocanale02:free}.

That the additive fragment of linear logic is particularly relevant to the
perspective of proofs as processes began to be suspected when Joyal started
his investigation of free bicompletions of categories \cite{joyal95:free,
joyal95:freeen} and their relationship to linear logic. Cockett and Seely,
in their efforts to understand Joyal's work, then investigated the logic of
finite sums and products~\cite{cockett01:finite} and realized that this was
precisely the logic of communication along a single (two-way) channel. Their
logic is a fragment of additive linear logic, restricted to exactly one
formula on each side of the turnstile. By removing this restriction and
allowing an arbitrary number of formulas on each side of the turnstile
we get the logic used in this thesis, which we (also) call $\Sigma\Pi$. 
$\Sigma\Pi$ is a slightly untraditional presentation of the additive
fragment of linear logic in that the sums and products are indexed by
arbitrary finite sets, which is equivalent to binary sums and products
with the nullary cases.

In this thesis we develop four equivalent views of concurrent processes:
a process semantic view, a proof theoretic view, a term representation view,
and a categorical view. The proof theoretic view is our starting point. We
begin with the $\Sigma\Pi$ logic and consider proofs as processes. The
process semantics is used to show that proofs in $\Sigma\Pi$ do, in fact,
correspond to concurrent processes. This, of course, is what validates the
claim that our other views are of concurrent processes. The term
representation provides terms as processes. Since it is known that the terms
correspond to concurrent processes, they can now be viewed as a programming
language for these processes. In the categorical view, concurrent processes
are represented by morphisms in a polycategory. The proof that all four views
of concurrent processes are equivalent also uses categorical ideas.

A more detailed account is as follows. We begin with the logic $\Sigma\Pi$. 
Cut-elimination rewrites are presented for this logic, but in order for the
cut-elimination rewrite system to be Church-Rosser, it is also necessary to
give some equivalence rewrites for proofs which allows the interchange of
rules: the so-called ``permuting conversions''.

Motivated by the interpretation of proofs as processes, two term
representations for proofs in $\Sigma\Pi$ are developed: a term calculus
representation and a ``programming language'' representation. The first step
to view processes as terms is to view formulas as ``protocols'' assigned to
a channel. For example, to view the formula $X+Y$ as a protocol it is
assigned a channel, say $\alpha$, and each subformula is assigned an
``event'' resulting in $\alpha:(x:X+y:Y)$. Proofs may now be viewed as
``multi-channel processes''. If $f$ is a proof the sequent $\alpha:X \vd
\beta:Z$ viewed as a process and $g$ a proof of $\alpha:Y \vd \beta:Z$
viewed as a process, then a proof of $\alpha:(x:X+y:Y) \vd \beta:Z$ may be
viewed as the process
\[\alpha\left\{\begin{array}{l}x \mapsto f\\ y \mapsto g \end{array}\right\}
\]
This process can be read as, ``wait for an input event to occur on the
channel $\alpha$. If it is $x$ then do the process $f$, and if it is
$y$ then it will do the process $g$.''

The cut-elimination rewrites and permuting conversions may be represented
using these terms. The notation ``$f ;_\gamma g$'' will be used to indicate
cutting the proofs $f$ and $g$ together on the channel $\gamma$. If we
interpret cut as communication, cut-elimination then is the dynamics of
communication. As an example consider the terms
\[\ora{\gamma}[a](f) \qquad \text{and} \qquad
\gamma\left\{\begin{array}{l}a \mapsto f\\ b \mapsto g \end{array}\right\}
\]
The first term is interpreted as, ``output $a$ on $\gamma$ and then do the
process $f$.'' The second term as, ``wait for an input event on $\gamma$. If
it is the event $a$ then do the process $g$ and if it is the event $b$ then
do the process $h$.'' If we tell these two processes to communicate on
$\gamma$ it should be obvious that after one step of communication they will
evolve to the process $f$ communicating with the process $g$ on $\gamma$,
i.e., 
\[\ora{\gamma}[a](f)\ ;_\gamma
\gamma\left\{\begin{array}{l}a \mapsto f\\ b \mapsto g \end{array}\right\}
\quad \Lra \quad f\ ;_\gamma g
\]
These processes look very sequential in nature, but the permuting
conversions are hiding the fact that these are concurrent processes!

As terms are much easier to manipulate than proofs, they are used to show
that the rewriting system induced by the cut-elimination rewrites is
Church-Rosser and terminating modulo the equivalences. A decision procedure
is presented which allows us to determine if two (cut-free) terms of the
same type (derivations of the same sequent) are equivalent, i.e., that they
are related by the equivalence rewrites.

A polycategory of processes may be formed by considering protocols as
objects, multi-channel processes (the terms) as morphisms, and communication
as composition. It is shown that this is the free polycategory with finite
sums and finite products. Following Joyal's and Cockett and Seely's lead, we
prove a ``Whitman theorem'' which gives a characterization of the free
polycategory with finite sums and products over a base polycategory in
terms of characterizations of the hom sets.

Finally, a process semantics for $\Sigma\Pi$-terms is presented. This
semantics illustrates precisely how the terms correspond to channel-based
concurrent processes. The rough idea is as follows. We define the notion of
a behaviour, which is an explicit representation of the history of a
process. These are typically represented using tables, e.g.,
\[\begin{array}{|c|c||c|c|}
& && b\\
d & e & & \ol{c} \\
\ol{c} & a & \ol{d} & a \\
\hline
\alpha & \beta & \gamma & \delta \\
\hline
\end{array}
\]
This behaviour represents outputting $c$ and then receiving $d$ on channel
$\alpha$, receiving $a$ and then $e$ on $\beta$, and so on. From behaviours
there is no way to tell the order that the events occurred, however, this is
necessary. To see this consider the process:
\[\alpha\left\{
\begin{array}{l}
a \mapsto \ora{\beta}[c](f) \\
b \mapsto \ora{\beta}[d](g)
\end{array} \right\}
\]
It is easily seen that this process is unable to produce an output on $\beta$
until it has first received an event on $\alpha$. That is, there is a
hierarchy on events. This motivates the use of ``entailments''. An
entailment is a behaviour with a distinguished output event. If all the
events in the behaviour have occurred then the output event may be performed.
The following set of entailments
\[\Bigg\{\ \begin{array}{|c||c|}
a & \\
\hline \alpha & \beta \\ \hline
\end{array} \vd \ol{\beta}[c],
\quad
\begin{array}{|c||c|}
& b \\
\hline \alpha & \beta \\ \hline
\end{array} \vd \ol{\beta}[d],
\quad
\begin{array}{|c||c|}
a & \ol{c} \\
\hline \alpha & \beta \\ \hline
\end{array} \vd f,
\quad
\begin{array}{|c||c|}
b & \ol{d} \\
\hline \alpha & \beta \\ \hline
\end{array} \vd g \Bigg\}
\]
encodes the hierarchy on events of the process above. Sets of entailments
(satisfying the seven rules defined in Chapter~\ref{chap-semantics}) are what
we call ``extensional processes''. These are what gives us our interpretation
of proofs as concurrent processes.

A polycategory of extensional processes is constructed and shown to be
equivalent to the polycategory of processes. This proves that this model
satisfies the property that every process is the denotation of a unique
cut-free proof, i.e., it is a full and faithfully complete model.

Some of the previous results are essentially an extension from the (ordinary)
categorical case to the polycategorical case of the results of Cockett and
Seely~\cite{cockett01:finite}. On occasion I have used the exposition found
therein when it is much clearer than I could hope to achieve. Errors in
these sections, as in the rest of the text, are of course solely my
responsibility.

\paragraph{Outline of this thesis}

This thesis is organized as follows: Chapter~\ref{chap-prelim} introduces
some standard concepts of category theory and the notion of a polycategory.
The concepts here will be needed throughout the thesis.
Chapter~\ref{chap-logic} presents the $\Sigma\Pi$-logic. In
Chapter~\ref{chap-syntax} we develop two term representations for
derivations in this logic: a term calculus and a ``programming language''
representation. A rewrite system for cut-elimination is presented
and shown to be Church-Rosser and terminating. The proof of decidability for
these terms is also presented in this chapter. Chapter~\ref{chap-catsem}
is devoted to proving that $\Sigma\Pi_\cat{A}$ is the free polycategory
built over an arbitrary polycategory \cat{A}. In Chapter~\ref{chap-semantics}
we present a process semantics for our logic. This semantics is shown to
provide a full and faithfully complete polycategorical model for our logic.

\paragraph{Contributions of this thesis}

In Chapter~\ref{chap-logic} the description of additive linear logic is a new
presentation of this fragment. Chapter~\ref{chap-syntax} contains two term
representations which are essentially new to this thesis, as is the
presentation of cut-elimination and the decision procedure for terms.
Chapter~\ref{chap-catsem} contains a new polycategorical presentation of
additive linear logic and the development of the Whitman theorem for this
setting. Chapter~\ref{chap-semantics} contains the description of extensional
processes, which is the main novel aspect of this thesis.

\chapter{Categorical Preliminaries} \label{chap-prelim}

This chapter is meant to give a brief introduction to category theory and
also to help accustom the reader with the notation that will be used
throughout this thesis. For a more complete introduction to category theory 
see, e.g., \cite{maclane98:categories} or \cite{barr99:category}.

\section{Categories}

A \textbf{category} \cat{C} consists of a class of objects, $\ob(\cat{C})$
(or $\cat{C}_0$), for each pair of objects $A$ and $B$, a class of morphisms
(or arrows), $\cat{C}(A,B)$, for each object $A$, an identity morphism
$\id_A \in \cat{C}(A,A)$, and for each triple of objects, $A$, $B$, and $C$,
a composition law
\[\cat{C}(A,B) \times \cat{C}(B,C) \To^{\ ;\ } \cat{C}(A,C)
\]
satisfying:
\begin{itemize}
\item if $f \in \cat{C}(A,B)$ then the domain (or source) of $f$ is 
$\mathrm{dom}(f) = A$ and the codomain (or target) of $f$ is 
$\mathrm{cod}(f) = B$.

\item if $f:A \ra B$ then $\id_A ; f = f ; \id_B = f$.

\item $f;(g;h) = (f;g);h$ whenever either side is defined, i.e., composition
is associative.

\end{itemize}

When it is clear from the context we will omit the label $\mathrm{ob}(\cat{C})$
and simply refer to an object $A \in \ob(\cat{C})$ as $A \in \cat{C}$. We will
sometimes write $\Hom(A,B)$ or just $(A,B)$ to mean $\cat{C}(A,B)$ when this
will not lead to confusion. A morphism $f \in \cat{C}(A,B)$ may also be
written as $f:A \rightarrow B$ or $A \to^f B$. For any composible pair of
arrows $f:A \rightarrow B$ and $g:B \rightarrow C$ we will sometimes denote
their composition in the usual manner with a $\circ$ as $g \circ f$ (as
opposed to the diagrammatic order $f;g$ that we typically use). We often omit
the $\circ$ and simply write $gf$ but we will never omit the ;\ .

Note that when we defined the notion of a category we did not say that the 
collection of objects constitute a set. Indeed, in the most famous of
categories, \Set, where the objects are sets and the morphisms are functions
between sets, the collection of all sets is not itself a set. This motivates
the following definitions: a category is called \textbf{small} if its objects
constitute a set, and \textbf{large} otherwise. If $\cat{C}(A,B)$ is a set for
all objects $A,B \in \cat{C}$ then \cat{C} is called \textbf{locally small}.

If \cat{C} is a category, then its \textbf{dual} $\cat{C}^{\op}$, is
defined by $\ob(\cat{C}^{\op}) = \ob(\cat{C})$ and $\cat{C}^{\op}(A,B) = 
\cat{C}(B,A)$. That is, the dual of a category is the category with all its
arrows reversed: if $f:A \ra B$ is an arrow in \cat{C} then $f:B \ra A$ 
is an arrow in $\cat{C}^{\op}$. It is clear that the dual of a category is 
also a category.  

\section{Properties of morphisms}

An arrow $f:A \ra B$ in a category \cat{C} is called \textbf{monic} (or a
\textbf{monomorphism}) if for any object $C \in \cat{C}$ and arrows
$x,y:C \ra A$ such that $x;f = y;f$ then $x=y$. The dual of a monomorphism is an
\textbf{epic} (or \textbf{epimorphism}), i.e., an arrow $f:A \ra B$ is epic
if for any object $C$ and arrows $x,y:B \ra C$ such that $f;x = f;y$ then
$x=y$. An arrow $f:A \ra B$ is called an \textbf{isomorphism} if there is an
arrow $g:B \ra A$ such that $f;g = \id_A$ and $g;f = \id_B$. If $f$ is an
isomorphism then the arrow $g$ is uniquely determined and is usually called
the inverse of $f$. If such an isomorphism exists, we say that $A$ is
isomorphic to $B$ and denote this as $A \iso B$.

\section{Sums and products} \label{prelim-sums}

If $A$ and $B$ are objects in a category \cat{C} then by the \textbf{product}
of $A$ and $B$ we mean an object $C$ together with arrows $p_1:C \ra A$ and
$p_2:C \ra B$ (called the first and second projections respectively) such
that for any object $D$ and arrows $f:D \ra A$ and $g:D \ra B$ there is a
unique arrow $q:D \ra C$ making the following diagram commute:
\[\xymatrix@=1.5cm{
  & D \ar[dl]_f \ar[dr]^g \ar@{-->}[d]^q & \\
A & C \ar[l]^{p_1} \ar[r]_{p_2} & B}
\]

The dual of a product is a \textbf{sum} (or \textbf{coproduct}). That is,
if $A$ and $B$ are objects in a category \cat{C} then the coproduct of $A$
and $B$ is an object $C$ together with arrows $b_1:A \ra C$ and $b_2:B \ra C$
(called the first and second injections respectively) such that for any 
object $D$ and arrows $f:A \ra D$ and $g:B \ra D$ there is a unique arrow
$q:C \ra D$ making the following diagram commute:
\[\xymatrix@=1.5cm{
& D & \\
A \ar[ur]^f \ar[r]_{b_1} & C \ar@{-->}[u]_q & B \ar[ul]_g \ar[l]^{b_2}} 
\]

The product of $A$ and $B$ is typically denoted as $A \times B$ and the unique
arrow from $D$ to $A \times B$ as $\langle f,g \rangle$. The coproduct of
$A$ and $B$ is typically denoted as $A+B$ and the unique arrow from $A+B$ to
$D$ as $\langle f|g \rangle$.

For products (and dually for sums), it can be shown that $A \x B$ is
isomorphic to $B \x A$. To see this consider the diagram
\[\xymatrix{&A \x B \ar[dl]_{p_1} \ar[dr]^{p_2} \\ A  && B \\
& B \x A \ar[ul]^{p_2} \ar[ur]_{p_1}}
\]
and the unique arrows from the definition that make this diagram commute.

\section{Equalizers and coequalizers}

Let \cat{C} be a category and $f,g:A \two B$ be a parallel pair of arrows.
An \textbf{equalizer} of $f$ and $g$ is an object $E$ together with an arrow
$e$ such that
\begin{enumerate}[{\upshape (i)}]
\item $e;f = e;g$, and
\item for any object $Q$ and arrow $q:Q \ra A$ such that $q;f = q;g$, there
is a unique arrow $h:Q \ra E$ such that $h;e = q$.
\end{enumerate}

Dually, a \textbf{coequalizer} of $f,g:A \two B$ is an object $C$ together
with an arrow $c:B \ra C$ such that
\begin{enumerate}[{\upshape (i)}]
\item $f;c = g;c$, and
\item for any object $Q$ and arrow $q:B \ra Q$ such that $f;q = g;q$, there
is a unique arrow $h:C \ra Q$ such that $c;h = q$.
\end{enumerate}

\section{Functors}

Let \cat{C} and \cat{D} be categories. A \textbf{functor} $F:\cat{C} \ra 
\cat{D}$ between categories is a pair of maps $F_0: \cat{C}_0 \ra \cat{D}_0$
and $F_1: \cat{C}(A,B) \ra \cat{D}(F_0(A),F_0(B))$, for all objects $A,B \in
\cat{C}$, satisfying
\begin{itemize}
\item $F_1(\id_A) = \id_{F_0(A)}$
\item $F_1(f;g) = F_1(f) ; F_1(g)$
\end{itemize}

It is standard practice to omit the subscripts from the functor when the 
context is clear. The brackets may also be omitted when they are not required
to disambiguate precedence. Thus, we will sometimes write $FA$ and $Ff$ to
mean $F_0(A)$ and $F_1(f)$ respectively.

For any category \cat{C} there is an identity functor, $\id_\cat{C}: \cat{C}
\ra \cat{C}$, defined in the obvious way.

Two categories \cat{C} and \cat{D} are said to be \textbf{isomorphic},
denoted $\cat{C} \iso \cat{D}$, if there are a pair of functors
$F:\cat{C} \ra \cat{D}$ and $G:\cat{D} \ra \cat{C}$ such that
$F;G = \id_\cat{C}$ and $G;F = \id_\cat{D}$.

\section{Natural transformations}

Given two functors $F,G:\cat{C} \rightarrow \cat{D}$, a \textbf{natural 
transformation} $\alpha:F \Rightarrow G$ consists of a family of morphisms 
$\alpha_A: FA \Rightarrow GA$, one for each object $A \in \cat{C}$, such 
that for any morphism $f:A \rightarrow B$ in \cat{C} the following diagram 
commutes.
\[\xymatrix@=2cm{
FA \ar[r]^{\alpha_A} \ar[d]_{Ff} & GA \ar[d]^{Gf} \\
FB \ar[r]_{\alpha_B} & GB}
\]
We call the $\alpha_A$ the \textbf{component} of the natural 
transformation at $A$.

For any functor $F:\cat{C} \ra \cat{D}$ there is an identity natural 
transformation $\id_F:F \Ra F$ defined by $(\id_F)_A = \id_{FA}$.

A natural transformation $\alpha:F \Ra G$ is called a \textbf{natural
isomorphism} if every component $\alpha_A$ is invertible in \cat{D}. In this
case we say that $F$ and $G$ are isomorphic and write $F \iso G$.
 
Two categories \cat{C} and \cat{D} are said to be \textbf{equivalent} if
there are functors $F:\cat{C} \ra \cat{D}$ and $G:\cat{D} \ra \cat{C}$ such
that $F;G \iso \id_\cat{C}$ and $G;F \iso \id_\cat{D}$.

\section{Polycategories}

A \textbf{planar polycategory} \cat{P} consists of the following data

\begin{itemize}
\item a class $\cat{P}_0$ of objects of \cat{P},

\item for each $m, n \in \mathbb{N}$ and $x_1,\ldots,x_m,y_1,\ldots,y_n \in
\cat{P}_0$, a set
\[\cat{P}(x_1,\ldots,x_m\ ;\ y_1,\ldots,y_n)\]
whose elements are called polymorphisms.
Using $\Gamma$ and $\Delta$ to represent strings of elements of $\cat{P}_0$,
the polymorphisms in $\cat{P}(\Gamma\ ;\ \Delta)$ may be denoted $f:\Gamma \ra
\Delta$ or $\Gamma \vd_f \Delta$ where $\dom(f) = \Gamma$ and $\cod(f) =
\Delta$. 
\end{itemize}
together with

\begin{itemize}
\item for each $x \in \cat{P}_0$, an identity morphism $1_x \in \cat{P}(x,x)$

\item an operation
\[\cat{P}(\Gamma;\Delta_1,x,\Delta_2) \x \cat{P}(\Gamma_1,x,\Gamma_2;\Delta) 
\ra \cat{P}(\Gamma_1,\Gamma,\Gamma_2;\Delta_1,\Delta,\Delta_2)
\]
called \textbf{cut}, restricted to the cases where either $\Gamma_1$ or
$\Delta_1$ is empty and either $\Gamma_2$ or $\Delta_2$ is empty. (This
restriction is called the \textbf{crossing} (or \textbf{planarity condition}),
cf. circuit diagrams below.) Explicitly, this gives four cut rules:
\begin{itemize}
\item $\cat{P}(\Gamma_1;\Delta_1,x) \x \cat{P}(x,\Gamma_2;\Delta_2)
\to^{\cut_r} \cat{P}(\Gamma_1,\Gamma_2;\Delta_1,\Delta_2)$

\item $\cat{P}(\Gamma_1;x,\Delta_1) \x \cat{P}(\Gamma_2,x;\Delta_2)
\to^{\cut_l} \cat{P}(\Gamma_2,\Gamma_1;\Delta_2,\Delta_1)$

\item $\cat{P}(\Gamma;x) \x \cat{P}(\Gamma_1,x,\Gamma_2;\Delta)
\to^{\text{mcut}} \cat{P}(\Gamma_1,\Gamma,\Gamma_2;\Delta)$

\item $\cat{P}(\Gamma;\Delta_1,x,\Delta_2) \x \cat{P}(x;\Delta)
\to^{\text{mcut}^{\op}} \cat{P}(\Gamma; \Delta_1,\Delta,\Delta_2)$

\end{itemize}
\end{itemize}
These data are subject to three axioms:

\begin{itemize}
\item cut has identities, i.e., the following diagrams commute

\[\xymatrix@C=10ex{
\cat{P}(\Gamma_1,x,\Gamma_2;\Delta) \ar[r]^-{1_x \x -} \ar@{=}[dr] &
\cat{P}(x,x) \x \cat{P}(\Gamma_1,x,\Gamma_2;\Delta) \ar[d]^-{\cut} \\
& \cat{P}(\Gamma_1,x,\Gamma_2;\Delta)}
\]

\[\xymatrix@C=10ex{
\cat{P}(\Gamma;\Delta_1,y,\Delta_2) \ar[r]^-{- \x 1_y} \ar@{=}[dr] &
\cat{P}(\Gamma;\Delta_1,y,\Delta_2) \x \cat{P}(y,y) \ar[d]^-{\cut} \\
& \cat{P}(\Gamma;\Delta_1,y,\Delta_2)}
\]

\item cut is associative, i.e., the following diagram commutes
{\scriptsize
\[\xymatrix@M=1ex@C=50ex@R=20ex@!0{
\txt{$\cat{P}(\Gamma_1;\Gamma_2,x,\Gamma_3)$ \\
$\x \cat{P}(\Delta_1,x,\Delta_2;\Delta_3,y,\Delta_4)$ \\
$\x \cat{P}(\Phi_1,y,\Phi_2;\Phi_3)$}
\ar[r]^-{1 \x \cut} \ar[d]_-{\cut \x 1}
& \txt{$\cat{P}(\Gamma_1;\Gamma_2,x,\Gamma_3)$ \\
$ \x \cat{P}(\Phi_1,\Delta_1,x,\Delta_2,\Phi_2;\Delta_3,\Phi_3,\Delta_4)$}
\ar[d]^-{\cut} \\
\txt{$\cat{P}(\Delta_1,\Gamma_1,\Delta_2;\Gamma_2,\Delta_3,x,\Delta_4,\Gamma_3)$
\\ $\x \cat{P}(\Phi_1,y,\Phi_2;\Phi_3)$} \ar[r]_-{\cut}
& \cat{P}(\Phi_1,\Delta_1,\Gamma_1,\Delta_2,\Phi_2;\Gamma_2,\Delta_3,\Phi_3,
\Delta_4,\Gamma_3)}
\]}
Recall that cut is subject to the crossing condition; writing this
restriction explicitly would result in nine separate commutative diagrams.

\item cut satisfies the interchange law (originally referred to as
``commutativity'' by Lambek~\cite{lambek69:deductive}), i.e., the following
diagrams commute
{\scriptsize
\[\xymatrix@M=1ex@C=55ex@R=15ex@!0{
\txt{$\cat{P}(\Gamma_1;\Gamma_2,x,\Gamma_3,y,\Gamma_4)
\x \cat{P}(\Delta_1,x,\Delta_2;\Delta_3)$ \\
$\x \cat{P}(\Phi_1,y,\Phi_2;\Phi_3)$}
\ar[d]_-{\cut \x 1} \ar[r]^-{1 \x \gamma} &
\txt{$\cat{P}(\Gamma_1;\Gamma_2,x,\Gamma_3,y,\Gamma_4)
\x \cat{P}(\Phi_1,y,\Phi_2;\Phi_3)$ \\
$\x \cat{P}(\Delta_1,x,\Delta_2;\Delta_3)$} \ar[d]^-{\cut \x 1} \\
\txt{$\cat{P}(\Delta_1,\Gamma_1,\Delta_2;\Gamma_2,\Delta_3,\Gamma_3,y,\Gamma_4)$
\\ $\x \cat{P}(\Phi_1,y,\Phi_2;\Phi_3)$} \ar[d]_-{\cut} &
\txt{$\cat{P}(\Phi_1,\Gamma_1,\Phi_2;\Gamma_2,x,\Gamma_3,\Phi_3,\Gamma_4)$ \\
$\x \cat{P}(\Delta_1,x,\Delta_2;\Delta_3)$} \ar[d]^-{\cut} \\
\cat{P}(\Phi_1,\Delta_1,\Gamma_1,\Delta_2,\Phi_2;\Gamma_2,\Delta_3,
\Gamma_3,\Phi_3,\Gamma_4) \ar@{=}[r] &
\cat{P}(\Delta_1,\Phi_1,\Gamma_1,\Phi_2,\Delta_2;
\Gamma_2,\Delta_3,\Gamma_3,\Phi_3,\Gamma_4)}
\]}

{\scriptsize
\[\xymatrix@M=1ex@C=55ex@R=15ex@!0{
\txt{$\cat{P}(\Gamma_1;\Gamma_2,x,\Gamma_3)$ \\
$\x \cat{P}(\Delta_1;\Delta_2,y,\Delta_3)
\x \cat{P}(\Phi_1,x,\Phi_2,y,\Phi_3;\Phi_4)$}
\ar[d]_-{1 \x \cut} \ar[r]^-{\gamma \x 1} &
\txt{$\cat{P}(\Delta_1;\Delta_2,y,\Delta_3)$ \\
$\x \cat{P}(\Gamma_1;\Gamma_2,x,\Gamma_3)
\x \cat{P}(\Phi_1,x,\Phi_2,y,\Phi_3;\Phi_4)$} \ar[d]^-{1 \x \cut} \\
\txt{$\cat{P}(\Gamma_1;\Gamma_2,x,\Gamma_3)$ \\ 
$\x \cat{P}(\Phi_1,x,\Phi_2,\Delta_1,\Phi_3;\Delta_2,\Phi_4,\Delta_3)$}
\ar[d]_-{\cut} &
\txt{$\cat{P}(\Delta_1;\Delta_2,y,\Delta_3)$ \\
$\x \cat{P}(\Phi_1,\Gamma_1,\Phi_2,y,\Phi_3;\Gamma_2,\Phi_4,\Gamma_3)$}
\ar[d]^-{\cut} \\
\cat{P}(\Phi_1,\Gamma_1,\Phi_2,\Delta_1,\Phi_3;
\Gamma_2,\Delta_2,\Phi_4,\Delta_3,\Gamma_3)
\ar@{=}[r] &
\cat{P}(\Phi_1,\Gamma_1,\Phi_2,\Delta_1,\Phi_3;
\Delta_2,\Gamma_2,\Phi_4,\Gamma_3,\Delta_3)
}
\]}
where $\gamma:A \x B \to^\sim B \x A$ is the isomorphism for products.
Similarly here writing out the crossing condition explicitly would give
four separate diagrams for each of the above diagrams.

\end{itemize}

A \textbf{symmetric polycategory} \cat{P} is a polycategory equipped with
a symmetric action, i.e., for permutations $\sigma \in \textbf{S}_m,\
\tau \in \textbf{S}_n$ (where $\textbf{S}_k$ is the group of permutations
on $k$ objects), a map
\[\cat{P}(\Gamma;\Delta) \to^{c_{\sigma,\tau}} \cat{P}(\sigma\Gamma;\tau\Delta)
\]
where $\sigma\Gamma = \sigma(x_1,\ldots,x_n) = (x_{\sigma(1)},\ldots
x_{\sigma(n)})$, satisfying the following coherence conditions.

\begin{itemize}
\item The symmetric actions may be composed.

\[\vcenter{\xymatrix@R=7ex@C=7ex{
\cat{P}(\Gamma;\Delta) \ar[r]^{c_{\sigma,\tau}} \ar[dr]_{c_{\sigma;\sigma',
\tau;\tau'}} & \cat{P}(\sigma\Gamma,\tau\Delta) \ar[d]^{c_{\sigma',\tau'}} \\
& \cat{P}(\sigma'(\sigma\Gamma),\tau'(\tau\Delta))}}
\]

\item That the next four diagrams commute assert that cutting and then
permuting the objects is equivalent to first permuting the objects and then
cutting.

\begin{equation*}
\vcenter{\xymatrix@R=7ex@C=7ex{
\cat{P}(\Gamma;\Delta,x) \x \cat{P}(x,\Gamma';\Delta')
\ar[r]^-{\cut_r} \ar[d]_{c_{\sigma,\tau|x} \x c_{x|\sigma',\tau'}} &
\cat{P}(\Gamma,\Gamma';\Delta,\Delta') \ar[d]^{c_{\sigma|\sigma',\tau|\tau'}} \\
\cat{P}(\sigma\Gamma;\tau\Delta,x) \x
\cat{P}(x,\sigma'\Gamma';\tau'\Delta') \ar[r]_-{\cut_r} &
\cat{P}(\sigma\Gamma,\sigma'\Gamma';\tau\Delta,\tau'\Delta')}}
\end{equation*}

\begin{equation*}
\vcenter{\xymatrix@R=7ex@C=7ex{
\cat{P}(\Gamma;x,\Delta) \x \cat{P}(\Gamma',x;\Delta)
\ar[r]^-{\cut_l} \ar[d]_{c_{\sigma,x|\tau} \x c_{\sigma'|x,\tau'}} &
\cat{P}(\Gamma',\Gamma;\Delta',\Delta) \ar[d]^{c_{\sigma'|\sigma,\tau'|\tau}}\\
\cat{P}(\sigma\Gamma;x,\tau\Delta) \x \cat{P}(\sigma'\Gamma',x;\tau'\Delta')
\ar[r]_-{\cut_l} &
\cat{P}(\sigma'\Gamma',\sigma\Gamma;\tau'\Delta',\tau\Delta)}}
\end{equation*}

\begin{equation*}
\vcenter{\xymatrix@R=7ex@C=7ex{
\cat{P}(\Gamma;x) \x \cat{P}(\Gamma_1,x,\Gamma_2;\Delta)
\ar[d]_c \ar[r]^-{\text{mcut}} &
\cat{P}(\Gamma_1,\Gamma,\Gamma_2;\Delta)
\ar[d]^c \\
\cat{P}(\sigma\Gamma;x) \x \cat{P}(\Gamma_1',x,\Gamma_2';\tau\Delta)
\ar[r]_-{\text{mcut}} &
\cat{P}(\Gamma_1',\sigma\Gamma,\Gamma_2';\tau\Delta)}}
\end{equation*}
where $\Gamma_1',-,\Gamma_2' = \sigma(\Gamma_1,-,\Gamma_2)$.

\begin{equation*}
\vcenter{\xymatrix@R=7ex@C=7ex{
\cat{P}(\Gamma;\Delta_1,x,\Delta_2) \x \cat{P}(x;\Delta)
\ar[r]^-{\text{mcut}^\op} \ar[d]_c &
\cat{P}(\Gamma;\Delta_1,\Delta,\Delta_2) \ar[d]^c \\
\cat{P}(\sigma\Gamma;\Delta_1',x,\Delta_2') \x \cat{P}(x;\tau'\Delta)
\ar[r]_-{\text{mcut}^\op} &
\cat{P}(\sigma\Gamma;\Delta_1',\tau'\Delta,\Delta_2')}}
\end{equation*}
where $\Delta_1',-,\Delta_2 = \tau(\Delta_1,-,\Delta_2)$
\end{itemize}

We will only be concerned with symmetric polycategories in this thesis and
so will refer to the symmetric version of polycategories simply as
polycategories.

\subsection{Circuit diagrams for polycategories} \label{sec-pol-cd}

\paragraph{Note} This section has been removed for the arXived version of
this thesis. The diagrams would not compile. For a complete version see
\senta{\texttt{<http://www.cpsc.ucalgary.ca/$\sim$pastroc/>}.}

\subsection{Morphisms of polycategories}

Let \cat{P} and \cat{Q} be polycategories. A \textbf{morphism} $\cat{P}
\to^F \cat{Q}$ between polycategories is a pair of maps $F_0:\cat{P}_0 \ra
\cat{Q}_0$ and $F_1:\cat{P}(\Gamma;\Delta) \ra
\cat{Q}(F^*(\Gamma);F^*(\Delta))$, where $F^*(X_1,\ldots,X_n) =
F_0(X_1),\ldots,F_0(X_n)$, such that the functor preserves identities,
\[
F_1(\id_A) = \id_{F_0(A)}
\]
preserves composition,
{\scriptsize
\begin{equation*}
\vcenter{\xymatrix@R=9ex{
\cat{P}(\Gamma_1;\Gamma_2,x,\Gamma_3) \x \cat{P}(\Delta_1,x,\Delta_2;\Delta_3)
\ar[r]^-{\cut} \ar[d]_F & 
\cat{P}(\Delta_1,\Gamma_1,\Delta_2;\Gamma_2,\Delta_1,\Gamma_3) \ar[d]^F \\
\txt{$\cat{Q}(F^*(\Gamma_1);F^*(\Gamma_2),F_0(x),F^*(\Gamma_3))$ \\
$\x \cat{Q}(F^*(\Delta_1),F_0(x),F^*(\Delta_2);F^*(\Delta_3))$} 
\ar[r]_-{\cut} &
\cat{Q}(F^*(\Delta_1),F^*(\Gamma_1),F^*(\Delta_2);
F^*(\Gamma_2),F^*(\Delta_3),F^*(\Gamma_3))}}
\end{equation*}}
and for symmetric polycategories, preserves the symmetric action,
\[
\cat{Q}(F^*(\sigma\Gamma);F^*(\tau\Delta))
= \cat{Q}(\sigma(F^*(\Gamma));\tau(F^*(\Delta)))
\]
Again here, we usually drop the subscripts and write $F$ for both $F_0$ and
$F_1$.

Notice that a morphism of polycategories is the polycategorical notion of a
functor between regular categories. There is a more elaborate notion of a
polyfunctor introduced in~\cite{cockett03:polybicat} which will not be
discussed here.

As with functors, for any polycategory \cat{P}, there is an obvious notion of
an identity morphism of polycategories $\id_\cat{P}:\cat{P} \ra \cat{P}$
defined in the obvious way.

Two polycategories \cat{P} and \cat{Q} are said to be \textbf{isomorphic},
denoted $\cat{P} \iso \cat{Q}$, if there are a pair of morphisms of
polycategories $F:\cat{P} \ra \cat{Q}$ and $G:\cat{Q} \ra \cat{P}$ such that
$F ; G = \id_\cat{P}$ and $G ; F = \id_\cat{Q}$.

For a much more detailed exposition on poly\-cate\-gories and
poly-bi\-cate\-gories see, e.g.,~\cite{cockett03:polybicat,cockett97:weakly}.

\begin{remark}[Notation for polycategories]
Until now we have avoid\-ed using any notation to represent the cut rule.
Typically, one represents the cut rule using a ``positional'' notation, i.e.,
the object to cut on is specified by an index. For example, given
\[\Gamma_1 \to^f \Gamma_2 \qquad \text{and} \qquad \Delta_1 \to^g \Delta_2
\]
the notation $f {}_i;_j g$ indicates that we are cutting the $i$-th
component of $\Gamma_2$ with the $j$-th component of $\Delta_1$, provided
the requirements for cut are satisfied.

In this thesis we take a somewhat different approach. Instead of using a
``positional'' notation as above, we will use a ``referential'' notation.
For each morphism in a polycategory we label the ``wires'' (or ``channels'')
and use these labels to indicate which object we are cutting on. To prevent
ambiguity, we require that each wire of a polymorphism receives a unique name,
and make the restriction that cuts may only occur between polymorphisms which
have no wire name in common (so that after the cut each wire still has a
unique name). (This can be accomplished via a renaming procedure for wire
names; we will not describe the details here.) For example, let
\[\alpha_1:X,\alpha_2:Y \to^f \alpha_3:Z,\alpha:W \qquad \text{and} \qquad
\beta:W,\beta_1:S \to^g \beta_2:T
\]
be polymorphisms. We may cut $f$ on $\alpha:W$ with $g$ on $\beta:W$, which
we will write $f {}_\alpha ;_\beta g$.
\end{remark}

\begin{example}[Polycategories] \quad
\begin{enumerate}
\item Any category is a polycategory with one input and one output.
\item Any (symmetric) multicategory (see, e.g., \cite{cheng02:weak}) is a
(symmetric) polycategory with one output.
\item Modules and multilinear maps form a multicategory (and hence a
polycategory).
\item Any (symmetric) linearly distributive category\footnote{Originally
referred to as weakly distributive categories.}~\cite{cockett97:weakly} is an
example of a (symmetric) polycategory. As symmetric linearly distributive
categories with negation are the same as $*$-autonomous categories
(see~\cite{cockett97:weakly} for the details), this implies that any
$*$-autonomous category is an example of a polycategory.

\item The primary example that this thesis is concerned with is: A Gentzen
style sequent calculus (with multiple formulas on either side of the
turnstile ``$\vd$'') with formulas as objects and (equivalence classes of)
derivations as the polymorphisms forms a polycategory.
\end{enumerate}
\end{example}

Below is a concrete example which is worked out in more detail. This example
provides the polycategorical version of an operad, and so is of some
independent interest.

\begin{example}[A one-object polycategory]

\item We may describe morphisms in a one-object polycategory (or
\textbf{polyad}) using pairs consisting of labelled cyclic graphs and their
related signatures, up to renaming of nodes. The nodes of the cyclic graphs
are labeled with ``channel'' names. The signature of a graph indicates (from
left-to-right) the ordering of its input and output wires. We will write
these pairs simply as $G:S$, where $G$ is a cyclic graph, and $S$ is the
signature, e.g.,
\begin{equation} \label{pcm} \tag{$\bigstar$}
\vcenter{\xymatrix@R=2ex@C=1ex{
              & \alpha_1 \ar@/^/[rr]  && \alpha_3 \ar@/^/[dr] \\
\alpha_2 \ar@/^/[ur] &             &&               & \alpha_5 \ar@/^/[dll] \\
              &                  & \alpha_4 \ar@/^/[ull]
}}:\<\alpha_1,\alpha_2,\alpha_3\> \ra \<\alpha_4,\alpha_5\>
\end{equation}
The above polymorphism may be thought of graphically as
\[\vcenter{\xymatrix@R=2ex@C=1ex{
\alpha_1 && \alpha_2 && \alpha_3 \\
&&&& \\
& \alpha_1  \ar@{.}[uurrr] \ar@/^/[rr] && \alpha_3 \ar@{.}[dddd] \ar@/^/[dr] \\
\alpha_2 \ar@{.}[uuu] \ar@/^/[ur] &&&& \alpha_5 \ar@{.}[dddlll] \ar@/^/[dll] \\
              &                  & \alpha_4 \ar@{.}[uuuu] \ar@/^/[ull] \\
&&&& \\
& \alpha_4 && \alpha_5
}}
\]
or, in circuit notation, as
\[\vcenter{\begin{xy} 
(2,0)*+{\alpha_1}="x1",
(8,0)*+{\alpha_2}="x2",
(14,0)*+{\alpha_3}="x3",
(5,-14)*+{}="f1",
(8,-14)*+{}="f2",
(10,-14)*+{}="f3",
(8,-14)*{\framebox[1cm][c]{}}="f",
(6,-14)*+{}="f4",
(10,-14)*+{}="f5",
(5,-28)*+{\alpha_4}="x4",
(11,-28)*+{\alpha_5}="x5",
\ar@{-}@(d,u) "x1";"f1"
\ar@{-}@(dl,dl) "x2";"f4"
\ar@{-}@(d,u) "x3";"f2"
\ar@{-}@(u,d) "x4";"f5"
\ar@{-}@(u,ur) "x5";"f3"
\end{xy}}
\]

Note that we are labeling the channels where one typically labels the objects,
however, this is natural as, in this example, there is only one object.

This polymorphism may also be denoted using a cyclic permutation presentation
as follows
\[
(\alpha_1,\alpha_3,\alpha_5,\alpha_4,\alpha_2):
\<\alpha_1,\alpha_2,\alpha_3\> \ra \<\alpha_4,\alpha_5\>
\]

Any cyclic permutation of a polymorphism represents the same polymorphism.
Thus, an equivalent presentation of the above polymorphism is
\[(\alpha_2,\alpha_1,\alpha_3,\alpha_5,\alpha_4):
\<\alpha_1,\alpha_2,\alpha_3\> \ra \<\alpha_4,\alpha_5\>
\]
which can be presented in circuit notation as
\[\vcenter{\begin{xy} 
(2,0)*+{\alpha_1}="x1",
(8,0)*+{\alpha_2}="x2",
(14,0)*+{\alpha_3}="x3",
(4,-14)*+{}="f1",
(8,-14)*+{}="f2",
(12,-14)*+{}="f3",
(8,-14)*{\framebox[1cm][c]{}}="f",
(6,-14)*+{}="f4",
(10,-14)*+{}="f5",
(5,-28)*+{\alpha_4}="x4",
(11,-28)*+{\alpha_5}="x5",
\ar@{-}@(d,u) "x1";"f2"
\ar@{-}@(d,u) "x2";"f1"
\ar@{-}@(d,u) "x3";"f3"
\ar@{-}@(u,d) "x4";"f4"
\ar@{-}@(u,d) "x5";"f5"
\end{xy}}
\]

The identity polymorphism in this setting is given by
$(\alpha,\beta):\<\alpha\> \ra \<\beta\>$.

Composition is given by ``gluing'' two composible cyclic graphs together
(to make larger cycles). For example, the composite of
\[\vcenter{\xymatrix@R=1ex@C=1.3ex{
              & \alpha_1 \ar[r]  & \alpha_2 \ar@/^/[dr] \\
\alpha_5 \ar@/^/[ur] &               &               & \alpha_3 \ar@/^/[dl] \\
              & \alpha_4 \ar@/^/[ul] & x \ar[l]
}}:\<\alpha_1,\alpha_2,\alpha_3\> \ra \<\alpha_4,\alpha_5,x\>
\]
and
\[\vcenter{\xymatrix@R=1ex@C=1.3ex{
                & \beta_1 \ar@/^/[dr] \\
x \ar@/^/[ur] &                 & \beta_2 \ar@/^/[dl] \\
                & \beta_3 \ar@/^/[ul] 
}}:\<x,\beta_1\> \ra \<\beta_2,\beta_3\>
\]
at the node $x$ is
\[\vcenter{\xymatrix@R=2ex@C=2ex{
              & \alpha_1 \ar[r]  & \alpha_2 \ar@/^/[dr] \\
\alpha_5 \ar@/^/[ur] &               &               & \alpha_3 \ar[d] \\
\alpha_4 \ar[u]  &               &               & \beta_1 \ar@/^/[dl] \\
              & \beta_3 \ar@/^/[ul] & \beta_2 \ar[l]
}}:\<\alpha_1,\alpha_2,\alpha_3,\beta_1\> \ra \<\alpha_4,\alpha_5,\beta_2,\beta_3\>
\]

This data can now be seen to satisfy the requirements of a symmetric
polycategory. The actions of the symmetric group on these polymorphisms
simply permute the domain and codomain of the nodes in the signature. That is
\[c_{\sigma,\tau}(G:\Gamma \ra \Delta) = G:\sigma\Gamma \ra \tau\Delta
\]
For example, if $\sigma(\<\alpha_1,\alpha_2,\alpha_3\>) =
\<\alpha_2,\alpha_1,\alpha_3\>$ and $\tau(\<\alpha_4,\alpha_5\>) =
\<\alpha_5, \alpha_4\>$ then 
\begin{align*}
c_{\sigma,\tau}((\alpha_1,\alpha_3,\alpha_5,\alpha_4,\alpha_2):
\<\alpha_1,\alpha_2,\alpha_3\> \ra \<\alpha_4,\alpha_5\>) \\
=(\alpha_1,\alpha_3,\alpha_5,\alpha_4,\alpha_2):
\<\alpha_2,\alpha_1,\alpha_3\> \ra \<\alpha_5,\alpha_4\>
\end{align*}
which may be presented graphically as
\[\vcenter{\xymatrix@R=2ex@C=1ex{
\alpha_2 && \alpha_1 && \alpha_3 \\
&&&& \\
& \alpha_1 \ar@{.}[uurrr] \ar@/^/[rr] && \alpha_3 \ar@{.}[ddddll]\ar@/^/[dr] \\
\alpha_2 \ar@{.}[uuurr] \ar@/^/[ur] &&&& \alpha_5 \ar@{.}[dddl] \ar@/^/[dll] \\
            &                & \alpha_4 \ar@{.}[uuuull] \ar@/^/[ull] \\
&&&& \\
& \alpha_5 && \alpha_4
}}
\]

This polycategory is, in fact, the polyad of a non-commutative cyclic monoid
(in any linearly distributive category).
\end{example}

\chapter{$\Sigma\Pi$-Poly Logic} \label{chap-logic}

In this chapter we introduce $\Sigma\Pi$-poly logic; this logic is an
extension of the logic $\Sigma\Pi$ of Cockett and Seely~\cite{cockett01:finite}
which allows an arbitrary number of formulas on each side of the turnstile
instead of exactly one. In what follows we will refer to this logic (the
$\Sigma\Pi$-poly logic) as $\Sigma\Pi$.
If there is a need to differentiate between the two logics it will be
explicitly mentioned which logic, the $\Sigma\Pi$-poly logic or the
Cockett-Seely $\Sigma\Pi$-logic, is under consideration.

In Section~\ref{sec-sc} the sequent calculus for $\Sigma\Pi$ is introduced.
In Section~\ref{sec-ce}, the cut-elimination rewrites for
$\Sigma\Pi$-derivations are presented. The proof that these rewrites terminate
is left until the next chapter. In $\Sigma\Pi$ the identity axiom applies only
in the atomic case; in Section~\ref{sec-id}, it is shown that the identity
holds for arbitrary formulas. Finally in Section~\ref{sec-pc}, we give some
equivalence schema, the so-called ``permuting conversions'', which are
necessary in order for our system to have the Church-Rosser property.

\section{The sequent calculus} \label{sec-sc}

The logic is presented in a Gentzen sequent style: a sequent takes the form
$\Gamma \vd \Delta$, where the antecedent $\Gamma$, and the succedent
$\Delta$, are comma separated strings of formulas. For convenience, we take
the strings of formulas to be unordered which then eliminates the need for the
exchange rules:
\[\vcenter{\infer{\Gamma,Y,X,\Gamma' \vd \Delta}{\Gamma,X,Y,\Gamma' \vd \Delta}}
\qqquad \text{and} \qqquad
\vcenter{\infer{\Gamma \vd \Delta,Y,X,\Delta'}{\Gamma \vd \Delta,X,Y,\Delta'}}
\]

The propositions are either atoms (which we write as $A,B,\ldots$) or compound
formulas (which we write as $X,Y,\ldots$). A compound formula is either an
$I$-ary \textbf{sum}, where $I$ is a finite set, denoted $\sum_{i \in I} X_i$,
or a \textbf{ product}, denoted $\prod_{i \in I} X_i$. The index set $I$ may be
empty which gives the empty sum and product, denoted $\sum_\emptyset = 0$ and
$\prod_\emptyset = 1$ respectively. For sufficiently small index sets we may
write out the sum or product explicitly, e.g., $X+Y$ or $X \times Y$.

A typical rule in the Gentzen sequent style looks like
\[
\infer{\Gamma \vd \Delta}
{\Gamma_1 \vd \Delta_1 & \Gamma_2 \vd \Delta_2 & \cdots & \Gamma_n \vd 
\Delta_n}
\]
where we read this as an inference from top-to-bottom, i.e., if all the 
sequents on the top (the $\Gamma_i \vd \Delta_i$'s) can be derived, then the 
sequent on the bottom ($\Gamma \vd \Delta$) may be inferred. The special case 
is when there are no sequents on top

\[\infer{\Gamma \vd \Delta}{}
\]
in which case we may simply infer the bottom with no assumptions. A proof
(or deduction) of a sequent is a finite tree with the given sequent at the
root, axioms at the leaves and internal nodes corresponding to inference
rules. A Gentzen style presentation together with its inference rules is
called a \textbf{sequent calculus}. The rules of inference for $\Sigma\Pi$
are as follows:

\begin{center}
\ovalbox{
\parbox{10cm}{\begin{center}
\medskip
$\infer[\text{(identity)}]{A \vdash A}{}$ \bigskip\\
$\infer[\text{(cotuple)}]
   {\hspace{-1ex} \Gamma,\sum\limits_{i \in I} X_i \vd \Delta \hspace{1ex}}
   {\{\Gamma,X_i \vd \Delta\}_{i \in I}}
\qquad
\infer[\text{(tuple)}]
   {\hspace{-1ex} \Gamma \vd \prod\limits_{i \in I} Y_i, \Delta \hspace{1ex}}
   {\{\Gamma \vd Y_i,\Delta\}_{i \in I}}$ \bigskip\\
$\infer[\text{(projection)}]
   {\Gamma,\prod\limits_{i \in I} X_i \vd \Delta}
   {\Gamma,X_k \vd \Delta}
\qquad
\infer[\text{(injection)}]
   {\Gamma \vd \sum\limits_{i \in I} Y_i,\Delta}
   {\Gamma \vd Y_k,\Delta}$
\centro{where $k \in I, I \not= \emptyset$} \smallskip
$\infer[\text{(cut)}]{\Gamma,\Gamma' \vd \Delta,\Delta'}
      {\Gamma \vd \Delta,Z & Z,\Gamma' \vd \Delta'}$
\medskip
\end{center}}}
\end{center}

Notice that in the cotuple and tuple rules the index set $I$ may be empty, 
though not in the injection and projection rules. 

Observe that the inference system for $\Sigma\Pi$ is self-dual, that is, it
has an obvious sum-product symmetry. Explicitly, we may swap the direction of
the sequents while turning sums into products and products into sums to
obtain the same system. This means that each proof has a dual interpretation
and can be ``reused'' to prove a dual theorem.

We shall consider various augmentations of this basic logic:

\begin{itemize}
\item The ``initial logic'' is the logic with no atoms. Notice that this is
still a non-trivial logic because of the symbols $\Sigma_\emptyset$ and 
$\Pi_\emptyset$ which we may use to construct more complex formulas.  We 
shall denote this logic as $\Sigma\Pi_\emptyset$.

\item The ``pure logic'' is the logic with an arbitrary set of atoms $A$:
we shall denote this logic as $\Sigma\Pi_A$. 

\item The ``free logic'' is the logic with an arbitrary set of atoms and 
an arbitrary set of non-logical axioms relating lists of atoms. If $f$ is
a non-logical axiom from $A,B$ to $C,D$, this may be denoted as $f:A,B \ra
C,D$ or as an inference $\infer{A,B \vd_f C,D}{}$. The atoms will be
regarded as objects in a polycategory and the axioms as maps in 
that polycategory (with the ``essential cuts'' being provided by composition).
If the polycategory is denoted \cat{A}, the resulting logic will be denoted as
$\Sigma\Pi_\cat{A}$. 
\end{itemize}

If we think of the atoms of a pure logic as forming a discrete category, the
free logic on this discrete category is then just the ``pure'' logic. Each
variant above therefore includes the previous variant, and as it is more
general, we shall tend to consider only this last variant.

\begin{example}[$\Sigma\Pi$-derivations] \label{exam-logic}
Some typical proofs in this logic.
\begin{enumerate}
\item 
\[\infer{A + (B \times C) \vd (A+B) \times (A+C)}{
   \infer{A \vd (A+B) \times (A+C)}{ 
      \infer{A \vd A+B}{
         \infer{A \vd A}{}} &
      \infer{A \vd A+C}{
         \infer{A \vd A}{}}} &
   \infer{B \times C \vd (A+B) \times (A+C)}{ 
      \infer{B \times C \vd A+B}{
         \infer{B \vd A+B}{
            \infer{B \vd B}{}}} &
      \infer{B \times C \vd A+C}{
         \infer{C \vd A+C}{
            \infer{C \vd C}{}}}}}
\]

This is one direction of the proof that sums distribute over products. The
other direction cannot be proved in our system.

\item 
\[\infer{(A \times B) + C \vd (A \times B) + C}{
   \infer{A \times B \vd (A \times B) + C}{
      \infer{A \times B \vd A \times B}{
            \infer{A \times B \vd A}{
               \infer{A \vd A}{}} &
            \infer{A \times B \vd B}{
               \infer{B \vd B}{}}}} &
   \infer{C \vd (A \times B) + C}{
         \infer{C \vd C}{}}}
\]

This proves the identity inference. 

\item The above examples have only one proposition on each side of the 
turnstile. Here is an example in which each side has more than one
proposition. In this case we will need some non-logical axioms.
\[\begin{array}{ll}
i:A,E \ra G,I \quad & k:C,F \ra H,I\\
j:B,E \ra G,J       & l:D,F \ra H,K
\end{array}
\]

{\small
\[\infer{(A \times B)+(C \times D),\ E \times F \vd G+H,\ I \times (J+K)}{
\infer{A \times B,\ E \times F \vd G+H,\ I \times (J+K)}{
  \infer{A \times B,\ E \vd G+H,\ I \times (J+K)}{
    \infer{A \times B,\ E \vd G,\ I \times (J+K)}{
        \infer{A \times B,E \vd G,I}{
          \infer{A,\ E \vd_i G,\ I}{}} &
        \infer{A \times B,E \vd G,J+K}{
          \infer{B,\ E \vd G,\ J+K}{
            \infer{B,\ E \vd_j G,\ J}{}}}}}} &
\infer{C \times D,\ E \times F \vd G+H,\  I \times (J+K)}{
  \infer{C \times D,\ F \vd G+H,\  I \times (J+K)}{
    \infer{C \times D,\ F \vd H,\  I \times (J+K)}{
      \infer{C \times D,F \vd H,I}{
          \infer{C,\ F \vd_k H,\ I}{}} &
      \infer{C \times D,F \vd H,J+K}{
        \infer{D,\ F \vd H,\ J+K}{
          \infer{D,\ F \vd_l H,\ K}{}}}}}}}
\]}

\end{enumerate}
\end{example}

\section{Cut-elimination} \label{sec-ce}

In this section we show that cut-elimination holds for the free logic
$\Sigma\Pi_\cat{A}$, i.e., any proof can be rewritten so that it does
not contain any applications of the cut rule. Of course, this process will
get stuck on the introduced atomic polymorphisms. A cut between atomic axioms
is called an \textbf{essential cut}: 
\[\infer{\Gamma_1,\Delta_1 \vd_{f ; g} \Gamma_2,\Delta_2}
{\Gamma_1 \vd_f \Gamma_2,A & A,\Delta_1 \vd_g \Delta_2}\]

\begin{proposition}[Cut-elimination]
Any proof in the free logic $\Sigma\Pi_\cat{A}$ can be transformed to a
proof in which the only cuts are essential.
\end{proposition}

We shall provide a family of rewrites for $\Sigma\Pi$-derivations and show
that they terminate. As is typical in cut-elimination proofs, the rewrites will
either replace a cut by cuts involving simpler formulas, or ``push up'' a cut
into the surrounding proof. A proof that cannot be further rewritten using this
set of rewrites will be a ``cut-eliminated'' proof in the sense of having no
inessential cuts. 

The rewrites are as follows. It will always be assumed that $i \in I$ and $j
\in J$ for index sets $I$ and $J$. Duality will be used to reduce the number
of rewrites presented.

\begin{itemize}

\item Sequent-identity (identity-sequent): This rewrite removes the cut
below an identity axiom on the right.
\[
\vcenter{\infer{\Gamma \vd \Delta,A}{
   \infer{\Gamma \vd \Delta,A}{\pi} &
   \infer{A \vd A}{\iota}}}
\quad \Lra \quad \vcenter{
\infer{\Gamma \vd \Delta,A}{\pi}}
\]
The dual of this rewrite removes the cut below an identity axiom on the left.

\item Cotuple-sequent (sequent-tuple): This rewrite moves a cut which is
below a cotupling and an arbitrary sequent above the cotupling. 
\[
\vcenter{\infer{\Gamma,\sum X_i,\Gamma' \vd \Delta,\Delta'}{
   \infer{\Gamma,\sum X_i \vd \Delta,Z}{ 
      \left\{ \raisebox{-1.5ex}{ 
      \infer{\Gamma,X_i \vd \Delta,Z}{\pi_i}} \right\}_i} &
   \infer{Z,\Gamma' \vd \Delta'}{\pi}}} 
\quad \Lra \quad \vcenter{
\infer{\Gamma,\sum X_i,\Gamma' \vd \Delta,\Delta'}{ & \hspace{-2ex}
   \left\{\raisebox{-3ex}{
   \infer{\Gamma,X_i,\Gamma' \vd \Delta,\Delta'}{
     \infer{\Gamma,X_i \vd \Delta,Z}{\pi_i} &
   \infer{Z,\Gamma' \vd \Delta'}{\pi}}} \right\}_i & \hspace{-2ex}}}
\]
The dual of this rewrite moves the cut above a tupling on the right.

\item Injection-sequent (sequent-projection): This rewrite moves a cut which
is below an injection and an arbitrary sequent above the injection. There are
two cases to consider: the cut is on the injection, or it is on an arbitrary
formula. The rewrite for the former case is below (the injection-cotuple
rewrite); the rewrite for the latter case is as follows. 
\[
\vcenter{\infer{\Gamma,\Gamma' \vd \sum X_i,\Delta,\Delta'}{
   \infer{\Gamma \vd \sum X_i,\Delta,Z}{
      \infer{\Gamma \vd X_k,\Delta,Z}{\pi}} &
   \infer{Z,\Gamma' \vd \Delta'}{\pi'}}} 
\quad \Lra \quad \vcenter{
\infer{\Gamma,\Gamma' \vd \sum X_i,\Delta,\Delta'}{
   \infer{\Gamma,\Gamma' \vd X_k,\Delta,\Delta'}{
      \infer{\Gamma \vd X_k,\Delta,Z}{\pi} &
   \infer{Z,\Gamma' \vd \Delta'}{\pi'}}}}
\]
The dual of this rewrite moves the cut above a projection on the right.

\item Projection-sequent (sequent-injection): This rewrite moves a cut
which is below a projection and an arbitrary sequent above the projection.
\[
\vcenter{\infer{\Gamma,\prod X_i,\Gamma' \vd \Delta,\Delta'}{
   \infer{\Gamma,\prod X_i \vd \Delta,Z}{
      \infer{\Gamma,X_k \vd \Delta,Z}{\pi}} &
   \infer{Z,\Gamma' \vd \Delta'}{\pi'}}} 
\quad \Lra \quad \vcenter{
\infer{\Gamma,\prod X_i,\Gamma' \vd \Delta,\Delta'}{
   \infer{\Gamma,X_k,\Gamma' \vd \Delta,\Delta'}{
      \infer{\Gamma,X_k \vd \Delta,Z}{\pi} &
   \infer{Z,\Gamma' \vd \Delta'}{\pi'}}}}
\]
The dual of this rewrite moves the cut above an injection.

\item Tuple-sequent (sequent-cotuple): This rewrite moves a cut which is
below a tuple and an arbitrary sequent above the tuple. There are two cases
to consider: the cut is on the tupling, or the cut is on an arbitrary
formula. The rewrite for the former case is dual to the injection-cotuple
rewrite below; the rewrite for the latter case is as follows.
\[
\vcenter{\infer{\Gamma,\Gamma' \vd \prod X_i,\Delta,\Delta'}{
   \infer{\Gamma \vd  \prod X_i,\Delta,Z}{
      \left\{ \raisebox{-1.5ex}{ 
      \infer{\Gamma \vd  X_i,\Delta,Z}{\pi_i}} \right\}_i} &
   \infer{Z,\Gamma' \vd \Delta'}{\pi}}} 
\quad \Lra \quad \vcenter{
\infer{\Gamma,\Gamma' \vd \prod X_i,\Delta,\Delta'}{& \hspace{-2ex}
   \left\{\raisebox{-3ex}{
    \infer{\Gamma,\Gamma' \vd X_i,\Delta,\Delta'}{
              \infer{\Gamma \vd X_i,\Delta,Z}{\pi_i} &
              \infer{Z,\Gamma' \vd \Delta'}{\pi}
           }} \right\}_i & \hspace{-2ex}}}
\]
The dual of this rewrite moves the cut above a cotupling on the right.

\item Injection-cotuple (tuple-projection): This rewrite moves the cut
above an injection and cotupling.
\[
\vcenter{\infer{\Gamma,\Gamma' \vd \Delta,\Delta'}{
  \infer{\Gamma \vd \Delta,\sum X_i}{
    \infer{\Gamma \vd \Delta,X_k}{\pi}} &
  \infer{\sum X_i,\Gamma' \vd \Delta'}{
      \left\{\raisebox{-1.5ex}{\infer{X_i,\Gamma' \vd \Delta'}{\pi_i}}
      \right\}_i}}}
\quad \Lra \quad \vcenter{
\infer{\Gamma,\Gamma' \vd \Delta,\Delta'}{
   \infer{\Gamma \vd \Delta,X_k}{\pi} &
   \infer{X_k,\Gamma' \vd \Delta'}{\pi_k}}}
\]
The dual of this rewrite moves the cut above a tupling and projection.

\end{itemize}

This accounts for all the ways in which compound formulas are introduced
either on the left or on the right above a cut, and we have shown how to move
the cut above these rules. Thus, a proof which cannot be rewritten further
must have an axiom above the cut on each side; this is an essential cut.

It remains to show that this rewriting procedure terminates. For this we will
need to define a measure on proofs which is reduced by each rewrite. The
technical details and proof that the rewritings terminate will be presented
in Chapter~\ref{chap-syntax}.

\section{Identity derivations} \label{sec-id}

Our goal is to view this proof system as a polycategory where cut is the
composition. The cut-elimination process therefore provides part of the
dynamics of composition: the activity which takes place when two proofs
are plugged together. 

Part of proving that cut acts as a composition is showing that there are
identity derivations which behave in the correct manner. The identity
derivations are defined inductively as follows.

\begin{itemize}

\item The identity atomic sequent: \[\infer[\text{(identity)}]{A \vd A}{}\]

\item The identity derivation on sums is given by
\[
\infer{\sum X_i \vd \sum X_i}{
   \left\{ \raisebox{-3ex}{ 
      \infer{X_i \vd \sum X_i}{
         \infer{X_i \vd X_i}{\iota_{X_i}}}} 
   \right\}_{i}}
\]
where the identity derivation $\iota_{X_i}$ of $X_i \vd X_i$ is given
by induction on the structure of $X_i$.

\item The identity on products is given by the dual of the proof above.
Explicitly,
\[
\infer{\prod X_i \vd \prod X_i}{
   \left\{ \raisebox{-3ex}{ 
      \infer{\prod X_i \vd X_i}{
         \infer{X_i \vd X_i}{\iota_{X_i}}}} 
   \right\}_{i}}
\]
where the identity derivation $\iota_{X_i}$ of $X_i \vd X_i$ is given
by induction on the structure of $X_i$.

\end{itemize}

The following lemma now proves that cutting (composing) any sequent
derivation together with the identity derivation (and vise versa)
results in the same sequent derivation. 

\begin{lemma} \label{lem-id}
The sequent-identity and identity-sequent cut-elimination reductions are
derivable reductions for the general identity derivations as defined above.
That is,
\[
\vcenter{\infer{\Gamma \vd \Delta,X}{
   \infer{\Gamma \vd \Delta,X}{\pi} &
   \infer{X \vd X}{\iota_X}}}
\quad \Lra \quad \vcenter{
\infer{\Gamma \vd \Delta,X}{\pi}}
\]
and similarly for the dual rule.
\end{lemma}

\begin{proof}
We shall suppose that the identity derivation is on the right; duality covers
the other case. The proof is by structural induction on the derivation $\pi$.

\begin{itemize}

\item The base case is a cut with an atomic sequent: here the cut-elimination
step removes the atomic identity and the result is immediate.

\item Next suppose the identity is on a sum type:
\[
\infer{\Gamma \vd \Delta,\sum X_i}{
   \infer{\Gamma \vd \Delta,\sum X_i}{\pi} 
  &\infer{\sum X_i \vd \sum X_i}{\iota}}
\]

There are four possibilities for the root inference of $\pi$.

(1) If the root inference is a cotupling the cut-elimination step moves the
cut onto smaller proofs. We may now apply the inductive hypothesis to each of
these smaller proofs which yields the required result.
{\small\[
\vcenter{\infer{\Gamma,\sum Y_j \vd \Delta,\sum X_i}{
   \infer{\Gamma,\sum Y_j \vd \Delta,\sum X_i}{ 
      \left\{ \raisebox{-1.5ex}{ 
      \infer{\Gamma, Y_j \vd \Delta,\sum X_i}{\pi_j}} \right\}_{j}} &
   \infer{\sum X_i \vd \sum X_i}{\iota}}} 
\ \Lra \ 
\vcenter{\infer{\Gamma,\sum Y_j \vd \Delta,\sum X_i}{ & \hspace{-2ex}
   \left\{\raisebox{-4ex}{
   \infer{\Gamma,Y_j \vd \Delta,\sum X_i}{
     \infer{\Gamma,Y_j \vd \Delta,\sum X_i}{\pi_j} &
   \infer{\sum X_i \vd \sum X_i}{\iota}}} \right\}_{j} & \hspace{-2ex}}}
\]}

(2) If the root inference is an injection there are two possibilities:
the cut is on $\sum X_i$ or it is not. In the first case we use the
injection-cotuple rewrite and the second case the injection-sequent rewrite.
The injection-cotuple rewrite is as follows.
\[\vcenter{\infer{\Gamma \vd \Delta,\sum X_i}{
   \infer{\Gamma \vd \Delta,\sum X_i}{
      \infer{\Gamma \vd \Delta,X_k}{\pi}} &
   \infer{\sum X_i \vd \sum X_i}{\left\{\raisebox{-3ex}{
     \infer{X_i \vd \sum X_i}{
       \infer{X_i \vd X_i}{\iota_i}}} \right\}_i}}}
\quad \Lra \quad \vcenter{
\infer{\Gamma \vd \Delta,\sum X_i}{
   \infer{\Gamma \vd \Delta,X_k}{
      \infer{\Gamma \vd \Delta,X_k}{\pi} &
      \infer{X_k \vd X_k}{\iota_k}}}}
\]
which moves the cut onto a smaller proof and so the inductive hypothesis
applies and we are done. The injection-sequent rewrite is
{\small\[
\vcenter{\infer{\Gamma \vd \sum Y_j,\Delta,\sum X_i}{
   \infer{\Gamma \vd \sum Y_j,\Delta,\sum X_i}{
      \infer{\Gamma \vd Y_k,\Delta,\sum X_i}{\pi}} &
   \infer{\sum X_i \vd \sum X_i}{\iota}}} 
\ \Lra\ \vcenter{
\infer{\Gamma \vd \sum Y_j,\Delta,\sum X_i}{
   \infer{\Gamma \vd Y_k,\Delta,\sum X_i}{
      \infer{\Gamma \vd Y_k,\Delta,\sum X_i}{\pi} &
   \infer{\sum X_i \vd \sum X_i}{\iota}}}}
\]}
which again moves the cut onto a smaller proof

(3) If the root inference is a projection the cut-elimination step moves the
cut onto a smaller proof and so we are done. 
{\small\[
\vcenter{\infer{\Gamma,\prod Y_j \vd \Delta,\sum X_i}{
   \infer{\Gamma,\prod Y_j \vd \Delta,\sum X_i}{
      \infer{\Gamma,Y_k \vd \Delta,\sum X_i}{\pi}} &
   \infer{\sum X_i \vd \sum X_i}{\iota}}} 
\quad \Lra \quad \vcenter{
\infer{\Gamma,\prod Y_j \vd \Delta,\sum X_i}{
   \infer{\Gamma,Y_k \vd \Delta,\sum X_i}{
      \infer{\Gamma,Y_k \vd \Delta,\sum X_i}{\pi} &
   \infer{\sum X_i \vd \sum X_i}{\iota}}}}
\]}

(4) If the root inference is a tupling the cut-elimination step moves the
cut onto a smaller proof and so we are done. 
{\small\[
\vcenter{\infer{\Gamma \vd \prod Y_j,\Delta,\sum X_i}{
   \infer{\Gamma \vd \prod Y_j,\Delta,\sum X_i}{
      \left\{ \raisebox{-1.5ex}{ 
      \infer{\Gamma \vd Y_j,\Delta,\sum X_i}{\pi_j}} \right\}_{j}} &
   \infer{\sum X_i \vd \sum X_i}{\iota}}}
\ \Lra \
\vcenter{\infer{\Gamma \vd \prod Y_j,\Delta,\sum X_i}{& \hspace{-2ex}
   \left\{\raisebox{-4ex}{
    \infer{\Gamma \vd Y_j,\Delta,\sum X_i}{
              \infer{\Gamma \vd Y_j,\Delta,\sum X_i}{\pi_j} &
              \infer{\sum X_i \vd \sum X_i}{\iota}
           }} \right\}_{j}  & \hspace{-2ex}}}
\]}

\item Finally, suppose the identity is on a product type:
\[\infer{\Gamma \vd \Delta,\prod X_i}{
   \infer{\Gamma \vd \Delta,\prod X_i}{\pi} 
  &\infer{\prod X_i \vd \prod X_i}{\iota}}
\]

This case is very similar to the one above. In all the cases for the root
inference of $\pi$, other than tupling, the appropriate cut-elimination
rewrite (e.g., when the root inference is a cotupling use the cotuple-sequent
rewrite) will move the cut onto a smaller proof (or smaller proofs).

For the case where the root inference is a tupling there are two subcases
to consider: the cut is on $\prod X_i$ or it is not. In the latter case the
tuple-sequent cut-elimination rewrite suffices to move the cut onto smaller
proofs. The former case is unique in that we must apply two rewrites: a
sequent-tuple followed by a tuple-projection.

{\small\allowdisplaybreaks\begin{align*}
\vcenter{\infer{\Gamma \vd \Delta,\prod X_i}{
   \infer{\Gamma \vd \Delta,\prod X_i}{\left\{\raisebox{-1.5ex}{
      \infer{\Gamma \vd \Delta,X_i}{\pi_i}} \right\}} &
   \infer{\prod X_i \vd \prod X_i}{\left\{\raisebox{-3ex}{
     \infer{\prod X_i \vd X_i}{
       \infer{X_i \vd X_i}{\iota_i}}} \right\}_{i}}}}
&\quad \Lra \quad \vcenter{
\infer{\Gamma \vd \Delta,\prod X_i}{& \hspace{-2ex} \left\{\raisebox{-6ex}{
  \infer{\Gamma \vd \Delta,X_i}{
    \infer{\Gamma \vd \Delta,\prod X_i}{\left\{\raisebox{-1.5ex}{
      \infer{\Gamma \vd \Delta,X_i}{\pi_i}} \right\}_{i}} &
    \infer{\prod X_i \vd X_i}{
      \infer{X_i \vd X_i}{\iota_i}}}} \right\}_{i} & \hspace{-2ex}}} 
\bigskip\medskip\\
&\quad \Lra \quad \vcenter{
\infer{\Gamma \vd \Delta,\prod X_i}{& \hspace{-2ex} \left\{\raisebox{-3ex}{
  \infer{\Gamma \vd \Delta,X_i}{
    \infer{\Gamma \vd \Delta,X_i}{\pi_i} &
    \infer{X_i \vd X_i}{\iota_i}}} \right\}_{i} & \hspace{-2ex}}}
\end{align*}}
\qed
\end{itemize}
\end{proof}

\section{Permuting conversions} \label{sec-pc}

In order to obtain a normal form for sequent derivations, we would like to
show that the cut-elimination rewrites are Church-Rosser. Currently this 
is not the case; for example, consider a derivation with a 
cotupling and tupling immediately above the cut:
\[
\infer{\Gamma,\sum_i X_i,\Gamma' \vd \Delta,\prod_j Y_j,\Delta'}{
\infer{\Gamma,\sum_i X_i \vd \Delta,Z}{\left\{ \raisebox{-1.5ex}{
  \infer{\Gamma,X_i \vd \Delta,Z}{\pi_i}} \right\}_i} &
\infer{Z,\Gamma' \vd \prod_j Y_j,\Delta'}{\left\{\raisebox{-1.5ex}{
  \infer{Z,\Gamma' \vd Y_j,\Delta'}{\pi'_j}} \right\}_j}}
\]

In this case one may apply the cotuple-sequent rewrite or the sequent-tuple
rewrite to reduce the derivation, but there seems to be no way in which to
resolve these rewrites. This motivates the use of additional rewrites which
will allow us to interchange these two rules. Similar considerations for the
nine other possible critical pairs (cotuple vs. cotuple, cotuple vs. injection,
etc.) leads us to the following ten conversions (which we denote by $\pc$).
Again, duality is used to reduce the number of conversions presented.

\begin{itemize}

\item Cotuple-cotuple (tuple-tuple) interchange:
\[\vcenter{\infer{\sum X_i,\Gamma,\sum Y_j \vd \Delta}{
  \left\{ \raisebox{-4ex}{
  \infer{X_i,\Gamma,\sum Y_j \vd \Delta}{
    \left\{ \raisebox{-1.5ex}{
    \infer{X_i,\Gamma,Y_j \vd \Delta}{\pi_{ij}}}
    \right\}_j }} \right\}_i}} 
\quad \pc \quad \vcenter{
\infer{\sum X_i,\Gamma,\sum Y_j \vd \Delta}{
  \left\{ \raisebox{-4ex}{
  \infer{\sum X_i,\Gamma,Y_j \vd \Delta}{
    \left\{ \raisebox{-1.5ex}{
    \infer{X_i,\Gamma,Y_j \vd \Delta}{\pi_{ij}}}
    \right\}_i }} \right\}_j}}
\]

\item Cotuple-injection (projection-tuple) interchange:
\[\vcenter{\infer{\Gamma,\sum X_i \vd \sum Y_j,\Delta}{
  \left\{ \raisebox{-3ex}{
  \infer{\Gamma,X_i \vd \sum Y_j,\Delta}{
    \infer{\Gamma,X_i \vd Y_k,\Delta}{\pi_i}}
  } \right\}_i}} 
\quad \pc \quad \vcenter{
\infer{\Gamma,\sum X_i \vd \sum Y_j,\Delta}{
  \infer{\Gamma,\sum X_i \vd Y_k,\Delta}{
    \left\{ \raisebox{-1.5ex}{
    \infer{\Gamma,X_i \vd Y_k,\Delta}{\pi_i}
    } \right\}_i}}}
\]

\item Cotuple-projection (injection-tuple) interchange:
\[\vcenter{
\infer{\sum X_i,\Gamma,\prod Y_j \vd \Delta}{
  \left\{ \raisebox{-3ex}{
  \infer{X_i,\Gamma,\prod Y_j \vd \Delta}{
    \infer{X_i,\Gamma,Y_k \vd \Delta}{\pi_i}}} 
  \right\}_i}} 
\quad \pc \quad \vcenter{
\infer{\sum X_i,\Gamma,\prod Y_j \vd \Delta}{
  \infer{\sum X_i,\Gamma,Y_k \vd \Delta}{
    \left\{ \raisebox{-1.5ex}{   
    \infer{X_i,\Gamma,Y_k \vd \Delta}{\pi_i}} 
    \right\}_i}}}
\]

\item Cotuple-tuple interchange:
\[\vcenter{\infer{\Gamma,\sum X_i \vd \prod Y_j,\Delta}{
  \left\{ \raisebox{-4ex}{
  \infer{\Gamma,X_i \vd \prod Y_j,\Delta}{
    \left\{ \raisebox{-1.5ex}{
    \infer{\Gamma,X_i \vd Y_j,\Delta}{\pi_{ij}}}
    \right\}_j}} \right\}_i}} 
\quad \pc \quad \vcenter{
\infer{\Gamma,\sum X_i \vd \prod Y_j,\Delta}{
  \left\{ \raisebox{-4ex}{
  \infer{\Gamma,\sum X_i \vd Y_j,\Delta}{
    \left\{ \raisebox{-1.5ex}{
    \infer{\Gamma,X_i \vd Y_j,\Delta}{\pi_{ij}}}
    \right\}_i}} \right\}_j}}
\]

\item Injection-injection (projection-projection) interchange:
\[\vcenter{\infer{\Gamma \vd \sum X_i,\Delta,\sum Y_j}{
  \infer{\Gamma \vd X_k,\Delta,\sum Y_j}{
    \infer{\Gamma \vd X_k,\Delta,Y_{k'}}{\pi}}}} 
\quad \pc \quad \vcenter{
\infer{\Gamma \vd \sum X_i,\Delta,\sum Y_j}{
  \infer{\Gamma \vd \sum X_i,\Delta,Y_{k'}}{
    \infer{\Gamma \vd X_k,\Delta,Y_{k'}}{\pi}}}}
\]

\item Projection-injection interchange:
\[\vcenter{\infer{\Gamma,\prod X_i \vd \sum Y_j,\Delta}{
  \infer{\Gamma,X_k \vd \sum Y_j,\Delta}{
    \infer{\Gamma,X_k \vd Y_{k'},\Delta}{\pi}}}} 
\quad \pc \quad \vcenter{
\infer{\Gamma,\prod X_i \vd \sum Y_j,\Delta}{
  \infer{\Gamma,\prod X_i \vd Y_{k'},\Delta}{
    \infer{\Gamma,X_k \vd Y_{k'},\Delta}{\pi}}}}
\]

\end{itemize}

It can now be shown that the cut-elimination procedure for
$\Sigma\Pi$-deri\-va\-tions is confluent modulo these permuting conversions.
It is, however, more convenient to do so after providing a term
representation for the $\Sigma\Pi$-derivations. As this is the purpose of
the next chapter, the proof of confluence will also be presented there in
Section~\ref{sec-cr}.

\chapter{Term Logic for $\Sigma\Pi$-Poly Maps} \label{chap-syntax}

In this chapter we prove that the cut-elimination rewrites are Church-Rosser
and terminating modulo the permuting conversions. It will, therefore, be
convenient to have a more compact notation for sequent derivations; this
leads us to introduce a system of terms typed by sequents. 

The terms for the $\Sigma\Pi$-logic introduced in this chapter will reflect
the view of the formulas as protocols and the proofs as processes. That is,
that a proof of a sequent $\Gamma \vd \Delta$ may be regarded as a process
between the protocols (the formulas) in $\Gamma$ and $\Delta$. As we shall
see in the next chapter, these terms provide a categorical semantics for the
logic: the free polycategory with sums and products.

In order to motivate the term logic we shall start by introducing the
view of this system as protocols and processes. In Section~\ref{sec-pap},
we then introduce two term representations for $\Sigma\Pi$-derivations:
a compact term calculus representation and a ``programming language''
representation. A rewriting system for cut-elimination using these
representations is developed. In Section~\ref{sec-au}, it is made explicit
how the term calculus handles the additive units (the objects
$\Sigma_\emptyset$ and $\Pi_\emptyset$). Then, with our term calculus
representation in hand, in Sections \ref{sec-cr} and \ref{sec-cm}, we prove
that our cut-elimination rewriting system for the $\Sigma\Pi$-derivations
has the Church-Rosser property and is terminating. Finally, in
Section~\ref{sec-dp}, we give a decision procedure for $\Sigma\Pi$-derivations.

\section{Protocols and processes}

If we adopt the system/environment perspective it is natural to view
polymorphisms as processes and the objects as protocols: the system is          connected to the environment via channels along which the environment may       input to the system and the system may output to the environment, where the     channels here allow two-way communication, i.e., on any given channel one       may receive input and/or produce output.

Protocols are assigned to channels along which processes can communicate
(in both input and output mode). A protocol on a channel may be thought of
as a convention between the system and the environment that dictates who may
place the next event on that channel. We will adopt the convention of using
Greek letters for channel names. Let $S$ represent a system and $E$ the
environment such that they are connected via a channel $\alpha$. Suppose $P$
is a protocol which is only going to allow the system to output either the
event $a$ or the event $b$ to the environment. Placing $P$ on $\alpha$          between $S$ and $E$, $\xymatrix@1{S \ar@{<->}[r]^{\alpha:P} & E}$, allows $S$   to output either the event $a$ or the event $b$, and moreover, prevents $E$
from inputing anything to $S$.

A process should be regarded as an agent which acts on any number of channels
by either accepting input events or producing output events in accordance
with the protocols associated with each channel. Thus, it is a system
communicating with the environment along a number of channels. For reasons
that will become apparent, the channels a process interacts along are divided
into (unordered) domain (or ``input'') channels and (unordered) codomain (or
``output'') channels. A process $f$ with domain channels $\alpha_1, \ldots,
\alpha_n$ and codomain channels $\beta_1,\ldots,\beta_m$ may be represented
graphically as
\[
\begin{xy} (2,0)*+{\alpha_1}="a", (8,2)*+{}="b", (14,0)*+{\alpha_n}="c",
(5,-6)*+{}="d", (11,-6)*+{}="e",
(8,-8)*+{\fbox{$\quad f \quad$}}="f",
(2,-16)*+{\beta_1}="a'", (8,-19)*+{}="b'", (14,-16)*+{\beta_m}="c'",
(5,-10)*+{}="d'", (11,-10)*+{}="e'",
\ar@{-}@/_0.5ex/ "a";"d"
\ar@{-}@/^0.5ex/ "c";"e"
\ar@{} "b";"f" |\cdots
\ar@{-}@/^0.5ex/ "a'";"d'"
\ar@{-}@/_0.5ex/ "c'";"e'"
\ar@{} "b'";"f" |\cdots
\end{xy}
\]

Our convention will be to have the domain channels (the $\alpha_i$'s) on the
top and the codomain channels (the $\beta_j$'s) on the bottom. Recall that
we are considering these channels unordered. If a codomain channel
$\alpha$ of a process $f$, and a domain channel $\beta$ of a process $g$,
share a common protocol then $f$ and $g$ may be ``plugged'' (or composed)
together on $\alpha$ and $\beta$ to form a new process. The domain channels
of this composed process are the domain channels of $f$ and the domain
channels of $g$ (minus $\beta$) and the codomain channels are the codomain
channels of $f$ (minus $\alpha$) and the codomain channels of $g$. This may
be pictured graphically as
\[
\begin{xy} 
(0,3)*+{}="gamma", (14,3)*+{}="gamma'",
(0,-5)*{\framebox[1cm][c]{$f$}}="f",
(12,-17)*{\framebox[1cm][c]{$g$}}="g",
(-2,-25)*+{}="delta", (12,-25)*+{}="delta'"
\ar@{=} "gamma";"f"
\ar@{=}@(d,u) "gamma'";"g"
\ar@{-}@(dr,ul) "f";"g" |{\alpha:\beta}
\ar@{=} "g";"delta'"
\ar@{=}@(d,u) "f";"delta"
\end{xy}
\]
where the double lines represent strings of channels. (This plugging together
of processes is exactly the cut rule. This is explained in detail in
Section~\ref{ssec-cer}).

Note that any two processes may be plugged together on at most one channel.
That is, situations like the following are not allowed:
\[
\begin{xy} 
(0,3)*+{}="gamma", (14,3)*+{}="gamma'",
(0,-5)*{\framebox[1cm][c]{$f$}}="f",
(0,-17)*{\framebox[1cm][c]{$g$}}="g",
(0,-25)*+{}="delta", (14,-25)*+{}="delta'"
\ar@{=} "gamma";"f"
\ar@{-}@/^3ex/ "f";"g"
\ar@{-}@/_3ex/ "f";"g"
\ar@{=} "g";"delta"
\end{xy}
\]

Physically, one could think of the channels as ``wires'', and the processes
as black boxes which send events along these wires, where the protocols
determine which events are allowed to be passed though the wires at any given
moment. 

\section{Formulas as protocols} \label{sec-fap}

``Protocols'' in our system will simply be $\Sigma\Pi$-terms annotated with
channel names and ``events''. Formally, a \textbf{protocol} may be either:

\begin{itemize}
\item an \textbf{atomic} protocol: $A,B,C,\ldots$ (the objects of a
polycategory)

\item a \textbf{coproduct} protocol: 
$\sum_{i \in I} a_i:X_i = \{a_i:X_i \mid i \in I\}$
where each $X_i$ is a protocol. If we suppose $I = \{1,\ldots,n\}$, we could
represent this by the following tree:
\[\xymatrix@R=4ex@C=2ex@M=0.2ex{& \circ \ar@{-}[dl]_{a_1}  
\ar@{-}[dr]^{a_n} \\ X_1 &\dots& X_n}
\]

\item a \textbf{product} protocol: 
$\prod_{j \in J} b_j:Y_j = (b_j:Y_j \mid j \in J)$ where each $Y_j$ is a
protocol. Again, we could represent this by the following tree:
\[\xymatrix@R=4ex@C=2ex@M=0.2ex{& \bullet \ar@{-}[dl]_{b_1}  
\ar@{-}[dr]^{b_m} \\ Y_1 & \dots & Y_m}
\]
\end{itemize}
The $a_i$'s and $b_j$'s are thought of as the events. An event, as
introduced here, is not a formal object; for our purposes they may be
thought of as names (or tags, constructors, etc.).

Note that in specifying protocols we allow the index sets $I$ and $J$ to
be empty; this gives two atomic protocols
\[\sum\nolimits_\emptyset = \{\ \} = \circ \qquad \text{and} \qquad
\prod\nolimits_\emptyset = (\ ) = \bullet
\]

\begin{example}[Protocols] \label{exam-protocols} \quad
\begin{enumerate}
\item $X = \{a:\{d:A\}, b:(e:B,f:C), c:\{g:D,h:E\}\}$
\[\xymatrix@R=5ex@C=3ex@M=0.2ex{
&&&\circ \ar@{-}[dlll]_a \ar@{-}[d]_b \ar@{-}[drrr]^c \\ 
\circ \ar@{-}[d]_d &&& \bullet \ar@{-}[dl]_e \ar@{-}[dr]^f &&&
\circ \ar@{-}[dl]_g \ar@{-}[dr]^h \\ 
A && B && C & D && E}
\]

\item $Y = (a:\{c:A,d:B\},b:\{e:A,f:C\})$
\[\xymatrix@R=4ex@C=3ex@M=0.2ex{
&&& \bullet \ar@{-}[dll]_a \ar@{-}[drr]^b \\
& \circ \ar@{-}[dl]_c \ar@{-}[dr]^d &&&& \circ \ar@{-}[dl]_e
\ar@{-}[dr]^f \\
A && B && A && C}
\]

\item $Z = (a:(c:A,d:B,e:\{\}),b:C)$
\[\xymatrix@R=5ex@C=5ex@M=0.2ex{
&& \bullet \ar@{-}[dl]_a \ar@{-}[dr]^b \\ 
& \bullet \ar@{-}[dl]_c \ar@{-}[d]^d \ar@{-}[dr]^e && C \\
A & B & \circ}
\]
\end{enumerate}
where $A$, $B$, $C$, $D$, and $E$ are protocols. 
\end{example}

Given a protocol $X$, we may form its \textbf{dual} protocol $X^*$ by flipping
the product and coproduct structure and leaving the channel names and events
alone. For example, the duals of our protocols above are:

\begin{enumerate}
\item $\dual{X} = (a:(d:\dual{A}),b:\{e:\dual{B},f:\dual{C}\}, c:(g:\dual{D},
h:\dual{E}))$
\[\xymatrix@R=5ex@C=3ex@M=0.2ex{
&&& \bullet \ar@{-}[dlll]_a \ar@{-}[d]_b \ar@{-}[drrr]^c \\ \bullet
\ar@{-}[d]_d &&& \circ \ar@{-}[dl]_e \ar@{-}[dr]^f &&& \bullet
\ar@{-}[dl]_g \ar@{-}[dr]^h\\ 
A^* && B^* && C^* & D^* && E^*}
\]

\item $\dual{Y} = \{a:(c:\dual{A},d:\dual{B}), b:(e:\dual{A},f:\dual{C})\}$
\[\xymatrix@R=4ex@C=3ex@M=0.2ex{
&&& \circ \ar@{-}[dll]_a \ar@{-}[drr]^b \\ & \bullet \ar@{-}[dl]_c
\ar@{-}[dr]^d &&&& \bullet \ar@{-}[dl]_e \ar@{-}[dr]^f \\
A^* && B^* && A^* && C^*}
\]

\item $\dual{Z} = \{a:\{c:\dual{A},d:\dual{B},e:()\},b:\dual{C}\}$
\[\xymatrix@R=5ex@C=5ex@M=0.2ex{
&& \circ \ar@{-}[dl]_a \ar@{-}[dr]^b \\ 
& \circ \ar@{-}[dl]_c \ar@{-}[d]^d \ar@{-}[dr]^e && C^* \\
A^* & B^* & \bullet}
\]
\end{enumerate}

Notice that taking the dual of some protocol implies that the atoms in that
protocol must have a dual in the underlying polycategory. In the ``initial''
$\Sigma\Pi$-logic (the logic with no atoms) taking the dual of a dual gets
us back to the original formula, i.e., $\dual{(\dual{X})} = X$. This is true
in the ``pure'' and ``free'' logic if and only if it is true for the atoms.

What we now want to do is to assign protocols to channels. Channels, as
mentioned earlier, will be denoted with Greek letters. The notation $\alpha:X$
will denote the assignment of the protocol $X$ to channel $\alpha$.

\subsection{Protocol transitions} 

A protocol may make a transition to one of its subprotocols (subformulas) by
following one of the edges from the ``root'' of the protocol. For example,
given the protocol
\[
\xymatrix@R=5ex@C=3ex@M=0.2ex{
&&& \circ \ar@{-}[dlll]_a \ar@{-}[d]_b \ar@{-}[drrr]^c \\ \circ \ar@{-}[d]_d
&&& \bullet \ar@{-}[dl]_e \ar@{-}[dr]^f &&& \circ \ar@{-}[dl]_g \ar@{-}[dr]^h\\
A && B && C & D && E}
\]
the transitions $a$, $b$, and $c$ (respectively) lead to the following three
subprotocols:
\[
\xymatrix@M=0.2ex{\circ \ar@{-}[d]_d \\ A} \qqquad
\xymatrix@M=0.2ex{& \bullet \ar@{-}[dl]_e \ar@{-}[dr]^f \\ B && C} \qqquad
\xymatrix@M=0.2ex{& \circ \ar@{-}[dl]_g \ar@{-}[dr]^h \\ D && E}
\]

In our system, in order for a protocol to make a transition it needs be
triggered by an event, of which there are two kinds: input events and output
events. Given a channel $\alpha$ and an event $a$ we use $\alpha[a]$ to
denote that the event $a$ has been input on channel $\alpha$ and
$\ol{\alpha}[a]$ to denote that the event $a$ has been output on channel
$\alpha$ (following, for example, the convention of the $\pi$-calculus).
Output events are thought of as being generated by the system and input events
as being generated by the environment. There is an obvious analogy here with
the game theoretic view of player and opponent moves. 

A process has two types of channels: domain channels and codomain channels.
A protocol will have a different (actually dual) behaviour depending
on which type of channel we associate it with. When a process is using a 
protocol $X$ on an domain channel $\alpha$, it will use the
\textbf{domain protocol} of $X$, denoted $\dom(X)$. When a process is using
$X$ on a codomain channel $\beta$, it will use the \textbf{codomain protocol}
of $X$, denoted $\cod(X)$.

The protocols $\dom(X)$ and $\cod(X)$ have as states all the subprotocols of
$X$. We are then able to associate with each protocol $X$ and channel $\alpha$,
two finite labeled trees (or labeled transition systems): one for $\dom(X)$
and one for $\cod(X)$. In a domain transition system a coproduct protocol may
only respond to input events while a product protocol may only respond to
output events. Dually, in a codomain transition system a coproduct protocol
may only respond to output events while a product protocol may only respond
to input events. We summarize these expected transitions in Table~\ref{inout}.

\begin{table}
\[\begin{array}{|c|c|c|} \hline
X & \dom(X) \text{ expects:} & \cod(X) \text{ expects:} \\
\hline
\{a_i:X_i\}_i & \text{input} & \text{output} \\
\hline
(a_i:X_i)_i & \text{output} & \text{input} \\
\hline
\end{array}
\]
\caption{Expected transitions for a protocol $X$ on a channel $\alpha$} 
\label{inout}
\end{table}

\begin{example}[Protocol transitions]
\quad
\begin{enumerate}

\item The transitions of \[\dom(\alpha:X) = 
\alpha:\{a:\{d:A\}, b:(e:B,f:C), c:\{g:D,h:E\}\}
\]
may be displayed as
\[\xymatrix@R=5ex@C=3ex@M=0.2ex{
&&& \circ \ar@{-}[dlll]_{\alpha[a]} \ar@{-}[d]_{\alpha[b]}
\ar@{-}[drrr]^{\alpha[c]} \\ \circ \ar@{-}[d]_{\alpha[d]} &&& \bullet
\ar@{-}[dl]_{\ol{\alpha}[e]} \ar@{-}[dr]^{\ol{\alpha}[f]} &&& \circ
\ar@{-}[dl]_{\alpha[g]} \ar@{-}[dr]^{\alpha[h]} \\ 
A && B && C & D && E}
\]

This protocol says that on $\alpha$ in the domain the system must wait
for one of the events $a$, $b$, or $c$ to be input from the environment. If
$a$ is received the system must to wait for the event $d$ to be input, after
which it will continue on with the protocol $A$. Alternatively, if a $b$ is
received the system enters an output state at which point it may output either
$e$ or $f$, and then continue with protocol $B$ or $C$ respectively. Lastly,
if $c$ is received it must then wait for more input, either $g$ or $h$ after
which it will continue with $D$ or $E$ respectively.

\item The transitions of 
\[\cod(\alpha:X) = 
\alpha:\{a:\{d:A\}, b:(e:B,f:C), c:\{g:D,h:E\}\}
\]
may be displayed as
\[\xymatrix@R=5ex@C=3ex@M=0.2ex{
&&&\circ \ar@{-}[dlll]_{\ol{\alpha}[a]} \ar@{-}[d]^{\ol{\alpha}[b]}
\ar@{-}[drrr]^{\ol{\alpha}[c]}\\ 
\circ \ar@{-}[d]_{\ol{\alpha}[d]} &&& \bullet \ar@{-}[dl]_{\alpha[e]}
\ar@{-}[dr]^{\alpha[f]} &&& \circ \ar@{-}[dl]_{\ol{\alpha}[g]}
\ar@{-}[dr]^{\ol{\alpha}[h]} \\ 
A && B && C & D && E}
\]
Notice that we could turn input events into output events and output events
into input events to obtain the same transition system as $(1)$ above.

\item The transitions of 
\[\dom(\alpha:Y) = 
\alpha:(a:\{c:A,d:B\}, b:\{e:A,f:C\})
\]
may be displayed as
\[\xymatrix@R=4ex@C=3ex@M=0.2ex{ &&& \bullet \ar@{-}[dll]_{\ol{\alpha}[a]} 
\ar@{-}[drr]^{\ol{\alpha}[b]} \\ & \circ \ar@{-}[dl]_{\alpha[c]} 
\ar@{-}[dr]^{\alpha[d]} &&&& \circ \ar@{-}[dl]_{\alpha[e]}
\ar@{-}[dr]^{\alpha[f]} \\
A && B && A && C}
\]

\item The transitions of 
\[\dom(\alpha:Z) = 
\alpha:(a:(c:A,d:B,e:\{\}),b:C)
\]
may be displayed as
\[\xymatrix@R=5ex@C=5ex@M=0.2ex{
&& \bullet \ar@{-}[dl]_{\ol{\alpha}[a]} \ar@{-}[dr]^{\ol{\alpha}[b]} \\ 
& \circ \ar@{-}[dl]_{\alpha[c]} \ar@{-}[d]|{\alpha[d]}
\ar@{-}[dr]^{\ol{\alpha}[e]}
&& C \\ A & B & \circ}
\]
\end{enumerate}
\end{example}

\section{Proofs as processes} \label{sec-pap}

A \textbf{process} in our system is a representation of a
$\Sigma\Pi$-derivation. In the previous section we saw that on a domain
channel, a coproduct protocol allows the environment to input certain events
to a process, and a product protocol allows a process to output certain events
to the environment. Protocols on a codomain channel behave dually, i.e., on an
codomain channel, a coproduct protocol allows a process to output certain
events to the environment, and a product protocol allows the environment to
input certain events to a process. Additionally, there are atomic protocols
between which there can be an ``atomic process'' (i.e., an atomic map).

This means that the kind of events a process may perform depend on the ``type''
of the process. The \textbf{type} of a process $f$ is defined to be two
strings of channel names (the domain channels and the codomain channels),
and an assignment of protocols to those channels. This will be denoted
\[\alpha_1:X_1,\ldots,\alpha_n:X_n \to_f \beta_1:Y_1,\ldots,\beta_m:Y_m
\]

The channels to the left of the arrow are the domain channels, and the
channels to the right of the arrow are codomain channels. In the following,
$\Gamma$ and $\Delta$ will be used to represent comma separated strings of
channels and their assigned protocols, e.g., $\Gamma = \alpha_1:X_1,\ldots,
\alpha_n:X_n$. 

In this section we introduce two term representations for
$\Sigma\Pi$-derivations (processes). The first is a compact term calculus
which will allow for easy manipulations of the derivations (i.e., cut
reduction steps and permuting conversions). This representation will be
used for proving the Church-Rosser property and cut-elimination. The second
term representation uses a ``programming language like'' syntax which will
be more representative of the ``proofs as processes''.

\subsection{A term calculus representation} \label{ssec-tc}

The formation rules for this representation are given in
Table~\ref{compactreptable}. The notation $\alpha\{\}$ will be used for the
map from the empty sum on a domain channel $\alpha$ and $\beta()$ for the map
to the empty product on a codomain channel $\beta$.

\begin{table}
\begin{center}
\ovalbox{
\parbox{65ex}{\begin{center}
\medskip
{\small
$\infer[\text{(identity)}]{\alpha:A \vd_{1_A} \beta:A}{}$
\bigskip\\

$\infer{\Gamma, \alpha:\sum\limits_{i \in I} a_i:X_i \vd_{\alpha\{a_i
\mapsto f_i\}_{i \in I}} \Delta}
{\{\Gamma,\alpha:X_i \vd_{f_i} \Delta\}_{i \in I}}
\qquad
\infer{\Gamma \vd_{\alpha(a_i \mapsto f_i)_{i \in I}} \alpha:\prod\limits_{i
\in I} a_i:X_i,\Delta}
{\{\Gamma \vd_{f_i} \alpha:X_i, \Delta\}_{i \in I}}$
\bigskip\\

$\infer
  {\Gamma, \alpha:\prod\limits_{i \in I} a_i:X_i \vd_{\ola{\alpha}[a_k](f)}
  \Delta}{\Gamma,\alpha:X_k \vd_f \Delta}
\qquad
\infer
  {\Gamma \vd_{\ora{\alpha}[a_k](f)} \alpha:\sum\limits_{i \in I} a_i:X_i,
  \Gamma}{\Gamma \vd_f \alpha:X_k,\Gamma}$
\centro{where $k \in I, I \neq \emptyset$} \smallskip

$\infer{\Gamma,\Gamma' \vd_{f ;_\gamma g} \Delta,\Delta'}
  {\Gamma \vd_{f} \Delta,\gamma:Z & \gamma:Z,\Gamma' \vd_{g} \Delta'}$
}\medskip
\end{center}}}
\end{center}
\caption{$\Sigma\Pi$ term calculus formation rules} \label{compactreptable}
\end{table}

If the domain and the codomain (i.e., both sides of the turnstile) consist of
the same single atom $A$ (e.g., $A \vd A$), the identity rule asserts that the
process may output the atomic identity map $1_A:A \ra A$ (from the underlying
polycategory) and then end. The axiom rule asserts that if a process is
sitting between atoms $\Gamma$ in the domain and $\Delta$ in the codomain,
and $f:\Gamma \ra \Delta$ is an atomic map, then the process may output $f$
and then end.

Given a process that is interacting on a channel $\alpha$ in an input state
(the protocol assigned to $\alpha$ is in a state in which input is allowed:
a coproduct protocol in the domain or a product protocol in the codomain),
the cotuple and tuple rules assert that the process must be prepared to accept
any one of the possible input events that the protocol will allow on $\alpha$,
in this case any of the $a_i$'s, for $i \in I$. After receiving one of the
input events $a_i$ the process will then continue on with its subprocess $f_i$
(where, in the cotuple case, the type of $f_i$ is $\Gamma, \alpha:X_i \ra
\Delta$, and in the tuple case, $\Gamma \ra \alpha:X_i, \Delta$).

Given a process that is interacting on a channel $\alpha$ in an output state,
the projection and injection rules assert that the process may output one of
the specified events, in this case, one of the $a_i$'s, for $i \in I$, and then
continue on with its subprocess $f$ of type $\Gamma,\alpha:X_i \ra \Delta$
or $\Gamma \ra \alpha:X_i,\Delta$ respectively.

\begin{example}[Term calculus representation]  \label{exam-termcal} \quad
\begin{enumerate}
\item Suppose that the following are atomic maps
\[i:A,E \ra G,I \qqquad
j:B,E \ra G,J
\]\[
k:C,F \ra H,I \qqquad
l:D,F \ra H,K
\]
and consider the assignment of protocols to channels:

\[\begin{array}{r}
\alpha:\left(\vcenter{
\xymatrix@R=4ex@C=1ex@M=0.2ex{ &&& \circ \ar@{-}[dll]_{a} \ar@{-}[drr]^{b} \\
& \bullet \ar@{-}[dl]_{c} \ar@{-}[dr]^{d} &&&& \bullet \ar@{-}[dl]_{e}
\ar@{-}[dr]^{f}
\\ A && B && C && D}}\right)
\medskip \\
\beta:\left(\vcenter{
\xymatrix@R=4ex@C=1ex@M=0.2ex{& \bullet \ar@{-}[dl]_{g} \ar@{-}[dr]^{h} \\
E && F}}\right)
\end{array}
\to^{\qqqquad}
\begin{array}{l}
\gamma:\left(\vcenter{
\xymatrix@R=4ex@C=1ex@M=0.2ex{ & \circ \ar@{-}[dl]_{a'} \ar@{-}[dr]^{b'} \\
G && H}}\right)
\bigskip \\
\delta:\left(\vcenter{
\xymatrix@R=4ex@C=1ex@M=0.2ex{& \bullet \ar@{-}[dl]_{c'} \ar@{-}[dr]^{d'} \\
I && \circ \ar@{-}[dl]_{e'} \ar@{-}[dr]^{f'} \\
& J && K}}\right)
\end{array}
\]

A process between $\alpha$ and $\beta$ in the domain and $\gamma$ and $\delta$
in the codomain is:

\[\alpha\left\{
\begin{array}{l}
a \mapsto \ola{\beta}[g](\ora{\gamma}[a'](\delta\left(
    \begin{array}{l}
    c' \mapsto \ola{\alpha}[c](i) \\
    d' \mapsto \ola{\alpha}[d](\ora{\delta}[e'](j))
    \end{array}\right)))
\\
b \mapsto \ola{\beta}[h](\ora{\gamma}[b'](\delta\left(
    \begin{array}{l}
    c' \mapsto \ola{\alpha}[e](k) \\
    d' \mapsto \ola{\alpha}[f](\ora{\delta}[f'](l))
    \end{array}\right)))
\end{array}\right\}
\]

This process must first wait for the environment to input either an the event
$a$ or the event $b$ on $\alpha$; say it receives the event $a$. The process
will then proceed to output first $g$ on $\beta$, and then $a'$ on $\gamma$,
after which it will wait for either a $c'$ or a $d'$ to be input on
$\delta$. Suppose $c'$ is input, then it will output $c$ on $\alpha$ and 
finish by outputting the atomic map $i$.

This, of course, is just describing one possible interaction that could
happen. The environment could have first input $b$ on $\alpha$, in which
case a different interaction would have occurred.

This term represents the derivation of the sequent in the third example
of Example~\ref{exam-logic} annotated in the following way:
\begin{align*}
\alpha\mathbf{:}\{a:(c:A,d:B),b:(e:C,f:D)\},\ \beta:(g:E,h:F) \qqquad \\
\vd \gamma:\{a':G,b':H\},\ \delta:(c':I,d':\{e':J,f':K\})
\end{align*}

\item Consider the following annotated version Example~\ref{exam-logic}(2):
{\small
\[\infer{\alpha:\{a:(e:A,f:B),b:C\} \vd \beta:\{c:(g:A,h:B),d:C\}}{
   \infer{\alpha:(e:A,f:B) \vd \beta:\{c:(g:A,h:B),d:C\}}{
      \infer{\alpha:(e:A,f:B) \vd \beta:(g:A,h:B)}{
            \infer{\alpha:(e:A,f:B) \vd \beta:A}{
               \infer{\alpha:A \vd \beta:A}{}} &
            \infer{\alpha:(e:A,f:B) \vd \beta:B}{
               \infer{\alpha:B \vd \beta:B}{}}}} &
   \infer{\alpha:C \vd \beta:\{c:(A,B),d:C\}}{
         \infer{\alpha:C \vd \beta:C}{}}}
\]}
The term corresponding to this sequent derivation is

\[\alpha\left\{
\begin{array}{l}
a \mapsto \ora{\beta}[c](\beta\left(
    \begin{array}{l}
    g \mapsto \ola{\alpha}[e](1_A) \\
    h \mapsto \ola{\alpha}[f](1_B)
    \end{array}\right))) \\
b \mapsto \ora{\beta}[d](1_C)
\end{array}\right\}
\]

\end{enumerate}
\end{example}

\subsection{Cut-elimination rewrites} \label{ssec-cer}

If a domain channel of a process $f$, and a codomain channel of a process $g$,
share a common protocol then these processes may be ``plugged'' (or composed)
together on these channels to form a new process. For example, the two
processes
\[
\alpha:X \to_f \gamma:Z \qquad \text{and} \qquad \gamma:Z \to_g \beta:Y
\]
share the protocol $Z$ in common where $Z$ is in the codomain of $f$
and in the domain of $g$, and thus, $f$ and $g$ may be ``plugged''
together on $\gamma$ to form a new process
\[
\alpha:X \to_{f ;_\gamma g} \beta:Y
\]

In this new process, any output event on $\gamma$ from $f$ becomes an input
event on $\gamma$ for $g$, and vise versa. This means that the environment
may no longer input to either $\gamma$ as it is hidden (in the same sense as
in the $\pi$-calculus).

In general, processes have many domain and codomain channels, e.g.,
\[
\Gamma \to_f \Delta,\gamma:Z \qquad \text{and} \qquad
\gamma:Z,\Gamma' \to_g \Delta'
\]
Plugging these processes together on $\gamma$ results in
\[
\Gamma,\Gamma' \to_{f ;_\gamma g} \Delta,\Delta'
\]
where $\Gamma$ and $\Gamma'$ are the domain channels and $\Delta$ and $\Delta'$
are the codomain channels for this new process.

Notice that plugging processes together can only occur when they have exactly
one channel name in common: after this plugging process all the channel names
must be distinct. This means that channels may have to be renamed in order to
compose. We shall use (simultaneous) channel name substitution to indicate this
renaming process
\[f[\alpha_1'/\alpha_1,\ldots,\alpha_n'/\alpha_n]
\]
where $\alpha_1,\ldots,\alpha_n$ must be distinct channel names, as must
$\alpha_1',\ldots,\alpha_n'$, and
\[(\mathrm{channels}(f) \bs \{\alpha_1,\ldots,\alpha_n\}) \cap
\{\alpha_1',\ldots,\alpha_n'\} = \emptyset
\]
so that the new channel names are distinct from the old channel names which
are not replaced.

\begin{example}[Renaming channels] \quad
\begin{enumerate}
\item In order to compose
\[
\alpha:W,\beta:X \vd_f \gamma:Y,\delta:Z \qquad \text{and} \qquad
\delta:Z,\alpha:W' \vd_g \beta:X', \gamma:Y'
\]
on $\delta$ we must first rename $\alpha,\ \beta$, and $\gamma$ in
$f$ or $g$. Suppose we rename in $g$:
\[
g[\epsilon/\alpha,\ \eta/\beta,\  \theta/\gamma] =
\delta:Z,\ \epsilon:W' \vd_g \eta:X',\ \theta:Y'
\]
The result of composing $f$ and the renamed $g$ on $\delta$ is then
\[
\alpha:W,\ \beta:X,\ \epsilon:W'\vd_{f ;_\delta g}
    \gamma:Y,\ \eta:X',\ \theta:Y'
\]

\item Suppose we wish to compose the following two processes
\[
\alpha:X \vd_f \beta:Z \qquad \text{and} \qquad \gamma:Z \vd_g \delta:Y
\]
using $\beta:Z$ from $f$ and $\gamma:Z$ from $g$. In this case, we may
substitute either $\beta$ for $\gamma$, $\gamma$ for $\beta$, or a fresh
channel name (a channel name which does not occur in either $f$ or $g$)
for both $\beta$ and $\gamma$. Suppose we substitute $\gamma$ for $\beta$
in $f$ which yields $\alpha:X \vd_f \gamma:Z$. Now $f$ and $g$ may be composed
on $\gamma$:
\[\alpha:X \vd_{f ;_\gamma g} \delta:Y
\]
\end{enumerate}

In general we may end up renaming the channels on which we want to compose
to have the same ``name'' and the rest of the channels to be distinct from
one another.
\end{example}

To simplify the exposition, in the following we will simply assume that
the channel names of processes are distinct unless otherwise specified.
For example, given the following two processes

\[\Gamma \to_f \Delta,\gamma:X \qquad \text{and} \qquad
\gamma:X,\Gamma' \to_g \Delta'
\]
it will be assumed that the only channel name $f$ and $g$ have in common is
$\gamma$.

It should now be obvious that plugging together two processes is exactly
an application of the cut rule. The dynamics of cut-elimination then can be
seen as the way in which two processes communicate.

The cut-elimination reductions and the permuting conversions are summarized
in Table~\ref{compactrewrites}. (Typing information has been omitted as it
can be inferred from the terms, and in any case these have been displayed as
sequent derivations in the previous chapter.) Note that apart from (19) and
(22), these come in dual pairs, so there are six rewrites, eight conversions
and their duals, and two other conversions (which are self-dual): essentially
16 rewrites.

\begin{table}
\begin{center}
\ovalbox{
\parbox{65ex}{\begin{center}
\[\begin{array}{rrcl}
(1) & f ;_\gamma 1 &\Lra& f \medskip\\
(2) & 1 ;_\gamma f &\Lra& f \medskip\\
(3) & \alpha\{a_i \mapsto f_i\}_i ;_\gamma g &\Lra&
\alpha\{a_i \mapsto f_i ;_\gamma g\}_i \medskip\\
(4) & f ;_\gamma \beta(b_i \mapsto g_i)_i &\Lra&
\beta(b_i \mapsto f ;_\gamma g_i)_i \medskip\\
(5) & \ora{\alpha}[a_k](f) ;_\gamma g &\Lra&
\ora{\alpha}[a_k](f ;_\gamma g) \medskip\\
(6) & f ;_\gamma \ola{\beta}[b_k](g) &\Lra&
\ola{\beta}[b_k](f ;_\gamma g) \medskip\\
(7) & \ola{\alpha}[a_k](f) ;_\gamma g &\Lra&
\ola{\alpha}[a_k](f ;_\gamma g) \medskip\\
(8) & f ;_\gamma \ora{\beta}[b_k](g) &\Lra&
\ora{\beta}[b_k](f ;_\gamma g) \medskip\\
(9) & \alpha(a_i \mapsto f_i)_i ;_\gamma g &\Lra&
\alpha(a_i \mapsto f_i ;_\gamma g)_i \medskip\\
(10) & f ;_\gamma \beta\{b_i \mapsto g_i\}_i &\Lra&
\beta\{b_i \mapsto f ;_\gamma g_i\}_i \medskip\\
(11) & \ora{\gamma}[a_k](f) ;_\gamma \gamma\{a_i \mapsto g_i\}_i &\Lra&
f ;_\gamma g_k \medskip\\
(12) & \gamma(a_i \mapsto f_i)_i ;_\gamma \ola{\gamma}[a_k](g)  &\Lra&
f_k ;_\gamma g \medskip\\
%
(13) & \alpha\{a_i \mapsto \beta\{b_j \mapsto f_{ij}\}_j\}_i &\pc&
\beta\{b_j \mapsto \alpha\{a_i \mapsto f_{ij}\}_i\}_j \medskip\\
(14) & \alpha(a_i \mapsto \beta(b_j \mapsto f_{ij})_j)_i &\pc&
\beta(b_j \mapsto \alpha(a_i \mapsto f_{ij})_i)_j \medskip\\
(15) & \alpha\{a_i \mapsto \ora{\beta}[b_k](f_i)\}_i &\pc&
\ora{\beta}[b_k](\alpha \{a_i \mapsto f_{i}\}_i) \medskip\\
(16) & \ola{\alpha}[a_k](\beta(b_i \mapsto g_i)_i) &\pc&
\beta(b_i \mapsto \ola{\alpha}[a_k](g_i))_i \medskip\\
(17) & \alpha\{a_i \mapsto \ola{\beta}[b_k](f_i)\}_i &\pc&
\ola{\beta}[b_k](\alpha\{a_i \mapsto f_{i}\}_i) \medskip\\
(18) &\ora{\alpha}[a_k](\beta(b_i \mapsto g_i)_i) &\pc&
\beta(b_i \mapsto \ora{\alpha}[a_k](g_i))_i \medskip\\
(19) & \alpha\{a_i \mapsto \beta(b_j \mapsto f_{ij})_j\}_i &\pc&
\beta(b_j \mapsto \alpha\{a_i \mapsto f_{ij}\}_i)_j \medskip\\
(20) & \ora{\alpha}[a_k](\ora{\beta}[b_{l}](f)) &\pc&
\ora{\beta}[b_{l}](\ora{\alpha}[a_k](f)) \medskip\\
(21) & \ola{\alpha}[a_k](\ola{\beta}[b_{l}](f)) &\pc&
\ola{\beta}[b_{l}](\ola{\alpha}[a_k](f)) \medskip\\
(22) &\ora{\alpha}[a_k](\ola{\beta}[b_{l}](f)) &\pc&
\ola{\beta}[b_{l}](\ora{\alpha}[a_k](f))
\end{array}
\]
\end{center}}}
\end{center}
\caption{$\Sigma\Pi$ conversion rules} \label{compactrewrites}
\end{table}

Recall that we allow the index sets $I$ and $J$ to be empty, except where the
projection and injection rules are involved (the rules (5), (7), (11), (15),
(17), (20), (22) and their duals); in these cases, since reference is made
to an element $k$ or $l$, it does not make sense for the corresponding index
set to be empty. In (15) and (17) (and their duals) the index set $J$ for the
named element $k$ must not be empty, but the other index set $I$ may be. In
(13), (14), and (19) either (or both or neither) index set may be empty. An
explicit treatment of these nullary cases may be found in Section~\ref{sec-au}.

\begin{example}[Process communication (composition)] \label{exam-comp}
Suppose that
\begin{itemize}
\item $X = \alpha:\{a:(c:A,d:B),b:(e:A,f:C)\}$
\item $Y=\beta:\{a:B,b:C\}$
\item $Z= \gamma:(a:A,b:\{c:B,d:C\})$
\end{itemize}
are protocols and consider the following situation:
{\scriptsize
\[\alpha:X \to_{
\alpha\left\{
\begin{array}{l}
a \mapsto \gamma\left(
    \begin{array}{l}
    a \mapsto \ola{\alpha}[c](1_A) \\
    b \mapsto \ora{\gamma}[c](\ola{\alpha}[d](1_B))
    \end{array}\right) \\
b \mapsto \gamma\left(
    \begin{array}{l}
    a \mapsto \ola{\alpha}[e](1_B) \\
    b \mapsto \ora{\gamma}[d](\ola{\alpha}[f](1_C))
    \end{array}\right)
\end{array}\right\}}
\gamma:Z
\quad
\gamma:Z \to_{
\ola{\gamma}[b](\gamma\left(
    \begin{array}{l}
    c \mapsto \ora{\beta}[a](1_B) \\
    d \mapsto \ora{\beta}[b](1_C)
    \end{array}\right)}
\beta:Y
\]}
Composing these two processes on $\gamma:Z$ gives

{\scriptsize
\[
\alpha\left\{
\begin{array}{l}
a \mapsto \gamma\left(
    \begin{array}{l}
    a \mapsto \ola{\alpha}[c](1_A) \\
    b \mapsto \ora{\gamma}[c](\ola{\alpha}[d](1_B))
    \end{array}\right) \\
b \mapsto \gamma\left(
    \begin{array}{l}
    a \mapsto \ola{\alpha}[e](1_B) \\
    b \mapsto \ora{\gamma}[d](\ola{\alpha}[f](1_C))
    \end{array}\right)
\end{array}\right\}
;_\gamma
\ola{\gamma}[b](\gamma\left(
    \begin{array}{l}
    c \mapsto \ora{\beta}[a](1_B) \\
    d \mapsto \ora{\beta}[b](1_C)
    \end{array}\right)
\] \medskip
\[
\Lra \alpha\left\{
\begin{array}{l}
a \mapsto \gamma\left(
    \begin{array}{l}
    a \mapsto \ola{\alpha}[c](1_A) \\
    b \mapsto \ora{\gamma}[c](\ola{\alpha}[d](1_B))
    \end{array}\right)
;_\gamma
\ola{\gamma}[b](\gamma\left(
    \begin{array}{l}
    c \mapsto \ora{\beta}[a](1_B) \\
    d \mapsto \ora{\beta}[b](1_C)
    \end{array}\right) \\
b \mapsto \gamma\left(
    \begin{array}{l}
    a \mapsto \ola{\alpha}[e](1_B) \\
    b \mapsto \ora{\gamma}[d](\ola{\alpha}[f](1_C))
    \end{array}\right)
;_\gamma
\ola{\gamma}[b](\gamma\left(
    \begin{array}{l}
    c \mapsto \ora{\beta}[a](1_B) \\
    d \mapsto \ora{\beta}[b](1_C)
    \end{array}\right)
\end{array}\right\} \quad \text{(by (3))}
\] \medskip
\[
\Lra \alpha\left\{
\begin{array}{l}
a \mapsto \ora{\gamma}[c](\ola{\alpha}[d](1_B))
;_\gamma
\gamma\left(
    \begin{array}{l}
    c \mapsto \ora{\beta}[a](1_B) \\
    d \mapsto \ora{\beta}[b](1_C)
    \end{array}\right) \\
b \mapsto \ora{\gamma}[d](\ola{\alpha}[f](1_C))
;_\gamma
\gamma\left(
    \begin{array}{l}
    c \mapsto \ora{\beta}[a](1_B) \\
    d \mapsto \ora{\beta}[b](1_C)
    \end{array}\right)
\end{array}\right\} \quad \text{(by (12))}
\] \medskip
\[
\Lra \alpha\left\{
\begin{array}{l}
a \mapsto \ola{\alpha}[d](1_B)
;_\gamma
\gamma(\ora{\beta}[a](1_B)) \\
b \mapsto \ola{\alpha}[f](1_C)
;_\gamma
\gamma(\ora{\beta}[b](1_C))
\end{array}\right\} \quad \text{(by (11))}
\] \medskip
\[
\Lra \alpha\left\{
\begin{array}{l}
a \mapsto \ola{\alpha}[d](1_B ;_\gamma \gamma(\ora{\beta}[a](1_B))) \\
b \mapsto \ola{\alpha}[f](1_C ;_\gamma \gamma(\ora{\beta}[b](1_C)))
\end{array}\right\} \quad \text{(by (7))}
\] \medskip
\[
\Lra \alpha\left\{
\begin{array}{l}
a \mapsto \ola{\alpha}[d](\ora{\beta}[a](1_B)) \\
b \mapsto \ola{\alpha}[f](\ora{\beta}[b](1_C))
\end{array}\right\} \quad \text{(by (2))}
\]}

It is easy to see the resulting process is indeed a process between
$\alpha:X$ and $\beta:Y$.
\end{example}

\subsection{A ``programming language'' representation} \label{ssec-plr}

In this section we introduce the second of our term representations for
derivations in $\Sigma\Pi$: the ``programming language'' representation. This
representation was suggested to Robin Cockett and the author by Robert Seely. 

The term formation rules for this representation are given in
Table~\ref{plreptable}. Notice that the cut rule in this representation does
not require that the channel names be the same, however, channels may still
have to be renamed (substituted) after cutting two processes together to keep
the channel names distinct. It is also worth noting that this representation
does not differentiate between input/output on a domain channel and
input/output on a codomain channel; the typing disambiguates between the two.

\begin{table}
\begin{center}
\ovalbox{\parbox[c]{81ex}{\begin{center}
\medskip
{\small
$\infer{\alpha:A \vd_{1_A} \beta:A}{}$
\bigskip\\

$\infer{\Gamma, \alpha:\sum\limits_{i \in I} a_i:X_i
\vd_{\mathtt{input\ on}\ \alpha\ \mathtt{of}\ \underset{i \in I}{|} a_i
\mapsto f_i} \Delta}
{\{\Gamma,\alpha:X_i \vd_{f_i} \Delta\}_{i \in I}}
\qquad
\infer{\Gamma \vd_{\mathtt{input\ on}\ \alpha\ \mathtt{of}\ \underset{i \in I}{|}
a_i \mapsto f_i} \alpha:\prod\limits_{i \in I} a_i:X_i,\Delta}
{\{\Gamma \vd_{f_i} \alpha:X_i, \Delta\}_{i \in I}}$
\bigskip\\

$\infer{\Gamma, \alpha:\prod\limits_{i \in I} a_i:X_i
\vd_{\mathtt{output}\ a_k\ \mathtt{on}\ \alpha\ \mathtt{then}\ f} \Delta}
{\Gamma,\alpha:X_k \vd_f \Delta}
\qquad
\infer{\Gamma \vd_{\mathtt{output}\ a_k\ \mathtt{on}\ \alpha\ \mathtt{then}\ f}
\alpha:\sum\limits_{i \in I} a_i:X_i,\Gamma}{\Gamma \vd_f \alpha:X_k,\Gamma}$
\centro{where $k \in I, I \neq \emptyset$} \smallskip

$\infer{\Gamma,\Gamma' \vd_{\mathtt{plug}\ \alpha\
\mathtt{in}\ f\ \mathtt{to}\ \beta\ \mathtt{in}\ g} \Delta,\Delta}
  {\Gamma \vd_{f} \Delta,\alpha:Z & \beta:Z,\Gamma' \vd_{g} \Delta'}$
}\medskip
\end{center}}}
\end{center}
\caption{$\Sigma\Pi$ programming language term formation rules}
\label{plreptable}
\end{table}

\begin{example}[Programming language representation] \quad
\begin{enumerate}

\item The following is the programming language representation of the first
example from Example~\ref{exam-termcal}.
\texttt{\begin{tabbing}
inp\=ut on $\alpha$ of \\
\> $\mid a$ \= $\mapsto$ output $g$ on $\beta$ then \\
\>\> out\=put $a'$ on $\gamma$ then  \\
\>\>\> inp\=ut on $\delta$ of \\
\>\>\>\> $\mid c' \mapsto$ out\=put $c$ on $\alpha$ then $i$ \\
\>\>\>\> $\mid d' \mapsto$ output $d$ on $\alpha$ then $j$ \\
\>\>\>\>\>                 output $e'$ on $\delta$ then $j$ \\
\> $\mid b$ \= $\mapsto$ output $h$ on $\beta$ then \\
\>\> output $b'$ on $\gamma$ then  \\
\>\>\> input on $\delta$ of \\
\>\>\>\> $\mid c' \mapsto$ output $e$ on $\alpha$ then $k$ \\
\>\>\>\> $\mid d' \mapsto$ output $f$ on $\alpha$ then $l$ \\
\>\>\>\>\>                 output $f'$ on $\delta$ then $l$
\end{tabbing}}

\item This is the second example from Example~\ref{exam-termcal}.
\texttt{\begin{tabbing}
inp\=ut on $\alpha$ of \\
\> $\mid a$ \= $\mapsto$ output $c$ on $\beta$ then \\
\>\> inp\=ut on $\beta$ of \\
\>\>\> $\mid g \mapsto$ output $e$ on $\alpha$ then $1_A$ \\
\>\>\> $\mid h \mapsto$ output $f$ on $\alpha$ then $1_B$ \\
\> $\mid b$ \= $\mapsto$ output $d$ on $\beta$ then $1_C$
\end{tabbing}}

\end{enumerate}
\end{example}

As this representation is self-dual, many of the cut-elimination reductions and
permuting conversions are identical. Therefore, only the unique rewrites are
presented; the rewrites which are identical will be indicated in brackets,
$[\_]$ (the dual rewrites are left out).

\[
\begin{array}{rl}
(1) & \mathtt{plug}\ \gamma\ \mathtt{in}\ f\ \mathtt{to}\ \delta\
\mathtt{in}\ 1 \Lra f \medskip\\

(3),[(9)]\ &
 \mathtt{plug}\ \gamma\ \mathtt{in}\ (\texttt{input on}\ \alpha\ \mathtt{of}\
\mid_i a_i \mapsto f_i)\ \mathtt{to}\ \delta\ \mathtt{in}\ g
\smallskip\\ &\Lra 
\texttt{input on}\ \alpha\ \mathtt{of}\ \mid_i a_i \mapsto
(\mathtt{plug}\ \gamma\ \mathtt{in}\ f_i\ \mathtt{to}\ \delta\ \mathtt{in}\ g)
\medskip\\

(5), [(7)]\ &
\mathtt{plug}\ \gamma\ \mathtt{in}\ (\texttt{output}\ a_k\ \mathtt{on}\
\alpha\ \mathtt{then}\ f)\ \mathtt{to}\ \delta\ \mathtt{in}\ g
\smallskip\\ &\Lra
\texttt{output}\ a_k\ \mathtt{on}\ \alpha\ \mathtt{then}\ (\mathtt{plug}\
\gamma\ \mathtt{in}\  f\ \mathtt{to}\ \delta\ \mathtt{in}\ g)
\medskip\\

(11)\ & \mathtt{plug}\ \gamma\ \mathtt{in}\
(\mathtt{output}\ a_k\ \mathtt{on}\ \gamma\ \mathtt{then}\ f)\
\smallskip\\ & \mathtt{to}\ \gamma\ \mathtt{in}\
(\texttt{input on}\ \gamma\ \mathtt{of}\ \mid_i a_i \mapsto g_i) 
\smallskip\\ & \Lra
\mathtt{plug}\ \gamma\ \mathtt{in}\ f\ \mathtt{to}\ \gamma\ \mathtt{in}\ g_k
\medskip\\

(13),[(19)]\ &
\texttt{input on}\ \alpha\ \mathtt{of}\ \mid_i a_i \mapsto
(\texttt{input on}\ \beta\ \mathtt{of}\ \mid_j b_j \mapsto f_{ij})
\smallskip\\ &\pc
\texttt{input on}\ \beta\ \mathtt{of}\ \mid_j b_j \mapsto
(\texttt{input on}\ \alpha\ \mathtt{of}\ \mid_i a_i \mapsto f_{ij})
\medskip\\

(15),[(17)] &
\texttt{input on}\ \alpha\ \mathtt{of}\ \mid_i a_i \mapsto
(\texttt{output}\ b_k\ \mathtt{on}\ \beta\ \mathtt{then}\ f_i)
\smallskip\\ &\pc
\texttt{output}\ b_k\ \mathtt{on}\ \beta\ \mathtt{then}\
(\texttt{input on}\ \alpha\ \mathtt{of}\ \mid_i a_i \mapsto f_i)
\medskip\\

(20),[(22)] &
\texttt{output}\ a_k\ \texttt{on}\ \alpha\ \mathtt{then}\
(\texttt{output}\ b_k\ \mathtt{on}\ \beta\ \mathtt{then}\ f_i)
\smallskip\\ &\pc
\texttt{output}\ b_k\ \texttt{on}\ \beta\ \mathtt{then}\
(\texttt{output}\ a_k\ \mathtt{on}\ \alpha\ \mathtt{then}\ f_i)
\end{array}
\]

The following example of cut-elimination is the ``programming language
representation'' of Example~\ref{exam-comp}.
\begin{itemize}
\item We begin with the following process.
\texttt{\begin{tabbing}
plug $\gamma$ \= in \\
\>inp\=ut on $\alpha$ of \\
\>\> $\mid a$ \= $\mapsto$ input on $\gamma$ of \\
\>\>\> $\mid a \mapsto$ \= output $c$ on $\alpha$ then $1_A$ \\
\>\>\> $\mid b \mapsto$ output $c$ on $\gamma$ then \\
\>\>\>\>                output $d$ on $\alpha$ then $1_B$ \\
\>\> $\mid b$ \= input on $\gamma$ of \\
\>\>\> $\mid a \mapsto$ output $e$ on $\alpha$ then $1_B$ \\
\>\>\> $\mid b \mapsto$ output $d$ on $\gamma$ then \\
\>\>\>\>                output $f$ on $\alpha$ then $1_C$ \\
to\ $\gamma$ in \\
\> output $b$ on $\gamma$ then \\
\>\> input on $\gamma$ of \\
\>\>\> $\mid c \mapsto$ output $a$ on $\beta$ then $1_B$ \\
\>\>\> $\mid d \mapsto$ output $b$ on $\beta$ then $1_C$
\end{tabbing}}

\item Applying rewrite (3) yields.
\texttt{\begin{tabbing}
inp\=ut on $\alpha$ of \\
\> $\mid a \mapsto$ \= plu\=g $\gamma$ in \\
\>\>\> inp\=ut on $\gamma$ of \\
\>\>\>\> $\mid a$ \= $\mapsto$ \= output $c$ on $\alpha$ then $1_A$ \\
\>\>\>\> $\mid b$ $\mapsto$ output $c$ on $\gamma$ then \\
\>\>\>\>\>\>                output $d$ on $\alpha$ then $1_B$ \\
\>\> to $\gamma$ in \\
\>\>\> output $b$ on $\gamma$ then \\
\>\>\>\> input on $\gamma$ of \\
\>\>\>\>\> $\mid c \mapsto$ output $a$ on $\beta$ then $1_B$ \\
\>\>\>\>\> $\mid d \mapsto$ output $b$ on $\beta$ then $1_C$ \\
\> $\mid b \mapsto$  plug $\gamma$ in \\
\>\>\> input on $\gamma$ of \\
\>\>\>\> $\mid a$ $\mapsto$ output $e$ on $\alpha$ then $1_B$ \\
\>\>\>\> $\mid b$ $\mapsto$ output $d$ on $\gamma$ then \\
\>\>\>\>\>                output $f$ on $\alpha$ then $1_C$ \\
\>\> to $\gamma$ in \\
\>\>\> output $b$ on $\gamma$ then \\
\>\>\>\> input on $\gamma$ of \\
\>\>\>\>\> $\mid c \mapsto$ output $a$ on $\beta$ then $1_B$ \\
\>\>\>\>\> $\mid d \mapsto$ output $b$ on $\beta$ then $1_C$
\end{tabbing}}

\item To this result, applying rewrite (12) (the dual of (11)) yields.
\texttt{\begin{tabbing}
inp\=ut on $\alpha$ of \\
\> $\mid a \mapsto$ \= plu\=g $\gamma$ in \\
\>\>\> out\=put $c$ on $\gamma$ then \\
\>\>\>\> output $d$ on $\alpha$ then $1_B$ \\
\>\> to $\gamma$ in \\
\>\>\> input on $\gamma$ of \\
\>\>\>\> $\mid c \mapsto$ output $a$ on $\beta$ then $1_B$ \\
\>\>\>\> $\mid d \mapsto$ output $b$ on $\beta$ then $1_C$ \\
\> $\mid b \mapsto$  plug $\gamma$ in \\
\>\>\> output $d$ on $\gamma$ then \\
\>\>\>\> output $f$ on $\alpha$ then $1_C$ \\
\>\> to $\gamma$ in \\
\>\>\> input on $\gamma$ of \\
\>\>\>\> $\mid c \mapsto$ output $a$ on $\beta$ then $1_B$ \\
\>\>\>\> $\mid d \mapsto$ output $b$ on $\beta$ then $1_C$
\end{tabbing}}

\item Rewrite (11).
\texttt{\begin{tabbing}
inp\=ut on $\alpha$ of \\
\> $\mid a \mapsto$ \= plu\=g $\gamma$ in \\
\>\>\> output $d$ on $\alpha$ then $1_B$ \\
\>\> to $\gamma$ in \\
\>\>\> output $a$ on $\beta$ then $1_B$ \\
\> $\mid b \mapsto$  plug $\gamma$ in \\
\>\>\> output $f$ on $\alpha$ then $1_C$ \\
\>\> to $\gamma$ in \\
\>\>\> output $b$ on $\beta$ then $1_C$
\end{tabbing}}

\item Rewrite (7).
\texttt{\begin{tabbing}
inp\=ut on $\alpha$ of \\
\> $\mid a$ \= $\mapsto$ output $d$ on $\alpha$ then \\
\>\> plu\=g $\gamma$ in $1_B$\\
\>\> to $\gamma$ in \\
\>\>\> output $a$ on $\beta$ then $1_B$ \\
\> $\mid b$ $\mapsto$ output $f$ on $\alpha$ then \\
\>\> plug $\gamma$ in $1_C$ \\
\>\> to $\gamma$ in \\
\>\>\> output $b$ on $\beta$ then $1_C$
\end{tabbing}}

\item And finally, by rewrite (2).
\texttt{\begin{tabbing}
inp\=ut on $\alpha$ of \\
\> $\mid a$ \= $\mapsto$ output $d$ on $\alpha$ then \\
\>\> output $a$ on $\beta$ then $1_B$ \\
\> $\mid b$ $\mapsto$ output $f$ on $\alpha$ then \\
\>\> output $b$ on $\beta$ then $1_C$
\end{tabbing}}
\end{itemize}

Although this syntax gives a clear intuition into processes, the reader
will agree that it is quite verbose. Thus, in what follows, we shall favor
the compact term calculus introduced in Section~\ref{ssec-tc}.

\section{The additive units} \label{sec-au}

The way in which the cut-elimination process handles the reductions and
permuting conversions when the index set $I = \emptyset$ is quite subtle.
To clarify this, in this section we make these special cases explicit.
In the following the abbreviations $\sum_\emptyset = 0$ and $\prod_\emptyset
= 1$ will be used. The nullary versions of the cotuple and tuple rules are:
\[\infer[\text{(cotuple)}]{\Gamma,\alpha:0 \vd_{\alpha\{\}} \Delta}{} \qqquad
\infer[\text{(tuple)}]{\Gamma \vd_{\beta()} \beta:1,\Delta}{}
\]

The notation here is ambiguous as one cannot derive the context from the
terms. To correct this we shall write the terms above as $\alpha\{\}_{\Gamma
\vd \Delta}$ and $\beta()_{\Gamma \vd \Delta}$ respectively.

There are four reductions that are relevant to this setting, corresponding to
the rewrites (3), (4), (9), and (10) where $I = \emptyset$.
Given terms $f:\Gamma \vd \Delta,\gamma:X$ and $g:\gamma:X,\Gamma' \vd \Delta'$
we have the following reductions.

\[\begin{array}{rrcl}
(3) & \alpha\{\}_{\Gamma \vd \Delta,\gamma:X}\ ;_\gamma g & \Lra  &
\alpha\{\}_{\Gamma, \Gamma' \vd \Delta,\Delta'} \medskip\\
(4) & f\ ;_\gamma \beta()_{\gamma:X,\Gamma' \vd \Delta'} & \Lra  &
\beta()_{\Gamma, \Gamma' \vd \Delta,\Delta'} \medskip\\
(9) & \alpha()_{\Gamma \vd \Delta,\gamma:X}\ ;_\gamma g & \Lra  &
\alpha()_{\Gamma, \Gamma' \vd \Delta,\Delta'} \medskip\\
(10) & f\ ;_\gamma \beta\{\}_{\gamma:X,\Gamma' \vd \Delta'} & \Lra  &
\beta\{\}_{\Gamma, \Gamma' \vd \Delta,\Delta'}
\end{array}
\]

In addition, there are 13 permuting conversions, corresponding to the
cases (15) and (17) (and their duals), three variants of (13)
(corresponding to the cases when only $I=\emptyset$, only $J=\emptyset$,
and both $I = J = \emptyset$) (and their duals), and three variants of (19).
Fortunately, they are all very similar, so we present only (15) and the
three variants of (19). Here we drop the typing on the term and indicate
it in the brackets.

\[\begin{array}{rrcll}
(15) & \alpha\{\} &\pc& \ora{\beta}[b_k](\alpha\{\})
& (\Gamma,\alpha:0 \vd \beta:\sum Y_i,\Delta) \medskip\\
(19) & \alpha\{a_i \mapsto \beta()\}_i & \pc & \beta()
& (\Gamma,\alpha:\sum X_i \vd \beta:1,\Delta) \medskip\\
(19) & \alpha\{\} & \pc & \beta(b_j \mapsto \alpha\{\})_j
& (\Gamma,\alpha:0 \vd \beta:\prod Y_j,\Delta) \medskip\\
(19) & \alpha\{\} & \pc & \beta()
& (\Gamma,\alpha:0 \vd \beta:1,\Delta) \medskip\\
\end{array}
\]

\section{Proof of the Church-Rosser property} \label{sec-cr}

In this section we present a proof of the Church-Rosser property for
$\Sigma\Pi$-morphisms. The proof presented here follows very closely the
proof in~\cite{cockett01:finite}, extended to the ``poly'' case. 

We wish to show that given any two $\Sigma\Pi$-morphisms related by a series
of reductions and permuting conversions

\[\xymatrix{t_1 \ar@{<=}[r] & t_2 \ar@{|=|}[r] & t_3 \ar@{=>}[r] &
\quad \cdots \quad &
\ar@{=>}[l] t_{n-2} \ar@{|=|}[r] & t_{n-1} \ar@{=>}[r] & t_n}
\]
there is an alternative way of arranging the reductions and permuting
conversions so that $t_1$ and $t_n$ can be reduced to terms which are
related by the permuting conversions alone. That is, we wish to show that
there is a convergence of the following form:

\[\xymatrix{t_1 \ar@{=>}[dr]_{*} & && & t_n \ar@{=>}[dl]^{*} \\
& t_1' \ar@{|=|}[rr]_{*} && t_n'}
\]

When the rewriting system terminates (in the appropriate sense) this allows
the decision procedure for the equality of $\Sigma\Pi$-terms to be reduced to
the decision procedure for the permuting conversions (see
Section~\ref{sec-dp}). In order to test the equality of two terms,
one can rewrite both terms into a reduced form (one from which there are no
further reductions), and these will be equal if and only if the two reduced
forms are equivalent through the permuting conversions alone. In the current
situation the reduction process is the cut-elimination procedure. In this
section we will show that this is a terminating procedure.

We begin with a couple of definitions. The first is the multiset ordering of
Dershowitz and Manna~\cite{dershowitz79:proving}. Let $(S,\succ)$ be a
partially ordered set, and let $\script{M}(S)$ denote the multisets (or bags)
over $S$. For $M,N \in \script{M}(S)$, $M > N$ (``$>$'' is called
the \textbf{multiset} (or \textbf{bag}) \textbf{ordering}), if there are
multisets $X,Y \in \script{M}(S)$, where $\emptyset \neq X \subseteq M$, such
that

\[N = (M \bs X) \cup Y \quad \text{and} \quad (\forall y \in Y)
(\exists x \in X)\ x \succ y
\]
where $\cup$ here is the multiset union.

For example,
\[[3] > [2,2,1,1], \quad [4,3] > [4], \quad [3,2] > [3,1]
\]

Recall from~\cite{dershowitz79:proving} that if $(S,\succ)$ is a total order
(linear order) then $\script{M}(S)$ is a total order. To see this consider
$M,N \in \script{M}(S)$. To determine whether $M > N$ sort the elements of
both $M$ and $N$ and then compare the two sorted sequences lexicographically. 

Following~\cite{cockett01:finite} we say a rewrite system is \textbf{locally
confluent modulo equations} if any (one step) divergence of the following
form

\[\vcenter{\xymatrix{& t_0 \ar@{=>}[dl] \ar@{=>}[dr] \\ t_1 && t_2}} 
\qquad \text{or} \qquad
\vcenter{\xymatrix{& t_0 \ar@{|=|}[dr] \ar@{=>}[dl] \\ t_1 && t_2}}
\]
(where ``$\Lra$'' denotes a reduction and ``$\pc$'' an equation) has a
convergence, respectively, of the form

\[\vcenter{\xymatrix{t_1 \ar@{|=>}[dr]_{*} && t_2 \ar@{|=>}[dl]^{*} \\ & t'}}
\qquad \text{or} \qquad
\vcenter{\xymatrix{t_1 \ar@{|=>}[ddr]_{*} && t_2 \ar@{=>}[d]  \\
&& t_2' \ar@{|=>}[dl]^{*} \\ & t'}}
\]
where the new arrow ``$\xymatrix@1{\ar@{|=>}[r]&}$'' indicates either
an equality or a reduction in the indicated direction.

\begin{proposition} \label{prop-conf}
Suppose $(N,\mathcal{R},\mathcal{E})$ is a rewriting system with the
equations equ\-ip\-ped with a well-ordered measure on the rewrite arrows such
that the measure of the divergences is strictly greater than the measure of
the convergences then the system is confluent modulo equations if and only
if it is locally confluent modulo equations.
\end{proposition}

\begin{proof}
If the system is confluent modulo equations it is certainly locally
confluent modulo equations. Conversely suppose we have a chain of reductions,
equations, and expansions. We may associate with it the bag of measures of
the arrows of the sequence.

The idea will be to show that replacing any local divergence in this chain
by a local confluence will result in a new chain whose bag measure is
strictly smaller. However, this can be seen by inspection as we are removing
the arrows associated with the divergence and replacing them with the arrows
associated with the convergence. The measure on the arrows associated with
the divergence is strictly greater then that of the measure on the arrows
associated with the convergence.

Thus, each rewriting reduces the measure and, therefore, any sequence
of rewriting on such a chain must terminate. However, it can only terminate
when there are no local divergences to resolve. This then implies that the
end result must be a confluence modulo equations.
\end{proof}

\subsection{Resolving critical pairs locally}

The proof of the Church-Rosser property involves examining all the possible
critical pairs involving reductions or reductions and conversions, and
showing that they are all of the form shown above and that they may be
resolved in the way shown above. It then must be shown that there is some
measure on the arrows which decreases when replacing a divergence with a
convergences. This will then suffice to show that our system is locally
confluent modulo equations, so that by Proposition~\ref{prop-conf}, it is
confluent modulo equations. The rewrites (1)-(12) are the ``reductions''
and the permuting conversions (13)-(22) are the ``equations''.

The resolutions of the critical pairs will be presented as reduction
diagrams. We begin with the rewrites involving the empty cotuple and empty
tuple rules.

\begin{itemize}
\item The resolution of the critical pair (1)-(3) (dually (2)-(4)) is
indicated by the following reduction diagram.

\[\xymatrix{\alpha\{\,\} ;_\gamma 1 \ar@{=>}[r]^{(1)}_{(3)} & \alpha\{\,\}}
\]
The critical pair (1)-(9) is handled similarly.

\item There are three cases for the resolution of the critical pair (3)-(4):
only $I = \emptyset$, only $J = \emptyset$, both $I = J = \emptyset$. The
first two cases correspond respectively to the following reductions diagrams:

\[\vcenter{\xymatrix@M=1ex@C=8ex@R=9ex@!0{
& \alpha\{\,\} ;_\gamma \beta(b_j \mapsto g_j)_j \ar@{=>}[dl]_{(3)}
\ar@{=>}[dr]^{(4)} \\
\alpha\{\,\} \ar@{|=|}[dr]_-{(19)}
&& \beta(b_j \mapsto \alpha\{\,\} ;_\gamma g_j)_j \ar@{=>}[dl]^{\beta((3))} \\
& \beta(b_j \mapsto \alpha\{\,\})_j}}
\qquad
\vcenter{\xymatrix@M=1ex@C=8ex@R=9ex@!0{
& \alpha\{a_i \mapsto f_i\} ;_\gamma \beta(\,) \ar@{=>}[dl]_{(3)}
\ar@{=>}[dr]^{(4)} \\
\alpha\{a_i \mapsto f_i ;_\gamma \beta(\,)\} \ar@{=>}[dr]_-{\alpha((4))}
&& \beta(\,) \ar@{=>}[dl]^{(19)} \\
& \alpha\{a_i \mapsto \beta(\,)\}}}
\]
while the reduction diagram for the third case is:

\[\xymatrix@M=1ex@C=8ex@R=9ex@!0{
& \alpha\{\,\} ;_\gamma \beta(\,) \ar@{=>}[dl]_{(3)} \ar@{=>}[dr]^{(4)} \\
\alpha\{\,\} \ar@{|=|}[rr]_-{(19)} && \beta(\,)}
\]
The rest of the critical pairs involving reductions from the empty cotuple
and empty tuple rules will have similar reductions. So, we now look at the
cases where there is a critical pair involving a reduction and a conversion.

\item There are three cases for the resolution of the critical pair (3)-(13)
(dually (4)-(14): only $I = \emptyset$, only $J = \emptyset$, both $I = J =
\emptyset$. The first two are dual so we present one the first and third
case. In the first case there are two subcases. There correspond to whether
the apex (of the reduction diagram) starts with $\alpha\{\,\}$ or
$\beta\{b_j \mapsto \alpha\{\,\}\}_j$. The reduction diagrams for these
subcases are as follows:

\[\xymatrix@M=1ex@C=12ex@R=10ex@!0{
& \alpha\{\,\} ;_\gamma g \ar@{=>}[dl]_{(3)} \ar@{|=|}[dr]^{(13)} \\
\alpha\{\,\} \ar@{|=|}[d]_{(13)} 
&& \beta\{b_j \mapsto \alpha\{\,\}\}_j ;_\gamma g 
\ar@{=>}[d]^{\beta\{(3)\}}  \\
\beta\{b_j \mapsto \alpha\{\,\}\}_j &&
\beta\{b_j \mapsto \alpha\{\,\} ;_\gamma g\}_j \ar@{=>}[ll]^-{\beta\{(3)\}}}
\quad
\xymatrix@M=1ex@C=9ex@R=10ex@!0{
& \beta\{b_j \mapsto \alpha\{\,\}\}_j ;_\gamma g
\ar@{=>}[dl]_{(3)} \ar@{|=|}[dr]^{(13)} \\
\beta\{b_j \mapsto \alpha\{\,\} ;_\gamma g\}_j \ar@{=>}[d]_{\beta\{(3)\}} 
&& \alpha\{\,\} ;_\gamma g \ar@{=>}[d]^{(3)}  \\
\beta\{b_j \mapsto \alpha\{\,\}\}_j
&& \alpha\{\,\} \ar@{|=|}[ll]^-{(13)}}
\]
In the third case we have

\[\xymatrix@M=1ex@C=10ex@R=9ex@!0{
& \alpha\{\,\} ;_\gamma g \ar@{=>}[dl]_{(3)} \ar@{|=|}[dr]^{(13)} \\
\alpha\{\,\} \ar@{|=|}[dr]_{(13)}
&& \beta\{\,\} ;_\gamma g \ar@{=>}[dl]^{(3)} \\
& \beta\{\,\}}
\]
The rest of the critical pairs involving reductions and critical pairs from
the empty cotuple and empty tuple rules where neither operate on the cut
channel will have similar reductions. The case where the terms operate on
the cut channel must now be looked at.

\item The resolution of the critical pair (11)-(13) (dually (12)-(14)) is
indicated by the following reduction diagram:

\[\xymatrix@M=1ex@C=10ex@R=10ex@!0{
& \ora{\gamma}[a_k](f) ;_\gamma \gamma\{a_i \mapsto \beta\{\,\}\}
\ar@{=>}[dl]_{(11)} \ar@{|=|}[dr]^{1;(13)} \\
f ;_\gamma \beta\{\,\} \ar@{=>}[dr]_{(10)}
&& \ora{\gamma}[a_k](f) ;_\gamma \beta\{\,\} \ar@{=>}[dl]^{(10)} \\
& \beta\{\,\}}
\]
The rest of the critical pairs involving reductions and critical pairs from
the empty cotuple and empty tuple rules where they operate on the cut
channel will have similar reductions.
\end{itemize}

We now move on to examining the cases for non-empty index sets. In this case
each of the reduction diagrams has one of five shapes. Due to the number of
critical pairs, we limit ourselves to presenting one of each shape here, and
the rest in Appendix~\ref{chap-rcp}.

\begin{itemize}
\item The resolution of the critical pair (1)-(3) (dually (2)-(4))
is indicated by the following reduction diagram. 

\[\xymatrix@R=5ex@C=0ex{
& \alpha\{a_i \mapsto f_i\}_i ;_\gamma 1 \ar@{=>}[dl]_-{(1)}
\ar@{=>}[dr]^-{(3)} \\
\alpha\{a_i \mapsto f_i\}_i && \alpha\{a_i \mapsto f_i ;_\gamma 1\}_i
\ar@{=>}[ll]^-{\alpha\{(1)\}} }
\]

\item The resolution of the critical pair (3)-(6) (dually (4)-(5))
is indicated by the following reduction diagram. 

\[\xymatrix@M=1ex@C=20ex@R=10ex@!0{
& \alpha\{a_i \mapsto f_i\}_i ;_\gamma \ola{\beta}[b](g)
\ar@{=>}[dl]_{(3)} \ar@{=>}[dr]^{(6)} \\
\alpha\{a_i \mapsto f_i ;_\gamma \ola{\beta}[b](g)\}_i 
\ar@{=>}[d]_{\alpha\{(6)\}} 
&& \ola{\beta}[b](\alpha\{a_i \mapsto f_i\}_i ;_\gamma g)
\ar@{=>}[d]^{\ola{\beta}((3))} \\
\alpha\{a_i \mapsto \ola{\beta}[b](f_i ;_\gamma g)\}_i 
\ar@{|=|}[rr]_{(17)} 
&& \ola{\beta}[b](\alpha\{a_i \mapsto f_i ;_\gamma g\}_i)}
\]

\item The resolution of the critical pair (5)-(18) (dually (6)-(17)),
where $\beta \neq \gamma$, is indicated by the following reduction 
diagram. 

\[\xymatrix@M=1ex@C=20ex@R=10ex@!0{
& \ora{\alpha}[a_k](\beta(b_i \mapsto f_i)_i) ;_\gamma g 
\ar@{=>}[dl]_{(5)} \ar@{|=|}[dr]^{(18);1} \\
\ora{\alpha}[a_k](\beta(b_i \mapsto f_i)_i ;_\gamma g) 
\ar@{=>}[dd]_{\ora{\alpha}((9))} 
&& \beta(b_i \mapsto \ora{\alpha}[a_k](f_i))_i ;_\gamma g
\ar@{=>}[d]^{(9)} \\
&& \beta(b_i \mapsto \ora{\alpha}[a_k](f_i) ;_\gamma g)_i
\ar@{=>}[d]^{\beta((5))} \\
\ora{\alpha}[a_k](\beta(b_i \mapsto f_i ;_\gamma g)_i) 
\ar@{|=|}[rr]_{(18)}
&& \beta(b_i \mapsto \ora{\alpha}[a_k](f_i ;_\gamma g))_i}
\]

\item The resolution of the critical pair (7)-(16) (dually (8)-(15)),
where we start with the morphism $\ola{\alpha}[a](\gamma(b_i \mapsto f_i)_i)
;_\gamma \ora{\beta}[c](g)$, is indicated by the following reduction 
diagram. 

\[\xymatrix@M=1ex@C=20ex@R=10ex@!0{
& \ola{\alpha}[a](\gamma(b_i \mapsto f_i)_i) ;_\gamma \ora{\beta}[c](g)
\ar@{=>}[dl]_{(7)} \ar@{|=|}[dr]^{(16);1} \\
\ola{\alpha}[a](\gamma(b_i \mapsto f_i)_i ;_\gamma \ora{\beta}[c](g))
\ar@{=>}[ddd]_{\ola{\alpha}((8))} 
&& \gamma(b_i \mapsto \ola{\alpha}[a](f_i))_i ;_\gamma \ora{\beta}[c](g)
\ar@{=>}[d]^{(8)} \\
&& \ora{\beta}[c](\gamma(b_i \mapsto \ola{\alpha}[a](f_i))_i ;_\gamma g)
\ar@{|=|}[d]^{\ora{\beta}((16);1)} \\
&& \ora{\beta}[c](\ola{\alpha}[a](\gamma(b_i \mapsto f_i)_i) ;_\gamma g)
\ar@{=>}[d]^{\ora{\beta}((7))} \\
\ola{\alpha}[a](\ora{\beta}[c](\gamma(b_i \mapsto f_i)_i ;_\gamma g))
\ar@{|=|}[rr]_{(22)}
&& \ora{\beta}[c](\ola{\alpha}[a_k](\gamma(b_i \mapsto f_i)_i ;_\gamma g))}
\]

\item The resolution of the critical pair (11)-(13) (dually (12)-(13)) is
indicated by the following reduction diagram. 

\[\xymatrix@M=1ex@C=20ex@R=10ex@!0{
& \ora{\gamma}[a_k](f) ;_\gamma \gamma\{a_i \mapsto \alpha\{b_j \mapsto
g_{ij}\}_j\}_i \ar@{=>}[dl]_{(11)} \ar@{|=|}[dr]^{1;(13)} \\
f ;_\gamma \alpha\{b_j \mapsto g_{kj}\}_j \ar@{=>}[d]_{(10)}
&& \ora{\gamma}[a_k](f) ;_\gamma \alpha\{b_j \mapsto \gamma\{a_i \mapsto
g_{ij}\}_i\}_j \ar@{=>}[d]^{(10)} \\
\alpha\{b_j \mapsto f ;_\gamma g_{kj}\}_j
&& \alpha\{b_j \mapsto \ora{\gamma}[a_k](f) ;_\gamma \gamma\{a_i \mapsto
g_{ij}\}_i\}_j
\ar@{=>}[ll]^-{\alpha\{(11)\}}}
\]

\end{itemize}

Notice that each of the five shapes of reduction diagram fits the required
form to show local confluence modulo equations, and hence, confluence
modulo equations. 

The next step in proving Church-Rosser is to show that the reduction
steps terminate. This is done by associating a bag of cut costs with
a sequent derivation and showing that each reduction strictly reduces the 
bag while each equality leaves it stationary. This then will imply that
the reductions terminate. The construction of this cost criterion is our
next task.

\section{The cut measure on $\Sigma\Pi$-morphisms} \label{sec-cm}

The purpose of this section is to show that the cut elimination procedure
(defined in Section~\ref{sec-ce}) terminates. To this end we define a bag of
cut heights and show that the bag is strictly reduced on each of the cut
elimination rewrites.

We begin by defining the \textbf{height} of a term as:

\begin{itemize}
\item $\hgt[a] = 1$ when $a$ is an atomic map (or an identity)
\item $\hgt[\alpha\{a_i \mapsto f_i\}_{i \in I}] = 1+ \max\{\hgt[f_i] \mid
i \in I\}$
\item $\hgt[\ol{\alpha}[a_k] \cdot f] = 1+\hgt[f]$
\item $\hgt[f ; g] = \hgt[f] + \hgt[g]$
\end{itemize}

The \textbf{height of a cut} is defined simply as its height, e.g.,
$\cuthgt[f;g] = \hgt[f;g]$. Define a function $\Lambda: T \ra \bag(\mathbb{N})$
which takes a term to its bag of cut heights.

\begin{proposition} \quad
\begin{enumerate}[{\upshape (i)}]
\item If $\xymatrix{t_1 \ar@{=>}[r]& t_2}$ then $\Lambda(t_1) >
\Lambda(t_2)$.
\item If $\xymatrix{t_1 \ar@{|=|}[r]^{(a)}& t_2}$ and $(a)$ is an
interchange which does not involve the nullary cotuple or tuple then
$\Lambda(t_1) = \Lambda(t_2)$.
\end{enumerate}
\end{proposition}

\begin{proof}
We begin with the proof of part (i). There are three properties that must be
shown: $\hgt[t_1] \geq \hgt[t_2]$, the height of each non-principal cut does
not increase, and the height of any cut produced from the principal cut is
strictly less than the height of the principal cut.

A simple examination of the rewrites will confirm that if $t_1 \Lra t_2$
then $\hgt[t_1] \geq \hgt[t_2]$:

\noindent (1) (and dually (2)):
\[\hgt[f;1] = \hgt[f] + \hgt[1] > \hgt[f]
\]

\noindent (3) (and similarly (4), (9), and (10)):
\begin{align*}
\hgt[\alpha\{a_i \mapsto f_i\}_{i \in I} ;_\gamma g]
&= \hgt[\alpha\{a_i \mapsto f_i\}_{i \in I}] + \hgt[g] \\
&= 1+ \max\{\hgt[f_i] \mid i \in I\} + \hgt[g] \\
&= 1+ \max\{\hgt[f_i] + \hgt[g] \mid i \in I\}  \\
&= \hgt[\alpha\{a_i \mapsto f_i ;_\gamma g\}_{i \in I}]
\end{align*}
If the index set $I$ is empty we have
\begin{align*}
\hgt[\alpha\{\,\} ;_\gamma g] &= \hgt[\alpha\{\,\}] + \hgt[g] \\
&= 1+ \hgt[g] \\
&> \hgt[\alpha\{\,\}]
\end{align*}

\noindent (5) (and similarly (6), (7), and (8)):
\begin{align*}
\hgt[\ora{\alpha}[a_k](f) ;_\gamma g]
&= \hgt[\ora{\alpha}[a_k](f)]+\hgt[g] \\
&= 1+ \hgt[f] + \hgt[g] \\
&= \hgt[\ora{\alpha}[a_k](f ;_\gamma g)]
\end{align*}

\noindent (11) (and dually (12)):
\begin{align*}
\hgt[\ora{\gamma}[a_k](f) ;_\gamma \gamma\{a_i \mapsto g_i\}_{i \in I}]
&= \hgt[\ora{\gamma}[a_k](f)] + \hgt[\gamma\{a_i \mapsto g_i\}_{i \in I}] \\
&= 1+ \hgt[f] + 1 + \max\{\hgt[g_i] \mid i \in I\} \\
&> \hgt[f] + \max\{\hgt[g_i] \mid i \in I\} \\
&\geq \hgt[f] + \hgt[g_k] \\
&= \hgt[f ;_\gamma g_k]
\end{align*}

Moreover, this implies that cuts below and cuts above the redex will not
increase their cut height on a rewriting. 

Finally, consider the principal cut of the reduction. Rewrite (1) (dually (2))
removes a cut and so strictly reduces the bag of cut heights. It is an easy
observation that (5), (7), and (11) (and their duals) each replace a cut with
one of lesser height, and that (3) (and its dual) replace a cut with zero or
more cuts of lesser height. Thus applying any of the rewrites strictly
reduces the bag.

We know prove part (ii). For the equations (13) through (19) we
assume that the index sets are non-empty. This then implies that the
permuting conversions are all of the form $\alpha(\beta(f))$ and thus
\[\hgt[\alpha(\beta(f))] = 1 + \hgt[\beta(f)] = 1 + 1+ \hgt[f] = 
\hgt[\beta(\alpha(f))]
\]
which proves that the height does not change across these (non-empty
tuple and cotuple) interchanges.
\end{proof}

To see that the height is not invariant across the empty cotuple (dually the
tuple) rule recall one of the nullary versions of the rewrite (13):
\[\alpha\{\,\} \pc \beta\{b_j \mapsto \alpha\{\,\}\}_j
\]
The height on the left-hand side is one, while on the right-hand side the
height is two.

\subsection{The measure on the rewriting arrows}

We define a measure $\lambda:A \ra \bag(\mathbb{N})$ on the rewriting arrows
as follows:
\begin{itemize}
\item if $\xymatrix{t_1 \ar@{=>}[r]^x & t_2}$ then $\lambda(x) =
\min\{\Lambda(t_1),\ \Lambda(t_2)\}$
\item if $\xymatrix{t_1 \ar@{|=|}[r]^x & t_2}$ then $\lambda(x) =
\max\{\Lambda(t_1),\ \Lambda(t_2)\}$
\end{itemize}
where $\Lambda(t)$ is the bag of cut heights of $t$.

A quick examination of the reduction diagrams now confirms that this measure
will decrease when we replace a divergence with a convergence.

This completes the proof of the proposition:

\begin{proposition}
$\Sigma\Pi_\cat{A}$ under the rewrites (1)-(12) is confluent modulo the
equations (13)-(22).
\end{proposition}

\section{Deciding the $\Sigma\Pi$-conversions} \label{sec-dp}

From the above, it is clear that given any two derivations, deciding their
equivalence reduces to deciding the equivalence of cut-free proofs. Of
course, this means that any atomic cuts (cuts involving atomic formula)
must be replaced with the atomic sequent given by the appropriate composition
in the generating category \cat{A}. Thus, the decision procedure is a relative
one depending on the decision procedure for \cat{A}.

The decision procedure presented for the $\Sigma\Pi$-conversions is
graphically-inspired. It operates on pairs of terms representing cut-free
derivations of a given sequent. One of the terms is used as a ``template''
for transforming the other term into one of the same shape. The idea is to
force the second term to start with the same proof rule as the template.
If this is possible, then proceed inductively with the subterms, and
otherwise, if it is not possible, the two terms must then not be equivalent.
By using one of the terms in this manner one provides an order to search for
the conversions which make the two terms the same. 

This can be described using the term calculus or even the derivations
themselves, but is clearer with a simple graphical representation of
the terms. With a term we can associate a term-graph, whose nodes represent
the subterms of the term. Tupling and cotupling will be denoted with a
triangle decorated with the channel name it operates on, which has ``output''
edges for each component of the (co)tuple. These output edges will be
decorated with the event associated with each component. The typing of the
term will indicate which rule (tuple or cotuple) each triangle represents.
Injections and projections will be denoted by boxes decorated with the
channel name it operates on. The output edge of an injection/projection
will be decorated by the event in which it outputs. Similarly here, the
typing of the term determines whether a box represents an injection or a
projection. Atomic sequents will be represented by circular nodes containing
the atomic term, as will identities on atomic formulas.

With these conventions the permuting conversions may be represented by
the following graph equivalences (where $\alpha \neq \beta$).

\begin{itemize}
\item Cotuple-cotuple, cotuple-tuple, tuple-tuple
\[
\vcenter{\xymatrix@R=5ex@C=0ex@M=0ex{
&&& \ar@{-}[d] \\
&&& \alphatri \ar@{-}[dll]!U  \ar@{-}[drr]!U \\
& \betatri \ar@{-}[dl] \ar@{-}[dr] \ar@{}[d]|\cdots && \cdots &&
  \betatri \ar@{-}[dl] \ar@{-}[dr]  \ar@{}[d]|\cdots \\
t_{11} && t_{1n} && t_{m1} && t_{mn}
}}
\quad \pc \quad 
\vcenter{\xymatrix@R=5ex@C=0ex@M=0ex{
&&& \ar@{-}[d] \\
&&& \betatri \ar@{-}[dll]!U  \ar@{-}[drr]!U \\
& \alphatri \ar@{-}[dl] \ar@{-}[dr] \ar@{}[d]|\cdots && \cdots &&
  \alphatri \ar@{-}[dl] \ar@{-}[dr]  \ar@{}[d]|\cdots \\
t_{11} && t_{m1} && t_{1n} && t_{mn}
}}
\]

\item Cotuple-injection, cotuple-projection, tuple-injection, tuple-projection
\[
\vcenter{\xymatrix@R=5ex@C=2ex@M=0ex{
& \ar@{-}[d] \\
& \alphatri \ar@{-}[dl]  \ar@{-}[dr] \\
\betabox \ar@{-}[d]_-{i} &\cdots& \betabox \ar@{-}[d]^-{i} \\
&&
}}
\quad \pc \quad 
\vcenter{\xymatrix@R=5ex@C=3ex@M=0ex{
& \ar@{-}[d] \\
& \betabox \ar@{-}[d]_-{i} \\
& \alphatri \ar@{-}[dl]  \ar@{}[d]|-\ldots \ar@{-}[dr] \\
&&
}}
\]

\item Injection-injection, injection-projection, projection-projection
\[
\vcenter{\xymatrix@R=5ex@C=1ex@M=0ex{
\ar@{-}[d] \\
\alphabox \ar@{-}[d]_-{i} \\
\betabox \ar@{-}[d]_-{j} \\ 
&
}}
\quad \pc \quad 
\vcenter{\xymatrix@R=5ex@C=1ex@M=0ex{
\ar@{-}[d] \\
\betabox \ar@{-}[d]_-{j} \\ 
\alphabox \ar@{-}[d]_-{i} \\
&
}}
\]

\end{itemize}

To illustrate the graphical representation, the third derivation in
Example~\ref{exam-logic} can be represented by the following graph.
Note that the graph is quite a direct representation of the derivation tree.

{\scriptsize
\[\xymatrix@R=4ex@C=0ex@M=0ex{
&&& \ar@{-}[d] \\
&&& \alphatri \ar@{-}[dll]!U_-a \ar@{-}[drr]!U^-b \\
& \betabox \ar@{-}[d]_-g &&&& \betabox \ar@{-}@{-}[d]^-h \\
& \gammabox \ar@{-}[d]_-{a'} &&&& \gammabox \ar@{-}[d]^-{b'} \\
& \deltatri \ar@{-}[dl]!U_-{c'} \ar@{-}[dr]!U^-{d'} &&&&
  \deltatri \ar@{-}[dl]!U_-{c'} \ar@{-}[dr]!U^-{d'} \\
\alphabox \ar@{-}[d]_-c && \alphabox \ar@{-}[d]^-d &&
  \alphabox \ar@{-}[d]_-e && \alphabox \ar@{-}[d]^-f \\
*+<10pt>[Fe]{i} && \deltabox \ar@{-}[d]^-{e'} && *+<8pt>[Fe]{k}
&& \deltabox \ar@{-}[d]^-{f'} \\
&& *+<8pt>[Fe]{j} &&&& *+<10pt>[Fe]{l}}
\]
}

An equivalent derivation is given as follows

{\tiny
\[\infer{(A \times B)+(C \times D),\ E \times F \vd G+H,\ I \times (J+K)}{
\infer{(A \times B)+(C \times D),\ E \times F \vd G+H,\ I}{
  \infer{A \!\times\! B,E \!\times\! F \vd G \!+\! H,I}{
    \infer{A,\ E \times F \vd G+H,\ I}{
      \infer{A,\ E \vd G+H,\ I}{
        \infer{A,\ E \vd_i G,\ I}{}}}} &
  \infer{C \!\times\! D,E \!\times\! F \vd G\!+\!H,I}{
    \infer{C,\ E \times F \vd G+H,\ I}{
      \infer{C,\ F \vd G+H,\ I}{
        \infer{C,\ F \vd_k H,\ I}{}}}}} &
\infer{(A \times B)+(C \times D),\ E \times F \vd G+H,\ J+K}{
  \infer{A \!\times\! B,E \!\times\! F \vd G\!+\!H,J\!+\!K}{
    \infer{B,\ E \times F \vd G+H,\ J+K}{
      \infer{B,\ F \vd G+H,\ J+K}{
        \infer{B,\ F \vd G,\ J+K}{
          \infer{B,\ F \vd_j G,\ J}{}}}}} &
  \infer{C \!\times\! D,E \!\times\! F \vd G\!+\!H,J\!+\!K}{
    \infer{D,\ E \times F \vd G+H,\ J+K}{
      \infer{D,\ F \vd G+H,\ J+K}{
        \infer{D,\ F \vd H,\ J+K}{
          \infer{D,\ F \vd_l H,\ K}{}}}}}}}
\]}
which is given (in its annotated version) by the graph below:

{\scriptsize
\[\xymatrix@R=4ex@C=0ex@M=0ex{
&&& \ar@{-}[d] \\
&&& \deltatri \ar@{-}[dll]!U_-{c'} \ar@{-}[drr]!U^-{d'} \\
& \alphatri \ar@{-}[dl]_-a \ar@{-}[dr]^-b &&&&
  \alphatri \ar@{-}[dl]_-a \ar@{-}[dr]^-b \\
\alphabox \ar@{-}[d]_-c && \alphabox \ar@{-}[d]^-e &&
  \alphabox \ar@{-}[d]_-d && \alphabox \ar@{-}[d]^-f \\
\betabox \ar@{-}[d]_-g && \betabox \ar@{-}[d]^-h &&
  \betabox \ar@{-}[d]_-g && \betabox \ar@{-}[d]^-h \\
\gammabox \ar@{-}[d]_-{a'} && \gammabox \ar@{-}[d]^-{b'} &&
  \gammabox \ar@{-}[d]_-{a'} && \gammabox \ar@{-}[d]^-{b'} \\
*+<10pt>[Fe]{i} && *+<8pt>[Fe]{k} && \deltabox \ar@{-}[d]_-{e'}
&& \deltabox \ar@{-}[d]^-{f'} \\
&&&& *+<8pt>[Fe]{j} && *+<10pt>[Fe]{l}}
\]
}

We shall illustrate the decision procedure with this example. Take the 
first graph as a template. The first step in the procedure is to see if the
second graph can start the same way as the first graph. This means we have
to move an $\alpha$-triangle up to the topmost level. This involves searching
through the second graph until an $\alpha$-triangle is found that can be
moved upwards in the necessary manner. In this case there is one at the
second level. Moving it up gives the graph on the left below. (The labels
have been removed as this will not cause any ambiguity here, but in any
case they may be inferred from the previous graph.)

{\scriptsize
\[\xymatrix@R=4ex@C=0ex@M=0ex{
&&& \ar@{-}[d] \\
&&& \alphatri \ar@{-}[dll]!U \ar@{-}[drr]!U \\
& \deltatri \ar@{-}[dl]!U \ar@{-}[dr]!U
&&&& \deltatri \ar@{-}[dl]!U \ar@{-}[dr]!U \\
\alphabox \ar@{-}[d] && \alphabox \ar@{-}[d] 
&& \alphabox \ar@{-}[d] && \alphabox \ar@{-}[d] \\
\betabox \ar@{-}[d] && \betabox \ar@{-}[d]
&& \betabox \ar@{-}[d] && \betabox \ar@{-}[d] \\
\gammabox \ar@{-}[d] && \gammabox \ar@{-}[d]
&& \gammabox \ar@{-}[d] && \gammabox \ar@{-}[d] \\
*+<10pt>[Fe]{i} && \deltabox \ar@{-}[d] && *+<8pt>[Fe]{k}
&& \deltabox \ar@{-}[d] \\
&& *+<8pt>[Fe]{j} &&&& *+<10pt>[Fe]{l}}
\qqqquad 
\qqqquad 
\xymatrix@R=4ex@C=0ex@M=0ex{
&&& \ar@{-}[d] \\
&&& \alphatri \ar@{-}[dll]!U \ar@{-}[drr]!U \\
& \betabox \ar@{-}[d] &&&& \betabox \ar@{-}@{-}[d] \\
& \deltatri \ar@{-}[dl]!U \ar@{-}[dr]!U
&&&& \deltatri \ar@{-}[dl]!U \ar@{-}[dr]!U \\
\alphabox \ar@{-}[d] && \alphabox \ar@{-}[d]
&& \alphabox \ar@{-}[d] && \alphabox \ar@{-}[d] \\
\gammabox \ar@{-}[d] && \gammabox \ar@{-}[d]
&& \gammabox \ar@{-}[d] && \gammabox \ar@{-}[d] \\
*+<10pt>[Fe]{i} && \deltabox \ar@{-}[d] && *+<8pt>[Fe]{k}
&& \deltabox \ar@{-}[d] \\
&& *+<8pt>[Fe]{j} &&&& *+<10pt>[Fe]{l}}
\]
}

Moving down a level we inductively repeat the process for all subterms at
the second level nodes: the next step is to pull up the $\beta$-boxes on the
left path and right path. This is done by interchanging the $\alpha$-boxes
and $\beta$-boxes in all the paths and then flipping the $\delta$-triangles
and $\beta$-boxes, which results in the graph on the right above.
The last step would be to bring the $\gamma$-boxes up to the third level.
This is done by interchanging the $\alpha$-boxes and the $\gamma$-boxes
in all the paths and then flipping the $\delta$-triangles and $\gamma$-boxes;
this produces the required graph (i.e., we were able to transform the second
graph into the ``template'' graph), and so completes the proof that the two
original derivations are equivalent. In general, the decision procedure will
proceed in this recursive manner.

\subsection{The details and proof of the decision procedure}

A term is \textbf{$\alpha$-inert} if it does not contain as subterms
$\ola{\alpha}[a](t),\ \ora{\alpha}[a](t),\ \alpha\{a_i \mapsto t_i\}_{i \in I}$,
or $\alpha(a_i \mapsto t_i)_{i \in I}$. Clearly if $t$ is $\alpha$-inert
then there is no equality involving $\alpha$ which applies to it.

Let $C_\alpha$ be the constructors $\ola{\alpha}[\cdot](\_)$,
$\ora{\alpha}[\cdot](\_)$, $\alpha\{\_\}_{i \in I}$, or $\alpha(\_)_{i \in I}$.
We shall say a term starts with constructor $C_\alpha$ in case the first
constructor in the term is $C_\alpha$.

The \textbf{$C_\alpha$-prefix} of a term $t$, denoted
$\mathrm{prefix}_{C_\alpha}[t]$, is defined as follows.

\begin{itemize}

\item If $t$ starts with $C_\alpha$ then $\mathrm{prefix}_{C_\alpha}[t] = *$
(where $*$ is the ``anonymous'' variable, i.e., a distinct variable which
has not been used before and will not be used again).

\item If $t$ does not start with $C_\alpha$ then
\begin{itemize}

\item if $t = \ola{\beta}[b](t')$ then $\mathrm{prefix}_{C_\alpha}
[\ola{\beta}[b](t')] = \ola{\beta}[b](\mathrm{prefix}_{C_\alpha}[t'])$.

\item if $t = \ora{\beta}[b](t')$ then $\mathrm{prefix}_{C_\alpha}
[\ora{\beta}[b](t')] = \ora{\beta}[b](\mathrm{prefix}_{C_\alpha}[t'])$.

\item if $t = \beta\{b_i \mapsto t_i\}_{i \in I}$ then
\[\mathrm{prefix}_{C_\alpha} [\beta\{b_i \mapsto t_i\}_{i \in I}] =
\beta\{b_i \mapsto \mathrm{prefix}_{C_\alpha}[t_i]\}_{i \in I}.
\]

\item if $t = \beta(b_i \mapsto t_i)_{i \in I}$ then
\[\mathrm{prefix}_{C_\alpha}[\beta(b_i \mapsto t_i)_{i \in I}] =
\beta(b_i \mapsto \mathrm{prefix}_{C_\alpha}[t_i])_{i \in I}.
\]
\end{itemize}
\end{itemize}

The \textbf{$C_\alpha$-frontier} of a term with a $C_\alpha$-prefix is
those first occurrences across the term of the constructor $C_\alpha$.

\begin{lemma}
Suppose a term $t$ starts with constructor $C_\alpha$. Then in any series of
equalities
\[\xymatrix{t \ar@{|=|}[r] & t_1 \ar@{|=|}[r] & t_2 \ar@{|=|}[r] &
\ \cdots\
\ar@{|=|}[r] & t_n}
\]
the $C_\alpha$-prefix of each $t_i$ is $\alpha$-inert.
\end{lemma}

\begin{proof}
Suppose that $t = \ola{\alpha}[a](t')$ and $\mathrm{prefix}_{\ola{\alpha}(\_)}
[t_i] = w$, where $w$ is $\alpha$-inert. Either $t_i \pc t_{i+1}$ is an 
application of an equality at the $\ola{\alpha}[\cdot](\_)$-frontier of the
inert term $w$ or it is not. If it is beyond the frontier then
$\mathrm{prefix}_{\ola{\alpha}[\cdot](\_)}[t_i] =
\mathrm{prefix}_{\ola{\alpha}[\cdot](\_)}[t_{i+1}]$ and if it is before the
frontier then it simple shuffles the $\alpha$-inert prefix.
If it is on the frontier either it moves structure out of the inert term by
shrinking the frontier (in which case
$\mathrm{prefix}_{\ola{\alpha}[\cdot](\_)}[t_{i+1}]$ is certainly still inert
(if smaller)), or it moves structure into the prefix by expanding the
frontier. However, only $\alpha$-inert structure can be moved over
$\ola{\alpha}[\cdot](\_)$, so again
$\mathrm{prefix}_{\ola{\alpha}[\cdot](\_)}[t_{i+1}]$ is $\alpha$-inert.
\end{proof}

In a series of equalities beginning with a term which starts with $C_\alpha$
we may distinguish the steps which increase the $C_\alpha$-inert prefix,
$\xymatrix{t_i \ar @{|=>} [r] & t_{i+1}}$, those which decrease the 
$C_\alpha$-inert prefix, $\xymatrix{t_i & \ar@{|=>}[l] t_{i+1}}$, and those
which do not affect the $C_\alpha$-inert prefix, 
$\xymatrix{t_i & \ar@{|=|}[l] t_{i+1}}$.

\begin{lemma}
\[\xymatrix{t_i \ar@{|=>}[r]^x & t_{i+1} \ar@{|=|}[r]^y & t_{i+2}}
\]
can be rearranged as
\[\xymatrix{t_i \ar@{|=|}[r]^y & t'_{i+1} \ar@{|=>}[r]^x & t_{i+2}}
\]
\end{lemma}

\begin{proof}
The redex of $y$ cannot be within the inert prefix, nor by assumption is it
on the frontier, and thus, it must be independent of $x$ (on the frontier),
and thus, the equalities can be rearranged.
\end{proof}

This means that we can rearrange the steps in any proof of equality so that
no $\alpha$-inert prefix-increasing step happens before a step which does
not affect the inert prefix. However, we are not able to move these
increasing steps past an inert prefix-decreasing step, but as a decreasing
step is only possible if there has already been the corresponding (reverse)
increasing step, this decreasing step may be cancelled with the increasing
step. Thus, we may conclude the following.
 
\begin{lemma}
In any series of steps
\[\xymatrix{t \ar@{|=|}[r]_{*} & t_1 \ar@{|=>}[r]_{*}^x & t_2 &
t_3 \ar@{|=>}[l]_y}
\]
the decreasing step $y$ can be cancelled with its corresponding increasing
step.
\end{lemma}

The following proposition now follows.

\begin{proposition}
Any proof of equality from $t$ to $t'$ can be rearranged as
\[\xymatrix{t \ar@{|=|}[r]_{*}^x & t_1 \ar@{|=>}[r]_{*}^y & t'}
\]
where the initial equalities do not touch the root constructor.
\end{proposition}

Note that the prefix-increasing steps $y$ are essentially unique. There may
be independent expansions of the inert frontier which can be reordered but
every equivalent proof must do the same expansions. The point is that reading
these steps in reverse gives a procedure which pulls the root structure of $t$
to the root of $t'$.

\begin{corollary}
The structure $C_\alpha$ may be pulled to the root of $t$ if and only if the
$C_\alpha$-prefix of $t$ is $\alpha$-inert.
\end{corollary}

\begin{proof}
It has already been observed that such a pulling up process results in an
inert prefix (of the appropriate sort). Conversely given an inert prefix
of the appropriate sort clearly means that the $C_\alpha$-frontier can be
contracted shrinking the $C_\alpha$-prefix.
\end{proof}

As the equalities in the first part of the proof, $x$, do not touch the
root constructor each equality must apply to one of the arguments of that
constructor. Thus, for each argument we then have an equality proof, but
each of these proofs can now be ``normalized'' into the form of the corollary.
This gives a normal form for (directed) equality proofs and whence an
algorithm for determining equality which amounts to matching the structure
of the first term $t$ starting from the root and pulling up that structure
to the root of the second term and then proceeding inductively with its
subterms.

\chapter{Polycategorical Semantics of Processes} \label{chap-catsem}

In this chapter we establish that $\Sigma\Pi_\cat{A}$ is the free
polycategory with sums and products generated from an arbitrary 
polycategory \cat{A}. The derivations, when cut is eliminated and treated
up to conversion, will be the morphisms of this polycategory.

In Section~\ref{sec-sums}, we define and show that $\Sigma\Pi_\cat{A}$ has
poly-sums and poly-products. Then, in Section~\ref{sec-softness}, it is shown
that the sums and products are ``soft''. Softness of the sum and products is
necessary to establish the ``Whitman Theorem'' which characterizes the free
polycategory with sums and products up to equivalence.

We begin by establishing that $\Sigma\Pi_\cat{A}$ is a polycategory.

\section{$\Sigma\Pi_\cat{A}$ is a polycategory}

The reduction rules and the permuting conversions together define an
equivalence relation (which we denote by $\sim$) on derivations of a
sequent. Our categorical semantics will have derivations modulo this
equivalence as morphisms. Thus, the goal of this section is to prove:

\begin{theorem} \label{thm-sppc}
$\Sigma\Pi_\cat{A}$ is a polycategory whose objects are the formulas of
the logic, and whose morphisms are $\sim$-equivalence classes of
derivations.
\end{theorem}

We shall start by presenting some technical results about the derivations
in $\Sigma\Pi_\cat{A}$.

A notion which will be useful throughout the rest of this thesis is that of
an output sequent. A sequent $\Gamma \vd \Delta$ is an \textbf{output
sequent} in case:
\begin{itemize}
\item $X \in \Gamma$ implies that $X$ is an atom or of the form
$X = \prod_{i \in I} X_i$,
\item $Y \in \Delta$ implies that $Y$ is an atom or of the from
$Y = \sum_{j \in J} Y_j$, and
\item $\Gamma \cup \Delta$ contains at least one compound formula.
\end{itemize}

\begin{proposition} \label{prop-insp}
In $\Sigma\Pi_\cat{A}$:
\begin{enumerate}[{\upshape (i)}]
\item Any cut-free derivation of a sequent $\Gamma,\sum_i X_i \vd \Delta$ is
equivalent to one whose principal rule is a cotupling applied to $\sum_i X_i$.

\item Any cut-free derivation of a sequent $\Gamma \vd \prod_i Y_i,\Delta$ is
equivalent to one whose principal rule is a tupling applied to $\prod_i Y_i$.

\item Any cut-free derivation of an output sequent $\Gamma \vd \Delta$ has as
its principal rule either an injection or a projection.

\item Any cut-free derivation of a sequent $A_1,\ldots,A_m \vd B_1,\ldots,B_n$,
where the $A_i$'s and $B_j$'s are atoms, must be an axiom (i.e., a morphism of 
\cat{A}).
\end{enumerate}
\end{proposition}

Notice that this result can be extended to arbitrary derivations (not only
cut-free ones) using the cut-elimination procedure. For example, any
derivation of a sequent $\Gamma,\sum_i X_i \vd \Delta$ can be transformed
to one whose principal rule is cotupling.

\begin{proof}
Notice that (iii) and (iv) are immediate given which (non-cut) sequent rules 
may be applied in the respective situations. (ii) is the dual of (i) and 
thus, it is sufficient to prove (i).

\begin{itemize}
\item If $\sum_i X_i$ is the only compound formula (i.e., all the other 
formulas are atoms) then the cotuple rule is the only one which applies, so 
this must be the principal rule.

\item If the principal rule is another cotupling 
\[\infer{\sum_j Y_j,\Gamma,\sum_i X_i \vd \Delta}
{\left\{Y_j,\Gamma,\sum_i X_i \vd \Delta\right\}_j}
\]
then, using our inductive hypothesis, the proofs above may be transformed to 
have cotupling applied to $\sum_i X_i$ as their principal rule. This allows 
us to use the cotuple-cotuple interchange to transform the derivation so that
cotupling applied to $\sum_i X_i$ is the principal rule for the original proof.

\item If the principal rule is a tupling
\[\infer{\Gamma,\sum_i X_i \vd \prod_j Y_j,\Delta'}{
\left\{\Gamma,\sum_i X_i \vd Y_j,\Delta'\right\}_j}
\]
then, using our inductive hypothesis, the proofs above may be transformed to 
have cotupling applied to $\sum_i X_i$ as their principal rule. This allows 
us to use the cotuple-tuple interchange to transform the derivation so that
cotupling applied to $\sum_i X_i$ is the principal rule for the original proof.

\item If the principal rule is a projection
\[\infer{\prod_j Y_j,\Gamma,\sum_i X_i \vd \Delta}
{Y_k,\Gamma,\sum_i X_i \vd \Delta}
\]
then, using our inductive hypothesis, the proof above may be transformed to 
have cotupling applied to $\sum_i X_i$ as its principal rule. This allows 
us to use the cotuple-projection interchange to transform the derivation so
that cotupling applied to $\sum_i X_i$ is the principal rule for the original
proof.

\item If the principal rule is an injection
\[\infer{\Gamma,\sum_i X_i \vd \sum_j Y_j,\Delta'}
{\Gamma,\sum_i X_i \vd Y_k,\Delta'}
\]
then, using our inductive hypothesis, the proof above may be transformed to 
have cotupling applied to $\sum_i X_i$ as its principal rule. This allows us
to use the cotuple-injection interchange to transform the derivation so that
cotupling applied to $\sum_i X_i$ is the principal rule for the original proof.

\end{itemize}
\end{proof}

Lemma~\ref{lem-id} has already shown that there are identity derivations
which act in the correct manner. It is left to show that the composition
given by cut satisfies the associativity and interchange laws. 

\begin{proposition} 
Cut satisfies the associative law. That is, given sequents of the form
\[\Gamma \vd_f \Delta,\gamma:X \qquad \gamma:X,\Gamma' \vd_g \Delta',
\delta:Y \qquad \delta:Y,\Gamma'' \vd_h \Delta''
\]
the composites $(f ;_\gamma g) ;_\delta h$ and $f ;_\gamma (g ;_\delta h)$
are $\sim$-equivalent.
\end{proposition}

\begin{proof}
By structural induction on $f$, $g$, and $h$; without loss of generality we
may assume that $f$, $g$, and $h$ are all cut-free.

\begin{enumerate}[{\upshape (i)}]
\item $f = 1_X$. In this case we have $(1_X ;_\gamma g) ;_\delta h
\Lra g ;_\delta h$ and $1_X ;_\gamma (g ;_\delta h) \Lra g ;_\delta h$ which 
are equivalent. The case where $h = 1_Y$ is dual.

\item $f = \alpha\{a_i \mapsto f_i\}_i$ and $\Gamma = \Phi,\alpha:\sum_i W_i$.
Composing on the left first gives
\[(\alpha\{a_i \mapsto f_i\}_i ;_\gamma g) ;_\delta h \Lra 
\alpha\{a_i \mapsto f_i ;_\gamma g\}_i ;_\delta h \Lra 
\alpha\{a_i \mapsto (f_i ;_\gamma g) ;_\delta h\}_i
\]
whereas composing on the right first gives
\[\alpha\{a_i \mapsto f_i\}_i ;_\gamma (g ;_\delta h) \Lra 
\alpha\{a_i \mapsto f_i ;_\gamma (g ;_\delta h)\}_i
\]
which by induction are $\sim$-equivalent. Since any derivation of a sequent 
$\Gamma,\alpha:\sum_i W_i \vd \Delta$ may be equivalently written to end with 
cotupling applied on $\alpha$, this suffices. The case where $h =
(b_j \mapsto h_j)_j$ and $\Delta'' = \prod_j W_j,\Psi$
is dual.

\item We now consider the case where $f$ is a map to a coproduct.
There are two cases to consider.
\\\smallskip

(a) $f = \ora{\alpha}[a_k](f')$, $\alpha \neq \gamma$, and 
$\Delta = \alpha:\sum_i W_i,\Phi$. Composing on the left first gives
\[(\ora{\alpha}[a_k](f') ;_\gamma g) ;_\delta h \Lra
\ora{\alpha}[a_k](f' ;_\gamma g) ;_\delta h \Lra
\ora{\alpha}[a_k]((f' ;_\gamma g) ;_\delta h)
\]
whereas composing on the right first gives
\[\ora{\alpha}[a_k](f') ;_\gamma (g ;_\delta h) \Lra
\ora{\alpha}[a_k](f' ;_\gamma (g ;_\delta h))
\]
which by induction are $\sim$-equivalent. The case where $h =
\ola{\beta}[b_k](h')$ with $\beta \neq \delta$ and $\Gamma'' = \Psi, 
\beta:\prod_j W_j$ is dual.
\\\smallskip

(b) $f = \ora{\gamma}[a_k](f')$ and $X = \gamma:\sum_i X_i$. Here
we may suppose that $g = \gamma\{a_i \mapsto g_i\}_{i}$ (as $X$ is a sum) and
so composing on the left first gives
\[(\ora{\gamma}[a_k](f') ;_\gamma \gamma\{a_i \mapsto g_i\}_{i}) 
;_\delta h \Lra (f' ;_\gamma g_k) ;_\delta h 
\]
whereas composing on the right first gives
\begin{align*}
\ora{\gamma}[a_k](f') ;_\gamma (\gamma\{a_i \mapsto g_i\}_{i} ;_\delta h)
&\Lra \ora{\gamma}[a_k](f') ;_\gamma \gamma\{a_i \mapsto g_i ;_\delta h\}_{i}\\ 
&\Lra f' ;_\gamma (g_k ;_\delta h)
\end{align*}
which by induction are $\sim$-equivalent. The case where $h =
\ola{\delta}[b_k](h')$ and $Y= \delta:\prod_j Y_j$ is dual.

\item $f = \ola{\alpha}[a_k](f')$ and $\Gamma = \Phi,\alpha: \prod_i W_i$.
Composing on the left first gives
\[(\ola{\alpha}[a_k](f') ;_\gamma g) ;_\delta h \Lra
\ola{\alpha}[a_k](f' ;_\gamma g) ;_\delta h \Lra
\ola{\alpha}[a_k]((f' ;_\gamma g) ;_\delta h)
\]
whereas composing on the right first gives
\[\ola{\alpha}[a_k](f') ;_\gamma (g ;_\delta h) \Lra
\ola{\alpha}[a_k](f' ;_\gamma (g ;_\delta h))
\]
which by induction are $\sim$-equivalent. The case where $h=\ora{\beta}(h')$ 
and $\Delta'' = \beta:\sum_j W_j,\Psi$ is dual.

\item We now consider the case where $f$ is a map to a product. There are a
number of subcases to consider.
\\\smallskip

(a) $f = \alpha(a_i \mapsto f_i)_{i}$, $\alpha \neq \gamma$, and $\Delta =
\alpha:\prod_i W_i,\Phi$. Composing on the left first gives
\[(\alpha(a_i \mapsto f_i)_i ;_\gamma g) ;_\delta h \Lra 
\alpha(a_i \mapsto f_i ;_\gamma g)_i ;_\delta h \Lra 
\alpha(a_i \mapsto (f_i ;_\gamma g) ;_\delta h)_i 
\]
whereas composing on the right first gives
\[\alpha(a_i \mapsto f_i)_i ;_\gamma (g ;_\delta h) \Lra \alpha(a_i 
\mapsto f_i ;_\gamma (g ;_\delta h))_i
\]
which by induction are $\sim$-equivalent. 
Since any derivation of a sequent $\Gamma \vd \alpha:\prod_i W_i,\Delta$ may 
be equivalently written to end with tupling applied on $\alpha$, this suffices. 
The case where $h=\beta\{b_j \mapsto h_j\}_j$ with $\beta \neq \delta$
and $\Gamma'' = \Psi,\beta:\sum_j W_j$ is dual.
\\\smallskip

(b) $f = \gamma(a_i \mapsto f_i)_i$, $X = \gamma:\prod_i X_i$, and
$g = \ola{\gamma}[a_k](g')$. Composing on the left first gives 
\[(\gamma(a_i \mapsto f_i)_i ;_\gamma \ola{\gamma}[a_k](g')) ;_\delta 
h \Lra (f_k ;_\gamma g') ;_\delta h
\]
whereas composing on the right first gives
\begin{align*}
\gamma(a_i \mapsto f_i)_i ;_\gamma (\ola{\gamma}[a_k](g') ;_\delta h)
&\Lra \gamma(a_i \mapsto f_i)_i ;_\gamma \ola{\gamma}[a_k](g' ;_\delta h)\\
&\Lra f_k ;_\gamma (g' ;_\delta h)
\end{align*}
which by induction are $\sim$-equivalent. 
Since any derivation of a sequent $\Gamma \vd \alpha:\prod_i W_i,\Delta$ may 
be equivalently written to end with tupling applied on $\alpha$, this suffices. 
The case where $h=\delta\{b_j \mapsto h_j\}_j$ and $Y= \delta:
\sum_j Y_j$ is dual.
\\\smallskip

(c) $f = \gamma(a_i \mapsto f_i)_i$, $X = \gamma:\prod_i X_i$, and $g$
operates on a channel $\beta \neq \gamma$. In this case we must explore
the structure of $g$. First, let us assume $g = \beta\{b_j \mapsto g_j\}_j$.
Composing on the left first gives
\begin{align*}
(\gamma(a_i \mapsto f_i)_i ;_\gamma \beta\{b_j \mapsto g_j\}_j) ;_\delta h
&\Lra \beta\{b_j \mapsto \gamma(a_i \mapsto f_i)_i ;_\gamma g_j\}_j
;_\delta h\\ &\Lra
\beta\{b_j \mapsto (\gamma(a_i \mapsto f_i)_i ;_\gamma g_j) ;_\delta h\}_j
\end{align*}
whereas composing on the right first gives
\begin{align*}
\gamma(a_i \mapsto f_i)_i ;_\gamma (\beta\{b_j \mapsto g_j\}_j ;_\delta h)
&\Lra \gamma(a_i \mapsto f_i)_i ;_\gamma \beta\{b_j \mapsto g_j ;_\delta h\}_j 
\\ &\Lra 
\beta\{\gamma(a_i \mapsto f_i)_i ;_\gamma (b_j \mapsto g_j ;_\delta h)\}_j
\end{align*}
The remaining cases where $g$ is either $\beta(b_j \mapsto g_j)_j$,
$\ola{\beta}[b_k](g')$, or $\ora{\beta}[b_k](g')$ are all handled similarly.

\end{enumerate}

This concludes all the essential cases. If all of $f$, $g$, and $h$ are atomic
then composition is associative because it is associative in the underlying
polycategory. If some of $f$, $g$, and $h$ are atomic a quick check of the
possibilities will show that one ends up with a case essentially like one of
the cases above.

\end{proof}

For polycategories one must also prove that composition satisfies the
interchange law.

\begin{proposition} 
Cut satisfies the interchange law. That is, given sequents of the form
\[
\Gamma \vd_f \gamma:X,\Delta,\delta:Y \qquad \gamma:X,\Gamma' \vd_g \Delta' 
\qquad \delta:Y,\Gamma'' \vd_h \Delta''
\]
the composites $(f ;_\gamma g) ;_\delta h$ and $(f ;_\gamma h) ;_\delta g$
are $\sim$-equivalent. Dually, given sequents of the form
\[
\Gamma \vd_f \gamma:X,\Delta \qquad \Gamma' \vd_g \delta:Y,\Delta' \qquad
\gamma:X,\Gamma'',\delta:Y \vd_h \Delta''
\]
then $f ;_\gamma (g ;_\delta h)$ and $g ;_\delta (f ;_\gamma h)$ are 
$\sim$-equivalent.
\end{proposition}

\begin{proof}
We prove the first statement by structural induction on $f$, $g$, and $h$,
duality handles the second statement. Without loss of generality we may
assume that $f$, $g$, and $h$ are all cut-free.

\begin{enumerate}[{\upshape (i)}]
\item $f = \alpha\{a_i \mapsto f_i\}_i$ and $\Gamma = \Phi,
\alpha:\sum_i W_i$. Composing with $g$ first gives
\[(\alpha\{a_i \mapsto f_i\}_i ;_\gamma g) ;_\delta h \Lra 
\alpha\{a_i \mapsto f_i ;_\gamma g\}_i ;_\delta  h \Lra 
\alpha\{a_i \mapsto (f_i ;_\gamma g) ;_\delta h\}_i
\]
whereas composing with $h$ first gives
\[(\alpha\{a_i \mapsto f_i\}_i ;_\delta h) ;_\gamma h \Lra 
\alpha\{a_i \mapsto f_i ;_\delta h\}_i ;_\gamma  g \Lra 
\alpha\{a_i \mapsto (f_i ;_\delta h) ;_\gamma g\}_i
\]
which, by induction, are $\sim$-equivalent. Since any derivation of a sequent
$\Gamma,\sum_i W_i$ $\vd \Delta$ may be equivalently written to end with
cotupling applied on $\alpha$, this suffices.

\item We now look at the case where $f$ is a map to a coproduct. There are
three subcases to consider.
\\\smallskip

(a) $f = \ora{\gamma}[a_k](f')$ and $X = \gamma:\sum_i X_i$. 
Here we may suppose that $g = \gamma\{a_i \mapsto g_i\}_i$ (as $X$ is 
a sum) and so, composing with $g$ first gives
\[(\ora{\gamma}[a_k](f') ;_\gamma \gamma\{a_i \mapsto g_i\}_i) 
;_\delta h \Lra (f' ;_\gamma g_k) ;_\delta h
\]
whereas composing with $h$ first gives
\begin{align*}
(\ora{\gamma}[a_k](f') ;_\delta h) ;_\gamma \gamma\{a_i \mapsto g_i\}_ i
&\Lra \ora{\gamma}[a_k](f' ;_\delta h) ;_\gamma \gamma\{a_i \mapsto g_i\}_i\\
&\Lra (f' ;_\delta h) ;_\gamma g_k
\end{align*}
which, by induction, are $\sim$-equivalent.
\\\smallskip

(b) $f = \ora{\alpha}[a_k](f'),\ \alpha \neq \gamma,\ \alpha \neq \delta$,
and $\Delta = \Phi, \alpha:\sum_i W_i, \Psi$. Composing with $g$ first gives
\[(\ora{\alpha}[a_k](f') ;_\gamma g) ;_\delta h \Lra 
\ora{\alpha}[a_k](f' ;_\gamma g) ;_\delta h \Lra 
\ora{\alpha}[a_k]((f' ;_\gamma g) ;_\delta h)
\]
whereas composing with $h$ first gives 
\[(\ora{\alpha}[a_k](f') ;_\delta h) ;_\gamma g \Lra 
\ora{\alpha}[a_k](f' ;_\delta h) ;_\gamma g \Lra 
\ora{\alpha}[a_k]((f' ;_\delta h) ;_\gamma g)
\]
which, by induction, are $\sim$-equivalent.
\\\smallskip

(c) $f = \ora{\delta}[a_k](f')$ and $Y = \delta:\sum_i Y_i$. Here we
may suppose that $h = \delta\{a_i \mapsto h_i\}_i$ (as $Y$ is a sum) and so,
composing with $g$ first gives
\begin{align*}
(\ora{\gamma}[a_k](f') ;_\gamma g) ;_\delta \delta\{a_i \mapsto h_i\}_i
&\Lra \ora{\gamma}[a_k](f' ;_\gamma g) ;_\delta \delta\{a_i \mapsto h_i\}_i\\
&\Lra (f' ;_\gamma g) ;_\delta h_k
\end{align*}
whereas composing with $h$ first gives
\[(\ora{\gamma}[a_k](f') ;_\delta \delta\{a_i \mapsto h_i\}_i) 
;_\gamma g \Lra (f' ;_\delta h_k) ;_\gamma g
\]
which, by induction, are $\sim$-equivalent.
\\\smallskip

\item $f = \ola{\alpha}[a_k](f')$ and $\Gamma = \Phi,\alpha:
\prod_i W_i$. Composing with $g$ first gives
\[(\ola{\alpha}[a_k](f') ;_\gamma g) ;_\delta h \Lra
\ola{\alpha}[a_k](f' ;_\gamma g) ;_\delta h \Lra
\ola{\alpha}[a_k]((f' ;_\gamma g) ;_\delta h)
\]
whereas composing with $h$ first gives
\[(\ola{\alpha}[a_k](f') ;_\delta h) ;_\gamma g \Lra
\ola{\alpha}[a_k](f' ;_\delta h) ;_\gamma g \Lra
\ola{\alpha}[a_k]((f' ;_\delta h) ;_\gamma g)
\]
which, by induction, are $\sim$-equivalent.

\item We now look at the case where $f$ is a map to a product. Again, there
are a number of subcases to consider. 
\\\smallskip

(a) $f = \alpha(a_i \mapsto f_i)_i,\ \alpha \neq \gamma,\ \alpha \neq \delta$,
and $\Delta = \alpha:\prod_i W_i,\Phi$. Composing with $g$ first gives
\[(\alpha(a_i \mapsto f_i)_i ;_\gamma g) ;_\delta h \Lra
\alpha(a_i \mapsto f_i ;_\gamma g)_i ;_\delta h \Lra
\alpha(a_i \mapsto (f_i ;_\gamma g) ;_\delta h)_i
\]
whereas composing with $h$ first gives
\[(\alpha(a_i \mapsto f_i)_i ;_\delta h) ;_\gamma g \Lra
\alpha(a_i \mapsto f_i ;_\delta h)_i ;_\gamma g \Lra
\alpha(a_i \mapsto (f_i ;_\delta h) ;_\gamma g)_i
\]
which, by induction, are $\sim$-equivalent.
\\\smallskip

(b) $f = \gamma(a_i \mapsto f_i)_i,\ X = \gamma:\prod_i X_i$, and
$g = \ola{\gamma}[a_k](g')$. Composing with $g$ first gives
\[(\gamma(a_i \mapsto f_i)_{} ;_\gamma \ola{\gamma}[a_k](g')) 
;_\delta h \Lra (f_k ;_\gamma g') ;_\delta h
\]
whereas composing with $h$ first gives
\begin{align*}
(\gamma(a_i \mapsto f_i)_i ;_\delta h) ;_\gamma \ola{\gamma}[a_k](g')
&\Lra \gamma(a_i \mapsto f_i ;_\delta h)_i ;_\gamma \ola{\gamma}[a_k](g') \\
&\Lra (f_k ;_\delta h) ;_\gamma g' 
\end{align*}
which, by induction, are $\sim$-equivalent.
\\\smallskip

(c) $f = \gamma(a_i \mapsto f_i)_i$, $X = \gamma:\prod_i X_i$, and $g$
operates on a channel $\beta \neq \gamma$. In this case we must explore the
structure of $g$. First, let us assume that $g = \beta\{b_j \mapsto g_j\}_j$.
Composing with $g$ first gives 
\begin{align*}
(\gamma(a_i \mapsto f_i)_i ;_\gamma \beta\{b_j \mapsto g_j\}_j) ;_\delta h
&\Lra \beta\{b_j \mapsto \gamma(a_i \mapsto f_i)_i ;_\gamma g_j\}_j ;_\delta h
\\ &\Lra
\beta\{b_j \mapsto (\gamma(a_i \mapsto f_i)_i ;_\gamma g_j) ;_\delta h\}_j
\end{align*}
whereas composing with $h$ first gives
\begin{align*}
(\gamma(a_i \mapsto f_i)_i ;_\delta h) ;_\gamma \beta\{b_j \mapsto g_j\}_j
&\Lra \gamma(a_i \mapsto f_i ;_\delta h)_i ;_\gamma \beta\{b_j \mapsto g_j\}_j 
\\ &\Lra 
\beta\{b_j \mapsto \gamma(a_i \mapsto f_i ;_\delta h)_i ;_\gamma g_j\}_j
\end{align*}
which, by induction, are $\sim$-equivalent. The remaining cases on the
structure of $g$ are handled similarly.
\\\smallskip

(d) $f = \delta(a_i \mapsto f_i)_i$, $Y=\delta:\prod_i Y_i$ and
$h=\ola{\delta}[a_k](h')$. Composing with $g$ first gives
\begin{align*}
(\delta(a_i \mapsto f_i)_i ;_\gamma g) ;_\delta \ola{\delta}[a_k](h')
&\Lra \delta(a_i \mapsto f_i ;_\gamma g)_i ;_\delta \ola{\delta}[a_k](h') \\
&\Lra (f_k ;_\gamma g) ;_\delta h' 
\end{align*}
whereas composing with $h$ first gives
\[(\delta(a_i \mapsto f_i)_i ;_\delta \ola{\delta}[a_k](h')) ;_\gamma g
\Lra (f_k ;_\delta h') ;_\gamma g
\]
which, by induction, are $\sim$-equivalent.
\\\smallskip

(e) $f = \delta(a_i \mapsto f_i)_i$, $Y=\delta:\prod_i Y_i$ and
$h$ operates on a channel $\beta \neq \delta$. In this case we
must explore the structure of $h$. First, let us assume that $h = \beta\{b_j
\mapsto h_j\}_j$. Composing with $g$ first gives 
\begin{align*}
(\gamma(a_i \mapsto f_i)_i ;_\gamma g) ;_\delta \beta\{b_j \mapsto h_j\}_j
&\Lra \gamma(a_i \mapsto f_i ;_\gamma g)_i ;_\delta \beta\{b_j \mapsto h_j\}_j 
\\ &\Lra 
\beta\{b_j \mapsto \gamma(a_i \mapsto f_i ;_\gamma g)_i ;_\delta h_j\}_j
\end{align*}
whereas composing with $h$ first gives
\begin{align*}
(\gamma(a_i \mapsto f_i)_i ;_\delta \beta\{b_j \mapsto h_j\}_j) ;_\gamma g
&\Lra \beta\{b_j \mapsto \gamma(a_i \mapsto f_i)_i ;_\delta h_j\}_j ;_\gamma g
\\ &\Lra 
\beta\{b_j \mapsto (\gamma(a_i \mapsto f_i)_i ;_\delta h_j) ;_\gamma g\}_j
\end{align*}
which, by induction, are $\sim$-equivalent. The remaining cases on the
structure of $h$ are handled similarly.

Since any derivation of a sequent $\Gamma \vd \alpha:\prod_i W_i,\Delta$ may be 
equivalently written to end with tupling applied on $\alpha$, this suffices. 

\end{enumerate}

This concludes all the essential cases. If all of $f$, $g$, and $h$ are atomic
then composition satisfies the interchange property because it satisfies the
interchange property in the underlying polycategory. If some of $f$, $g$, and 
$h$ are atomic a quick check of the possibilities will show that one ends up 
with a case essentially like one of the cases above.

\end{proof}

This now shows that $\Sigma\Pi_\cat{A}$ is a polycategory proving
Theorem~\ref{thm-sppc}.

\section{Poly-sums and poly-products} \label{sec-sums}

In this section we show that $\Sigma\Pi_\cat{A}$ is the free polycategory
generated from (the polycategory) \cat{A} under (finite) poly-sums and
poly-products. We begin by defining sums and products in a polycategory.

In a polycategory \cat{A}, an object $\sum_{i \in I} X_i \in \cat{A}$ is
said to be the \textbf{poly-sum} (or \textbf{poly-coproduct}) of a family of
objects $X_i \in \cat{A}$, for $i \in I$, in case there is a poly-natural
correspondence 
\begin{equation} \label{natcorsum} \tag{$*$}
\vcenter{\infer={\Gamma,\alpha:\sum_{i \in I} X_i \vd_{\alpha\{f_i\}_i} \Delta}
{\{\Gamma, \alpha:X_i \vd_{f_i} \Delta\}_{i \in I}}}
\end{equation}
where by poly-natural we mean that the following two equivalences
\[h ;_\gamma \alpha\{f_i\}_i = \alpha\{h ;_\gamma f_i\}_i
\qquad \text{and} \qquad
\alpha\{f_i\}_i ;_\gamma h = \alpha\{f_i ;_\gamma h\}_i
\]
hold (when $\alpha \neq \gamma$). These equivalences assert that cutting on
an object and then forming the coproduct is the same as first forming the
coproduct and then performing the cut.

Products in polycategories are (as we expect) dual to coproducts. Explicitly,
an object $\prod_{i \in I} X_i \in \cat{A}$ is said to be the
\textbf{poly-product} of a family of objects $X_i \in \cat{A}$, for $i \in I$,
in case there is a poly-natural correspondence
\begin{equation} \label{natcorprod} \tag{$\star$}
\vcenter{\infer={\Gamma \vd_{\beta(f_i)_i} \beta:\prod_{i \in I} X_i,\Delta}
{\{\Gamma \vd_{f_i} \beta:X_i,\Delta\}_{i \in I}}}
\end{equation}

It will now be useful to define more ``standard'' injection and projection
maps. The reader may have noticed that the injections and projections seem
to be a little unfamiliar in their presentation.
\[
\infer[\text{(injection)}]{\Gamma \vd \sum_i X_i,\Delta}{\Gamma \vd X_k,\Delta}
\qqqquad
\infer[\text{(projection)}]{\Gamma,\prod_i X_i \vd \Delta}
    {\Gamma,X_k \vd \Delta}
\]

This may be ``remedied'' as follows. There are injection derivations
$X_k \vd_{b_k} \sum_i X_i$ for $k \in I$ given by $b_k = b_k(\iota_{X_k})$.
With these more ``standard'' injections, the general injection terms may be
identified with $f ; b_k$. Note this is a valid identification, since there
is a reduction of derivations
\[\vcenter{
\infer{\Gamma \vd \sum X_i,\Delta}{
  \Gamma \vd X_k,\Delta 
 &\infer{X_k \vd \sum X_i}{X_k \vd X_k}}}
\quad \Lra \quad
\vcenter{\infer{\Gamma \vd \sum X_i,\Delta}{\Gamma \vd X_k,\Delta}}
\]
Dually, the general projection terms may be identified with $p_k ; f$ where
$\prod_i X_i \vd_{p_k} X_k$.

\begin{remark}
The typing will sometimes be left off the injection/projection maps as there
is only one formula on each side of the turnstile and hence, it may be inferred
from the annotated composition symbol ``$;_\gamma$''. For example, if we write
$b_k\ ;_\gamma \Gamma, \gamma:\sum_i X_i \vd \Delta$ it is clear that we are
cutting on $\gamma:\sum_i X_i$.
\end{remark}

We will now work with the poly-coproduct (the dual observations hold for
the poly-product). The first thing to establish is the connection between the
poly-coproduct defined here and the standard way of viewing a coproduct. 

\begin{proposition}
For a polycategory \cat{A} the following are equivalent:
\begin{enumerate}[{\upshape (i)}]
\item \cat{A} has poly-coproducts (for a set $I$).
\item There is an object $\sum_i X_i$ for each family of objects $\{X_i\}_i$ 
which has injection maps $X_i \vd_{b_i} \sum_i X_i$ and a cotupling
operation producing a unique map $\Gamma,\gamma:\sum_i X_i
\vd_{\gamma\{f_i\}_i} \Delta$, where $\Gamma, \gamma:X_i \vd_{f_i} \Delta$,
such that
\begin{itemize}
\item $\{b_i\}_i = 1_{\sum_i X_i}$
\item $b_i\ ;_\gamma \gamma\{f_i\}_i = f_i$
\item $\alpha\{f_i\}_i\ ;_\gamma h = \alpha\{f_i\ ;_\gamma h\}_i$,
where $\alpha \neq \gamma$
\item $h\ ;_\gamma \alpha\{f_i\}_i = \alpha\{h\ ;_\gamma f_i\}_i$, 
where $\alpha \neq \gamma$
\end{itemize}
\end{enumerate}
\end{proposition}

\begin{proof}
To show that (i) implies (ii) consider the identity map,
\centro{$\sum\nolimits_i X_i \vd_{1_{\sum_i X_i}} \sum\nolimits_i X_i
$}
By the definition of poly-coproducts, this gives
\[
\infer={\sum_i X_i \vd_{1_{\sum_i X_i}} \sum_i X_i}
    {\{X_i \vd_{b_i} \sum_i X_i\}_{i \in I}}
\]
and hence, $\{b_i\}_i = 1_{\sum_i X_i}$.

Now consider the composite
\centro{$
X_k \vd_{b_k} \gamma:\sum\nolimits_i X_i \quad ;_\gamma \quad
\gamma:\sum\nolimits_i X_i, \Gamma  \vd_{\gamma:\{f_i\}} \Delta
$}
From the definition of poly-coproducts we have the following correspondence
\[
\infer={\sum_i X_i \vd_{1_{\sum_i X_i}} \gamma:\sum_i X_i \quad ;_\gamma \quad
\gamma:\sum_i X_i, \Gamma \vd_{\gamma\{f_i\}_{i \in I}} \Delta}
{\{X_i \vd_{b_i} \gamma:\sum_i X_i \quad ;_\gamma \quad
\gamma:\sum_i X_i, \Gamma \vd_{\gamma\{f_i\}_i} \Delta\}_{i \in I}}
\]
so that $b_i\ ;_\gamma \gamma\{f_i\}_i = f_i$.

Similarly, the following poly-natural correspondences
\[
\infer={\Gamma,\alpha:\sum_i X_i \vd_{\alpha\{f_i\}_{i \in I}} \Delta,\gamma:Z
\quad ;_\gamma \quad \gamma:Z,\Gamma' \vd_{h} \Delta'}
{\{\Gamma, \alpha:X_i \vd_{f_i} \Delta,\gamma:Z \quad ;_\gamma \quad
\gamma:Z,\Gamma' \vd_{h} \Delta'\}_{i \in I}}
\]
and
\[
\infer={\Gamma \vd_{h} \Delta,\gamma:Z \quad ;_\gamma \quad
\gamma:Z,\Gamma',\alpha:\sum_i X_i \vd_{\alpha\{f_i\}_{i \in I}} \Delta'}
{\{\Gamma \vd_h \Delta,\gamma:Z \quad ;_\gamma: \quad
\gamma:Z,\Gamma', \alpha:X_i \vd_{f_i} \Delta'\}_{i \in I}}
\]
show that $\alpha\{f_i\}_i\ ;_\gamma h = \alpha\{f_i\ ;_\gamma h\}_i$ and
$h\ ;_\gamma \alpha\{f_i\}_i = \alpha\{h\ ;_\gamma f_i\}_i$.

The argument for the uniqueness of the comparison map
\[
\xymatrix@=8ex{& \Delta \\
\Gamma, X_i \ar[ur]^{f_i} \ar[r]_{b_i} & \Gamma, \sum_i X_i
\ar@{-->}[u]_h}
\]
is $\{f_i\}_i = \{b_i\ ; h\}_i = \{b_i\}_i\ ; h = 1_{\sum_i X_i}\ ; h = h$.

For the implication from (ii) to (i) we may immediately conclude the
top-to-bottom direction of the correspondence from the definition of (ii).
The argument for the bottom-to-top direction is as follows:
\begin{align*}
\gamma\{f_i\}_i &= 1_{\sum_i X_i}\ ;_\gamma \gamma\{f_i\}_i \\
&= \alpha\{b_i\}_i\ ;_\gamma \gamma\{f_i\}_i \\
&= \alpha\{b_i\ ;_\gamma \gamma\{f_i\}_i\}_i \\
&= \alpha\{f_i\}_i
\end{align*}
Writing out the terms explicitly, this says that
\centro{$\Gamma, \gamma:\sum\nolimits_i X_i \vd_{\{f_i\}_i} \Delta =
\left\{\Gamma, \alpha:X_i \vd_{f_i} \Delta\right\}_i$}

This establishes both directions of the correspondence. That it is
poly-natural follows from the definition (the last two conditions).
\end{proof}

The next two propositions show that $\Sigma\Pi_\cat{A}$ is the free
polycategory generated from \cat{A} with finite sums and products.

\begin{proposition}
$\Sigma\Pi_\cat{A}$ has finite poly-sums and finite poly-products.
\end{proposition}

\begin{proof} 
In order to establish that $\Sigma\Pi_\cat{A}$ has finite poly-sums and 
finite poly-products we must show that the inferences
\[
\vcenter{\infer={\Gamma,\alpha:\sum\limits_{i \in I} i:X_i 
\vd_{\alpha\{i \mapsto f_i \}_{i \in I}} \Delta}
{\left\{\Gamma, \alpha:X_i \vd_{f_i} \Delta \right\}_{i \in I}}}
\qquad \text{and} \qquad
\vcenter{\infer={\Gamma \vd_{\beta(j \mapsto g_j)_{j \in J}} 
\beta:\prod\limits_{j \in J} j:Y_j, \Delta}
{\left\{\Gamma \vd_{g_i} \beta:Y_j, \Delta \right\}_{j \in J}}}
\]
are two-way and poly-natural.

We begin with coproducts.
Going from top-to-bottom is immediate via the cotupling derivation so we need
only prove the other direction. By Proposition~\ref{prop-insp} we know that
a sequent of the form $\Gamma,\sum_{i \in I} X_i \vd \Delta$ may be written in
a cut-free manner so that its principal rule is a cotupling which leaves us 
with the set of sequents $\left\{\Gamma, X_i \vd \Delta\right\}_{i \in I}$. 
That this correspondence respects the equivalence relation~$\sim$ is
immediate from our decision procedure as once the cotuple structure is made
principal, equality is determined by equality of the arguments. This proves
the bottom-to-top direction.

It remains to show that this correspondence is poly-natural. However, this
follows immediately from rewrite (3). Products are handled dually and
therefore this establishes that $\Sigma\Pi_\cat{A}$ has products and
coproducts.
\end{proof}

\begin{proposition}
$\Sigma\Pi_\cat{A}$ is the free polycategory generated from \cat{A} with
finite sums and products.
\end{proposition}

\begin{proof}
To show that $\Sigma\Pi_\cat{A}$ is the free polycategory generated from
\cat{A} with products and coproducts it suffices to show that all the 
identities (1)--(22) must hold in any polycategory with poly-products and 
poly-coproducts.

The identities (1) and (2) clearly hold. (3) and (9) (dually (4) and (10))
follow by poly-naturality, i.e.,
\[
\infer={\Gamma,\sum_i X_i \vd_{\{f_i\}_i} \Delta,Z\ ;\ Z,\Gamma' \vd_g \Delta'}
{\{\Gamma,X_i \vd_{f_i} \Delta,Z\ ;\ Z,\Gamma' \vd_g \Delta'\}_i}
\] 
(5) and (7) (dually (6) and (8)) follow from the interchange law, i.e.,
\[
\infer=
{(\Gamma \vd_f X_k,\Delta,Z\ ;\ Z,\Gamma' \vd_g \Delta')\ ;\ X_k \vd_{b_k}
\sum_i X_i}
{(\Gamma \vd_f X_k,\Delta,Z\ ;\ X_k \vd_{b_k} \sum_i X_i)\ ;\ Z,\Gamma' \vd_g
\Delta'}
\]
(11) and dually (12) follow from the associativity of cut
\[
\infer=
{\Gamma \vd_f \Delta,X_k,\ ;\ (X_k \vd_{b_k} \sum_i X_i\ ;\ \Gamma',\sum_i X_i
\vd_{\{g_i\}_i} \Delta')}
{(\Gamma \vd_f \Delta,X_k,\ ;\ X_k \vd_{b_k} \sum_i X_i)\ ;\ \Gamma',\sum_i X_i
\vd_{\{g_i\}_i} \Delta')}
\]
as $b_k\ ;\ \{g_i\}_i = {g_k}$.

(13) and dually (14) follow by the following argument (applications of
poly-naturality):
\begin{align*}
\{\{f_{ij}\}_i\}_j &= 1_{\sum_i X_i}\ ;\ \{\{f_{ij}\}_i\}_j \\
&= \{b_i\}_i\ ;\ \{\{f_{ij}\}_i\}_j \\
&= \{b_i\ ;\ \{\{f_{ij}\}_i\}_j\}_i \\
&= \{\{b_i\ ;\ \{f_{ij}\}_i\}_j\}_i \\
&= \{\{f_{ij}\}_j\}_i
\end{align*}

A similar proof shows that (19) holds. (15) and (17) (dually (16) and (18))
follow by poly-naturality. Explicitly, (15) is as follows:
\[
\infer={\Gamma,\sum_i X_i \vd_{\{f_i\}_i} \Delta,Y_k\ ;\
Y_k \vd_{b_k} \sum_j Y_j}
{\{\Gamma,X_i \vd_{f_i} \Delta,Y_k\ ;\ Y_k \vd_{b_k} \sum_j Y_j\}_i}
\] 
(20) and (21) follow by the interchange law $(f\ ;\ b_k)\ ;\ b_l = (f\ ;\
b_l)\ ;\ b_k$, while (22) follows by associativity $p_k\ ;\ (f\ ;\ b_k) =
(p_k\ ;\ f)\ ;\ b_k$.
\end{proof}

In the next section we show that the sums and products in $\Sigma\Pi_\cat{A}$
satisfy another important property.

\subsection{Softness of poly-sums and poly-products} \label{sec-softness}

The purpose of this section is to characterize the free sum and product
completion of a polycategory. We begin with a discussion of poly-hom-sets.

Given the types $\Gamma$ and $\Delta$, what does the poly-hom-set
$\Hom(\Gamma; \Delta)$ look like? If the domain $\Gamma$ contains a coproduct
$\sum_{i \in I} X_i$, or the codomain $\Delta$ contains a product
$\prod_{j \in J} Y_j$, then (by the inferences~\eqref{natcorsum}
and~\eqref{natcorprod}) we may ``break-up'' the poly-hom-sets
$\Hom(\Gamma,\sum_{i \in I} X_i\ ;\ \Delta)$ and
$\Hom(\Gamma\ ;\ \prod_{j \in J} Y_j,\Delta)$ into a product (in the
category of sets) of poly-hom-sets, respectively,
\[\prod_{i \in I} \Hom(\Gamma, X_i\ ;\ \Delta) \quad \text{and} \quad
\prod_{j \in J} \Hom(\Gamma\ ;\ Y_j,\Delta)
\]

In the case where there are only products or atoms in the domain and
coproducts or atoms in the codomain there is not in general a
description of the poly-hom-set. In the free case for (non-poly) categories
however there is a resolution due to Joyal~\cite[and other
references]{joyal95:free,joyal95:freeen}, which he derived from
Whitman's observation on lattices~\cite{whitman41:free}. Here we generalize
these ideas to polycategories.

Define the \textbf{output index} of a pair of types $(\Gamma, \Delta)$,
denoted $\oi(\Gamma ; \Delta)$, as:
\begin{itemize}
\item if $\alpha:\prod_{i \in I} X_i \in \Gamma$ then the
pairs $(\alpha,i)$, for $i \in I$, are in $\oi(\Gamma; \Delta)$.
\item if $\beta:\sum_{j \in J} Y_j \in \Delta$ then the pairs
$(\beta,j)$, for $j \in J$, are in $\oi(\Gamma; \Delta)$.
\end{itemize}
For example, given
\centro{$
\Gamma = \alpha:A,\ \beta:\prod_{i \in I} W_i,\ \gamma:\sum_{j \in J} X_j
\qquad \text{and} \qquad 
\Delta = \delta:\sum_{k \in K} Y_k,\ \epsilon:\prod_{l \in L} Z_l$}
where $A$ is atomic, the output index $\oi(\Gamma;\Delta) =
\{(\beta,i),(\delta,k) \mid i \in I, k \in K\}$.

Given a product $\alpha:\prod_i X_i$ on a domain channel and $(\alpha,k)$,
$k \in I$, we will denote by
\centro{$[\Hom(\Gamma, \alpha:\prod\nolimits_i X_i\ ;\ \Delta)]_{(\alpha,k)}
$}
the poly-hom-set $\Hom(\Gamma, \alpha:X_k\ ;\ \Delta)$. The same notation will
be used for a sum on a codomain channel. The following
\[\sum_{(\alpha,i) \in \oi(\Gamma ; \Delta)}
\negthickspace\negthickspace\negthickspace
[\Hom(\Gamma\ ;\ \Delta)]_{(\alpha,i)}
\]
will then be used to indicate the disjoint union of poly-hom-sets where each
product in the domain and coproduct in the codomain has been ``broken down''
into its composite elements. For example, applying this construction to
$\Gamma$ and $\Delta$ as defined above results in the following set of
poly-hom-sets:
\centro{$\left\{\Hom(A,W_i,\sum_{j \in J} X_j\ ;\ Y_k,\prod_{l \in L} Z_l)
\mid i \in I, k \in K \right\}$}

Suppose \cat{A} is a polycategory and \cat{B} is a polycategory with sums
and products such that $\cat{A} \to^{\mathcal{I}} \cat{B}$ is an inclusion
morphism of polycategories. Consider a poly-hom-set $\Hom(\Gamma; \Delta)$
in \cat{B} such that the output index $\oi(\Gamma;\Delta)$ is non-empty and
any channel not in the output index is of the form $\mathcal{I}(A)$ where
$A \in \cat{A}$. That is, each object of $\Gamma$ is either ``atomic'' or a
product, and each object of $\Delta$ is either ``atomic'' or a sum, so that
this is the polycategorical analogue of an output sequent (which we call an
\textbf{output poly-hom-set}). The morphism $\mathcal{I}:\cat{A} \ra \cat{B}$
is called a \textbf{semi-soft extension} of polycategories if for any output
poly-hom-set $\Hom(\Gamma; \Delta) \in \cat{B}$ the map
\[\sum_{(\alpha,i) \in \oi(\Gamma;\Delta)}
\negthickspace\negthickspace\negthickspace
[\Hom(\Gamma; \Delta)]_{(\alpha,i)}
\to^{\{\ol{\alpha}[i](\_)\}_{(\alpha,i)}} \Hom(\Gamma; \Delta)
\]
is a surjection, where $\{\ol{\alpha}[i](\_)\}_{(\alpha,i)}$ is the cotupling
map (in \category{Set}) and the underscore in $\ol{\alpha}[i](\_)$ is used
to represent any map in $[\Hom(\Gamma; \Delta)]_{(\alpha,i)}$. (The overline
$\ol{\alpha}$ is used to indicate that the map $\ol{\alpha}[i](\_)$ may
represent an injection or a projection; it will depend upon the type of the
$\alpha$ channel.)

In terms of processes, a semi-soft extension guarantees that every process in
an output poly-hom-set $\Hom(\Gamma;\Delta)$ has as its next action an output
event.

\begin{example}[Semi-soft extension]
Consider the simple output poly-hom-set $\Hom(\alpha:A \times B\ ;\ \beta:C+D)$.
Dropping the channel names we have:

\[\sum_{(\gamma,k) \in \oi(A \times B ; C+D)}
\negthickspace\negthickspace\negthickspace\negthickspace
\negthickspace\negthickspace\negthickspace\negthickspace
\Hom(A \times B\ ;\ C+D)
= \left\{
\begin{array}{l}
\Hom(A\ ;\ C+D) \smallskip\\
\Hom(B\ ;\ C+D) \smallskip\\
\Hom(A \times B\ ;\ C) \smallskip\\
\Hom(A \times B\ ;\ D)
\end{array} \right\}
\]

The (injection and projection) maps are
\[\begin{array}{rcl}
\ola{\alpha}[1](\_) &:& \Hom(A\ ;\ C+D) \ra \Hom(A \times B\ ;\ C+D)
\medskip\\
\ola{\alpha}[2](\_) &:& \Hom(B\ ;\  C+D) \ra \Hom(A \times B\ ;\  C+D)
\medskip\\
\ora{\beta}[1](\_) &:& \Hom(A \times B\ ;\  C) \ra \Hom(A \times B\ ;\  C+D)
\medskip\\
\ora{\beta}[2](\_) &:& \Hom(A \times B\ ;\  D) \ra \Hom(A \times B\ ;\  C+D)
\end{array}
\]
so that considering them as a whole gives
\[
\{\ol{\gamma}[k](\_)\}_{(\gamma,k)}:
\sum_{(\gamma,k)}
\Hom(A \times B; C+D) \ra \Hom(A \times B; C+D)
\]
where $(\gamma,k) \in \oi(A \times B,C+D)$.
\end{example}

The inclusion $\mathcal{I}:\cat{A} \ra \cat{B}$ is called a \textbf{soft
extension} in case for any output poly-hom-set $\Hom(\Gamma;\Delta)$ the
following diagram 
\begin{equation} \label{dia-coeq} \tag{$\bigstar$}
\vcenter{
\xymatrix@R=8ex{{\sum\limits_{
((\alpha,i), (\beta,j)) \in \oi(\Gamma,\Delta)}
\negthickspace\negthickspace\negthickspace\negthickspace
\negthickspace\negthickspace\negthickspace\negthickspace
[\Hom(\Gamma ; \Delta)]_{((\alpha,i),(\beta,j))}}
\ar@<2ex>[d]^-{\sum\limits_{(\beta,j)}\{\ol{\alpha}[i](\_)\}_{((\alpha,i),(\beta,j))}}
\ar@<-2ex>[d]_-{\sum\limits_{(\alpha,i)}\{\ol{\beta}[j](\_)\}_{((\alpha,i),(\beta,j))}}
\\
\sum\limits_{(\gamma,k) \in \oi(\Gamma,\Delta)}
\negthickspace\negthickspace\negthickspace\negthickspace
[\Hom(\Gamma ; \Delta)]_{(\gamma,k)}
\ar[d]^-{\{\ol{\gamma}[k](\_)\}_{(\gamma,k)}}
\\
\Hom(\Gamma ; \Delta)
}}
\end{equation}
is a coequalizer diagram in the category of sets, where
$((\alpha,i),(\beta,j))$ for $(\alpha,i),\ (\beta,j) \in \oi(\Gamma; \Delta)$,
$\alpha \neq \beta$, represents choosing pairs from the members of
$\oi(\Gamma; \Delta)$ and the map $\sum_{(\beta,j)} \{\ol{\alpha}[i](\_)\}_{
((\alpha,i),(\beta,j))}$ represents the coproduct (of the $\beta$ components)
of the cotupling maps of the $\alpha$ components, and vise versa when
$\alpha$ and $\beta$ are switched.

In terms of processes a soft extension means that given a process in
$\Hom(\Gamma; \Delta)$ which is able to output on two distinct
channels, the order in which one chooses to do the outputs is irrelevant, e.g.,
the maps $\ola{\alpha}(\ora{\beta}(f))$ and $\ora{\beta}(\ola{\alpha}(f))$ 
will be $\sim$-equivalent.

\begin{remark} \label{remark-coeq}
In \Set, given two parallel arrows $f,g: A \two B$, the coequalizer of this
pair always exists (see, e.g.,~\cite{barr99:category}). Thus, in order
to show that the above diagram~\eqref{dia-coeq} is a coequalizer diagram we
may assume that the coequalizer is given by $(q,Q)$ and then show that the
unique map $h:Q \ra \Hom(\Gamma;\Delta)$ is an isomorphism, i.e., that the
following diagram commutes:
\[\xymatrix@R=14ex@C=15ex@!0{
{\sum\limits_{((\alpha,i), (\beta,j)) \in \oi(\Gamma,\Delta)}
\negthickspace\negthickspace\negthickspace\negthickspace
\negthickspace\negthickspace\negthickspace\negthickspace
[\Hom(\Gamma ; \Delta)]_{((\alpha,i),(\beta,j))}}
\ar@<2ex>[d]^-{\sum\limits_{(\beta,j)}\{\ol{\alpha}[i](\_)\}_{((\alpha,i),(\beta,j))}}
\ar@<-2ex>[d]_-{\sum\limits_{(\alpha,i)}\{\ol{\beta}[j](\_)\}_{((\alpha,i),(\beta,j))}}
\\
\sum\limits_{(\gamma,k) \in \oi(\Gamma,\Delta)}
\negthickspace\negthickspace\negthickspace\negthickspace
[\Hom(\Gamma ; \Delta)]_{(\gamma,k)}
\ar[d]_-{q} \ar[dr]^-{\{\ol{\gamma}[k](\_)\}_{(\gamma,k)}}
\\
Q \ar@{-->}[r]_-h^-\sim& \Hom(\Gamma ; \Delta)
}
\]
\end{remark}

\begin{example}[Soft extension]
Consider the poly-hom-set $\Hom(\alpha:A \times B\ ;\ \beta:C+D)$
from the example above and the following coequalizer diagram (again dropping
the channel names):

\[\xymatrix@R=5ex@C=1ex{
D_1 = \left\{\txt{$\Hom(A\ ;\ C)\ \ \Hom(B\ ;\ C)$ \\
$\Hom(A\ ;\ D)\ \ \Hom(B\ ;\ D)$} \right\} 
\ar@<-2ex>[d]_-{\sum_\beta \ola{\alpha}}
\ar@<2ex>[d]^-{\sum_\alpha \ora{\beta}} \\
D_2 = \left\{\txt{$\Hom(A \times B\ ;\ C)\ \ \Hom(A\ ;\ C+D)$ \\
            $\Hom(A \times B\ ;\ D)\ \ \Hom(B\ ;\ C+D)$} \right\}
\ar[d]^-{\ol{\gamma}} \\
\Hom(A \times B\ ;\ C+D)}
\]
where $\sum_\beta \ola{\alpha}$ and $\sum_\alpha \ora{\beta}$ will be
constructed in what follows and $\ol{\gamma}$ is the map from the example
above.

In $D_1$ all the ``splittings'' of pairs has been done. Fixing the first
component of $\beta$, the pairs $((\alpha,1),(\beta,1))$ and
$((\alpha,2),(\beta,1))$ pick out the poly-hom-sets $\Hom(A\ ;\ C)$ and
$\Hom(B\ ;\ C)$ respectively; a map from each of these is:
\[\begin{array}{rcl}
\ola{\alpha}[1](\_)_{(\beta,1)} &:& \Hom(A\ ;\ C) \ra \Hom(A \times B\ ;\ C)
\medskip\\
\ola{\alpha}[2](\_)_{(\beta,1)} &:& \Hom(B\ ;\ C) \ra \Hom(A \times B\ ;\ C)
\end{array}
\]
where the subscripted $(\beta,1)$ is used to indicate that the first component
of $\beta$ is fixed. The cotupling map of both of these maps gives a map
\[\left\{\ola{\alpha}[1](\_)_{(\beta,1)},\ \ola{\alpha}[2](\_)_{(\beta,1)}
\right\}: D_1 \ra D_2
\]

Similarly, fixing the second component of $\beta$ gives the maps
\[\begin{array}{rcl}
\ola{\alpha}[1](\_)_{(\beta,2)} &:& \Hom(A\ ;\ D) \ra \Hom(A \times B\ ;\ D)
\medskip\\
\ola{\alpha}[2](\_)_{(\beta,2)} &:& \Hom(B\ ;\ D) \ra \Hom(A \times B\ ;\ D)
\end{array}
\]
and thus the cotupling map,
\[\left\{\ola{\alpha}[1](\_)_{(\beta,2)},\ \ola{\alpha}[2](\_)_{(\beta,2)}
\right\}: D_1 \ra D_2
\]

Together these (cotupling) maps give the following map
\[\sum_{j \in \{1,2\}}
\left\{\ola{\alpha}[1](\_)_{(\beta,j)},\ \ola{\alpha}[2](\_)_{(\beta,j)}
\right\}: D_1 \ra D_2
\]

In the same way a second map may be constructed where $\alpha$ is fixed
instead:
\[\sum_{i \in \{1,2\}}
\left\{\ora{\beta}[1](\_)_{(\alpha,i)},\ \ora{\beta}[2](\_)_{(\alpha,i)}
\right\}:D_1 \ra D_2
\]
\end{example}

\begin{definition}
An inclusion $\mathcal{I}:\cat{A} \to \cat{B}$ which is both a semi-soft
and soft extension is said to be \textbf{soft}.
\end{definition}

\begin{lemma}
The inclusion morphism of polycategories $\cat{A} \to^{\mathcal{I}}
\Sigma\Pi_\cat{A}$ is soft.
\end{lemma}

\begin{proof}
That the inclusion $\mathcal{I}$ is semi-soft follows by
Proposition~\ref{prop-insp} and which (non-cut) rules we may apply. 

As mentioned in Remark~\ref{remark-coeq}, to show that it is a soft extension
we will assume that the diagram~\eqref{dia-coeq} has a coequalizer $(q,Q)$, and
show that $Q$ is isomorphic to $\Hom(\Gamma,\Delta)$. As we are in \Set, a
bijective correspondence between the two objects suffices to establish the
isomorphism.

As $Q$ is the coequalizer, there exists a unique map $h:Q \ra
\Hom(\Gamma;\Delta)$ such that $q;h = \{\ol{\gamma}[k](\_)\}_{(\gamma,k)}$.
But $\{\ol{\gamma}[k](\_)\}_{(\gamma,k)}$ is a surjection and hence, so is
$h$. It is left to establish that $h$ is an injection.

Suppose that $\ol{\alpha}[i](f)$ and $\ol{\beta}[j](g)$ are equivalent
morphisms in $\Hom(\Gamma; \Delta)$. From the decision procedure for
$\Sigma\Pi$-morphisms we know that the equality of morphisms is determined by
the equality of their arguments so that $\ol{\alpha}[i](f)$ and
$\ol{\beta}[j](g)$ are equivalent to morphisms of the form
$\ol{\alpha}[i](\ol{\beta}[j](f'))$ and $\ol{\beta}[j](\ol{\alpha}[i](g'))$
respectively. This implies that $f'$ and $g'$ are equivalent in
$\sum\limits_{((\alpha,i),(\beta,j))}
\negthickspace\negthickspace
[\Hom(\Gamma; \Delta)]_{((\alpha,i),(\beta,j))}$,
and therefore, must be coequalized in $Q$ establishing that $h$ is injective.

Thus, the map $h$ is a bijection and $Q \iso \Hom(\Gamma; \Delta)$,
establishing that
\[\left(\{\ol{\gamma}[k](\_)\}_{(\gamma,k)},\ \Hom(\Gamma; \Delta) \right)
\]
is the coequalizer.
\end{proof}

Let $\cat{A} \to^{\mathcal{I}} \Sigma\Pi_\cat{A}$ be an inclusion morphism
of polycategories and $\cat{A} \to^F \cat{B}$ be a full inclusion
morphism of polycategories (injective on objects and an isomorphism on each
poly-hom-set), where the objects of \cat{B} are generated from the objects
of \cat{A} under finite sums and products. That $\Sigma\Pi_\cat{A}$ is free
guarantees that there is a unique comparison morphism of polycategories
$\Sigma\Pi_\cat{A} \to^{F^*} \cat{B}$ such that the following diagram commutes
\[\xymatrix@=7ex{\cat{A} \ar[r]^{\mathcal{I}} \ar[dr]_F &
   \Sigma\Pi_\cat{A} \ar@{-->}[d]^{F^*} \\
& \cat{B}}
\]

Given this data, the ``Whitman theorem'' which characterizes the free
polycat\-e\-gory with sums and products is as follows.

\begin{theorem}[Whitman theorem] \label{thm-whitman}
For $F$ and $F^*$ as above:
\begin{enumerate}[{\upshape (i)}]
\item If $F$ is a semi-soft extension then $F^*$ is full.
\item If $F$ is soft then $F^*$ is an equivalence of polycategories.
\end{enumerate}
\end{theorem}

\begin{proof}
\begin{enumerate}[{\upshape (i)}]
\item That $F^*$ is full follows from the following induction. Suppose
$f:F^*(\Gamma) \ra F^*(\Delta) \in \cat{B}$ (where $F^*(X_1,\ldots, X_n) =
F^*(X_1),\ldots,F^*(X_n)$). If $f$ is atomic this implies that $f \in
\cat{A}$, and thus
also in $\Sigma\Pi_\cat{A}$. So suppose that $f$ is not atomic. The point will
be to show that $f$ may be ``decomposed'' into a ``$\Sigma\Pi$-word'' of
simpler (smaller) functions which, by induction, are in the image of $F^*$.
Then, as $F^*$ preserves sums and products the same $\Sigma\Pi$-word may be
used to place $f$ in the image of $F^*$.

To see that $f$ may be decomposed consider the form of $F^*(\Gamma)$ and
$F^*(\Delta)$. If $F^*(\Gamma)$ contains a sum or $F^*(\Delta)$ contains a
product, then by the equivalences~\eqref{natcorsum} and~\eqref{natcorprod},
$f$ may be decomposed on the sum or the product respectively; if neither of
these is the case then, as $F$ is a semi-soft extension, $f$ must have been
the result of an injection or a projection map. In either case, we are able to
decompose $f$ into smaller functions.

\item To show that $F^*$ is an equivalence of polycategories it suffices to
show that $F^*$ is full, faithful, and essentially surjective, where by
essentially surjective we mean that each object in \cat{B} is isomorphic
to $F^*(X)$ for some $X$ in $\Sigma\Pi_\cat{A}$. If $F$ is soft, then by
definition it is a semi-soft extension, and thus, $F^*$ is full.

To show that $F^*$ is faithful, consider a parallel pair of arrows $F^*(f),
F^*(g): F^*(\Gamma) \two F^*(\Delta)$ such that $F^*(f) = F^*(g)$; as above
we may decompose $F^*(f)$ and $F^*(g)$ into $\Sigma\Pi$-words. Again, the
decomposition is obvious if either $F^*(\Gamma)$ contains a sum or
$F^*(\Delta)$ contains a product, so suppose not, and moreover, suppose
that we are able to decompose $F^*(f)$ and $F^*(g)$ into different
$\Sigma\Pi$-words (via injection or projection maps). Then, by softness,
these words may again be decomposed so that they are equivalent. Thus, in
both cases, $F^*(f)$ and $F^*(g)$ may be decomposed into different
substitution instances of the same $\Sigma\Pi$-word. As $F^*$ is full,
these words involve simpler functions of the form $F^*(h)$. By induction the
corresponding subterms $h$ are equal, and thus $f=g$.

To see that $F^*$ is essentially surjective, notice that since $\cat{B}$ is
also generated from \cat{A} under finite sums and products, each object of
\cat{B} must be isomorphic to an object in the image of $F^*$.
\end{enumerate}
\end{proof}

\chapter{Process Semantics} \label{chap-semantics}

The term calculus for $\Sigma\Pi$ used in the previous chapters was
motivated by the interpretation of $\Sigma\Pi$-derivations as processes.
The purpose of this chapter is to explain precisely how these terms
correspond to processes. To this end a semantics for processes is
introduced. This semantics could alternatively have been described as a
game theoretic semantics, and indeed many of the ideas and terminology
are derived from that view (see, e.g., \cite{abramsky98:game,
abramsky00:full, hyland00:full}).

\section{Behaviours} \label{sec-behav}

\subsection{Legal transitions}

A protocol may be in an atomic state, an output state, or an input state.
These states may be assigned one of three different \textbf{roles} which
will be used to determine the legality of transitions:

\begin{itemize}
\item \textbf{Source roles}: these are denoted by superscripting a
$\mathbf{0}$.
\item \textbf{Sink roles}: these are denoted by superscripting a
$\mathbf{1}$.
\item \textbf{Flow roles}: these are denoted by superscripting a $+$.
\end{itemize}

Given a protocol on a domain channel we may calculate its role inductively 
as follows:

\begin{itemize}
\item any atomic formula $A$ has a flow role: $A^+$;
\item the empty coproduct has a source role: $\{\,\}^\mathbf{0}$;
\item the empty product has a sink role: $(\,)^\mathbf{1}$;
\item the role $x$ of a coproduct, $\{a_i:A_i^{x_i}\}_i^x$ is determined by 
the roles $x_i$ of its subformula as follows:
\[x =
\begin{cases}
\mathbf{0} & \text{if and only if } (\forall i)\; x_i = \mathbf{0} \\
\mathbf{1} & \text{if and only if } (\exists !k)\; x_k = \mathbf{1} 
\text{ and } (\forall i \neq k)\; x_i = \mathbf{0}  \\
+ & \text{otherwise}
\end{cases}
\]

\item the role $x$ of a product, $(a_i:A_i^{x_i})_i^x$ is determined by the 
roles $x_i$ of its subformula as follows:
\[x =
\begin{cases}
\mathbf{0} & \text{if and only if } (\exists !k)\; x_k = \mathbf{0} 
\text{ and } (\forall i \neq k)\; x_i = \mathbf{1}  \\
\mathbf{1} & \text{if and only if } (\forall i)\; x_i = \mathbf{1} \\
+ & \text{otherwise}
\end{cases}
\]

\end{itemize}

The role of a protocol on a codomain channel is calculated dually, i.e., swap 
the \textbf{0} and the \textbf{1} while leaving $+$ alone. It is easy to 
see that any protocol which contains an atomic formula will have a flow
role.

It should be mentioned that the purpose of the role of a protocol is to
identify the initial (source), final (sink), and other (flow) protocols,
in the usual categorical sense.

\begin{example}[Roles] \label{rolecalc-example}
Some example protocols (on domain channels) with their roles calculated.
\begin{enumerate}

\item
\[\xymatrix@R=4ex@C=3ex@M=0ex{ 
&&& \bullet^+ \ar@{-}[dll] \ar@{-}[drr] \\ & 
\bullet^{\bf 0} \ar@{-}[dl] \ar@{-}[dr] &&&& \circ^+ \ar@{-}[dl] \ar@{-}[dr] \\ 
\bullet^{\bf 1} && \circ^{\bf 0} && \bullet^{\bf 1} && \bullet^{\bf 1}}
\]

\item
\[\xymatrix@R=4ex@C=5ex@M=0ex{
&& \bullet^+ \ar@{-}[dl] \ar@{-}[dr] \\
& \circ^+ \ar@{-}[dl] \ar@{-}[d] \ar@{-}[dr] && \circ^{\bf 0} \\ 
\bullet^{\bf 1} & \bullet^{\bf 1} & \circ^{\bf 0}}
\]

\item 
\[\xymatrix@R=4ex@C=4ex@M=0ex{
&&& \circ^+ \ar@{-}[dlll] \ar@{-}[d] \ar@{-}[drrr] \\
\circ^\mathbf{1} \ar@{-}[d] &&& \bullet^\mathbf{0} \ar@{-}[dl] \ar@{-}[dr] &&&
    \circ^+ \ar@{-}[dl] \ar@{-}[dr] \\
\bullet^\mathbf{1} && \circ^\mathbf{0} && \bullet^\mathbf{1} & 
\bullet^\mathbf{1} && \bullet^\mathbf{1}}
\]

\end{enumerate}
\end{example}

A protocol transition is \textbf{legal} (is a \textbf{legal transition}) in
case it starts at a state in a flow role and is either
\begin{itemize}
\item an output transition which ends at a state which does not have a 
sink role, or
\item an input transition which ends at a state which does not have a 
source role.
\end{itemize}
Graphically, a transition in the domain is legal in case it is either
\[\vcenter{\xymatrix@R=6ex{s^+ \ar@{-}[d] \\ s^+}}
\text{\qquad or \qquad}
\vcenter{\xymatrix@R=6ex{\circ^+ \ar@{-}[d] \\ s^{\bf 1}}}
\text{\qquad or \qquad}
\vcenter{\xymatrix@R=6ex{\bullet^+ \ar@{-}[d] \\ s^{\bf 0}}}
\]
where $s$ represents a product or a coproduct state. In the codomain
a transition is legal in case it is either
\[\vcenter{\xymatrix@R=6ex{s^+ \ar@{-}[d] \\ s^+}}
\text{\qquad or \qquad}
\vcenter{\xymatrix@R=6ex{\circ^+ \ar@{-}[d] \\ s^{\bf 0}}}
\text{\qquad or \qquad}
\vcenter{\xymatrix@R=6ex{\bullet^+ \ar@{-}[d] \\ s^{\bf 1}}}
\]

\begin{example}[Legal transitions]
Using the protocols from the above example (Example~\ref{rolecalc-example}),
we indicate the legal transitions as solid lines and the non-legal
transitions as dashed lines. 

\begin{enumerate}
\item
\[\xymatrix@R=4ex@C=3ex@M=0ex{ 
&&& \bullet^+ \ar@{-}[dll] \ar@{-}[drr] \\ & 
\bullet^{\mathbf{0}} \ar@{--}[dl] \ar@{--}[dr] &&&& \circ^+ \ar@{-}[dl]
\ar@{-}[dr] \\ 
\bullet^{\mathbf{1}} && \circ^{\mathbf{0}} && \bullet^{\mathbf{1}} &&
\bullet^{\mathbf{1}}}
\]

\item
\[\xymatrix@R=4ex@C=5ex@M=0ex{
&& \bullet^+ \ar@{-}[dl] \ar@{-}[dr] \\
& \circ^+ \ar@{-}[dl] \ar@{-}[d] \ar@{--}[dr] && \circ^{\bf 0} \\ 
\bullet^{\bf 1} & \bullet^{\bf 1} & \circ^{\bf 0}}
\]

\item
\[\xymatrix@R=4ex@C=4ex@M=0ex{
&&& \circ^+ \ar@{-}[dlll] \ar@{--}[d] \ar@{-}[drrr] \\
\circ^{\bf 1} \ar@{--}[d] &&& \bullet^{\bf 0} \ar@{--}[dl] \ar@{--}[dr] &&& 
\circ^+ \ar@{-}[dl] \ar@{-}[dr] \\
\bullet^{\bf 1} && \circ^{\bf 0} && \bullet^{\bf 1} & \bullet^{\bf 1} && 
\bullet^{\bf 1}}
\]

\end{enumerate}
\end{example}

\begin{lemma} \label{lem-interaction}
For any protocol $X$ in an input flow state there is a sequence of legal input
events which drives the protocol into either an atomic, output or sink state.
\end{lemma}

\begin{proof}
We prove this for a protocol $X$ on a domain channel; duality covers the
case where $X$ is on a codomain channel. 

For $X$ on a domain channel, to be in an input flow state it must be of the
form $X = \{a_i:X_i\}_{i \in I}$ where $I \neq \emptyset$. As $X$ has a flow
role, this means that not every $X_i$ has a source role, i.e., there is at
least one $k \in I$ such that $X_k^{\textbf{1}}$ or $X_k^+$. Thus, input $a_k$
to $X$. In the former case we are done as $X_k^{\textbf{1}}$ has a sink role.
In the latter case $X_k^+$ has a flow role which may be either an atomic,
output or input state; if it is an atomic or output state then the desired
conclusion has be reached, and otherwise, it is an input state which we may
inductively assume has the desired property.
\end{proof}

\subsection{Legal behaviours}

Suppose a channel $\alpha$ is assigned a protocol $X$. A \textbf{legal
channel behaviour} (or simply \textbf{channel behaviour}) for $\alpha$,
denoted $\gp_\alpha$, is a finite sequence of legal events on $\alpha$
satisfying its assigned protocol $X$ (see Section~\ref{sec-fap}). For
example, the behaviour $\<a,\ol{b}\>$ for a channel $\alpha$ represents
that on $\alpha$ the input event $a$ was received, after which, the event
$\ol{b}$ was output.

An (input or output) event $e$ may be appended to a channel behaviour
$\gp_\alpha$, denoted $\gp_\alpha * e$. Similarly, an (input or output)
event $e'$ may be prepended to $\gp_\alpha$, which is denoted
$e' * \gp_\alpha$. For example, $\<a,\ol{b}\> * c = \<a,\ol{b},c\>$
and $c * \<a,\ol{b}\> = \<c,a,\ol{b}\>$. In this manner (finite) sequences
of events may be appended or prepended (or both) to a channel behaviour.

If $\gp_\alpha$ and $\gq_\alpha$ are two channel behaviours on $\alpha$, we 
say that $\gp_\alpha$ is a \textbf{prefix} of $\gq_\alpha$, denoted 
$\gp_\alpha \prefix \gq_\alpha$, if there exists a sequence of events 
$\gp'_\alpha$ such that $\gp_\alpha * \gp'_\alpha = \gq_\alpha$. It is a 
\textbf{proper prefix}, denoted $\gp_\alpha \sqsubset \gq_\alpha$, if 
$\gp_\alpha \prefix \gq_\alpha$ and $\gp_\alpha \neq \gq_\alpha$. The
notation $\gp \prefix_i \gq$ and $\gp \prefix_o \gq$ is used to indicate
that $\gp$ is a prefix of $\gq$ separated only by input or output events
respectively.

A \textbf{legal behaviour} (or simply \textbf{behaviour}) $\gp =
(\gp_{\alpha_1},\ldots,\gp_{\alpha_n})$ is a tuple of channel
behaviours, one for each channel along which the process interacts. We will
typically make use of \textbf{tables} to denote behaviours. For example, if
$\<a,\ol{b}\>_\alpha$, $\<c,\ol{d},\ol{e}\>_\beta$, $\<f\>_\gamma$, and
$\<g,h\>_\delta$ are channel behaviours for $\alpha$, $\beta$, $\gamma$, and
$\delta$ respectively, where the domain consists of $\alpha$ and $\beta$ and
the codomain consists of $\gamma$ and $\delta$, the behaviour consisting of
these channel behaviours may be denoted

\[\begin{array}{|c|c||c|c|} 
       & \ol{e} &        &   \\
\ol{b} & \ol{d} &        & h \\
a      & c      & f      & g \\
\hline
\alpha & \beta  & \gamma & \delta \\
\hline
\end{array}
\]

The \textbf{dual behaviour} $\gp^*$ of $\gp$ is the behaviour $\gp$
in which input events are considered as output events and output events are
considered as input events. That is, any input event in $\gp$ is an output
event in $\gp^*$, and similarly any output event in $\gp$ is an input event in
$\gp^*$. For example, the dual of the behaviour above is:

\[\begin{array}{|c|c||c|c|} 
       & e &    & \\
b      & d &    & \ol{h} \\
\ol{a} & \ol{c} & \ol{f} & \ol{g}  \\
\hline
\alpha & \beta & \gamma & \delta \\
\hline
\end{array}
\]

The last states of a behaviour is called its \textbf{frontier}. For example,
given the protocols
\[\alpha:\left(\vcenter{
\xymatrix@R=4ex@C=0.8ex@M=0.2ex{ &&& \circ \ar@{-}[dll]_{a} \ar@{-}[drr]^{b} \\
& \bullet \ar@{-}[dl]_{c} \ar@{-}[dr]^{d} &&&& \bullet \ar@{-}[dl]_{e}
\ar@{-}[dr]^{f}
\\ A && B && C && D}}\right)
\
\beta:\left(\vcenter{
\xymatrix@R=4ex@C=0.8ex@M=0.2ex{& \bullet \ar@{-}[dl]_{g} \ar@{-}[dr]^{h} \\
E && F}}\right)
\
\gamma:\left(\vcenter{
\xymatrix@R=4ex@C=0.8ex@M=0.2ex{ & \circ \ar@{-}[dl]_{i} \ar@{-}[dr]^{j} \\
G && H}}\right)
\
\delta:\left(\vcenter{
\xymatrix@R=4ex@C=0.8ex@M=0.2ex{& \bullet \ar@{-}[dl]_{k} \ar@{-}[dr]^{m} \\
I && \circ \ar@{-}[dl]_{n} \ar@{-}[dr]^{f'} \\
& J && K}}\right)
\]
and the behaviour
\[\begin{array}{|c|c||c|c|}
&&& \\
\ol{c} &&& \\
a & & \ol{j} & k \\
\hline
\alpha & \beta & \gamma & \delta \\
\hline
\end{array}
\]
The frontier consists of $\alpha:A,\ \beta:(g:E,h:F),\ \gamma:H$, and
$\delta:I$. Notice that the frontier is actually representing a sequent: it
is what is yet to be proved. In the above example, the sequent the frontier
represents is:
\[
\alpha:A,\ \beta:(g:E,h:F) \vd \gamma:H,\ \delta:I
\]

Given any behaviour $\gp$ in which $\alpha$ is a channel of this behaviour,
we may form the \textbf{restriction} of $\gp$ to $\alpha$, denoted
$\gp_\alpha$. For example,

\[\gp =
\begin{array}{|c|c||c|c|} 
       & \ol{e} &        &   \\
\ol{b} & \ol{d} &        & h \\
a      & c      & f      & g \\
\hline
\alpha & \beta  & \gamma & \delta \\
\hline
\end{array}
\qqqquad 
\gp_\alpha =
\begin{array}{|c|}
\ \\
\ol{b} \\
a \\
\hline
\alpha \\
\hline
\end{array}
\]

The prefix and proper prefix of behaviours are defined in the
obvious way, i.e., $\gp \prefix \gq$ if for each channel $\alpha$ of $\gp$
we have $\gp_\alpha \prefix \gq_\alpha$, and $\gp \sqsubset \gq$ if $\gp
\prefix \gq$ and $\gp \neq \gq$. Appending and prepending events is similar
to the case for channel behaviours except that the channel to append/prepend on
must be specified, e.g., $\gp * \beta[a]$ is interpreted as
$(\gp_{\alpha_1},\ldots,\gp_\beta * \beta[a],\ldots,\gp_{\alpha_n})$ and
$\beta[b] * \gp$ is interpreted as
$(\gp_{\alpha_1},\ldots,\beta[a] * \gp_\beta,\ldots,\gp_{\alpha_n})$.

A behaviour in which all the states at the frontier have either flow ($+$) or
sink ($\mathbf{1}$) roles is called an \textbf{antecedent behaviour}. An
antecedent behaviour is called \textbf{saturated} if for any state at the
frontier, it is either atomic, has a sink role, or has an output flow role.
If all the states at the frontier are atomic, we may refer to this behaviour
as being \textbf{atomic saturated}.

Two behaviours $\gp$ and $\gq$ are \textbf{compatible}, denoted $\gp \smile
\gq$, if for each channel $\alpha$ either $\gp_\alpha \prefix \gq_\alpha$
or $\gq_\alpha \prefix \gp_\alpha$.
Given two compatible behaviours $\gp$ and $\gq$, their \textbf{join},
denoted $\gp \join \gq$, is defined as:
\[
(\gp \join \gq)_\alpha = 
\begin{cases}
\gq_\alpha & \text{if } \gp_\alpha \prefix \gq_\alpha \\
\gp_\alpha & \text{if } \gq_\alpha \prefix \gp_\alpha
\end{cases}
\text{\ \ for all } \alpha
\]

\begin{example}[Compatable behaviours]
The following behaviours
\[\gp =
\begin{array}{|c|c||c|c|} 
&&&\\
\ol{b} &   & \ol{a} &  \\
a & \ol{c} & d &  \\
\hline
\alpha & \beta & \gamma & \delta \\
\hline
\end{array}
\qquad \text{and} \qquad 
\gq = 
\begin{array}{|c|c||c|c|} 
& e &&\\
\ol{b} & d & & \ol{c} \\
a & \ol{c} & d & a \\
\hline
\alpha & \beta & \gamma & \delta \\
\hline
\end{array}
\]
are compatible with their join given by:

\[\gp \join \gq = \begin{array}{|c|c||c|c|} 
& e &&\\
\ol{b} & d & \ol{a} & \ol{c} \\
a & \ol{c} & d & a \\
\hline
\alpha & \beta & \gamma & \delta \\
\hline
\end{array}
\]
\end{example}

\begin{definition}
A \textbf{behavioural entailment} (or simply \textbf{entailment}) is either
of the form
\[\gp \vd \ol{\alpha}[a]
\]
where $\gp$ is an antecedent behaviour and $a$ is a legal output event on
channel $\alpha$ so that $\gp * \ol{\alpha}[a]$ is itself a behaviour, or
\[\gp \vd f
\]
where $f:\Gamma \ra \Delta$ is an atomic morphism on the frontier of $\gp$.
We call $\gp$ the \textbf{antecedent} of the entailment and the
$\ol{\alpha}[a]$ or $f$ the \textbf{conclusion} of the entailment. If the
conclusion of an entailment is an atomic morphism it is called an
\textbf{atomic entailment}. The notation $\gp \vd \times$ will be used to
denote that conclusion may be either an output event or an atomic morphism.
\end{definition}

Notice that no events may be added (or removed) to the antecedent of an
atomic entailment, effectively ``ending'' the entailment. A \textbf{hanging
entailment} is an entailment in which the conclusion is an output event which
gives a transition to a state with a source ($\mathbf{0}$) role. Observe
that after a hanging entailment $\gp \vd \ol{\alpha}[a]$, the ensuing
behaviour $\gp * \ol{\alpha}[a]$ will no longer be an antecedent behaviour.

We will usually talk about a set of entailments $\script{Q}:\Gamma \ra \Delta$,
where every entailment in \script{Q} has as domain channels $\Gamma$ and
codomain channels $\Delta$. An antecedent behaviour which occurs in \script{Q}
(i.e., an antecedent of some entailment in \script{Q}) will be called a
\script{Q}-\textbf{antecedent}. The notation $\alpha[a] * \script{Q}$ will be
used to denote the set of entailments $\script{Q}$ with the event $\alpha[a]$
prepended to each of them, e.g., if $\script{Q} = \{\gp_i \vd \times_i \mid i
\in I\}$ then $\alpha[a] * \script{Q} = \{\alpha[a] * \gp_i \vd \times_i \mid
i \in I\}$.

An output event $\ol{\alpha}[a] \in \gp$ is \script{Q}-\textbf{justified} if 
there is an entailment $\gu \vd \ol{\alpha}[a] \in \script{Q}$ such that
$\gu * \ol{\alpha}[a] \prefix \gp$. A \script{Q}-\textbf{justified behaviour}
is a behaviour in which all the output events are \script{Q}-justified.
A \script{Q}-justified behaviour which is also an antecedent behaviour (not
necessarily a \script{Q}-antecedent) will be called a
\script{Q}-\textbf{preantecedent}. Note that this notion of justification is
not related to the Hyland-Ong~\cite[and other references]{hyland00:full}
notion of justification, and is actually a reachability condition.

\begin{lemma} \label{lem-forceinteraction}
Let \script{Q} be a set of entailments. Every \script{Q}-preantecedent
can, via a sequence (possible empty) of legal inputs, evolve to a saturated
\script{Q}-preantecedent.
\end{lemma}

\begin{proof}
Let $\gp$ be a \script{Q}-justified antecedent behaviour. Then all the
states at the frontier of $\gp$ must have either flow or sink roles.
Consider the states with flow roles. By Lemma~\ref{lem-interaction} we are
able, via a sequence of input events, to drive each of these into either an
atomic, output, or sink state. Thus, we are able to produce a behaviour $\gp'$
such that every state in its frontier is either an atomic state, an output
state, or a sink state, i.e., a saturated behaviour. Moreover, since input
events do not affect justification we have that $\gp'$ is \script{Q}-justified.
\end{proof}

\subsection{Extensional processes} \label{ssec-ep}

\begin{definition}
A set of entailments $\script{Q}:\Gamma \ra \Delta$ is an
\textbf{extensional process} if it satisfies the following conditions:

\begin{enumerate}[EP-1 ]
\item All \script{Q}-antecedents are \script{Q}-justified.

\item If $\gp \vd \ol{\alpha}[a],\ \gq \vd \ol{\alpha}[b] \in \script{Q}$ 
such that $\gp_\alpha = \gq_\alpha$ and $\gp \smile \gq$, then $a=b$.

\item For any \script{Q}-preantecedent $\gp$ and saturated
\script{Q}-preantecedent $\gq$ where $\gp \prefix_i \gq$, there is a
\script{Q}-antecedent $\gp'$ such that $\gp \prefix_i \gp' \prefix_i \gq$.

\item If $\alpha$ and $\beta$ are distinct channels and $\gp \vd 
\ol{\alpha}[a],\ \gq \vd \ol{\beta}[b] \in \script{Q}$ such that $\gp_\alpha
= \gq_\alpha$ and $\gp \prefix \gq$, then $\gq * \ol{\alpha}[a] \vd
\ol{\beta}[b] \in \script{Q}$.

\item If $\alpha$ and $\beta$ are distinct channels such that $\gp \vd 
\ol{\alpha}[a] \in \script{Q}$ and $\beta[b]$ is a legal input event given
$\gp$, then $\gp * \beta[b] \vd \ol{\alpha}[a] \in \script{Q}$.

\item If $\alpha$ and $\beta$ are distinct channels such that $\gp *
\ol{\beta}[b] \vd \ol{\alpha}[a] \in \script{Q}$, then $\gp \vd \ol{\alpha}[a]
\in \script{Q}$.

\item If $\alpha$ and $\beta$ are distinct channels and $\gp$ is a
\script{Q}-preantecedent, such that $\gp * \beta[b] \vd \ol{\alpha}[a] \in
\script{Q}$ for every possible legal input $b$ on $\beta$, then $\gp \vd
\ol{\alpha}[a] \in \script{Q}$.
\end{enumerate}
\end{definition}

\begin{remark}
The definition for an extensional process is rather dense, and so it may
help to have some intuition behind the rules, which we give here. First
of all, behaviours may be thought of as partial processes. The behaviour
is explicit representation of the history of the process. Entailments may
be thought of as a partial processes which may perform an output event
(the conclusion of the entailment).

In this system, justification is a reachability condition for processes,
i.e., to output certain events certain other events must have already been
performed. So, EP-1 asserts that any output event in the antecedent of each
entailment in an extensional process is reachable given what events have
been previously performed. EP-2 ensures that processes are deterministic.
That is, two compatible behaviours (two processes with compatible histories)
which are able to output at a given state on a particular channel, must
output the same event. EP-3 states that if a process is able to perform an
output event, then it must eventually perform one. In game semantic lingo,
this says that the player may not give up. The rule EP-4 states that if a
process can perform an output event at a certain stage, then it can certainly
perform that output event at any later stage (say after first receiving some
other input events and performing some other output events). The rule EP-5
allows the environment to input events without affecting the the output event
the process is about to produce. EP-6 states that output events
do not affect the justification (reachability) of other output events.
Finally, EP-7 states that if for every possible input events on a given
channel, a process is able to output the same event, then this process
must then be able to perform this output event immediately without first
receiving any of these inputs.
\end{remark}

\begin{remark}
We have found it useful to think of behaviours and entailments using
``landscape'' diagrams. For example, the diagram on the left below is a
representation of a behaviour and the diagram on the right below is a
representation of an entailment, the black square representing the conclusion
of the entailment.
\begin{center}
\epsfig{file=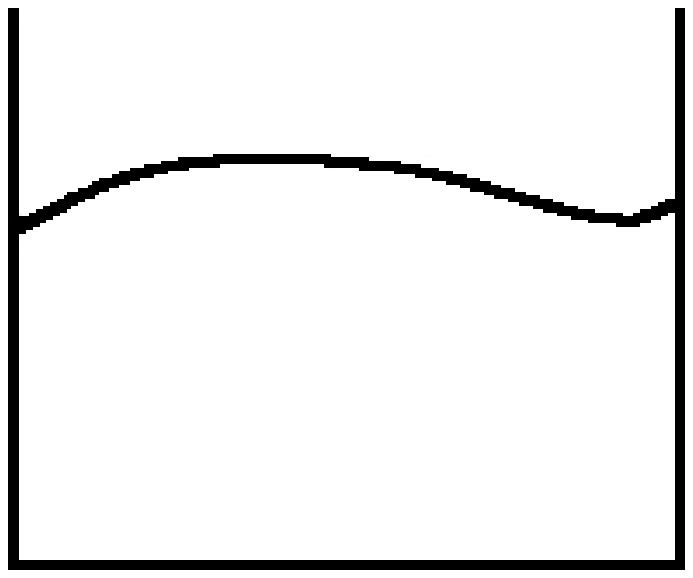, width=5cm, height=4cm} \qqqquad
\epsfig{file=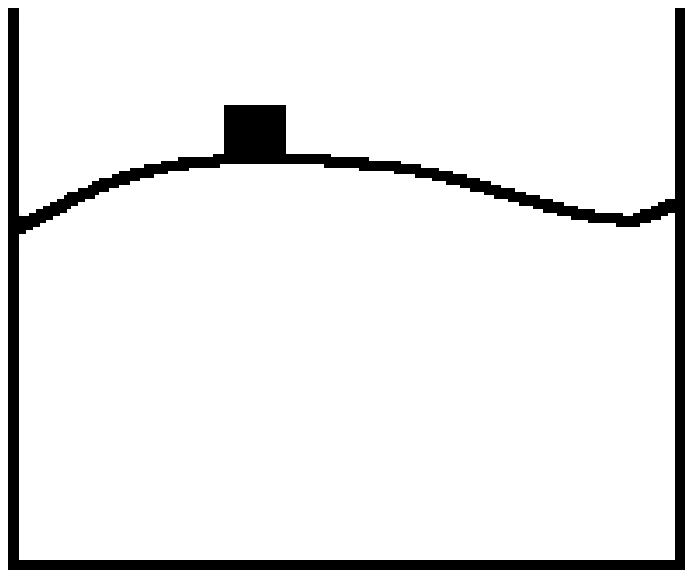, width=5cm, height=4cm}
\end{center}
In this way we may think of a justified output event in a behaviour as:
\begin{center}
\epsfig{file=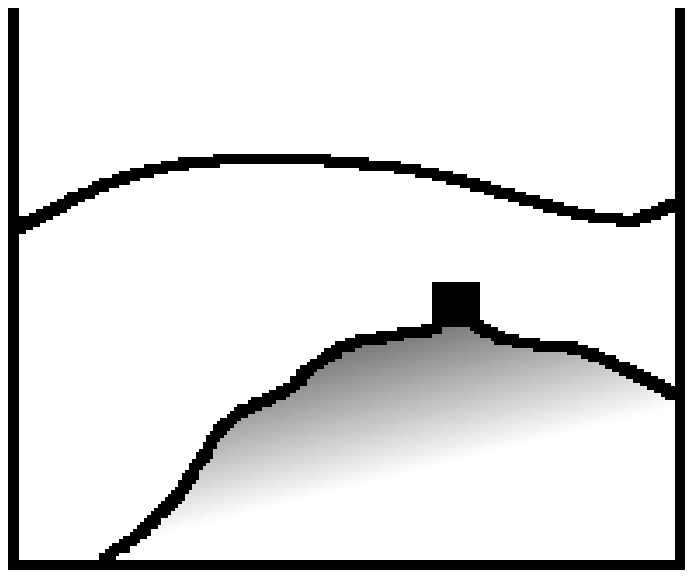, width=5cm, height=4cm}
\end{center}
Using these diagrams, the premise and conclusion of the rule EP-4 are
respectively the left and right diagrams below.
\centro{\epsfig{file=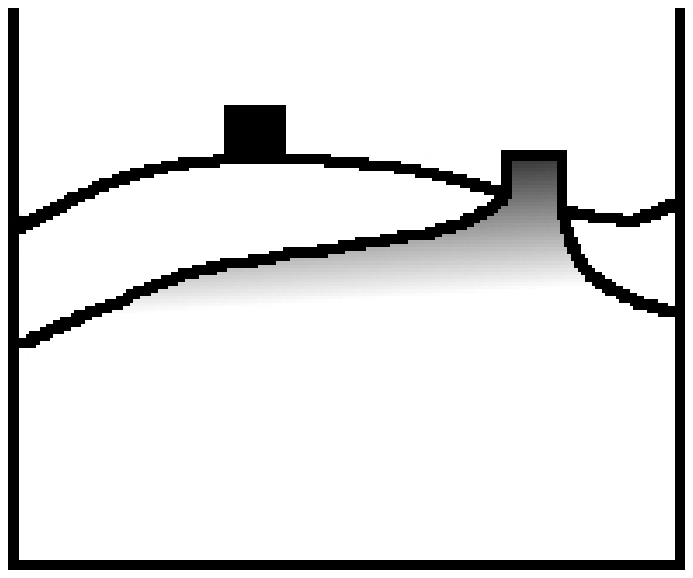, width=5cm, height=4cm} \qqqquad
\epsfig{file=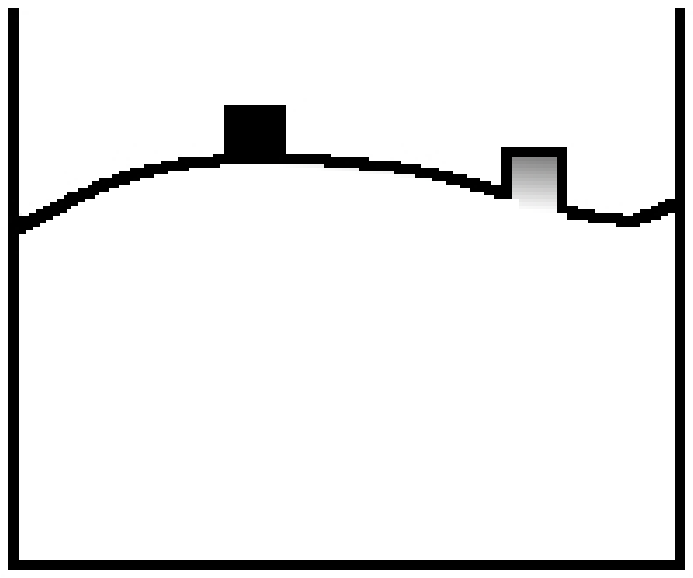, width=5cm, height=4cm}}
\end{remark}

\begin{lemma} \label{lem-joinpre}
Let $\script{Q}:\Gamma \ra \Delta$ be a set of entailments satisfying EP-1
(i.e., a set of \script{Q}-justified entailments). If $\gp$ and $\gq$ are
\script{Q}-preantecedents such that $\gp \smile \gq$, then $\gp \join \gq$
is a \script{Q}-preantecedent.
\end{lemma}

\begin{proof}
This is obvious as any output event in $\gp$ is justified in $\gp \prefix
\gp \join \gq$, and similarly for output events in $\gq$.
\end{proof}

\begin{lemma} \label{lem-pfant}
Let $\script{Q}:\Gamma \ra \Delta$ be a set of entailments satisfying EP-1.
Suppose that $\gp$ and $\gq$ are \script{Q}-preantecedents such that $\gp
\sqsubset \gq$ and for all channels $\alpha$ either $\gp_\alpha =
\gq_\alpha$, or there is an output event $\ol{\alpha}[a] \in \gq$ such that
$\gp_\alpha * \ol{\alpha}[a] \prefix \gq_\alpha$. Then there is an output
event $\ol{\beta}[b] \in \gq$ such that its justification $\gu$ is a
prefix of $\gp$.
\end{lemma}

\begin{proof}
Let $\alpha_i$ for $i \in \{1,\ldots,n\}$ where $n \geq 1$, be all the channels
in which $\gp_{\alpha_i} \sqsubset \gq_{\alpha_i}$ and consider
output events $\ol{\alpha}_i[a_i] \in \gq$ just past the frontier of $\gp$,
i.e., $\gp * \ol{\alpha}_i[a_i] \prefix \gq$.

Consider $\ol{\alpha}_{k}[a_k] \in \gq$ where $k \in \{1,\ldots,n\}$. Its
justification $\gu_k \in \script{Q}$ must have its channel behaviour for
$\alpha_k$ contained in $\gp$. If $\gu_k \not\prefix \gp$ then there exists
a $k' \in \{1,\ldots,n\} \bs \{k\}$ such that $\ol{\alpha}_{k'}[a_{k'}] \in
\gu_k$ but $\ol{\alpha}_{k'}[a_{k'}] \not\in \gp$. The justification of
$\alpha_{k'}[a_{k'}]$ in $\gu_k$, say $\gu_{k'}$, must be a prefix of
$\gu_k$ and thus, must have channel behaviours for $\alpha_k$ and
$\alpha_{k'}$ contained in $\gp$. Continuing this process means that
eventually we will find an output $\ol{\beta}[b]$ with justification
$\gu \in \script{Q}$ such that $\gu \prefix \gp$.
\end{proof}

\begin{lemma} \label{lem-doitbefore}
Let \script{Q} be a set of entailments satisfying EP-1, EP-4, and EP-5. If $\gp
\vd \x \in \script{Q}$, and $\gq$ is a \script{Q}-preantecedent such that $\gp
\prefix \gq$ and $\gp_\alpha = \gq_\alpha$, then $\gq \vd \x \in \script{Q}$.
\end{lemma}

\begin{proof}
We proceed by induction on the number of events between $\gp$ and $\gq$.
If there are no events between $\gp$ and $\gq$ then $\gp = \gq$ and
therefore, $\gp \vd \x \in \script{Q}$.

So suppose that there are $N$ events between $\gp$ and $\gq$ and consider
the channels $\beta_i$, for $i \in I$, in which $\gp_{\beta_i} \sqsubset
\gq_{\beta_i}$. There are two cases to consider: on some $\beta_i$ we can add
input events and on no $\beta_i$ can we add an input event. 

In the first case, use EP-5 to add an input event to $\gp \vd x$. This new
entailment and $\gq$ are now separated by $N-1$ events so that by the
inductive hypothesis $\gq \vd \x \in \script{Q}$.

In the second case we may conclude by Lemma~\ref{lem-pfant} that there is
an output event $\ol{\beta_k}[b] \in \gq$ (but not in $\gp$) for some 
$k \in I$ such that its justification is a prefix of $\gp$. By EP-4 we
may add this event to $\gp \vd \x$, i.e., $\gp * \ol{\beta}[b] \vd \x \in
\script{Q}$. Applying the inductive hypothesis to $\gp * \ol{\beta}[b] \vd
\x$ and $\gq$ now yields $\gq \vd \x \in \script{Q}$.
\end{proof}

The following useful fact is a consequence of the previous lemma.

\begin{corollary} \label{cor-joinentail}
Let \script{Q} be a set of entailments satisfying EP-1, EP-4, and EP-5.
If $\gp \vd \x \in \script{Q}$ and $\gq$ is a \script{Q}-preantecedent such
that $\gp \smile \gq$ and $\gp_\alpha = \gq_\alpha$, then $\gp \join \gq \vd
\x \in \script{Q}$.
\end{corollary} 

\begin{proof} 
It is clear that $\gp \join \gq$ is a \script{Q}-preantecedent with $\gp
\prefix (\gp \join \gq)$ and $\gp_\alpha = (\gp \join \gq)_\alpha$ so that
by Lemma~\ref{lem-doitbefore}, $\gp \join \gq \vd \x \in \script{Q}$.
\end{proof}

The next proposition shows that an extensional process is maximal.

\begin{proposition} \label{prop-subthensame}
Let $\script{Q}, \script{Q'}: \Gamma \ra \Delta$ be extensional processes.
If $\script{Q} \subseteq \script{Q'}$ then $\script{Q} = \script{Q'}$.
\end{proposition}

\begin{proof}
To show equality we must show that $\script{Q}' \subseteq \script{Q}$. To this
end, suppose that $\gp' \vd \x \in \script{Q'}$ and let $\gp$ be a
\script{Q}-preantecedent such that $\gp \prefix \gp'$ ($\gp$ may be the empty
preantecedent). We proceed by induction on the number of events between $\gp$
and $\gq$.  If there are no events between $\gp$ and $\gq$ then $\gp = \gq$
and therefore, $\gq \vd \x \in \script{Q}$.

So, suppose that there are $N$ events between $\gp$ and $\gq$. There are two
cases: we can add an input event on some channel $\alpha$ or on no channel
can we add an input event.

In the first case use EP-5 to add an input event to $\gp$ which corresponds to
an event in $\gq$. This preantecedent and $\gq$ are now separated by $N-1$
events so that by the inductive hypothesis $\gq \vd \x \in \script{Q}$.

In the second case, consider $\gp$ to be a \script{Q'}-preantecedent. By
Lemma~\ref{lem-pfant}, there is an output event $\ol{\alpha}[a] \in \gq$
such that its justification is a prefix of $\gp$. This implies that $\gp *
\ol{\alpha}[a]$ is a \script{Q'}-preantecedent, and therefore, a
\script{Q}-preantecedent. $\gp * \ol{\alpha}[a]$ and $\gq$ are now
separated by $N-1$ events so that by our inductive hypothesis $\gq \vd \x
\in \script{Q}$.

We have now shown that $\script{Q}' \subseteq \script{Q}$ and consequently
$\script{Q} = \script{Q}'$.
\end{proof}

\section{Proto-processes} \label{sec-pp}

A set of entailments \script{P} satisfying the conditions EP-1 through EP-3
is said to be a \textbf{proto-process}. Explicitly, a set of entailments
\script{P} is a proto-process if it satisfies the following conditions:

\begin{enumerate}[EP-1 ]
\item All \script{P}-antecedents are \script{P}-justified.

\item If $\gp \vd \ol{\alpha}[a],\ \gq \vd \ol{\alpha}[b] \in \script{P}$ 
such that $\gp_\alpha = \gq_\alpha$ and $\gp \smile \gq$, then $a=b$.

\item For any \script{P}-preantecedent $\gp$ and saturated
\script{P}-preantecedent $\gq$ where $\gp \prefix_i \gq$, there is a
\script{P}-antecedent $\gp'$ such that $\gp \prefix_i \gp' \prefix_i \gq$.
\end{enumerate}

In this section we show how to translate a (syntactic) process $f$, into a
set of entailments \script{F}, such that \script{F} forms a proto-process.
The main result of this section is that a proto-process \script{P} may be
enlarged in such a way that it additionally satisfies the rules EP-4 $\sim$
EP-7, i.e., a proto-process may be enlarged to an extensional process. This
gives a way to translate syntactic processes into extensional processes.

\subsection{From (syntactic) processes to proto-processes}

Given a (syntactic) process $f: \Gamma \ra \Delta$, the translation of $f$
into a set of entailments, denoted \script{F}, is defined inductively as
follows:

\begin{itemize}
\item If $\Gamma$ or $\Delta$ (or both) contain a protocol with a source role,
then set $\script{F} = \emptyset$.

\item If there are no protocols in $\Gamma$ and $\Delta$ which contain a
source role and

\begin{itemize}
\item $f$ is atomic, then set $\script{F} = \{\vd f\}$.

\item $f = \alpha\{a_i \mapsto f_i\}_{i \in I}$ or $f = \alpha(a_i \mapsto
f_i)_{i \in I}$, then set $\script{F} = \bigcup\limits_{i \in I} \{\alpha[a_i]
* \script{F}_i\}$, where $\script{F}_i$ is the translation of $f_i$.

\item $f = \ola{\alpha}[a](f')$ or $f = \ora{\alpha}[a](f')$, then set
$\script{F} = \{\vd \ol{\alpha}[a], \ \ol{\alpha}[a] * \script{F'}\}$, where
\script{F'} is the translation of $f'$.
\end{itemize}
\end{itemize}

\begin{example}[Translating (syntactic) processes] \quad
\begin{enumerate}
\item Given the following process (no source roles; assume that the atomic
morphisms are well-typed)
\[
\alpha \left\{
\begin{array}{l}
a \mapsto \ora{\beta}[a](f) \medskip\\
b \mapsto \alpha \left\{
  \begin{array}{l}
  d \mapsto \ora{\beta}[b](g) \\
  e \mapsto \ora{\beta}[b](h)
  \end{array}\right\} \medskip \\
c \mapsto \ora{\beta}[c](\beta \left(
  \begin{array}{l}
  d \mapsto \ola{\alpha}[k](i) \\
  e \mapsto \ola{\alpha}[l](j)
  \end{array}\right)) \\
\end{array}\right\}
\]
the translation produces the following set of entailments:

\[\begin{array}{llll}
\begin{array}{|c||c|} 
& \\
a & \\
\hline
\alpha & \beta \\
\hline
\end{array}
\vd \ol{\beta}[a],
& \quad
\begin{array}{|c||c|} 
& \\
a & \ol{a}\\
\hline
\alpha & \beta \\
\hline
\end{array}
\vd f,
& \quad
\begin{array}{|c||c|} 
d & \\
b & \\
\hline
\alpha & \beta \\
\hline
\end{array}
\vd \ol{\beta}[b],
& \quad
\begin{array}{|c||c|} 
d & \\
b & \ol{b}\\
\hline
\alpha & \beta \\
\hline
\end{array}
\vd g,
\bigskip \\ 
\begin{array}{|c||c|} 
e & \\
b & \\
\hline
\alpha & \beta \\
\hline
\end{array}
\vd \ol{\beta}[b],
& \quad
\begin{array}{|c||c|} 
e & \\
b & \ol{b} \\
\hline
\alpha & \beta \\
\hline
\end{array}
\vd h,
& \quad
\begin{array}{|c||c|} 
  & \\
c & \\
\hline
\alpha & \beta \\
\hline
\end{array}
\vd \ol{\beta}[c],
& \quad
\begin{array}{|c||c|} 
  & d \\
c & \ol{c}\\
\hline
\alpha & \beta \\
\hline
\end{array}
\vd \ol{\alpha}[k],
\bigskip \\ 
\begin{array}{|c||c|} 
\ol{k} & d \\
c      & \ol{c} \\
\hline
\alpha & \beta \\
\hline
\end{array}
\vd i, 
& \quad
\begin{array}{|c||c|} 
  & e \\
c & \ol{c}\\
\hline
\alpha & \beta \\
\hline
\end{array}
\vd \ol{\alpha}[l],
& \quad
\begin{array}{|c||c|} 
\ol{l}  & e \\
c & \ol{c}\\
\hline
\alpha & \beta \\
\hline
\end{array}
\vd j 
\end{array}
\]

\item This next example illustrates how the translation handles a
(sub)protocol with a source role. Consider the following protocols and
process (where $f:A \ra B$ is atomic):
\[
\alpha:
\vcenter{
\xymatrix@R=4ex@C=3ex@M=0ex{
& \circ^+ \ar@{-}[dl]_a \ar@{-}[dr]^b \\ \circ^{\mathbf{0}} && A^+}}
\to_{
\alpha\left\{\begin{array}{l}
a \mapsto \ora{\beta}[b](\alpha\{\,\}) \medskip \\
b \mapsto \ora{\beta}[a](f)
\end{array} \right\}}
\beta:
\vcenter{
\xymatrix@R=4ex@C=3ex@M=0ex{
& \circ^+ \ar@{-}[dl]_a \ar@{-}[dr]^b \\ B^+ && A^+}}
\]
Here $a$ on $\alpha$ is an input transition to a protocol with a source role.
The resulting set of entailments is:

\[\begin{array}{|c||c|} 
& \\
b & \\
\hline
\alpha & \beta \\
\hline
\end{array}
\vd \ol{\beta}[a], 
\qquad
\begin{array}{|c||c|} 
  &  \\
b & \ol{a} \\
\hline
\alpha & \beta \\
\hline
\end{array}
\vd f,
\]

\end{enumerate}
\end{example}
We now wish to establish that this translation procedure results in a
proto-process. First however, there is an observation which is required.

\begin{proposition} \label{prop-emptyep}
The empty set $\emptyset: \Gamma \ra \Delta$ is an extensional process if
and only if there is a protocol with a source role in the domain or in the
codomain (or both).
\end{proposition}

\begin{proof}
We first prove the ``if'' direction. All the premises of the EP-rules,
except for EP-3, involve having some entailments, and thus $\emptyset$
vacuously satisfies these rules. To see that $\emptyset$ satisfies EP-3 as
well recall the definition of an antecedent behaviour: a behaviour is an
antecedent behaviour in case each state at the frontier has either a flow
or a sink role. If $\Gamma$ or $\Delta$ contain a protocol with a source role,
there can be no antecedent (or preantecedent) behaviours. Thus, $\emptyset$
satisfies EP-3, and is therefore an extensional process.

For the ``only if'' direction suppose neither $\Gamma$ nor $\Delta$ contains
a protocol with a source role. This implies that the empty behaviour is an
preantecedent behaviour and, therefore, EP-3 would require $\emptyset$ to
contain some entailment. This is a contradiction and therefore $\Gamma$ or
$\Delta$ must contain a protocol with a source role. 
\end{proof}

\begin{proposition} \label{transtopp}
If a set of entailments \script{F} is the result of translating a
(syntactic) process $f:\Gamma \ra \Delta$, then \script{F} is a proto-process.
\end{proposition}

\begin{proof}
We proceed by structural induction on the process $f$. If $f$ is atomic or
one of $\Gamma$ or $\Delta$ contains a source role, then the translation of
$f$ is clearly a proto-process.

So assume now that the translation fits one of the other cases. It must
produce either

\[\mathrm{(i)}\
\bigcup_{i \in I} \{\alpha[a_i] * \script{F}_i\} \qquad \text{or} \qquad
\mathrm{(ii)}\ \{\vd \ol{\alpha}[a],\ \ol{\alpha}[a] * \script{F'}\}
\]
where we inductively assume that the $\script{F}_i$'s and \script{F'} are
proto-processes. The goal is to show that this new set of entailments is a
proto-process.

\begin{enumerate}[EP-1 ]
\item (i) Let $\gp$ be a
$(\bigcup_{i \in I} \{\alpha[a_i] * \script{F}_i\})$-antecedent and suppose
that $\ol{\beta}[b] \in \gp$.
That $\gp$ is a $(\bigcup_{i \in I} \{\alpha[a_i] * \script{F}_i\})$-antecedent
implies that $\gp$ must be of the form $\alpha[a_k] * \gp'$, for some $k \in
I$, where $\gp' \in \script{F}_k$. So, $\ol{\beta}[b]$ must be in $\gp'$, and
as $\script{F}_k$ is a proto-process, it must be justified
in $\script{F}_k$, say by $\gu$, but this implies that $\alpha[a_k] * \gu$
justifies it in $\bigcup_{i \in I} \{\alpha[a_i] * \script{F}_i\}$.

(ii) There are no output events in $\vd \ol{\alpha}[a]$ so it justified,
and any antecedent in $\ol{\alpha}[a] * \script{F'}$ is justified using an
argument similar to (i) above.

\item (i) Suppose that for some $k \in I$, $\alpha[c] * \gp \vd \ol{\beta}[a],\
\alpha[c] * \gq \vd \ol{\beta}[b] \in \alpha[c] * \script{F}_k$ such that
$(\alpha[c] * \gp) \smile (\alpha[c] * \gq)$ and $(\alpha[c] * \gp)_\beta =
(\alpha[c] * \gq)_\beta$. Then $\gp \vd \ol{\beta}[a]$, $\gq \vd \ol{\beta}[b]
\in \script{F}_k$ such that $\gp \smile \gq$ and $\gp_\beta = \gq_\beta$, and
so $a=b$ by EP-2 in $\script{F}_k$.

(ii) The entailment $\vd \ol{\alpha}[c]$ and any entailment in $\ol{\alpha}[c]
* \script{F'}$ will never have the same channel behaviour on $\alpha$, and
therefore, the only case we need to consider is if $\ol{\alpha}[c] * \gp \vd
\ol{\beta}[a]$, $\ol{\alpha}[c] * \gq \vd \ol{\beta}[b] \in \ol{\alpha}[c] *
\script{F}'$ such that $(\ol{\alpha}[c] * \gp) \smile (\ol{\alpha}[c] * \gq)$
and $(\ol{\alpha}[c] * \gp)_\beta = (\ol{\alpha}[c] * \gq)_\beta$. Using an
argument similar to (i) above, this implies $a=b$.

\item (i) Let $\gq$ be a saturated $(\bigcup_{i \in I} \{\alpha[a_i] *
\script{F}_i\})$-preantecedent. It must be of the form $\alpha[a_k] * \gq'$,
for some $k \in I$, where $\gq'$ is a saturated $\script{F}_k$-preantecedent.
By EP-3, $\gq' \vd \x \in \script{F}_k$, and thus, $\alpha[a_k] * \gq' \vd \x
\in \bigcup_{i \in I} \{\alpha[a_i] * \script{F}_i\}$ satisfying EP-3.

(ii) Here, similarly to (i) above, any saturated $\{\vd \ol{\alpha}[a],\
\ol{\alpha}[a] * \script{F'}\}$-preantecedent $\gq$ is either of the form
$\ol{\alpha}[a] * \gq'$ or $\gq'$, where $\gq' \in \script{F}'$. Thus, an
argument similar to (i) above suffices.

\end{enumerate}

Thus, since each step produces a proto-process, the complete translation
\script{F} is a proto-process.
\end{proof}

An important example of a proto-process is the \textbf{identity proto-process}.
As it is known how to form the identity process (from the identity derivation;
see Section~\ref{sec-id}), Proposition~\ref{transtopp} tells us how to
form the identity proto-process.

\begin{example}[Identity proto-process]
The identity process $\alpha:X \to^{1_X} \beta:X$ where $X$ is the protocol
\[X = \{a:A,b:\{d:B, e:(\,)\}, c:(f:C, g:\{\,\})\}
\]
is defined to be:
\[1_X = \alpha \left\{
\begin{array}{l}
a \mapsto \ora{\beta}[a](1_A) \medskip\\
b \mapsto \ora{\beta}[b](
  \alpha\left\{\begin{array}{l}
  d \mapsto \ora{\beta}[d](1_B) \\
  e \mapsto \ora{\beta}[e](1_{(\,)})
  \end{array}\right\}) \medskip\\
c \mapsto \ora{\beta}[c](
  \beta\left(\begin{array}{l}
  f \mapsto \ora{\beta}[f](1_C) \\
  g \mapsto \ora{\beta}[g](1_{\{\,\}})
  \end{array}\right))
\end{array}\right\}
\]

The proto-process resulting from the translations is as follows:

\[\begin{array}{llll}
\begin{array}{|c||c|} 
& \\
a & \\
\hline
\alpha & \beta \\
\hline
\end{array}
\vd \ol{\beta}[a],
& \quad
\begin{array}{|c||c|} 
& \\
a & \ol{a} \\
\hline
\alpha & \beta \\
\hline
\end{array}
\vd 1_A
& \quad
\begin{array}{|c||c|} 
& \\
b & \\
\hline
\alpha & \beta \\
\hline
\end{array}
\vd \ol{\beta}[b],
& \quad
\begin{array}{|c||c|} 
d & \\
b & \ol{b} \\
\hline
\alpha & \beta \\
\hline
\end{array}
\vd \ol{\beta}[d],
\bigskip \\ 
\begin{array}{|c||c|} 
d & \ol{d} \\
b & \ol{b} \\
\hline
\alpha & \beta \\
\hline
\end{array}
\vd 1_B
& \quad
\begin{array}{|c||c|} 
e & \\
b & \ol{b} \\
\hline
\alpha & \beta \\
\hline
\end{array}
\vd \ol{\beta}[e],
& \quad
\begin{array}{|c||c|} 
e & \ol{e} \\
b & \ol{b} \\
\hline
\alpha & \beta \\
\hline
\end{array}
\vd 1_{(\,)}
& \quad
\begin{array}{|c||c|} 
& \\
c & \\
\hline
\alpha & \beta \\
\hline
\end{array}
\vd \ol{\beta}[c],
\bigskip \\ 
\begin{array}{|c||c|} 
& f \\
c & \ol{c} \\
\hline
\alpha & \beta \\
\hline
\end{array}
\vd \ol{\alpha}[f],
& \quad
\begin{array}{|c||c|} 
\ol{f} & f \\
c & \ol{c} \\
\hline
\alpha & \beta \\
\hline
\end{array}
\vd 1_C
& \quad
\begin{array}{|c||c|} 
& g \\
c & \ol{c} \\
\hline
\alpha & \beta \\
\hline
\end{array}
\vd \ol{\alpha}[g],
& \quad
\begin{array}{|c||c|} 
\ol{g} & g \\
c & \ol{c} \\
\hline
\alpha & \beta \\
\hline
\end{array}
\vd 1_{\{\,\}}
\end{array}
\]
\end{example}

\subsection{Proto-processes to extensional processes} \label{sec-pp2ep}

The purpose of this section is to prove:

\begin{theorem} \label{thm-exproex}
If \script{P} is a proto-process then there exists a unique extensional
process \script{Q} with $\script{P} \subseteq \script{Q}$.
\end{theorem}

The idea will be to enlarge a proto-process with additional entailments
in such a way that this new set of entailments will satisfy the rules EP-4
$\sim$ EP-7. If this can be accomplished, the new set of entailments will
satisfy all the rules of an extensional process. To this end, we prove a
series of lemmas (one for each additional EP rule) which allows us to add
entailments to a proto-process such that the resulting set of entailments is
a proto-process.

\begin{lemma}
Suppose that $\script{P}_i$ is a set of entailments satisfying EP-$i$
for $i \in \{1,2,3\}$. An application of EP-4 to the entailments of
$\script{P}_i$ produces a larger set of entailments $\script{P}'_i$ which
satisfies EP-$i$. That is, if $\gp \vd \ol{\alpha}[a],\ \gq \vd \ol{\beta}[b]
\in \script{P}_i$, such that $\gp \prefix \gq$ and $\gp_\alpha =
\gq_\alpha$, then $\script{P}'_i = \script{P}_i \cup \{\gq * \ol{\alpha}[a]
\vd \ol{\beta}[b]\}$ satisfies EP-$i$.
\end{lemma}

\begin{proof}
For each $i \in \{1,2,3\}$, we must show that $\script{P}'_i$ (defined above)
satisfies EP-$i$.

\begin{enumerate}[EP-1 ]

\item Our new antecedent $\gq * \ol{\alpha}[a]$ is $\script{P}_1$-justified (as
$\gp \vd \ol{\alpha}[a] \in \script{P}$) and so it is $\script{P}'_1$-justified.

\item Suppose that $\gr \vd \ol{\beta}[c] \in \script{P}_2$ such that
$\gr_\beta = (\gq * \ol{\alpha}[a])_\beta$ and $\gr \smile (\gq *
\ol{\alpha}[a])$. Then $\gr_\beta = \gq_\beta$ and $\gr \smile \gq$ in
$\script{P}_2$ so that $b=c$.

\item Adding an output event to an entailment does not increase the number
of $\script{P}'_3$-preantecedents and thus, $\script{P}'_3$ satisfies EP-3.
\qed
\end{enumerate}
\end{proof}

This allows us to conclude:

\begin{corollary}
If \script{P} is a proto-process, the closure of \script{P} with respect to
EP-4 is a proto-process.
\end{corollary}

\begin{lemma}
Suppose that $\script{P}_i$ is a set of entailments satisfying EP-$i$
for $i \in \{1,2,3\}$. An application of EP-5 to an entailment in
$\script{P}_i$ produces a larger set of entailments $\script{P}'_i$ which
satisfies EP-$i$. That is, if $\gp \vd \ol{\alpha}[a] \in \script{P}_i$ and
$\beta[b]$ is a legal input event given $\gp$, then $\script{P}'_i =
\script{P}_i \cup \{\gp * \beta[b] \vd \ol{\alpha}[a]\}$ satisfies EP-$i$.
\end{lemma}

\begin{proof}
For each $i \in \{1,2,3\}$, we must show that $\script{P}'_i$ (defined above)
satisfies EP-$i$.

\begin{enumerate}[EP-1 ]

\item Adding an input event to a $\script{P}_1$-antecedents will not affect
its $\script{P}_1$-justification.

\item Suppose $\gp \vd \ol{\alpha}[a] \in \script{P}_2$ and that $\beta[b]$
is a legal input event given $\gp$. Adding this input to $\gp$ gives an
entailment $\gp * \beta[b] \vd \ol{\alpha}[a] \in \script{P}'_2$. Now
suppose that there is a $\gq \vd \ol{\alpha}[b] \in \script{P}_2$ such that
$(\gp * \beta[b])_\alpha = \gq_\alpha$ and $\gp * \beta[b] \smile \gq$.
Then clearly, $\gp_\alpha = \gq_\alpha$ and $\gp \smile \gq$ in $\script{P}_2$
so that $a=b$.

\item Adding input events does not increase the number of
$\script{P}_3$-preantecedents and thus, $\script{P}'_3$ satisfies EP-3.
\qed
\end{enumerate}
\end{proof}

Similarly here we may conclude:

\begin{corollary} \label{cor-closewith5}
If \script{P} is a proto-process, the closure of \script{P} with respect to
EP-5 is a proto-process.
\end{corollary}

Enlarging a proto-process with respect to EP-6 has a slightly different flavor
from the rest in that it cannot be done one step at a time. We must take
the closure of the set with respect to EP-6 in order for it to work. Another
difference is that the set we wish to enlarge must satisfy all the rules of
a proto-process, unlike the other lemmas.

\begin{lemma}
Let \script{P} be a proto-process. The closure of \script{P} with respect to
EP-6 produces a larger set of entailments \script{P'} which is itself a
proto-process.
\end{lemma}

\begin{proof}
We must show that \script{P'} satisfies EP-1 $\sim$ EP-3.

\begin{enumerate}[EP-1 ]

\item Suppose $\gp$ is a \script{P'}-antecedent and consider an output event
$\ol{\alpha}[a] \in \gp$. That $\gp$ is a \script{P'}-antecedent implies that
there is an \script{P}-antecedent $\gp'$ such that $\gp \prefix_o \gp'$. As
\script{P}-antecedents are justified $\ol{\alpha}[a]$ must have some
justifying antecedent in \script{P}, say $\gu$, i.e., $\gu \vd
\ol{\alpha}[a] \in \script{P}$ such that $\gu \prefix \gp'$. In $\gu$
there may be some channels $\gamma_j$ such that $\gp_{\gamma_j} \sqsubset_o
\gu_{\gamma_j}$, but simply applying EP-6 to the output events on these
channels gives an entailment $\gu' \vd \ol{\alpha}[a] \in \script{P'}$ such
that $\gu' \prefix \gp$. Thus, $\ol{\alpha}[a]$ is justified in \script{P'}.

\item Suppose that $\gp' \vd \ol{\alpha}[a],\ \gq' \vd \ol{\alpha}[b] \in
\script{P'}$ with $\gq'_\alpha = \gp'_\alpha$ and $\gq' \smile \gp'$. As EP-6
removes output events this implies that there are entailments $\gp \vd
\ol{\alpha}[a],\ \gq \vd \ol{\alpha}[b] \in \script{P}$ such that 
$\gp' \prefix_o \gp$ and $\gq' \prefix_o \gq$. Clearly $\gp_\alpha =
\gq_\alpha$, and by induction on the number of ``incompatible'' output events,
we will show that $\gp \smile \gq$ so that $a=b$ follows.

If there are no incompatible output events then $\gp \smile \gq$, so assume
that there are $N$ incompatible output events. Consider a particular one, say
$\ol{\beta}[c] \in \gp$ and $\ol{\beta}[d] \in \gq$, and their respective
justifications in \script{P}: $\gu \vd \ol{\beta}[c]$ and $\gothic{s} \vd
\ol{\beta}[d]$ where $\gu \sqsubset \gp$ and $\gothic{s} \sqsubset \gq$.
Now $\gu_\beta = \gothic{s}_\beta$ and $\gu$ and $\gothic{s}$ have less than
$N$ incompatible events so that by the inductive hypothesis $\gr \smile
\gothic{s}$ and so $c=d$. This now implies that $\gp$ and $\gp$ have $N-1$
incompatible events, so that again applying the inductive hypothesis $\gp
\smile \gq$.

\item Removing output events does not increase the number of
\script{P}-preantecedents and so \script{P'} satisfies EP-3.
\qed
\end{enumerate}
\end{proof}

\begin{lemma}
Suppose that $\script{P}_i$ is a set of entailments satisfying EP-$i$
for $i \in \{1,2,3\}$. An application of EP-7 to an entailment in
$\script{P}_i$ produces a larger set of entailments $\script{P}'_i$ which
satisfies EP-$i$. That is, if $\gp$ is a \script{P}-antecedent and
$\gp * \beta[b_i] \vd \ol{\alpha}[a] \in \script{P}_i$ for every possible
legal input $b_i$ on $\beta$, then $\script{P}'_i = \script{P}_i \cup \{\gp
\vd \ol{\alpha}[a]\}$ satisfies EP-$i$.
\end{lemma}

\begin{proof}
For each $i \in \{1,2,3\}$, we must show that $\script{P}'_i$ (defined above)
satisfies EP-$i$.

\begin{enumerate}[EP-1 ]

\item The new antecedent $\gp$ is justified in $\script{P}_1$, and so is
justified in $\script{P}'_1$. Therefore, $\script{P}'_1$ satisfies EP-1.

\item Suppose that $\gp \vd \ol{\alpha}[a] \in \script{P}'_2$ and  $\gq \vd
\ol{\alpha}[b] \in \script{P}_2$ with $\gp_\alpha = \gq_\alpha$ and $\gp
\smile \gq$. Then as $\gp * \beta[b] \vd \ol{\alpha}[a]$ for all possible
inputs $b_i$ on $\beta$, there is some event $b_k$ such that in $\script{P}_2$,
$(\gp * \beta[b_k])_\alpha = \gq_\alpha$ and $\gp * \beta[b_k] \smile \gq$.
It now follows that $a=b$.

\item It suffices to show that any saturated $\script{P}_3'$-preantecedent is
itself a $\script{P}_3'$-antecedent. This follows since given any saturated
$\script{P}'_3$-preantecedent $\gq$, it must also be a saturated
$\script{P}_3$-preantecedent. If input events were removed from $\gq$ then
the resulting behaviour would no longer be saturated. So, applying EP-3 in
$\script{P}_3$ to $\gq$ implies that $\gq$ is a $\script{P}_3$-antecedent and
therefore, a $\script{P}'_3$-antecedent.
\qed
\end{enumerate}
\end{proof}

We conclude:

\begin{corollary} \label{cor-closewith7}
If \script{P} is a proto-process, the closure of \script{P} with respect to
EP-7 is a proto-process.
\end{corollary}

The previous four lemmas show that a proto-process \script{P} may be enlarged
to a set of entailments \script{Q}, so that \script{Q} is itself a
proto-process which additionally satisfies EP-4 through EP-7, i.e., an
extensional process.

To complete the proof of Theorem~\ref{thm-exproex}, it is necessary to show
that that this enlarging procedure results in a unique extensional process.
To this end, suppose that $\script{P}$ is a proto-process and \script{Q} is
any extensional process such that $\script{P} \subseteq \script{Q}$ and denote
by $\script{E}(\script{P})$ the extensional process which results from
enlarging $\script{P}$. If $\gp \vd \times \in \script{E(P)}$ is the result
of closing \script{P} with respect to some EP-rule, then, as $\script{P}
\subseteq \script{Q}$, it must be in \script{Q} as well. Therefore,
$\script{E(P)} \subseteq \script{Q}$ and by
Proposition~\ref{prop-subthensame}, $\script{E(P)} = \script{Q}$.

This now completes the proof of Theorem~\ref{thm-exproex}. 

\section{The polycategory of extensional processes} \label{sec-pcep}

In this section we establish that protocols and extensional processes form
a polycategory which we will denote by $\script{E}\script{P}$. We begin by
introducing some new notation. Let $\gp$ and $\gq$ be behaviours. If
$\gamma$ is a channel in $\gp$, the notation $\gp \bs \gamma$ is used to
represent the behaviour $\gp$ with all the events on $\gamma$ (including the
channel name) removed. For example,
\[\gp =
\begin{array}{|c||c|c|}
&& \\
\ol{a} & c      & e      \\
b      & \ol{d} & \ol{f} \\
\hline
\alpha & \beta  & \gamma \\
\hline
\end{array}
\qqquad\qquad
\gp \bs \gamma =
\begin{array}{|c||c|}
& \\
\ol{a} & c      \\
b      & \ol{d} \\
\hline
\alpha & \beta  \\
\hline
\end{array}
\]

The notation $\gp || \gq$ is used to denote the behaviour whose input
channel behaviours are the input channel behaviours of $\gp$ and $\gq$ and
whose output channel behaviours are the output channel behaviours of $\gp$
and $\gq$ (assuming $\gp$ and $\gq$ have distinct channel names). For
example,

\[
\gp =
\begin{array}{|c||c|}
& \\
\ol{a} & c \\
b      & \ol{d} \\
\hline
\alpha & \beta \\
\hline
\end{array}
\qquad
\gq =
\begin{array}{|c||c|}
& \\
e      & \ol{g} \\
\ol{f} & h \\
\hline
\gamma & \delta \\
\hline
\end{array}
\qqquad\qquad
\gp || \gq =
\begin{array}{|c|c||c|c|}
&&&\\
\ol{a} & e      & c      & \ol{g} \\
b      & \ol{f} & \ol{d} & h \\
\hline
\alpha & \gamma & \beta & \delta \\
\hline
\end{array}
\]

We are now ready to show how to compose extensional processes.

\subsection{Composition of extensional processes}

Here, as in the case with syntactic processes (see Section~\ref{ssec-cer}),
in order to compose two extensional processes \script{F} and \script{G} they
must have exactly one channel name in common. This means that we must rename
channels so that the channel names in \script{F} and \script{G} are distinct
except for the channel on which they will be composed. This may be accomplished
through channel name substitution, however, as we did with syntactic processes,
we will assume that the channel names of any two extensional processes are
distinct unless otherwise specified.

\begin{definition} [Composition of extensional processes]
Let
\[\script{F}: \Gamma \ra \gamma:X,\Delta \qquad \text{and} \qquad
\script{G}: \Gamma',\gamma:X \ra \Delta'
\]
be extensional processes. The set of entailments
representing the composite of \script{F} and \script{G} on $\gamma$, denoted
$\script{F} ;_\gamma \script{G}: \Gamma,\Gamma' \ra \Delta,\Delta'$, is
defined to be:
\begin{align*}
\script{F} ;_\gamma \script{G} =
& \left\{(\gp\bs\gamma || \gq\bs\gamma) \vd f\ ;_\gamma g\ \mid\ \gp \vd f \in
\script{F} \meet \gq \vd g \in \script{G} \meet \gp_\gamma = \gq_\gamma^*
\right\} \\
& \cup 
\left\{(\gp\bs\gamma || \gq\bs\gamma) \vd \ol{\alpha}[a]\ \mid\
\gp \vd \ol{\alpha}[a] \in \script{F} \meet
\gq \text{ a \script{G}-preantecedent} \meet
\gp_\gamma = \gq^*_\gamma \right\} \\
& \cup
\left\{(\gp\bs\gamma || \gq\bs\gamma) \vd \ol{\beta}[a]\ \mid\
\gp \text{ a \script{F}-preantecedent} \meet
\gq \vd \ol{\beta}[a] \in \script{G} \meet
\gp_\gamma = \gq^*_\gamma \right\}
\end{align*}
\end{definition}

In order to simplify the notation, the shorthand $\gp ;_\gamma \gq$ will be
used to denote $\gp\bs\gamma || \gq\bs\gamma$ provided that $\gp_\gamma =
\gq_\gamma^*$.

Composition of extensional processes may become clearer with a couple of
examples.

\begin{example}[Composing extensional processes]
\quad
\begin{enumerate}
\item The following two atomic entailments may be composed on $\gamma$,
\[
\begin{array}{|c|c||c|c|}
&&& \\
  & \ol{e} &        & \ol{h} \\
a & b      & \ol{c} & d      \\
\hline
\alpha & \beta & \delta & \gamma \\
\hline
\end{array}
\vd f \in \script{F}
\qqquad
\begin{array}{|c|c||c|}
&& \\
\ol{h} &   & \ol{j} \\
d      & i & k      \\
\hline
\gamma & \epsilon & \eta \\
\hline
\end{array}
\vd g \in \script{G}
\]
with their composite given by:
\[
\begin{array}{|c|c|c||c|c|}
&&&& \\
  & \ol{e} &   &        & \ol{j} \\
a & b      & i & \ol{c} & k      \\
\hline
\alpha & \beta & \epsilon & \delta & \eta \\
\hline
\end{array}
\vd f;g \in \script{F} ;_\gamma \script{G}
\]

\item The following entailment in \script{F} and \script{G}-preantecedent may
be composed on $\gamma$,
\[
\begin{array}{|c|c||c|c|c|} 
&&& c &\\
\ol{b} &   & \ol{a} & \ol{b} & \\
a & \ol{a} & d & a &  \\
\hline
\alpha & \beta & \delta & \gamma & \epsilon \\
\hline
\end{array}
\vd \ol{\alpha}_3[c] \in \script{F}
\qqquad
\begin{array}{|c|c|c||c|} 
& \ol{c} & \ol{a} & \\
& b & c & \\
\ol{e} & \ol{a} & e & \ol{c} \\
\hline
\eta & \gamma & \theta & \iota \\
\hline
\end{array}
\quad \text{a \script{G}-preantecedent}
\]
with their composite given by:
\[
\begin{array}{|c|c|c|c||c|c|c|} 
       &        &        & \ol{a} &&& \\
\ol{b} &        &        & c      & \ol{a} & & \\
a      & \ol{a} & \ol{e} & e      & d      & & \ol{c}  \\
\hline
\alpha & \beta & \eta & \theta & \delta & \epsilon & \iota \\
\hline
\end{array}
\vd \ol{\alpha}_3[c] \in \script{F} ;_\gamma \script{G}
\]

\end{enumerate}
\end{example}

\begin{lemma} \label{lem-sameingamma}
If $\gp\ ;_\gamma \gq \vd \times,\ \gp'\ ;_\gamma \gq' \vd \times' \in
\script{F} ;_\gamma \script{G}$ such that $\gp\ ;_\gamma \gq \smile \gp'\
;_\gamma \gq'$, then $\gp \smile \gp'$ and $\gq \smile \gq'$.
\end{lemma}

\begin{proof}
We will show that $\gp$ and $\gp'$ are compatible; a similar argument can
be used to show that $\gq$ and $\gq'$ are compatible.

The behaviours $\gp$ and $\gp'$ are compatible on channels other than
$\gamma$ by definition, so what needs to be shown is that they are
compatible on $\gamma$. To this end, suppose they are not compatible on
$\gamma$ and consider the first position in which they differ, say by the
events $a \in \gp$ and $b \in \gp'$ in \script{F} (dually $\ol{a} \in \gq$
and $\ol{b} \in \gq'$ in \script{G}). There are two cases to consider:

\begin{enumerate}[(i)]
\item $a$ and $b$ are output events. Consider their justifications $\gu \vd
\ol{\gamma}[a] \in \gp$ and $\gu' \vd \ol{\gamma}[b] \in \gp'$. Clearly,
$\gu_\gamma = \gu'_\gamma$ and $\gu \smile \gu'$ so that
by EP-2, $a=b$.

\item $a$ and $b$ are input events. This implies that they are output events
in $\gq$ and $\gq'$ and moreover, this must be the first place in which
they differ on $\gamma$. Therefore, by repeating the argument used in (i)
in \script{G}, $a=b$
\end{enumerate}
Thus, $\gp$ and $\gp'$ are compatible.
\end{proof}

\begin{proposition}
Let $\script{F}: \Gamma \ra \gamma:X,\Delta$ and $\script{G}: \Gamma',
\gamma:X \ra \Delta'$ be extensional processes. The set of entailments
$\script{F} ;_\gamma \script{G}: \Gamma,\Gamma' \ra \Delta, \Delta'$ forms an
extensional process.
\end{proposition}

\begin{proof}
We must show that $\script{F} ;_\gamma \script{G}$ satisfies EP-1 through EP-7.

\begin{enumerate}[EP-1 ]

\item Suppose $\gp;_\gamma \gq$ is a $\script{F} ;_\gamma
\script{G}$-antecedent and $\ol{\alpha}[a]$ is an output event in $\gp;_\gamma
\gq$; without loss of generality assume that $\alpha \in \script{F}$. It
follows from the definition of composition that $\gp$ is justified in
\script{F} and $\gq$ is justified in \script{G} so that $\ol{\alpha}[a]$ is
justified in $\gp$, say by $\gu \vd \ol{\alpha}[a]$.

By EP-6 (in \script{F}), output events on $\gamma$ may be removed from
$\gu$ so that either $\gu_\gamma = \emptyset$ or the last event in
$\gu_\gamma$ is an input event. 

If $\gu_\gamma = \emptyset$ then $\gu \vd \ol{\alpha}[a]$ justifies
$\ol{\alpha}[a]$ in $\script{F};_\gamma \script{G}$. If $\gu_\gamma
\neq \emptyset$, then in \script{F} the last event on $\gamma$ is an input
event $\gamma[c]$, which means that it is an output event in $\gq$.
Moreover, it is justified in $\gq$, say by $\gu' \vd \ol{\gamma}[c]$, and so
$\gu' * \ol{\gamma}[c]$ is a justified behaviour in \script{G}. Therefore,
$\gu ;_\gamma (\gu'* \gamma[c]) \vd \ol{\alpha}[a] =  \gu ;_\gamma \gu' \vd
\ol{\alpha}[a] \in \script{F} ;_\gamma \script{G}$ and $\gu ;_\gamma \gu'
\prefix \gp ;_\gamma \gq$ so that $\ol{\alpha}[a]$ is justified in $\script{F}
;_\gamma \script{G}$.

If $\alpha \in \script{G}$ a similar argument in which we reverse the
components shows that $\ol{\alpha}[a]$ is justified in $\script{F} ;_\gamma
\script{G}$.

\item Suppose that $\gp ;_\gamma \gq \vd \ol{\alpha}[a],\ \gp' ;_\gamma
\gq' \vd \ol{\alpha}[b] \in \script{F} ;_\gamma \script{G}$ such that
$(\gp ;_\gamma \gq)_\alpha = (\gp' ;_\gamma \gq')_\alpha$ and $\gp
;_\gamma \gq \prefix \gp' ;_\gamma \gq'$. Without loss of generality
assume that $\alpha \in \script{F}$. This implies $\gp \vd \ol{\alpha}[a],\
\gp' \vd \ol{\alpha}[b] \in \script{F}$ such that $\gp_\alpha = \gp_\alpha$
and by Lemma~\ref{lem-sameingamma}, $\gp \smile \gp'$ so that by EP-2 in
\script{F}, $a=b$.

$\alpha \in \script{G}$ is handled similarly by considering
$\gq \vd \ol{\alpha}[a],\ \gq' \vd \ol{\alpha}[b] \in \script{G}$.

\item
Suppose that $\gp || \gp'$ is a $\script{F} ;_\gamma \script{G}$-preantecedent
and $\gq || \gq'$ is a saturated $\script{F}\ ;_\gamma
\script{G}$-preantecedent such that $\gp || \gp' \prefix_i \gq || \gq'$.

Consider (the possible unjustified) behaviours $\gp$ in \script{F} and $\gp'$
in \script{G}. Any output event $\ol{\alpha}[a_i] \in \gp$ must be justified
in $\script{F} ;_\gamma \script{G}$, say by $\gu_i ;_\gamma \gu'_i$ so that
$\gu_i \vd \ol{\alpha}[a_i]$ justifies $\ol{\alpha}[a_i]$ in \script{F}.
Similarly we may justify output events in $\gp'$ in \script{G}.

Choose the justification $\gu ;_\gamma \gu'$ such that $\gu^*_\gamma =
\gu'_\gamma$ is maximal. By Lemma~\ref{lem-sameingamma}, all the $\gu_i$'s are
compatible and all the $\gu'_i$'s are compatible so that $\gp \join \gu$ and
$\gp' \join \gu'$ are respectively \script{F} and \script{G}-preantecedents.
Notice that also $\gq \join \gu$ and $\gq' \join \gu'$ are also \script{F}
and \script{G}-preantecedents respectively (as they are separated from $\gp$
and $\gp'$ by input events). Consider the last state on $\gu_\gamma$: it may
be atomic, an input state, or an output state.

(i) Suppose it is atomic. In this case, at least one of $\gq \join \gu$ and
$\gq' \join \gu'$ is saturated; they both may be. Suppose both are saturated.
If they are both atomic saturated, by EP-3, this implies that $\gq \join \gu
\vd f \in \script{F}$ and $\gq' \join \gu' \vd g \in \script{G}$, and
therefore, $(\gq \join \gu);_\gamma (\gq' \join \gu') \vd f;g = \gq || \gq'
\vd f;g \in \script{F};_\gamma \script{G}$, which satisfies EP-3 in $\script{F}
;_\gamma \script{G}$. 

It may also be the case that one is atomic saturated and the other is
(non-atomic) saturated. Suppose that
$\gq' \join \gu'$ 
is atomic saturated and
$\gq \join \gu$
is (non-atomic) saturated. By EP-3, $\gq \join \gu \vd \ol{\alpha}[a] \in
\script{F}$ for some output event $\ol{\alpha}[a]$ where $\alpha \neq \gamma$.
Thus $(\gq \join \gu) ;_\gamma (\gq' \join \gu') \vd \ol{\alpha}[a] =
\gq || \gq' \vd \ol{\alpha}[a] \in \script{F};_\gamma \script{G}$ which
satisfies EP-3.

Now suppose only one is saturated and without loss of generality suppose that
it is $\gq \join \gu$ which is saturated. By EP-3 again, this implies that
$\gq \join \gu \vd \ol{\alpha}[a] \in \script{F}$ for some output event
$\ol{\alpha}[a]$, and thus, $(\gq \join \gu) ;_\gamma (\gq' \join \gu') \vd
\ol{\alpha}[a] = \gq || \gq' \vd \ol{\alpha}[a] \in \script{F} ;_\gamma
\script{G}$ satisfying EP-3.

(ii) Suppose the last event on $\gu_\gamma$ is not atomic. Then, one of $\gq
\join \gu$ and $\gq' \join \gu'$ is saturated. Without loss of generality
suppose that $\gq \join \gu$ is saturated. By EP-3, this implies that
$\gq \join \gu \vd \ol{\alpha}[a] \in \script{F}$ for some output event
$\ol{\alpha}[a]$. If $\alpha = \gamma$ then $\gu * \ol{\gamma}[a]$ and
$\gu * \gamma[a]$ are respectively \script{F} and \script{G}-preantecedents
such that $(\gu * \ol{\gamma}[a]) ;_\gamma (\gu * \gamma[a])$ is saturated
in $\script{F} ;_\gamma \script{G}$. Thus, we may apply this argument on
this preantecedent recursively to get an $\script{F} ;_\gamma
\script{G}$-antecedent which satisfies EP-3.

So suppose that $\alpha \neq \gamma$. Then $(\gq \join \gu) ;_\gamma (\gq'
\join \gu') \vd \ol{\alpha}[a] = \gq || \gq' \vd \ol{\alpha}[a] \in
\script{F};_\gamma \script{G}$ satisfying EP-3.

\item Suppose that $\gp ;_\gamma \gq \vd \ol{\alpha}[a],\ 
\gp';_\gamma \gq' \vd \ol{\beta}[b] \in \script{F} ;_\gamma \script{G}$ with
$(\gp ;_\gamma \gq)_\alpha = (\gp' ;_\gamma \gq')_\alpha$ and
$\gp ;_\gamma \gq \prefix \gp' ;_\gamma \gq'$. Without loss of generality
suppose that $\alpha \in \script{F}$. There are two cases depending on whether
$\beta \in \script{F}$ or $\beta \in \script{G}$.

(i) Suppose $\beta \in \script{F}$. This implies that $\gq$ and $\gq'$ are
compatible \script{G}-preantecedents so that by Lemma~\ref{lem-joinpre}, $\gq
\join \gq'$ is a \script{G}-preantecedent. Similarly, by
Corollary~\ref{cor-joinentail}, $\gp \join \gp' \vd \ol{\beta}[b] \in
\script{F}$. Notice that $(\gp \join \gp')_\alpha = (\gq \join \gq')_\alpha$.
Then, $\gp \vd \ol{\alpha}[a],\ \gp
\join \gp' \vd \ol{\beta}[b] \in \script{F}$ such that $\gp_\alpha =
\gp'_\alpha$ and $\gp \prefix \gp \join \gp'$ so that by EP-4, $\gp \join \gp'
* \ol{\alpha}[a] \vd \ol{\beta}[b] \in \script{F}$. Composing with $\gq \join
\gq'$ yields $(\gp \join \gp' * \ol{\alpha}[a]) ;_\gamma (\gq \join \gq') \vd
\ol{\beta}[b] = (\gp' ;_\gamma \gq') * \ol{\alpha}[a] \vd \ol{\beta}[b] \in
\script{F};_\gamma \script{G}$.

(ii) Suppose $\beta \in \script{G}$. By Corollary~\ref{cor-joinentail},
$\gp \join \gp' \vd \ol{\alpha}[a] \in \script{F}$ and $\gq \join \gq' \vd
\ol{\beta}[b] \in \script{G}$. This implies that $(\gp \join \gp') *
\ol{\alpha}[a]$ is a \script{F}-preantecedent, and as
$(\gp \join \gp')_\alpha =  (\gq \join \gq')_\alpha$ we may form
$(\gp \join \gp' * \ol{\alpha}[a]) ;_\gamma (\gq \join \gq') \vd
\ol{\beta}[b] = (\gp' ;_\gamma \gq') * \ol{\alpha}[a] \vd \ol{\beta}[b]
\in \script{F};_\gamma \script{G}$.

\item Let $\gp ;_\gamma \gq \vd \ol{\alpha}[a] \in \script{F} ;_\gamma 
\script{G}$. If $\beta[b]$ is a legal input event for this entailment, then
it is certainly a legal input event in its corresponding component so that
$(\gp ;_\gamma \gq) * \beta[b] \vd \ol{\alpha}[a] \in \script{F}
;_\gamma \script{G}$.

\item Suppose $(\gp ;_\gamma \gq) * \ol{\beta}[b] \vd \ol{\alpha}[a] \in
\script{F} ;_\gamma \script{G}$. Without loss of generality assume that
$\alpha \in \script{F}$. There are two cases depending on whether $\beta \in
\script{F}$ or $\beta \in \script{G}$.

(i) If $\beta \in \script{F}$ then $\gp * \ol{\beta}[b] \vd \ol{\alpha}[a] \in
\script{F}$, which, by EP-6, implies that $\gp \vd \ol{\alpha}[a] \in
\script{F}$. Therefore, $\gp ;_\gamma \gq \vd \ol{\alpha}[a] \in
\script{F} ;_\gamma \script{G}$.

(ii) If $\beta \in \script{G}$, by the definition of composition, $\gq *
\ol{\beta}[b]$ must be a \script{G}-preantecedent and therefore, $\gq$ must
also be a \script{G}-preantecedent (if $\gq * \ol{\beta}[b]$ is any
\script{G}-antecedent, by EP-6, $\gr$ must also be a \script{G}-antecedent).
Thus, $\gp ;_\gamma \gq \vd \ol{\alpha}[a] \in \script{F} ;_\gamma \script{G}$.

The case where $\alpha \in \script{G}$ is handled similarly.

\item Suppose that $\{\gr * \beta[b_i] \vd \ol{\alpha}[a] \in \script{F}
;_\gamma \script{G} \mid i \in I\}$, where $\{b_i \mid i \in I\}$ is the set
of all possible input events on $\beta$ given $\gr$. Without loss of
generality suppose that $\alpha \in \script{F}$. There are two cases depending
on whether $\beta \in \script{F}$ or $\beta \in \script{G}$.

(i) If $\beta \in \script{F}$ then $\gr * \beta[b_i] \vd \ol{\alpha}[a]$
came from $\gp_i * \beta[b_i] \vd \ol{\alpha}[a] \in \script{F}$ and a
\script{G}-preantecedent $\gq_i$, for each $i \in I$. We know that all the
$\gp_i$ and $\gq_i$ are justified and all the $\gp_i$ are compatible and all
the $\gq$ are compatible (Lemma~\ref{lem-sameingamma}). Thus, we may form
$\bigjoin_i \gp_i$ and $\bigjoin_i \gq_i$ which are \script{F} and
\script{G}-preantecedents respectively (Corollary~\ref{cor-joinentail}).
Clearly then, $(\bigjoin_i \gp_i) * \beta[b_i] \vd \ol{\alpha}[a] \in
\script{F}$ for $i \in I$, and therefore, $\bigjoin_i \gp_i \vd \ol{\alpha}[a]
\in \script{F}$ by EP-7. This implies that $(\bigjoin_i \gp_i) ;_\gamma
(\bigjoin_i \gq_i) \vd \ol{\alpha}[a] = \gr \vd \ol{\alpha}[a] \in \script{F}
;_\gamma \script{G}$.

(ii) If $\beta \in \script{F}$ then $\gr * \beta[b_i] \vd \ol{\alpha}[a]$
came from $\gp_i \vd \ol{\alpha}[a] \in \script{F}$ and a
\script{G}-preantecedent $\gq_i * \beta[b_i]$, for each $i \in I$. By
definition, all the $\gq_i$'s are justified so that each $\gq_i$ is a
\script{G}-preantecedent, and therefore, $\gp_i ;_\gamma \gq_i \vd
\ol{\alpha}[a] = \gr \vd \ol{\alpha}[a] \in \script{F} ;_\gamma \script{G}$,
for any $i \in I$. 

\end{enumerate}

We have shown that the set of entailments $\script{F} ;_\gamma \script{G}$
satisfies EP-1 $\sim$ EP-7 and is, therefore, an extensional process.
\end{proof}

Next it is shown that the identity extensional process behaves correctly.
We begin with a remark on the structure of the identity extensional process.

\begin{remark} \label{rem-idep}
Let $\alpha:X \to^{1_X} \beta:X$ be the identity extensional process.
If $\gp$ is a legal sequence of events for $X$ in the domain that ends at:

\begin{itemize}

\item the last state (atomic) Z, then $\gp || \gp^* \vd 1_Z \in 1_X$.

\item an output event $\ol{\alpha}[a]$, then $(\gp \bs \ol{a}) || \gp^* \vd
\ol{\alpha}[a] \in 1_X$.

\item an input event $\alpha[a]$, then $\gp || (\gp \bs a)^* \vd
\ol{\beta}[a] \in 1_X$.

\end{itemize}
This follows directly from the definition of the identity process and the
translation from a process to an extensional process (see
Section~\ref{sec-pp2ep}).
\end{remark}

\begin{proposition}
$\alpha:X \to^{1_X} \beta:X$ is the identity extensional process on $X$. That
is, given an extensional process $\script{F}: \Gamma \ra \alpha:X,\Delta$, the
composite $\script{F} ;_\alpha 1_X = \script{F}$ (up to a renaming of
channels), and given an extensional process $\script{G}: \Gamma,\beta:X \ra
\Delta$, the composite $1 ;_\beta \script{G} = \script{G}$ (up to a renaming
of channels).
\end{proposition}

\begin{proof}
We will prove that $\script{F} = \script{F} ;_\alpha 1_X$\ ; the statement
$\script{G} = 1_X ;_\beta \script{G}$ is dual.  
The idea will be to prove that $\script{F} \subseteq \script{F} ;_\alpha 1_X$,
which, by Proposition~\ref{prop-subthensame}, then implies that $\script{F} =
\script{F} ;_\alpha 1_X$. To this end, let $\gp \vd \times \in \script{F}$. The
goal is to show that there is a entailment $\gp' \vd \times' \in \script{F}
;_\alpha 1_X$ such that if $\alpha$ is substituted for $\beta$ in this
entailment we get back $\gp \vd \times \in \script{F}$.

There are three cases to consider.

\begin{enumerate}[(i)]

\item $\times = f$, where $f$ is atomic. This implies that the last state
of $\gp_\alpha$ is atomic, say $Z$. As $\gp_\alpha$ is a legal sequence of
events for $X$ (in the codomain) that ends at atomic $Z$, by
Remark~\ref{rem-idep}, $\gp^*_\alpha || \gp_\alpha \vd 1_Z \in 1_X$, and
therefore $\gp ;_\alpha (\gp^*_\alpha || \gp_\alpha) \vd f;1_Z = 
(\gp \bs \alpha) || \gp_\alpha \vd f \in \script{F};_\alpha 1_X$, which is
equivalent to $\gp \vd f$ up to a renaming of $\beta$ to $\alpha$.

\item $\times = \ol{\alpha}[a]$. This implies that $\gp * \ol{\alpha}[a]$ is
a \script{F}-preantecedent and that $\gp^*_\alpha * \alpha[a]$ is a legal
sequence of events for $X$ (in the domain). By our remark above, this implies
that $(\gp^*_\alpha * \alpha[a]) || \gp_\alpha
\vd \ol{\beta}[a] \in 1_X$. Therefore, $(\gp *\ol{\alpha}[a]) ;_\alpha
((\gp^*_\alpha * \alpha[a]) || \gp_\alpha) \vd \ol{\beta}[a] = (\gp \bs \alpha)
|| \gp_\alpha \vd \ol{\beta}[a] \in \script{F} ;_\alpha 1_X$, which is
equivalent to $\gp \vd \ol{\alpha}[a]$ up to a renaming of $\beta$ to $\alpha$.

\item $\times = \ol{\gamma}[c]$, where $\gamma \neq \alpha$. In this case
we must consider events in $\gp$ on $\alpha$. 

\begin{itemize}
\item If there are no events on $\alpha$, then $\gp \vd \ol{\gamma}[c] \in
\script{F} ;_\alpha 1_X$. 

\item Suppose that the last event on $\alpha$ is an output $\ol{\alpha}[a]$.
This implies that it is an input event in $1_X$ and that $\gp^*_\alpha
* \alpha[a]$ is a legal sequence of events for $X$ (in the domain)
so that, by our remark above, $(\gp^*_\alpha * \alpha[a]) || \gp_\alpha \vd
\ol{\beta}[a] \in 1_X$. Hence, $(\gp^*_\alpha * \alpha[a]) || (\gp_\alpha *
\ol{\beta}[a])$ is a $1_X$-preantecedent, and $\gp ;_\alpha ((\gp^*_\alpha *
\alpha[a]) || (\gp_\alpha * \ol{\beta}[a])) \vd \ol{\gamma}[c] = 
(\gp \bs \alpha) || (\gp_\alpha * \ol{\beta}[a]) \vd \ol{\gamma}[c]
\in \script{F} ;_\alpha 1_X$, which is equivalent to $\gp \vd \ol{\gamma}[c]$
up to a renaming of $\beta$ to $\alpha$.

\item Suppose that the last event on $\alpha$ is an input $\alpha[a]$. This
implies that it is an output event in $1_X$ and that $\gp^*_\alpha *
\ol{\alpha}[a]$ is a legal sequence of events for $X$ (in the domain) so that,
by our remark above, $\gp^*_\alpha || (\gp_\alpha * \beta[a]) \vd
\ol{\alpha}[a] \in 1_X$. Then, $(\gp^*_\alpha * \ol{\alpha}[a]) || (\gp_\alpha
* \beta[a])$ is a $1_X$-preantecedent, and therefore, $\gp ;_\alpha
((\gp^*_\alpha  * \ol{\alpha}[a]) || (\gp_\alpha * \beta[a])) \vd
\ol{\gamma}[c] = (\gp \bs \alpha) || (\gp_\alpha * \beta[a]) \vd \ol{\gamma}[c]
\in \script{F} ;_\alpha 1_X$, which is equivalent to $\gp \vd \ol{\gamma}[c]$
up to a renaming of $\beta$ to $\alpha$.
\qed
\end{itemize}
\end{enumerate}
\end{proof}

It is left to show that composition of extensional processes is associative
and satisfies the interchange law.

\begin{proposition}
Composition of extensional processes is associative. That is, if
\[\script{F}:\Gamma \ra \gamma:X,\Delta
\qquad \script{G}:\Gamma',\gamma:X \ra \delta:Y,\Delta'
\qquad \script{H}:\Gamma'',\delta:Y \ra \Delta''
\]
are extensional processes, then
$(\script{F}\ ;_\gamma \script{G})\ ;_\delta \script{H}
= \script{F}\ ;_\gamma (\script{G}\ ;_\delta \script{H})$.
\end{proposition}

\begin{proof}
Composition of atomic entailments is associative as composition is associative
in the underlying polycategory. That is, if $\gp \vd (f ;_\gamma g) ;_\delta h
\in (\script{F}\ ;_\gamma \script{G})\ ;_\delta \script{H}$, then
$\gp \vd f ;_\gamma (g ;_\delta h) \in 
\script{F}\ ;_\gamma (\script{G}\ ;_\delta \script{H})$.

To show that the composition of non-atomic entailments is associative 
suppose that $\alpha \in \script{F}$ and consider the entailment
\[(\gp ;_\gamma \gq) ;_\delta \gr \vd \ol{\alpha}[a] \in
(\script{F}\ ;_\gamma \script{G})\ ;_\delta \script{H}
\]
This implies that $\gr$ is a \script{H}-preantecedent and $\gp ;_\gamma
\gq \vd \ol{\alpha}[a] \in \script{F}\ ;_\gamma \script{G}$ such that
$\gq_\delta = \gr_\delta^*$. This in turn implies that $\gq$ is a
\script{G}-preantecedent and $\gp \vd \ol{\alpha}[a] \in \script{F}$
such that $\gp_\gamma = \gq_\gamma^*$. From these observations it should
now be clear that $\gq ;_\delta \gr$ is a $\script{G} ;_\delta
\script{H}$-preantecedent so that
\[\gp ;_\gamma (\gq ;_\delta \gr) \vd \ol{\alpha}[a] \in
\script{F}\ ;_\gamma (\script{G}\ ;_\delta \script{H})
\]

The argument is similar if $\alpha \in \script{H}$, and slightly different if
$\alpha \in \script{G}$. We describe the latter case and so suppose $\alpha
\in \script{G}$ and consider
\[(\gp ;_\gamma \gq) ;_\delta \gr \vd \ol{\alpha}[a] \in
(\script{F}\ ;_\gamma \script{G})\ ;_\delta \script{H}
\]
This implies that $\gr$ is a \script{H}-preantecedent and $\gp ;_\gamma \gq
\vd \ol{\alpha}[a] \in \script{F}\ ;_\gamma \script{G}$ such that $\gq_\delta
= \gr_\delta^*$. This in turn implies that $\gp$ is a \script{F}-preantecedent
and $\gq \vd \ol{\alpha}[a] \in \script{G}$ such that $\gp_\gamma =
\gq_\gamma^*$. From these observations it should now be clear that
$\gq ;_\delta \gr \vd \ol{\alpha}[a] \in \script{G}\ ;_\delta \script{H}$
and therefore,
\[\gp ;_\gamma (\gq ;_\delta \gr) \vd \ol{\alpha}[a]
\in \script{F}\ ;_\gamma (\script{G}\ ;_\delta \script{H})
\]

This argument may also be used to show that any entailment in $\script{F}\
;_\gamma (\script{G}\ ;_\delta \script{H})$ is an entailment in $(\script{F}\
;_\gamma \script{G})\ ;_\delta \script{H}$. Thus, we may conclude that both
inclusions hold and composition of extensional processes is associative.
\end{proof}

With some small changes the proof of associativity may be used to prove that
composition satisfies the interchange law and so we conclude:

\begin{proposition}
Composition of extensional processes satisfies the interchange law.
\end{proposition}

All the requirements for polycategories are now satisfied proving:

\begin{theorem}
Protocols and extensional processes form a polycategory.
\end{theorem}

The polycategory of protocols and extensional process built over an 
arbitrary polycategory \cat{A} will be denoted by $\script{EP}_\cat{A}$.

\subsection{Poly-sums and poly-products}

The purpose of this section is to show that the polycategory
$\script{EP}_\cat{A}$ has sums and products. We begin by defining sums and
products of extensional processes.

Let $\script{F}_i: \Gamma,\gamma:X_i \ra \Delta$, for $i \in I$, be
extensional processes. The \textbf{sum} (or \textbf{coproduct}) on $\gamma$
of the $\script{F}_i$'s, denoted $\gamma\<\script{F}_i\>_{i \in I}: \Gamma,
\gamma:\sum_{i \in I} X_i \ra \Delta$, is defined to be the set $\{\gamma[i] *
\script{F}_i \mid i \in I\}$. For example, if $\gp \vd \times \in \script{F}_k$
then $\gamma[k] * \gp \vd \times \in \gamma\<\script{F}_i\>_{i \in I}$.

The \textbf{product} is constructed dually, i.e., if $\script{G}_j:\Gamma \ra
\delta:Y_j,\Delta$, for $j \in J$, are extensional processes, then the
product of the $\script{G}_j$'s on $\delta$, denoted
$\delta\<\script{G}_j\>_{j \in J}:\Gamma \ra \delta:\prod_{j \in J} Y_j,
\Delta$, is defined to be the set $\{\delta[i] * \script{G}_j \mid j \in J\}$.

The extensional process for the $k^{th}$ injection on $\gamma$ is the
set $\script{B}_k = \ol{\gamma}[k] * 1_{X_k}$ where $\gamma$ is the codomain
channel of $1_{X_k}$. The extensional process for the $k^{th}$
projection on $\gamma$ is $\script{P}_k = \ol{\gamma}[k] * 1_{X_k}$ (the same
set as the injection), where $\gamma$ is the domain channel of $1_{X_k}$.

We will now work with the sum; the dual properties hold for the product.
Note that it still must be shown that what we have defined as the
``sum'' is actually a sum. First, however, it is necessary to show
that the ``sum'' satisfies the requirements of an extensional process.

\begin{proposition}
$\gamma\<\script{F}_i\>_{i \in I}$ is an extensional process.
\end{proposition}

\begin{proof}
We must show that $\gamma\<\script{F}_i\>_{i \in I}$ satisfies EP-1 $\sim$
EP-7.

\begin{enumerate}[EP-1 ]

\item Suppose that $\gamma[k] * \gp$, for some $k \in I$, is a $\gamma\<
\script{F}_i\>_{i \in I}$-antecedent and $\ol{\alpha}[a]$ is an output event
in $\gamma[k] * \gp$. Then $\ol{\alpha}[a]$ must be justified in
$\script{F}_k$, say by $\gu$, so that $\gamma[k] * \gu$ justifies it
in $\gamma\<\script{F}_i\>_{i \in I}$. Thus, all 
$\gamma\<\script{F}_i\>_{i \in I}$-antecedents are
$\gamma\<\script{F}_i\>_{i \in I}$-justified.

\item Suppose that $\gamma[k] * \gp \vd \ol{\alpha}[a],\ \gamma[l] * \gq 
\vd \ol{\alpha}[b] \in \gamma\<\script{F}_i\>_{i \in I}$, for some $k,l \in
I$, with $(\gamma[k] * \gp)_\alpha = (\gamma[l] * \gq)_\alpha$  and
$\gamma[k] * \gp \smile \gamma[l] * \gq$. Since they are compatible, $k=l$,
and therefore, $\gp \vd \ol{\alpha}[a]$, $\gq \vd \ol{\alpha}[b] \in
\script{F}_k$ with $\gp_\alpha = \gq_\alpha$ and $\gp \smile \gq$, and
as $\script{F}_k$ satisfies EP-2, $a=b$.

\item Let $\gq$ be a saturated $\gamma\<\script{F}_i\>_{i \in
I}$-preantecedent. It must be of the form $\gamma[k] * \gq'$, for some $k \in
I$, where $\gq'$ is a saturated $\script{F}_k$-preantecedent. By EP-3, $\gq'
\vd \times \in \script{F}_k$, and thus, $\gamma[k] * \gq \vd \times \in
\gamma\<\script{F}_i\>_{i \in I}$ satisfying EP-3.

\item Suppose that $\gamma[k] * \gp \vd \ol{\alpha}[a],\ \gamma[l] * \gq
\vd \ol{\beta}[b] \in \gamma\<\script{F}_i\>_{i \in I}$, for some $k,l \in I$,
with $(\gamma[k] * \gp)_\alpha = (\gamma[l] * \gq)_\alpha$ and
$\gamma[k] * \gp \prefix \gamma[l] * \gq$. It is clear then that $k=l$, and
therefore, $\gp \vd \ol{\alpha}[a]$, $\gq \vd \ol{\beta}[b] \in
\script{F}_k$ with $\gp_\alpha = \gq_\alpha$ and $\gp \prefix \gq$. As
$\script{F}_k$ satisfies EP-4, $\gq * \ol{\alpha}[a] \vd \ol{\beta}[a] \in
\script{F}_k$ and thus, $\gamma[k] * \gq * \ol{\alpha}[a] \vd \ol{\beta}[b]
\in \gamma\<\script{F}_i\>_{i \in I}$.

\item Let $\gamma[k]*\gp \vd \ol{\alpha}[a] \in \gamma\<\script{F}_i\>_{i \in
I}$, for some $k \in I$. If $\beta[b]$ is a legal input event for this
entailment, then it is certainly a legal input event for $\gp \vd
\ol{\alpha}[a] \in \script{F}_k$. Thus, $\gp * \beta[b] \vd \ol{\alpha}[a] \in
\script{F}_k$ and so, $\gamma[k] * \gp * \beta[b] \vd \ol{\alpha}[a] \in
\gamma\<\script{F}_i\>_{i \in I}$.

\item Let $\gamma[k] * \gp * \ol{\beta}[b] \vd \ol{\alpha}[a] \in
\gamma\<\script{F}_i\>_{i \in I}$, for some $k \in I$. Then, $\gp *
\ol{\beta}[b] \vd \ol{\alpha}[a] \in \script{F}_k$, and as $\script{F}_k$
satisfies EP-6, $\gp \vd \ol{\alpha}[a] \in \script{F}_k$, and thus,
$\gamma[k] * \gp \vd \ol{\alpha}[a] \in \gamma\<\script{F}_i\>_{i \in I}$.

\item Suppose $\{\gamma[k] * \gp * \beta[b_i] \vd \ol{\alpha}[a] \mid i \in
I\} \in \gamma\<\script{F}_i\>_{i \in I}$, for some $k \in I$, where $\{b_i
\mid i \in I\}$ is the set of all possible input events on $\beta$ given
$\gamma[k] * \gp$. This implies that $\{\gp * \beta[b_i] \vd \ol{\alpha}[a]
\mid i \in I\} \in \script{F}_k$, and as $\script{F}_k$ satisfies EP-7,
$\gp \vd \ol{\alpha}[a] \in \script{F}_k$, and therefore, $\gamma[k] * \gp \vd
\ol{\alpha}[a] \in \gamma\<\script{F}_i\>_{i \in I}$. 
\qed
\end{enumerate}
\end{proof}

We now show that $\gamma\<\script{F}_i\>_{i \in I}$ is a coproduct of
extensional processes.

\begin{proposition} 
Let $\script{F}:\Gamma,\gamma:X_i \ra \Delta$, for $i \in I$, be
extensional processes. Then $\gamma\<\script{F}_i\>_{i \in I}:\Gamma,
\gamma:\sum_{i \in I} X_i \ra \Delta$ is the poly-coproduct of
the $\script{F}_i$'s, $i \in I$.
\end{proposition}

\begin{proof}
We must show that the following four properties hold:
\begin{enumerate}[{\upshape (i)}]
\item $\gamma\<\script{B}_i\>_{i \in I} = 1_{\sum\limits_{i \in I} X_i}$
\item $\script{B}_k\ ;_\gamma \gamma\<\script{F}_i\>_{i \in I} = \script{F}_k$
\item $\gamma\<\script{F}_i\>_{i \in I}\ ;_\alpha \script{H} =
\gamma\<\script{F}_i\ ;_\alpha \script{H}\>_{i \in I}$, where $\gamma \neq
\alpha$
\item $\script{H}\ ;_\alpha \gamma\<\script{F}_i\>_{i \in I} =
\gamma\<\script{H}\ ;_\alpha \script{F}_i\>_{i \in I}$, where $\alpha \neq
\gamma$
\end{enumerate}
We take each one in turn.

\begin{enumerate}[{\upshape (i)}]

\item The goal will be to show that for any $\script{F}_i:\Gamma, \gamma: X_i
\ra \Delta$, where $i \in I$, and for $\script{B}_k: \alpha:X_k \ra
\gamma:\sum_{i \in I} X_i$
\[
\alpha\<\script{B}_i\>_{i \in I}\ ;_\gamma \gamma\<\script{F}_i\>_{i \in I} =
\alpha\<\script{F}_i\>_{i \in I}
\]
which shows that $\alpha\<\script{B}_i\>_{i \in I}$ acts as the identity on
$\sum_{i \in I} X_i$. Consider an entailment in $\gq \vd \times
\in \alpha\<\script{B}_i\>_{i \in I}\ ;_\gamma
\gamma\<\script{F}_i\>_{i \in I}$. There are three cases to consider
depending on whether $\x$ is atomic, an output event on $\alpha$, or an
output event on a channel other than $\alpha$. However, these cases are
subsumed if, in the following, we simply consider $\gr$ to be a
$1_{X_k}$-preantecedent and $\gp$ to be a $\script{F}_k$-preantecedent, such
that one is an antecedent (it does not matter which). Then, by the definition
of $\script{B}_k$, composition, and coproducts, $\gq \vd \times$ must be of the
form:
\begin{align*}
\gq \vd \times 
&= (\alpha[k] * \ol{\gamma}[k] * \gr) ;_\gamma (\gamma[k] * \gp) \vd \times \\
&= (\alpha[k] * \gr) ;_\gamma \gp \vd \times \\
&= \alpha[k] * (\gr ;_\gamma \gp) \vd \times 
\end{align*}
By the definition of the identity extensional process, $\gr ;_\gamma \gp \in
\script{F}_k$, and so $\alpha[k] * (\gr ;_\gamma \gp) \vd \times \in
\gamma\<\script{F}_i\>_{i \in I}$ (up to renaming $\alpha$ to $\gamma$).
Therefore, $\alpha\<\script{B}_i\>_{i \in I} ;_\gamma
\gamma\<\script{F}_i\>_{i \in I} \subseteq \gamma\<\script{F}_i\>_{i \in I}$,
and, by Proposition~\ref{prop-subthensame}, $\alpha\<\script{B}_i\>_{i \in I}
;_\gamma \gamma\<\script{F}_i\>_{i \in I} = \gamma\<\script{F}_i\>_{i \in I}$.
Thus, $\alpha\<\script{B}_i\>_{i \in I}$ is the identity on
$\sum_{i \in I} X_i$.

\item Let $\gq \vd \times \in \script{B}_k ;_\gamma \gamma\<\script{F}_i\>_{i
\in I}$. Again, in the following, assume that $\gr$ is a
$1_{X_k}$-preantecedent and that $\gp$ is a $\script{F}_k$-preantecedent. By
the definition of $\script{B}_k$, composition, and coproducts, $\gq =
(\ol{\gamma}[k] * \gr) ;_\gamma (\gamma[k] * \gp) = \gr ;_\gamma \gp$, and
hence, $\gq \vd \x \in \script{F}_k$ (up to renaming of channels). Therefore,
by Proposition~\ref{prop-subthensame}, $\script{B}_k ;_\gamma
\gamma\<\script{F}_i\>_{i \in I} = \script{F}_k$.

\item Given a $\gamma\<\script{F}_i\>_{i \in I}$-antecedent $\gamma[k] * \gp$,
for some $k \in I$, and a \script{H}-preantecedent $\gq$, the following 
equality holds
\[
(\gamma[k] * \gp) ;_\alpha \gq = \gamma[k] *(\gp ;_\alpha \gq)
\qquad \text{(when $\alpha \neq \gamma$)} 
\] 
so that any entailment $(\gamma[k] * \gp) ;_\alpha \gq \vd \times \in
\gamma\<\script{F}_i\>_{i \in I} ;_\alpha \script{H}$, may be written as
$\gamma[k] *(\gp ;_\alpha \gq) \vd \times \in \gamma\<\script{F}_i ;_\alpha
\script{H}\>_{i \in I}$ and vise versa, which shows that
$\gamma\<\script{F}_i\>_{i \in I} ;_\alpha \script{H} = \gamma\<\script{F}_i
;_\alpha \script{H}\>_{i \in I}$.

\item That $\script{H} ;_\alpha \gamma\<\script{F}_i\>_{i \in I} =
\gamma\<\script{H} ;_\alpha \script{F}_i\>_{i \in I}$ is true follows from the
equality
\[
\gq ;_\alpha (\gamma[k] * \gp) = \gamma[k] * (\gq ;_\alpha \gp) 
\qquad \text{(when $\alpha \neq \gamma$)}
\]
and an argument similar to (iii) above.
\qed
\end{enumerate}
\end{proof}

With duality this now proves:

\begin{proposition}
$\script{EP}_\cat{A}$ has finite poly-sums and poly-products.
\end{proposition}

\subsection{Extensional processes are soft}

The goal of this section is to prove the following theorem:

\begin{theorem}
$\script{EP}_\cat{A}$ is the free polycategory generated from \cat{A} under
(finite) sums and products.
\end{theorem}

Notice that the inclusion functor $\cat{A} \to^{\mathcal{I}} \script{EP}_\cat{A}$
is a full inclusion and the objects of $\script{EP}_\cat{A}$ are the objects
of $\Sigma\Pi_\cat{A}$ (i.e., they are generated under sums and products from
the objects of \cat{A}). Thus, if it can be shown that the inclusion
$\mathcal{I}$ is soft, the Whitman Theorem, Theorem~\ref{thm-whitman},
may be used to show that $\script{EP}_\cat{A}$ is equivalent to
$\Sigma\Pi_\cat{A}$, and hence, is the free polycategory generated from
\cat{A} under finite sums and products. Thus:

\begin{proposition}
The inclusion morphism of polycategories $\cat{A} \to^{\mathcal{I}}
\script{EP}_\cat{A}$ is soft. 
\end{proposition}

Recall that in order to show that $\mathcal{I}$ is soft, we need to show that
$\mathcal{I}$ is a semi-soft extension and a soft extension. That is, that the
map
\[\sum_{(\alpha,i) \in \oi(\Gamma;\Delta)}
\negthickspace\negthickspace
[\Hom(\Gamma;\Delta)]_{(\alpha,i)} \to^{\{\ol{\alpha}[i] *
\script{F}\}_{(\alpha,i)}} \Hom(\Gamma;\Delta)
\]
is a surjection, and that
\[
\xymatrix@R=8ex{{\sum\limits_{
((\alpha,i), (\beta,j)) \in \oi(\Gamma,\Delta)}
\negthickspace\negthickspace\negthickspace
[\Hom(\Gamma ; \Delta)]_{((\alpha,i),(\beta,j))}}
\ar@<2ex>[d]^-{\sum\limits_{(\beta,j)} \{\ol{\alpha}[i] *
\script{F}\}_{((\alpha,i),(\beta,j))}}
\ar@<-2ex>[d]_-{\sum\limits_{(\alpha,i)}\{\ol{\beta}[j] *
\script{F}\}_{((\alpha,i),(\beta,j))}}
\\
\sum\limits_{(\gamma,k) \in \oi(\Gamma,\Delta)}
[\Hom(\Gamma ; \Delta)]_{(\gamma,k)}
\ar[d]^-{\{\ol{\gamma}[k] * \script{F}'\}_{(\gamma,k)}}
\\
\Hom(\Gamma ; \Delta)}
\]
is a coequalizer diagram.

\begin{proof}
To see that $\mathcal{I}$ is a semi-soft extension consider the output
poly-hom-set $\Hom(\Gamma;\Delta)$. Any behaviour in $\Hom(\Gamma;\Delta)$
started out saturated and therefore, must be of the form $\ol{\alpha}[a] *
\script{F}$ where $\script{F} \in \sum_{(\alpha,i)}
[\Hom(\Gamma;\Delta)]_{(\alpha,i)}$. Thus, the map
\[\sum_{(\alpha,i) \in \oi(\Gamma;\Delta)}
\negthickspace\negthickspace
[\Hom(\Gamma;\Delta)]_{(\alpha,i)} \to^\{\ol{\alpha}[i] *
\script{F}\}_{(\alpha,i)} \Hom(\Gamma;\Delta)
\]
is a surjection and $\mathcal{I}$ is a semi-soft extension.

That $\mathcal{I}$ is soft follows follows from an argument similar to
the syntactic case. Explicitly, suppose the diagram has a coequalizer $(q,Q)$
and let $h:Q \ra \Hom(\Gamma;\Delta)$ be the unique map such that $q;h = 
\{\ol{\gamma}[k] * \script{F}_{(\gamma,k)}\}$. But $\{\ol{\gamma}[k] *
\script{F}_{(\gamma,k)}\}$ is a surjection, and therefore, so is $h$. 

Now consider the two following equivalent extensional processes
$\ol{\alpha}[i] * \script{F}$ and $\ol{\beta}[j] * \script{G}$
in $\Hom(\Gamma;\Delta)$. That they are equivalent means that they must be of
the form
\[\ol{\alpha}[i] * (\ol{\beta}[j] * \script{F'}) \qquad \text{and} \qquad
\ol{\beta}[j] * (\ol{\alpha}[i] * \script{G'})
\]
respectively. However, this implies that the extensional processes \script{F'}
and \script{G'} must be equivalent in 
$\sum\limits_{((\alpha,i), (\beta,j))}
\negthickspace\negthickspace
[\Hom(\Gamma ; \Delta)]_{((\alpha,i),(\beta,j))}$, and therefore, must be
coequalized in $Q$ establishing that $h$ is an injection. 

Therefore, the map $h$ is a bijection and $Q \iso \Hom(\Gamma;\Delta)$
establishing that
\[\left(\{\ol{\gamma}[k] * \script{F}_{(\gamma,k)}\}_{(\gamma,k)},\
\Hom(\Gamma ; \Delta) \right)
\]
is the coequalizer.
\end{proof}

\subsection{The additive units}

In this section we explain how the additive units are handled in this
system as they are handled in a somewhat subtle manner. (Indeed, when
Dominique Hughes first asked about the units, and we tried to explain it to
him, we instantly got lost in the subtleties.)

Recall that we are using the notation $\sum_\emptyset = 0$ and
$\prod_\emptyset = 1$ to denote the initial and final objects respectively.
In what follows it will be shown that any map in which the domain contains
the initial object and the codomain contains the final object is the empty
set of entailments.

Notice that if, for any map, the domain contains the initial object or the
codomain contains the final object, i.e.,
\[\Gamma,0 \vd \Delta  \qquad \text{or} \qquad \Gamma \vd 1,\Delta
\]
then the frontier of the empty behaviour has a source role (as the
initial object in the domain and the final object in the codomain have a
source role), and so cannot be the antecedent of any entailment.
Thus, the empty set of entailments results. Recall that in this situation
the empty set does constitute an extensional process
(Proposition~\ref{prop-emptyep}). Moreover, notice that if we were to
compose an extensional process with either
\[\Gamma,0 \vd \Delta  \qquad \text{or} \qquad \Gamma \vd 1,\Delta
\]
the resulting frontier would also contain a protocol with a source role and
no entailments could result. 

The roles are also used to determine where a given protocol is an initial or
a final object. That is, the initial (resp. final) object may be
``hidden'' within some structure and the roles are used to determine this.
For example, $0 \iso 0+0$ is the initial object, in the domain, its role
is calculated as
\[\xymatrix@R=4ex@C=3ex@M=0ex{ 
& \bullet^\mathbf{0} \ar@{--}[dl] \ar@{--}[dr] \\
\bullet^\mathbf{0} && \bullet^\mathbf{0}}
\]
Similarly for $1 \iso 1 \times 1$ in the codomain or any other ``hidden''
initial or final objects.

Consider a protocol in the domain of the form
\[\xymatrix@R=4ex@C=3ex@M=0ex{ 
&& \circ^+ \ar@{--}[dl] \ar@{-}[dr] \\ & 
\circ^{\bf 0} \ar@{--}[dl] \ar@{--}[dr] && A^+ \\
\circ^{\bf 0} && \circ^{\bf 0}}
\]
which contains an initial object and so is really just the protocol $A$. As
we can see from the tree, this is indeed the case as the only legal transition
is the one which leads to $A$.

In the free polycategory with sums and products one does not expect a map
$1 \ra 0$. This map is prevented from occurring by the rule EP-3.
As both $1$ in the domain and $0$ in the codomain are ``output'' protocols
the empty behaviour is saturated in $1 \ra 0$. As the empty behaviour is both
a preantecedent and a saturated preantecedent it fits the premise of EP-3.
However, it is clear that the empty behaviour cannot be made into an entailment
as there is nothing to output, and so, this ``set'' of entailments is not an
extensional process. 

\chapter{Conclusion and Further Directions} \label{chap-conclusion}

In this thesis we began by introducing a logic $\Sigma\Pi$, and constructing
a (syntactic) polycategory of protocols and processes from this logic. This
polycategory was shown to be the free polycategory with finite sums and
finite products. To characterize polycategories of this type we chose to
use Joyal's notion of softness, and so it was necessary to extend this notion
to the polycategorical case. Following the lead of Joyal and Cockett and
Seely, we prove a ``Whitman theorem'' which gives the characterization.
Next, it is shown how proofs in our logic can be interpreted as concurrent
channel-based processes by providing a process semantics. These process are
organized into a polycategory of protocols and extensional processes. Using
softness, we show that this polycategory is equivalent to the syntactic
polycategory. This then establishes that every extensional process is the
denotation of a unique cut-free proof in $\Sigma\Pi$, proving that this
model is full and faithfully complete.

\paragraph{Further directions}

The next obvious step is to add the multiplicative connectives of linear
logic into our interpretation. We have a fairly detailed idea of how this
can be accomplished, but it still needs to be fully formalized. Additionally,
from the viewpoint of process semantics, one would also like to model infinite
processes. This leads into fixed point logics and circular
proofs~\cite{santocanale02:circular} where there is considerable gaps to be
filled.

\bibliographystyle{alpha}
\bibliography{craig}

\appendix

\chapter{Resolving Critical Pairs} \label{chap-rcp}

The following notation will be used to reduce the number of cases of 
certain rewrites. Let $\delta(g)$ denote any of the following morphism,
\[\delta\{a_i \mapsto g_i\}, \quad \delta(a_i \mapsto g_i), \quad
\ola{\delta}[a](g), \quad \ora{\delta}[a](g)
\]

Then, $f ;_\gamma \delta(g) \Lra \delta(f ;_\gamma g)$ denotes
\[\begin{array}{cc}
f ;_\gamma \delta\{a_i \mapsto g_i\} \Lra \delta\{a_i \mapsto f ;_\gamma g_i\},
\quad &
f ;_\gamma \ola{\delta}[a](g) \Lra \ola{\delta}[a](f ;_\gamma g), \medskip\\
f ;_\gamma \delta(a_i \mapsto g_i) \Lra \delta(a_i \mapsto f ;_\gamma g_i),
\quad &
f ;_\gamma \ora{\delta}[a](g) \Lra \ora{\delta}[a](f ;_\gamma g)
\end{array}
\]
and similarly for $\delta(f) ;_\gamma g \Lra \delta(f ;_\gamma g)$.

The resolutions are as follows. The dual rewrite, if there is one, is
indicated in $[-]$.

\begin{description}
\item [(1)-(2)] obvious.

\item [{(1)-(3) [(2)-(4)]}] 
\[
\xymatrix{& \alpha\{a_i \mapsto f_i\}_i ;_\gamma 1 \ar@{=>}[dl]_-{(1)}
\ar@{=>}[dr]^-{(3)} \\
\alpha\{a_i \mapsto f_i\}_i && \alpha\{a_i \mapsto f_i ;_\gamma 1\}_i
\ar@{=>}[ll]^-{\alpha\{(1)\}} }
\]

\item [{(1)-(5) [(2)-(6)]}]
\[
\xymatrix{& \ora{\alpha}[a](f) ;_\gamma 1 \ar@{=>}[dl]_-{(1)}
\ar@{=>}[dr]^-{(5)} \\
\ora{\alpha}[a](f) && \ora{\alpha}(f ;_\gamma 1)
\ar@{=>}[ll]^-{\ora{\alpha}((1))}}
\]

\item [{(1)-(7) [(2)-(8)]}] 
\[
\xymatrix{& \ola{\alpha}[a](f) ;_\gamma 1 \ar@{=>}[dl]_-{(1)}
\ar@{=>}[dr]^-{(7)} \\
\ola{\alpha}[a](f) && \ola{\alpha}(f ;_\gamma 1)
\ar@{=>}[ll]^-{\ola{\alpha}((1))}}
\]

\item [{(1)-(9) [(2)-(10)]}]
\[
\xymatrix{& \alpha(a_i \mapsto f_i)_i ;_\gamma 1 \ar@{=>}[dl]_-{(1)}
\ar@{=>}[dr]^-{(9)} \\
\alpha(a_i \mapsto f_i)_i && \alpha(a_i \mapsto f_i ;_\gamma 1)_i
\ar@{=>}[ll]^-{\alpha((1))} }
\]

\end{description}
This handles all critical pairs involving (1) and (2). We now look at any
critical pairs involving (3) and (4).

\begin{description}
\item [(3)-(4)]
\[\xymatrix@M=1ex@C=20ex@R=10ex@!0{
& \alpha\{a_i \mapsto f_i\}_i ;_\gamma \beta(b_j \mapsto g_j)_j
\ar@{=>}[dl]_{(3)} \ar@{=>}[dr]^{(4)} \\
\alpha\{a_i \mapsto f_i ;_\gamma \beta(b_j \mapsto g_j)_j\}_i
\ar@{=>}[d]_{\alpha\{(4)\}} 
&& \beta(b_j \mapsto \alpha\{a_i \mapsto f_i\}_i ;_\gamma g_j)_j 
\ar@{=>}[d]^{\beta((3))} \\
\alpha\{a_i \mapsto \beta(b_j \mapsto f_i ;_\gamma g_j)_j\}_i
\ar@{|=|}[rr]_{(19)} 
&& \beta(b_j \mapsto \alpha\{a_i \mapsto f_i ;_\gamma g_j\}_i)_j}
\]

\item [{(3)-(6) [(4)-(5)]}] 
\[\xymatrix@M=1ex@C=20ex@R=10ex@!0{
& \alpha\{a_i \mapsto f_i\}_i ;_\gamma \ola{\beta}[b](g) 
\ar@{=>}[dl]_{(3)} \ar@{=>}[dr]^{(6)} \\
\alpha\{a_i \mapsto f_i ;_\gamma \ola{\beta}[b](g)\}_i
\ar@{=>}[d]_{\alpha\{(6)\}} 
&& \ola{\beta}[b](\alpha\{a_i \mapsto f_i\}_i ;_\gamma g)
\ar@{=>}[d]^{\ola{\beta}((3))} \\
\alpha\{a_i \mapsto \ola{\beta}[b](f_i ;_\gamma g)\}_i
\ar@{|=|}[rr]_{(17)} 
&& \ola{\beta}[b](\alpha\{a_i \mapsto f_i ;_\gamma g\}_i)}
\]

\item [{(3)-(8) [(4)-(7)]}] 
\[\xymatrix@M=1ex@C=20ex@R=10ex@!0{
& \alpha\{a_i \mapsto f_i\}_i ;_\gamma \ora{\beta}[b](g) 
\ar@{=>}[dl]_{(3)} \ar@{=>}[dr]^{(8)} \\
\alpha\{a_i \mapsto f_i ;_\gamma \ora{\beta}[b](g)\}_i 
\ar@{=>}[d]_{\alpha\{(8)\}} 
&& \ora{\beta}[b](\alpha\{a_i \mapsto f_i\}_i ;_\gamma g) 
\ar@{=>}[d]^{\ora{\beta}((3))} \\
\alpha\{a_i \mapsto \ora{\beta}[b](f_i ;_\gamma g)\}_i 
\ar@{|=|}[rr]_{(15)} 
&& \ora{\beta}[b](\alpha\{a_i \mapsto f_i ;_\gamma g\}_i)}
\]

\item [{(3)-(10) [(4)-(9)]}] 
\[\xymatrix@M=1ex@C=20ex@R=10ex@!0{
& \alpha\{a_i \mapsto f_i\}_i ;_\gamma \beta\{b_j \mapsto g_j\}_j
\ar@{=>}[dl]_{(3)} \ar@{=>}[dr]^{(10)} \\
\alpha\{a_i \mapsto f_i ;_\gamma \beta\{b_j \mapsto g_j\}_j\}_i 
\ar@{=>}[d]_{\alpha\{(10)\}} 
&& \beta\{b_j \mapsto \alpha\{a_i \mapsto f_i\}_i ;_\gamma g_j\}_j
\ar@{=>}[d]^{\beta\{(3)\}} \\
\alpha\{a_i \mapsto \beta\{b_j \mapsto f_i ;_\gamma g_j\}_j\}_i 
\ar@{|=|}[rr]_{(13)} 
&& \beta\{b_j \mapsto \alpha\{a_i \mapsto f_i ;_\gamma g_j\}_i\}_j}
\]

\item [{(3)-(13) [(4)-(14)]}]
\[\xymatrix@M=1ex@C=20ex@R=10ex@!0{
& \alpha\{a_i \mapsto \beta\{b_j \mapsto f_{ij}\}_j\}_i ;_\gamma g 
\ar@{=>}[dl]_{(3)} \ar@{|=|}[dr]^{(13);1} \\
\alpha\{a_i \mapsto \beta\{b_j \mapsto f_{ij}\}_j ;_\gamma g\}_i 
\ar@{=>}[dd]_{\alpha\{(3)\}} 
&& \beta\{b_j \mapsto \alpha\{a_i \mapsto f_{ij}\}_i\}_j ;_\gamma g 
\ar@{=>}[d]^{(3)} \\ 
&& \beta\{b_j \mapsto \alpha\{a_i \mapsto f_{ij}\}_i ;_\gamma g\}_j
\ar@{=>}[d]^{\beta\{(3)\}} \\ 
\alpha\{a_i \mapsto \beta\{b_j \mapsto f_{ij} ;_\gamma g\}_j\}_i 
\ar@{|=|}[rr]_{(13)} 
&& \beta\{b_j \mapsto \alpha\{a_i \mapsto f_{ij} ;_\gamma g\}_i\}_j}
\]

\item [{(3)-(15) [(4)-(16)]}] 
\[\alpha\{a_i \mapsto \ora{\beta}[b_k](f_i) ;_\gamma g\}_i \Lla 
\alpha\{a_i \mapsto \ora{\beta}[b_k](f_i)\}_i ;_\gamma g \pc
\ora{\beta}[b_k](\alpha\{a_i \mapsto f_i\}_i) ;_\gamma g 
\]
There are two subcases to consider.

(a) $\beta \neq \gamma$.
\[\xymatrix@M=1ex@C=20ex@R=10ex@!0{
& \alpha\{a_i \mapsto \ora{\beta}[b](f_i)\}_i ;_\gamma g 
\ar@{=>}[dl]_{(3)} \ar@{|=|}[dr]^{(15);1} \\
\alpha\{a_i \mapsto \ora{\beta}[b](f_i) ;_\gamma g\}_i 
\ar@{=>}[dd]_{\alpha\{(5)\}} 
&& \ora{\beta}[b](\alpha\{a_i \mapsto f_i\}_i) ;_\gamma g
\ar@{=>}[d]^{(5)} \\ 
&& \ora{\beta}[b](\alpha\{a_i \mapsto f_i\}_i ;_\gamma g)
\ar@{=>}[d]^{\ora{\beta}((3))} \\ 
\alpha\{a_i \mapsto \ora{\beta}[b](f_i ;_\gamma g)\}_i 
\ar@{|=|}[rr]_{(15)} 
&& \ora{\beta}[b](\alpha\{a_i \mapsto f_i ;_\gamma g\}_i)}
\]

(b) $\beta = \gamma$. Here we may assume that $g = \gamma\{b_j \mapsto g_j
\}_j$.
\[\xymatrix@M=1ex@C=20ex@R=10ex@!0{
& \alpha\{a_i \mapsto \ora{\gamma}[b_k](f_i)\}_i ;_\gamma \gamma\{b_j \mapsto
g_j\}_j \ar@{=>}[dl]_{(3)} \ar@{|=|}[dr]^{(15);1} \\
\alpha\{a_i \mapsto \ora{\gamma}[b_k](f_i) ;_\gamma \gamma\{b_j \mapsto g_j
\}_j\}_i \ar@{=>}[d]_{\alpha\{(11)\}} 
&& \ora{\gamma}[b_k](\alpha\{a_i \mapsto f_i\}_i) ;_\gamma \gamma\{b_j \mapsto
g_j\}_j \ar@{=>}[d]^{(11)} \\ 
\alpha\{a_i \mapsto f_i ;_\gamma g_k\}_i 
&& \alpha\{a_i \mapsto f_i\}_i ;_\gamma g_k
\ar@{=>}[ll]^{(3)}}
\]

\item [{(3)-(17) [(4)-(18)]}]
\[\xymatrix@M=1ex@C=20ex@R=10ex@!0{
& \alpha\{a_i \mapsto \ola{\beta}[b](f_i)\}_i ;_\gamma g 
\ar@{=>}[dl]_{(3)} \ar@{|=|}[dr]^{(17);1} \\
\alpha\{a_i \mapsto \ola{\beta}[b](f_i) ;_\gamma g\}_i 
\ar@{=>}[dd]_{\alpha\{(7)\}} 
&& \ola{\beta}[b](\alpha\{a_i \mapsto f_i\}_i) ;_\gamma g
\ar@{=>}[d]^{(7)} \\ 
&& \ola{\beta}[b](\alpha\{a_i \mapsto f_i\}_i ;_\gamma g)
\ar@{=>}[d]^{\ola{\beta}((3))} \\ 
\alpha\{a_i \mapsto \ola{\beta}[b](f_i ;_\gamma g)\}_i 
\ar@{|=|}[rr]_{(17)} && 
\ola{\beta}[b](\alpha\{a_i \mapsto f_i ;_\gamma g\}_i)}
\]

\item [{(3)-(19) [(4)-(19)]}] 
\[\xymatrix@M=1ex@C=20ex@R=9ex@!0{
& \alpha\{a_i \mapsto \beta(b_j \mapsto f_{ij})_j\}_i ;_\gamma \delta(g)
\ar@{=>}[dl] \ar@{|=|}[dr] \\
\alpha\{a_i \mapsto \beta(b_j \mapsto f_{ij})_j ;_\gamma \delta(g)\}_i &&
\beta(b_j \mapsto \alpha\{a_i \mapsto f_{ij}\}_i)_j ;_\gamma \delta(g)}
\]
There are three subcases to consider.

(a) $\beta \neq \gamma$
\[\xymatrix@M=1ex@C=20ex@R=10ex@!0{
& \alpha\{a_i \mapsto \beta(b_j \mapsto f_{ij})_j\}_i ;_\gamma \delta(g)
\ar@{=>}[dl]_{(3)} \ar@{|=|}[dr]^{(19);1} \\
\alpha\{a_i \mapsto \beta(b_j \mapsto f_{ij})_j ;_\gamma \delta(g)\}_i
\ar@{=>}[dd]_{\alpha\{(3)\}} 
&& \beta(b_j \mapsto \alpha\{a_i \mapsto f_{ij}\}_i)_j ;_\gamma \delta(g)
\ar@{=>}[d]^{(3)} \\ 
&& \beta(b_j \mapsto \alpha\{a_i \mapsto f_{ij}\}_i ;_\gamma \delta(g))_j
\ar@{=>}[d]^{\beta((3))} \\ 
\alpha\{a_i \mapsto \beta(b_j \mapsto f_{ij} ;_\gamma \delta(g))_j)\}_i 
\ar@{|=|}[rr]_{(13)} && 
\beta(b_j \mapsto \alpha\{a_i \mapsto f_{ij} ;_\gamma \delta(g)\}_i)_j}
\]

(b) $\beta = \gamma$ and $\delta \neq \gamma$
\[\xymatrix@M=1ex@C=22ex@R=10ex@!0{
& \alpha\{a_i \mapsto \gamma(b_j \mapsto f_{ij})_j\}_i ;_\gamma \delta(g)
\ar@{=>}[dl]_{(3)} \ar@{|=|}[dr]^{(19);1} \\
\alpha\{a_i \mapsto \gamma(b_j \mapsto f_{ij})_j ;_\gamma \delta(g)\}_i
\ar@{=>}[ddd]_{\alpha\{(4),(6),(8),(10)\}} 
&& \gamma(b_j \mapsto \alpha\{a_i \mapsto f_{ij}\}_i)_j ;_\gamma \delta(g)
\ar@{=>}[d]^{(4),(6),(8),(10)} \\ 
&& \delta(\gamma(b_j \mapsto \alpha\{a_i \mapsto f_{ij}\}_i)_j ;_\gamma g)
\ar@{|=|}[d]^{\delta((19);1)} \\ 
&& \delta(\alpha\{a_i \mapsto \gamma(b_j \mapsto f_{ij})_j\}_i ;_\gamma g)
\ar@{=>}[d]^{\delta((3))} \\ 
\alpha\{a_i \mapsto \delta(\gamma(b_j \mapsto f_{ij})_j ;_\gamma g)\}_i 
\ar@{|=|}[rr]_{(19),(17),(15),(13)} 
&& \delta(\alpha\{a_i \mapsto \gamma(b_j \mapsto f_{ij})_j ;_\gamma g\}_i)}
\]

(c) $\beta = \gamma$ and $\delta(g) = \ola{\gamma}[b_k](g)$
\[\xymatrix@M=1ex@C=20ex@R=10ex@!0{
& \alpha\{a_i \mapsto \gamma(b_j \mapsto f_{ij})_j\}_i ;_\gamma 
\ola{\gamma}[b_k](g) \ar@{=>}[dl]_{(3)} \ar@{|=|}[dr]^{(19);1} \\
\alpha\{a_i \mapsto \gamma(b_j \mapsto f_{ij})_j ;_\gamma 
\ola{\gamma}[b_k](g)\}_i \ar@{=>}[d]_{\alpha\{(12)\}} 
&& \gamma(b_j \mapsto \alpha\{a_i \mapsto f_{ij}\}_i)_j ;_\gamma 
\ola{\gamma}[b_k](g) \ar@{=>}[d]^{(12)} \\ 
\alpha\{a_i \mapsto f_{ik} ;_\gamma g\}_i 
&& \alpha\{a_i \mapsto f_{ik}\}_i ;_\gamma g
\ar@{=>}[ll]^{(3)}}
\]

\end{description}
This handles all critical pairs involving (3) and (4). We now look at any 
critical pairs involving (5) and (6). 

\begin{description}
\item [{(5)-(6)}] 
\[\xymatrix@M=1ex@C=20ex@R=10ex@!0{
& \ora{\alpha}[a](f) ;_\gamma \ola{\beta}[b](g) 
\ar@{=>}[dl]_{(5)} \ar@{=>}[dr]^{(6)} \\
\ora{\alpha}[a](f ;_\gamma \ola{\beta}[b](g)) 
\ar@{=>}[d]_{\ora{\alpha}((6))} 
&& \ola{\beta}[b](\ora{\alpha}[a](f) ;_\gamma g) 
\ar@{=>}[d]^{\ola{\beta}((5))} \\
\ora{\alpha}[a](\ola{\beta}[b](f ;_\gamma g)) 
\ar@{|=|}[rr]_{(22)}
&& \ola{\beta}[b](\ora{\alpha}[a](f ;_\gamma g))}
\]

\item [{(5)-(8) [(6)-(7)]}] 
\[\xymatrix@M=1ex@C=20ex@R=10ex@!0{
& \ora{\alpha}[a](f) ;_\gamma \ora{\beta}[b](g) 
\ar@{=>}[dl]_{(5)} \ar@{=>}[dr]^{(8)} \\
\ora{\alpha}[a](f ;_\gamma \ora{\beta}[b](g)) 
\ar@{=>}[d]_{\ora{\alpha}((8))} 
&& \ora{\beta}[b](\ora{\alpha}[a](f) ;_\gamma g) 
\ar@{=>}[d]^{\ora{\beta}((5))} \\
\ora{\alpha}[a](\ora{\beta}[b](f ;_\gamma g)) 
\ar@{|=|}[rr]_{(20)}
&& \ora{\beta}[b](\ora{\alpha}[a](f ;_\gamma g))}
\]

\item [{(5)-(10) [(6)-(9)]}]
\[\xymatrix@M=1ex@C=20ex@R=10ex@!0{
& \ora{\alpha}[a](f) ;_\gamma \beta\{b_i \mapsto g_i\}_i 
\ar@{=>}[dl]_{(5)} \ar@{=>}[dr]^{(10)} \\
\ora{\alpha}[a](f ;_\gamma \beta\{b_i \mapsto g_i\}_i) 
\ar@{=>}[d]_{\ora{\alpha}((10))} 
&& \beta\{b_i \mapsto \ora{\alpha}[a](f) ;_\gamma g_i\}_i
\ar@{=>}[d]^{\beta\{(5)\}} \\
\ora{\alpha}[a](\beta\{b_i \mapsto f ;_\gamma g_i\}_i) 
\ar@{|=|}[rr]_{(15)}
&& \beta\{b_i \mapsto \ora{\alpha}[a](f ;_\gamma g_i)\}_i}
\]

\item [{(5)-(15) [(6)-(16)]}] 
\[\xymatrix@M=1ex@C=20ex@R=10ex@!0{
& \ora{\alpha}[a](\beta\{b_i \mapsto f_i\}_i) ;_\gamma g 
\ar@{=>}[dl]_{(5)} \ar@{|=|}[dr]^{(15);1} \\
\ora{\alpha}[a](\beta\{b_i \mapsto f_i\}_i ;_\gamma g) 
\ar@{=>}[dd]_{\ora{\alpha}((3))} 
&& \beta\{b_i \mapsto \ora{\alpha}[a](f_i)\}_i ;_\gamma g
\ar@{=>}[d]^{(3)} \\
&& \beta\{b_i \mapsto \ora{\alpha}[a](f_i) ;_\gamma g\}_i
\ar@{=>}[d]^{\beta\{(5)\}} \\
\ora{\alpha}[a](\beta\{b_i \mapsto f_i ;_\gamma g\}_i) 
\ar@{|=|}[rr]_{(15)}
&& \beta\{b_i \mapsto \ora{\alpha}[a](f_i ;_\gamma g)\}_i}
\]

\item [{(5)-(18) [(6)-(17)]}] 
\[\xymatrix@M=1ex@C=20ex@R=10ex@!0{
& \ora{\alpha}[a](\beta(b_i \mapsto f_i)_i) ;_\gamma \delta(g)
\ar@{=>}[dl] \ar@{|=|}[dr] \\
\ora{\alpha}[a](\beta(b_i \mapsto f_i)_i ;_\gamma \delta(g)) &&
\beta(b_i \mapsto \ora{\alpha}[a](f_i))_i ;_\gamma \delta(g)}
\]
There are three subcases to consider.

(a) $\beta \neq \gamma$.
\[\xymatrix@M=1ex@C=20ex@R=10ex@!0{
& \ora{\alpha}[a](\beta(b_i \mapsto f_i)_i) ;_\gamma \delta(g)
\ar@{=>}[dl]_{(5)} \ar@{|=|}[dr]^{(18);1} \\
\ora{\alpha}[a](\beta(b_i \mapsto f_i)_i ;_\gamma \delta(g)) 
\ar@{=>}[dd]_{\ora{\alpha}((3))} 
&& \beta(b_i \mapsto \ora{\alpha}[a](f_i))_i ;_\gamma \delta(g)
\ar@{=>}[d]^{(9)} \\
&& \beta(b_i \mapsto \ora{\alpha}[a](f_i) ;_\gamma \delta(g))_i
\ar@{=>}[d]^{\beta((5))} \\
\ora{\alpha}[a](\beta(b_i \mapsto f_i ;_\gamma \delta(g))_i) 
\ar@{|=|}[rr]_{(18)}
&& \beta(b_i \mapsto \ora{\alpha}[a](f_i ;_\gamma \delta(g)))_i}
\]

(b) $\beta = \gamma$ and $\delta \neq \gamma$. 
\[\xymatrix@M=1ex@C=22ex@R=10ex@!0{
& \ora{\alpha}[a](\gamma(b_i \mapsto f_i)_i) ;_\gamma \delta(g)
\ar@{=>}[dl]_{(5)} \ar@{|=|}[dr]^{(18);1} \\
\ora{\alpha}[a](\gamma(b_i \mapsto f_i)_i ;_\gamma \delta(g))
\ar@{=>}[ddd]_{\ora{\alpha}((4),(6),(8),(10))} 
&& \gamma(b_i \mapsto \ora{\alpha}[a](f_i))_i ;_\gamma \delta(g)
\ar@{=>}[d]^{(4),(6),(8),(10)} \\ 
&& \delta(\gamma(b_i \mapsto \ora{\alpha}[a](f_i))_i ;_\gamma g)
\ar@{|=|}[d]^{\delta((18);1)} \\ 
&& \delta(\ora{\alpha}[a](\gamma(b_i \mapsto f_i)_i) ;_\gamma g)
\ar@{=>}[d]^{\delta((5))} \\ 
\ora{\alpha}[a](\delta(\gamma(b_i \mapsto f_i)_i ;_\gamma g))
\ar@{|=|}[rr]_{(18),(22),(20),(15)} 
&& \delta(\ora{\alpha}[a](\gamma(b_i \mapsto f_i)_i ;_\gamma g))}
\]

(c) $\beta = \gamma$ and $\delta(g) = \ola{\gamma}[b_k](g)$.
\[\xymatrix@M=1ex@C=20ex@R=10ex@!0{
& \ora{\alpha}[a](\gamma(b_i \mapsto f_i)_i) ;_\gamma \ola{\gamma}[b_k](g)
\ar@{=>}[dl]_{(5)} \ar@{|=|}[dr]^{(18);1} \\
\ora{\alpha}[a](\gamma(b_i \mapsto f_i)_i ;_\gamma \ola{\gamma}[b_k](g))
\ar@{=>}[d]_{\ora{\alpha}((12))} 
&& \gamma(b_i \mapsto \ora{\alpha}[a](f_i))_i ;_\gamma \ola{\gamma}[b_k](g)
\ar@{=>}[d]^{(12)} \\
\ora{\alpha}[a](f_k ;_\gamma g)
&& \ora{\alpha}[a](f_k) ;_\gamma g
\ar@{=>}[ll]^{(5)}}
\]

\item [{(5)-(20) [(6)-(21)]}] 
\[\ora{\alpha}[a](\ora{\beta}[b_k](f) ;_\gamma g) \Lla 
\ora{\alpha}[a](\ora{\beta}[b_k](f)) ;_\gamma g \pc
\ora{\beta}[b_k](\ora{\alpha}[a](f)) ;_\gamma g
\]
There are two subcases to consider.

(a) $\beta \neq \gamma$.
\[\xymatrix@M=1ex@C=20ex@R=10ex@!0{
& \ora{\alpha}[a](\ora{\beta}[b_k](f)) ;_\gamma g 
\ar@{=>}[dl]_{(5)} \ar@{|=|}[dr]^{(20);1} \\
\ora{\alpha}[a](\ora{\beta}[b_k](f) ;_\gamma g) 
\ar@{=>}[dd]_{\ora{\alpha}((5))} 
&& \ora{\beta}[b_k](\ora{\alpha}[a](f)) ;_\gamma g
\ar@{=>}[d]^{(5)} \\
&& \ora{\beta}[b_k](\ora{\alpha}[a](f) ;_\gamma g)
\ar@{=>}[d]^{\ora{\beta}((5))} \\
\ora{\alpha}[a](\ora{\beta}[b_k](f ;_\gamma g)) 
\ar@{|=|}[rr]_{(20)}
&& \ora{\beta}[b_k](\ora{\alpha}[a](f ;_\gamma g))}
\]

(b) $\beta=\gamma$. Here we may assume that $g=\gamma\{b_i \mapsto g_i\}_i$.
\[\xymatrix@M=1ex@C=20ex@R=10ex@!0{
&\ora{\alpha}[a](\ora{\gamma}[b_k](f)) ;_\gamma \gamma\{b_i \mapsto g_i\}_i
\ar@{=>}[dl]_{(5)} \ar@{|=|}[dr]^{(20);1} \\
\ora{\alpha}[a](\ora{\gamma}[b_k](f) ;_\gamma \gamma\{b_i \mapsto g_i\}_i) 
\ar@{=>}[d]_{\ora{\alpha}((11))} 
&& \ora{\gamma}[b_k](\ora{\alpha}[a](f)) ;_\gamma \gamma\{b_i \mapsto g_i
\}_i \ar@{=>}[d]^{(11)} \\
\ora{\alpha}[a](f ;_\gamma g_k) 
&& \ora{\alpha}[a](f) ;_\gamma g_k
\ar@{=>}[ll]^{(5)}}
\]

\item [{(5)-(22) [(6)-(21)]}] 
\[\xymatrix@M=1ex@C=20ex@R=10ex@!0{
& \ora{\alpha}[a](\ola{\beta}[b](f)) ;_\gamma g 
\ar@{=>}[dl]_{(5)} \ar@{|=|}[dr]^{(22);1} \\
\ora{\alpha}[a](\ola{\beta}[b](f) ;_\gamma g) 
\ar@{=>}[dd]_{\ora{\alpha}((7))} 
&& \ola{\beta}[b](\ora{\alpha}[a](f)) ;_\gamma g
\ar@{=>}[d]^{(7)} \\
&& \ola{\beta}[b](\ora{\alpha}[a](f) ;_\gamma g)
\ar@{=>}[d]^{\ola{\beta}((5))} \\
\ora{\alpha}[a](\ola{\beta}[b](f ;_\gamma g)) 
\ar@{|=|}[rr]_{(22)}
&& \ola{\beta}[b](\ora{\alpha}[a](f ;_\gamma g))}
\]

\end{description}
This handles all critical pairs involving (5) and (6). We now look at any 
critical pairs involving (7) and (8). 

\begin{description} 
\item [{(7)-(8)}] 
\[\xymatrix@M=1ex@C=20ex@R=10ex@!0{
& \ola{\alpha}[a](f) ;_\gamma \ora{\beta}[b](g) 
\ar@{=>}[dl]_{(7)} \ar@{=>}[dr]^{(8)} \\
\ola{\alpha}[a](f ;_\gamma \ora{\beta}[b](g)) 
\ar@{=>}[d]_{\ola{\alpha}((8))} 
&& \ora{\beta}[b](\ola{\alpha}[a](f) ;_\gamma g) 
\ar@{=>}[d]^{\ora{\beta}((7))} \\
\ola{\alpha}[a](\ora{\beta}[b](f ;_\gamma g)) 
\ar@{|=|}[rr]_{(22)}
&& \ora{\beta}[b](\ola{\alpha}[a](f ;_\gamma g))}
\]

\item [{(7)-(10) [(8)-(9)]}]
\[\xymatrix@M=1ex@C=20ex@R=10ex@!0{
& \ola{\alpha}[a](f) ;_\gamma \beta\{b_i \mapsto g_i\}_i 
\ar@{=>}[dl]_{(7)} \ar@{=>}[dr]^{(10)} \\
\ola{\alpha}[a](f ;_\gamma \beta\{b_i \mapsto g_i\}_i) 
\ar@{=>}[d]_{\ola{\alpha}((10))} 
&& \beta\{b_i \mapsto \ola{\alpha}[a](f) ;_\gamma g_i\}_i
\ar@{=>}[d]^{\beta\{(7)\}} \\
\ola{\alpha}[a](\beta\{b_i \mapsto f ;_\gamma g_i\}_i) 
\ar@{|=|}[rr]_{(17)}
&& \beta\{b_i \mapsto \ola{\alpha}[a](f ;_\gamma g)\}_i}
\]

\item [{(7)-(16) [(8)-(15)]}] 
\[\xymatrix@M=1ex@C=20ex@R=10ex@!0{
& \ola{\alpha}[a](\beta(b_i \mapsto f_i)_i) ;_\gamma \delta(g)
\ar@{=>}[dl] \ar@{|=|}[dr] \\
\ola{\alpha}[a](\beta(b_i \mapsto f_i)_i ;_\gamma \delta(g)) &&
\beta(b_i \mapsto \ola{\alpha}[a](f_i))_i ;_\gamma \delta(g)}
\]
There are three subcases to consider.

(a) $\beta \neq \gamma$
\[\xymatrix@M=1ex@C=20ex@R=10ex@!0{
& \ola{\alpha}[a](\beta(b_i \mapsto f_i)_i) ;_\gamma \delta(g)
\ar@{=>}[dl]_{(7)} \ar@{|=|}[dr]^{(16);1} \\
\ola{\alpha}[a](\beta(b_i \mapsto f_i)_i ;_\gamma \delta(g))
\ar@{=>}[dd]_{\ola{\alpha}((9))} 
&& \beta(b_i \mapsto \ola{\alpha}[a](f_i))_i ;_\gamma \delta(g)
\ar@{=>}[d]^{(9)} \\
&& \beta(b_i \mapsto \ola{\alpha}[a](f_i) ;_\gamma \delta(g))_i
\ar@{=>}[d]^{\beta((5))} \\
\ola{\alpha}[a](\beta(b_i \mapsto f_i ;_\gamma \delta(g))_i)
\ar@{|=|}[rr]_{(16)}
&& \beta(b_i \mapsto \ola{\alpha}[a](f_i ;_\gamma \delta(g)))_i}
\]

(b) $\beta = \gamma$ and $\delta \neq \gamma$.
\[\xymatrix@M=1ex@C=22ex@R=10ex@!0{
& \ola{\alpha}[a](\gamma(b_i \mapsto f_i)_i) ;_\gamma \delta(g)
\ar@{=>}[dl]_{(7)} \ar@{|=|}[dr]^{(16);1} \\
\ola{\alpha}[a](\gamma(b_i \mapsto f_i)_i ;_\gamma \delta(g))
\ar@{=>}[ddd]_{\ola{\alpha}((4),(6),(8),(10))} 
&& \gamma(b_i \mapsto \ola{\alpha}[a](f_i))_i ;_\gamma \delta(g)
\ar@{=>}[d]^{(4),(6),(8),(10)} \\
&& \delta(\gamma(b_i \mapsto \ola{\alpha}[a](f_i))_i ;_\gamma g)
\ar@{|=|}[d]^{\delta((16);1)} \\
&& \delta(\ola{\alpha}[a](\gamma(b_i \mapsto f_i)_i) ;_\gamma g)
\ar@{=>}[d]^{\delta((7))} \\
\ola{\alpha}[a](\delta(\gamma(b_i \mapsto f_i)_i ;_\gamma g))
\ar@{|=|}[rr]_{(16),(21),(22),(17)}
&& \delta(\ola{\alpha}[a](\gamma(b_i \mapsto f_i)_i ;_\gamma g))}
\]

(c) $\beta = \gamma$ and $\delta(g) = \ola{\gamma}[b_k](g)$.
\[\xymatrix@M=1ex@C=20ex@R=10ex@!0{
& \ola{\alpha}[a](\gamma(b_i \mapsto f_i)_i) ;_\gamma \ola{\gamma}[b_k](g)
\ar@{=>}[dl]_{(7)} \ar@{|=|}[dr]^{(16);1} \\
\ola{\alpha}[a](\gamma(b_i \mapsto f_i)_i ;_\gamma \ola{\gamma}[b_k](g))
\ar@{=>}[d]_{\ola{\alpha}((12))} 
&& \gamma(b_i \mapsto \ola{\alpha}[a](f_i))_i ;_\gamma \ola{\gamma}[b_k](g)
\ar@{=>}[d]^{(12)} \\
\ola{\alpha}[a](f_k ;_\gamma g)
&& \ola{\alpha}[a](f_k) ;_\gamma g
\ar@{=>}[ll]^{(7)}}
\]

\item [{(7)-(17) [(8)-(18)]}] 
\[\xymatrix@M=1ex@C=20ex@R=10ex@!0{
& \ola{\alpha}[a](\beta\{b_i \mapsto f_i\}_i) ;_\gamma g
\ar@{=>}[dl]_{(7)} \ar@{|=|}[dr]^{(17);1} \\
\ola{\alpha}[a](\beta\{b_i \mapsto f_i\}_i ;_\gamma g)
\ar@{=>}[dd]_{\ola{\alpha}((3))} 
&& \beta\{b_i \mapsto \ola{\alpha}[a](f_i)\}_i ;_\gamma g
\ar@{=>}[d]^{(3)} \\
&& \beta\{b_i \mapsto \ola{\alpha}[a](f_i) ;_\gamma g\}_i
\ar@{=>}[d]^{\beta\{(5)\}} \\
\ola{\alpha}[a](\beta\{b_i \mapsto f_i ;_\gamma g\}_i)
\ar@{|=|}[rr]_{(17)}
&& \beta\{b_i \mapsto \ola{\alpha}[a](f_i ;_\gamma g)\}_i}
\]

\item [{(7)-(21) [(8)-(20)]}] 
\[\xymatrix@M=1ex@C=20ex@R=10ex@!0{
& \ola{\alpha}[a](\ola{\beta}[b](f)) ;_\gamma g
\ar@{=>}[dl]_{(7)} \ar@{|=|}[dr]^{(21);1} \\
\ola{\alpha}[a](\ola{\beta}[b](f) ;_\gamma g)
\ar@{=>}[dd]_{\ola{\alpha}((3))} 
&& \ola{\beta}[b](\ola{\alpha}[a](f)) ;_\gamma g
\ar@{=>}[d]^{(3)} \\
&& \ola{\beta}[b](\ola{\alpha}[a](f) ;_\gamma g)
\ar@{=>}[d]^{\ola{\beta}((5))} \\
\ola{\alpha}[a](\ola{\beta}[b](f ;_\gamma g))
\ar@{|=|}[rr]_{(21)}
&& \ola{\beta}[b](\ola{\alpha}[a](f) ;_\gamma g))}
\]

\item [{(7)-(22) [(8)-(22)]}] 
\[\ola{\alpha}[a](\ora{\beta}[b_k](f) ;_\gamma g) \Lla
\ola{\alpha}[a](\ora{\beta}[b_k](f)) ;_\gamma g \pc
\ora{\beta}[b_k](\ola{\alpha}[a](f)) ;_\gamma g
\]
There are two subcases to consider.

(a) $\beta \neq \gamma$.
\[\xymatrix@M=1ex@C=20ex@R=10ex@!0{
& \ola{\alpha}[a](\ora{\beta}[b_k](f)) ;_\gamma g
\ar@{=>}[dl]_{(7)} \ar@{|=|}[dr]^{(22);1} \\
\ola{\alpha}[a](\ora{\beta}[b_k](f) ;_\gamma g)
\ar@{=>}[dd]_{\ola{\alpha}((3))} 
&& \ora{\beta}[b_k](\ola{\alpha}[a](f)) ;_\gamma g
\ar@{=>}[d]^{(3)} \\
&& \ora{\beta}[b_k](\ola{\alpha}[a](f) ;_\gamma g)
\ar@{=>}[d]^{\ora{\beta}((5))} \\
\ola{\alpha}[a](\ora{\beta}[b_k](f ;_\gamma g))
\ar@{|=|}[rr]_{(22)}
&& \ora{\beta}[b_k](\ola{\alpha}[a](f ;_\gamma g))}
\]

(b) $\beta = \gamma$. Here we may assume that $g = \gamma\{b_i \mapsto g_i\}_i$.
\[\xymatrix@M=1ex@C=20ex@R=10ex@!0{
& \ola{\alpha}[a](\ora{\gamma}[b_k](f)) ;_\gamma \gamma\{b_i \mapsto g_i\}_i
\ar@{=>}[dl]_{(7)} \ar@{|=|}[dr]^{(22);1} \\
\ola{\alpha}[a](\ora{\gamma}[b_k](f) ;_\gamma \gamma\{b_i \mapsto g_i\}_i)
\ar@{=>}[d]_{\ola{\alpha}((11))} 
&& \ora{\gamma}[b_k](\ola{\alpha}[a](f)) ;_\gamma \gamma\{b_i \mapsto g_i\}_i
\ar@{=>}[d]^{(11)} \\
\ola{\alpha}[a](f ;_\gamma g_k)
&& \ola{\alpha}[a](f) ;_\gamma g_k \ar@{=>}[ll]^{(7)}}
\]

\end{description}
This handles all critical pairs involving (7) and (8). We move on now to look
at any critical pairs involving (9) and (10).

\begin{description}
\item [{(9)-(10)}] 
\[\xymatrix@M=1ex@C=20ex@R=10ex@!0{
& \alpha(a_i \mapsto f_i)_i ;_\gamma \beta\{b_j \mapsto g_j\}_j
\ar@{=>}[dl]_{(9)} \ar@{=>}[dr]^{(10)} \\
\alpha(a_i \mapsto f_i ;_\gamma \beta\{b_j \mapsto g_j\}_j)_i 
\ar@{=>}[d]_{\alpha((10))} 
&& \beta\{b_j \mapsto \alpha(a_i \mapsto f_i)_i ;_\gamma g_j\}_j
\ar@{=>}[d]^{\beta\{(9)\}} \\
\alpha(a_i \mapsto \beta\{b_j \mapsto f_i ;_\gamma g_j\}_j)_i 
\ar@{|=|}[rr]_{(19)} 
&& \beta\{b_j \mapsto \alpha(a_i \mapsto f_i ;_\gamma g_j)_i\}_j}
\]

\item [{(9)-(14) [(10)-(13)]}] 
\[\xymatrix@M=1ex@C=20ex@R=10ex@!0{
& \alpha(a_i \mapsto \beta(b_j \mapsto f_{ij})_j)_i ;_\gamma \delta(g)
\ar@{=>}[dl] \ar@{|=|}[dr] \\
\alpha(a_i \mapsto \beta(b_j \mapsto f_{ij})_j ;_\gamma \delta(g))_i &&
\beta(b_j \mapsto \alpha(a_i \mapsto f_{ij})_i)_j ;_\gamma \delta(g)}
\]
There are three subcases to consider.

(a) $\beta \neq \gamma$.
\[\xymatrix@M=1ex@C=20ex@R=10ex@!0{
& \alpha(a_i \mapsto \beta(b_j \mapsto f_{ij})_j)_i ;_\gamma \delta(g)
\ar@{=>}[dl]_{(9)} \ar@{|=|}[dr]^{(14);1} \\
\alpha(a_i \mapsto \beta(b_j \mapsto f_{ij})_j ;_\gamma \delta(g))_i
\ar@{=>}[dd]_{\alpha((9))} 
&& \beta(b_j \mapsto \alpha(a_i \mapsto f_{ij})_i)_j ;_\gamma \delta(g)
\ar@{=>}[d]^{(9)} \\
&& \beta(b_j \mapsto \alpha(a_i \mapsto f_{ij})_i ;_\gamma \delta(g))_j
\ar@{=>}[d]^{\beta((9))} \\
\alpha(a_i \mapsto \beta(b_j \mapsto f_{ij} ;_\gamma \delta(g))_j)_i
\ar@{|=|}[rr]_{(14)} && 
\beta(b_j \mapsto \alpha(a_i \mapsto f_{ij} ;_\gamma \delta(g))_i)_j}
\]

(b) $\beta = \gamma$ and $\delta \neq \gamma$.
\[\xymatrix@M=1ex@C=22ex@R=10ex@!0{
& \alpha(a_i \mapsto \gamma(b_j \mapsto f_{ij})_j)_i ;_\gamma \delta(g)
\ar@{=>}[dl]_{(9)} \ar@{|=|}[dr]^{(14);1} \\
\alpha(a_i \mapsto \gamma(b_j \mapsto f_{ij})_j ;_\gamma \delta(g))_i
\ar@{=>}[ddd]_{\alpha\{(4),(6),(8),(10)\}} 
&& \gamma(b_j \mapsto \alpha(a_i \mapsto f_{ij})_i)_j ;_\gamma \delta(g)
\ar@{=>}[d]^{(4),(6),(8),(10)} \\ 
&& \delta(\gamma(b_j \mapsto \alpha(a_i \mapsto f_{ij})_i)_j ;_\gamma g)
\ar@{|=|}[d]^{\delta((14);1)} \\ 
&& \delta(\alpha(a_i \mapsto \gamma(b_j \mapsto f_{ij})_j)_i ;_\gamma g)
\ar@{=>}[d]^{\delta((9))} \\ 
\alpha(a_i \mapsto \delta(\gamma(b_j \mapsto f_{ij})_j ;_\gamma g))_i 
\ar@{|=|}[rr]_{(14),(16),(18),(19)} 
&& \delta(\alpha(a_i \mapsto \gamma(b_j \mapsto f_{ij})_j ;_\gamma g)_i)}
\]

(c) $\beta = \gamma$ and $\delta(g) = \ola{\gamma}[b_k](g)$
\[\xymatrix@M=1ex@C=20ex@R=10ex@!0{
& \alpha(a_i \mapsto \gamma(b_j \mapsto f_{ij})_j)_i ;_\gamma 
\ola{\gamma}[b_k](g) \ar@{=>}[dl]_{(9)} \ar@{|=|}[dr]^{(14);1} \\
\alpha(a_i \mapsto \gamma(b_j \mapsto f_{ij})_j ;_\gamma 
\ola{\gamma}[b_k](g))_i \ar@{=>}[d]_{\alpha((12))} 
&& \gamma(b_j \mapsto \alpha(a_i \mapsto f_{ij})_i)_j ;_\gamma 
\ola{\gamma}[b_k](g) \ar@{=>}[d]^{(12)} \\ 
\alpha(a_i \mapsto f_{ik} ;_\gamma g)_i 
&& \alpha(a_i \mapsto f_{ik})_i ;_\gamma g
\ar@{=>}[ll]^{(9)}}
\]

\item [{(9)-(16) [(10)-(15)]}] 
\[\xymatrix@M=1ex@C=20ex@R=10ex@!0{
& \alpha(a_i \mapsto \ola{\beta}[b](f_i))_i ;_\gamma g 
\ar@{=>}[dl]_{(9)} \ar@{|=|}[dr]^{(16);1} \\
\alpha(a_i \mapsto \ola{\beta}[b](f_i) ;_\gamma g)_i 
\ar@{=>}[dd]_{\alpha((7))} 
&& \ola{\beta}[b](\alpha(a_i \mapsto f_i)_i) ;_\gamma g
\ar@{=>}[d]^{(7)} \\ 
&& \ola{\beta}[b](\alpha(a_i \mapsto f_i)_i ;_\gamma g)
\ar@{=>}[d]^{\ola{\beta}((9))} \\ 
\alpha(a_i \mapsto \ola{\beta}[b](f_i ;_\gamma g))_i 
\ar@{|=|}[rr]_{(16)} && 
\ola{\beta}[b](\alpha(a_i \mapsto f_i ;_\gamma g)_i)}
\]

\item [{(9)-(18) [(10)-(17)]}] 
\[\alpha(a_i \mapsto \ora{\beta}[b_k](f_i) ;_\gamma g)_i \Lla 
\alpha(a_i \mapsto \ora{\beta}[b_k](f_i))_i ;_\gamma g \pc
\ora{\beta}[b_k](\alpha(a_i \mapsto f_i)_i) ;_\gamma g 
\]
There are two subcases to consider.

(a) $\beta \neq \gamma$.
\[\xymatrix@M=1ex@C=20ex@R=10ex@!0{
& \alpha(a_i \mapsto \ora{\beta}[b_k](f_i))_i ;_\gamma g 
\ar@{=>}[dl]_{(9)} \ar@{|=|}[dr]^{(18);1} \\
\alpha(a_i \mapsto \ora{\beta}[b_k](f_i) ;_\gamma g)_i 
\ar@{=>}[dd]_{\alpha((5))} 
&& \ora{\beta}[b_k](\alpha(a_i \mapsto f_i)_i) ;_\gamma g
\ar@{=>}[d]^{(5)} \\ 
&& \ora{\beta}[b_k](\alpha(a_i \mapsto f_i)_i ;_\gamma g)
\ar@{=>}[d]^{\ora{\beta}((9))} \\ 
\alpha(a_i \mapsto \ora{\beta}[b_k](f_i ;_\gamma g))_i 
\ar@{|=|}[rr]_{(18)} 
&& \ora{\beta}[b_k](\alpha(a_i \mapsto f_i ;_\gamma g)_i)}
\]

(b) $\beta = \gamma$. Here we may assume that $g = \gamma\{b_j \mapsto g_j
\}_j$.
\[\xymatrix@M=1ex@C=20ex@R=10ex@!0{
& \alpha(a_i \mapsto \ora{\gamma}[b_k](f_i))_i ;_\gamma \gamma\{b_j \mapsto
g_j\}_j \ar@{=>}[dl]_{(9)} \ar@{|=|}[dr]^{(18);1} \\
\alpha(a_i \mapsto \ora{\gamma}[b_k](f_i) ;_\gamma \gamma\{b_j \mapsto g_j
\}_j)_i \ar@{=>}[d]_{\alpha((11))} 
&& \ora{\gamma}[b_k](\alpha(a_i \mapsto f_i)_i) ;_\gamma \gamma\{b_j \mapsto
g_j\}_j \ar@{=>}[d]^{(11)} \\ 
\alpha(a_i \mapsto f_i ;_\gamma g_k)_i 
&& \alpha(a_i \mapsto f_i)_i ;_\gamma g_k
\ar@{=>}[ll]^{(9)}}
\]

\item [{(9)-(19) [(10)-(19)]}] 
\[\xymatrix@M=1ex@C=20ex@R=10ex@!0{
& \alpha(a_i \mapsto \beta\{b_j \mapsto f_{ij}\}_j)_i ;_\gamma g 
\ar@{=>}[dl]_{(9)} \ar@{|=|}[dr]^{(19);1} \\
\alpha(a_i \mapsto \beta\{b_j \mapsto f_{ij}\}_j ;_\gamma g)_i 
\ar@{=>}[dd]_{\alpha((3))} 
&& \beta\{b_j \mapsto \alpha(a_i \mapsto f_{ij})_i\}_j ;_\gamma g 
\ar@{=>}[d]^{(3)} \\ 
&& \beta\{b_j \mapsto \alpha(a_i \mapsto f_{ij})_i ;_\gamma g\}_j
\ar@{=>}[d]^{\beta\{(9)\}} \\ 
\alpha(a_i \mapsto \beta\{b_j \mapsto f_{ij} ;_\gamma g\}_j)_i 
\ar@{|=|}[rr]_{(19)} 
&& \beta\{b_j \mapsto \alpha(a_i \mapsto f_{ij} ;_\gamma g)_i\}_j}
\]

\end{description}
This handles all critical pairs involving (9) and (10). We move on now to look
at any critical pairs involving (11) and (12). 

\begin{description}
\item [{(11)-(13) [(12)-(14)]}]
\[\xymatrix@M=1ex@C=20ex@R=10ex@!0{
& \ora{\gamma}[a_k](f) ;_\gamma \gamma\{a_i \mapsto \alpha\{b_j \mapsto
g_{ij}\}_j\}_i \ar@{=>}[dl]_{(11)} \ar@{|=|}[dr]^{1;(13)} \\
f ;_\gamma \alpha\{b_j \mapsto g_{kj}\}_j \ar@{=>}[d]_{(10)}
&& \ora{\gamma}[a_k](f) ;_\gamma \alpha\{b_j \mapsto \gamma\{a_i \mapsto
g_{ij}\}_i\}_j \ar@{=>}[d]^{(10)} \\
\alpha\{b_j \mapsto f ;_\gamma g_{kj}\}_j
&& \alpha\{b_j \mapsto \ora{\gamma}[a_k](f) ;_\gamma \gamma\{a_i \mapsto
g_{ij}\}_i\}_j
\ar@{=>}[ll]^-{\alpha\{(11)\}}}
\]

\item [{(11)-(15) [(12)-(16)]}] There are two subcases on the structure.

(a) $\ora{\gamma}[a_k](\alpha\{b_j \mapsto f_j\}_j) ;_\gamma \gamma\{a_i 
\mapsto g_i\}_i$
\[\xymatrix@M=1ex@C=20ex@R=10ex@!0{
& \ora{\gamma}[a_k](\alpha\{b_j \mapsto f_j\}_j) ;_\gamma \gamma\{a_i \mapsto
g_i\}_i \ar@{=>}[dl]_{(11)} \ar@{|=|}[dr]^{(15);1} \\
\alpha\{b_j \mapsto f_j\}_j ;_\gamma g_k \ar@{=>}[d]_{(3)}
&& \alpha\{b_j \mapsto \ora{\gamma}[a_k](f_j)\}_j ;_\gamma \gamma\{a_i \mapsto
g_i\}_i \ar@{=>}[d]^{(3)} \\ 
\alpha\{b_j \mapsto f_j ;_\gamma g_k\}_j
&& \alpha\{b_j \mapsto \ora{\gamma}[a_k](f_j) ;_\gamma \gamma\{a_i \mapsto
g_i\}_i\}_j \ar@{=>}[ll]^-{\alpha\{(11)\}}}
\]

(b) $\ora{\gamma}[a_k](f) ;_\gamma \gamma\{a_i \mapsto 
\ora{\alpha}[b](g_i)\}_i$
\[\xymatrix@M=1ex@C=20ex@R=10ex@!0{
& \ora{\gamma}[a_k](f) ;_\gamma \gamma\{a_i \mapsto 
\ora{\alpha}[b](g_i)\}_i \ar@{=>}[dl]_{(11)} \ar@{|=|}[dr]^{1;(15)} \\
f ;_\gamma \ora{\alpha}[b](g_k) \ar@{=>}[d]_{(8)} 
&& \ora{\gamma}[a_k](f) ;_\gamma \ora{\alpha}[b](\gamma\{a_i \mapsto 
g_i\}_i) \ar@{=>}[d]^{(8)} \\
\ora{\alpha}[b](f ;_\gamma g_k) 
&& \ora{\alpha}[b](\ora{\gamma}[a_k](f) ;_\gamma \gamma\{a_i \mapsto 
g_i\}_i) \ar@{=>}[ll]^-{\ora{\alpha}((11))}}
\]

\item [{(11)-(17) [(12)-(18)]}]
\[\xymatrix@M=1ex@C=20ex@R=10ex@!0{
& \ora{\gamma}[a_k](f) ;_\gamma \gamma\{a_i \mapsto 
\ola{\alpha}[b](g_i)\}_i \ar@{=>}[dl]_{(11)} \ar@{|=|}[dr]^{1;(17)} \\
f ;_\gamma \ola{\alpha}[b](g_k) \ar@{=>}[d]_{(6)} 
&& \ora{\gamma}[a_k](f) ;_\gamma \ola{\alpha}[b](\gamma\{a_i \mapsto 
g_i\}_i) \ar@{=>}[d]^{(6)} \\
\ola{\alpha}[b](f ;_\gamma g_k) 
&& \ola{\alpha}[b](\ora{\gamma}[a_k](f) ;_\gamma \gamma\{a_i \mapsto 
g_i\}_i) \ar@{=>}[ll]^-{\ola{\alpha}((11))}}
\]

\item [{(11)-(18) [(12)-(17)]}]
\[\xymatrix@M=1ex@C=20ex@R=10ex@!0{
& \ora{\gamma}[a_k](\alpha(b_j \mapsto f_j)_j) ;_\gamma \gamma\{a_i \mapsto
g_i\}_i \ar@{=>}[dl]_{(11)} \ar@{|=|}[dr]^{(18);1} \\
\alpha(b_j \mapsto f_j)_j ;_\gamma g_k \ar@{=>}[d]_{(9)}
&& \alpha(b_j \mapsto \ora{\gamma}[a_k](f_j))_j ;_\gamma \gamma\{a_i \mapsto
g_i\}_i \ar@{=>}[d]^{(9)} \\ 
\alpha(b_j \mapsto f_j ;_\gamma g_k)_j
&& \alpha(b_j \mapsto \ora{\gamma}[a_k](f_j) ;_\gamma \gamma\{a_i \mapsto
g_i\}_i)_j \ar@{=>}[ll]^-{\alpha((11))}}
\]

\item [{(11)-(19) [(12)-(19)]}]
\[\xymatrix@M=1ex@C=20ex@R=10ex@!0{
& \ora{\gamma}[a_k](f) ;_\gamma \gamma\{a_i \mapsto \alpha(b_j \mapsto 
g_{ij})_j\}_i \ar@{=>}[dl]_{(11)} \ar@{|=|}[dr]^{1;(19)} \\
f ;_\gamma \alpha(b_j \mapsto g_{kj})_j \ar@{=>}[d]_{(4)} 
&& \ora{\gamma}[a_k](f) ;_\gamma \alpha(b_j \mapsto \gamma\{a_i \mapsto 
g_{ij}\}_i)_j \ar@{=>}[d]^{(4)} \\ 
\alpha(b_j \mapsto f ;_\gamma g_{kj})_j 
&& \alpha(b_j \mapsto \ora{\gamma}[a_k](f) ;_\gamma \gamma\{a_i \mapsto 
g_{ij}\}_i)_j 
\ar@{=>}[ll]^-{\alpha((11))}}
\]

\item [{(11)-(20) [(12)-(21)]}]
\[\xymatrix@M=1ex@C=20ex@R=10ex@!0{
& \ora{\gamma}[a_k](\ora{\alpha}[b](f)) ;_\gamma \gamma\{a_i \mapsto
g_i\}_i \ar@{=>}[dl]_{(11)} \ar@{|=|}[dr]^{(20);1} \\
\ora{\alpha}[b](f) ;_\gamma g_k \ar@{=>}[d]_{(5)}
&& \ora{\alpha}[b](\ora{\gamma}[a_k](f)) ;_\gamma \gamma\{a_i \mapsto
g_i\}_i \ar@{=>}[d]^{(5)} \\ 
\ora{\alpha}[b](f ;_\gamma g_k)
&& \ora{\alpha}[b](\ora{\gamma}[a_k](f) ;_\gamma \gamma\{a_i \mapsto
g_i\}_i) \ar@{=>}[ll]^-{\ora{\alpha}((11))}}
\]

\item [{(11)-(22) [(12)-(22)]}]
\[\xymatrix@M=1ex@C=20ex@R=10ex@!0{
& \ora{\gamma}[a_k](\ola{\alpha}[b](f)) ;_\gamma \gamma\{a_i \mapsto g_i\}_i
\ar@{=>}[dl]_{(11)} \ar@{|=|}[dr]^{(22);1} \\
\ola{\alpha}[b](f) ;_\gamma g_k \ar@{=>}[d]_{(7)}
&& \ola{\alpha}[b](\ora{\gamma}[a_k](f)) ;_\gamma \gamma\{a_i \mapsto g_i\}_i
\ar@{=>}[d]^{(7)} \\ 
\ola{\alpha}[b](f ;_\gamma g_k)
&& \ola{\alpha}[b](\ora{\gamma}[a_k](f) ;_\gamma \gamma\{a_i \mapsto g_i\}_i)
\ar@{=>}[ll]^-{\ola{\alpha}((11))}}
\]

\end{description}
This handles all critical pairs involving (11) and (12).

\end{spacing}
\end{document}